\documentclass[reqno,12pt,thmsa]{article}
\usepackage[left=1.5cm,right=1.5cm,top=2cm,bottom=2cm]{geometry}
\usepackage{amsmath, amssymb, graphics, enumerate}
\usepackage{algorithmicx, algpseudocode}
\usepackage{algorithm}
\usepackage{titlesec}
\usepackage[font={small}]{caption}
\usepackage{subcaption}
\usepackage{titlesec}
\usepackage{float}
\usepackage{tocloft} % for TOC tweaks (optional but useful)
\usepackage{xcolor}
\usepackage{graphicx} 
\usepackage{multirow}
\titleclass{\subsubsubsection}{straight}[\subsubsection]

\newcounter{subsubsubsection}[subsubsection]
\renewcommand\thesubsubsubsection{\thesubsubsection.\arabic{subsubsubsection}}

\titleformat{\subsubsubsection}
  {\normalfont\normalsize\bfseries}{\thesubsubsubsection}{1em}{}

\titlespacing*{\subsubsubsection}
{0pt}{3.25ex plus 1ex minus .2ex}{1.5ex plus .2ex}

\makeatletter

\newcommand\l@subsubsubsection{\@dottedtocline{4}{7em}{4em}} 
\makeatother

\usepackage[font={small}]{caption}
\newcommand{\mathsym}[1]{{}}

\usepackage{color}
\usepackage[utf8]{inputenc}
\usepackage{amsmath,amsthm,amsfonts}  %additional math symbols, extremely useful
\usepackage{sectsty}  %to change section title font size and type
\usepackage{graphicx} %to include images
\usepackage{hyperref} %to include url links
\usepackage{cleveref}

\author{Noura Al Helwani \thanks{Department of Physics, American University of Beirut, Lebanon (nmh82@mail.aub.edu)} \and Sophie Moufawad \thanks{ Department of Mathematics, American University of Beirut, Lebanon (sm101@aub.edu.lb)} \and Nabil Nassif \thanks{Department of Mathematics, American University of Beirut, Lebanon (nn12@aub.edu.lb)}  \thanks{This work was funded by the AUB AI Exploratory Research Grant (AI-ERG 2024-2025) and by the Mamdouha El-Sayed Bobst Deanship Fund, Faculty of Arts and Sciences, AUB (2026) }}
\title{A Two-Stage Learning PINN Approach for Solving the Inverse Problem of the 1D Porous Medium Equation}
\begin{document}
\maketitle

\begin{abstract}
\noindent The Porous Medium Equation (PME), given by 
\(
u_t = \Delta(u^m) \quad\text{for}\quad m > 1,
\)
is a degenerate nonlinear parabolic partial differential equation that arises in various physical applications such as fluid flow in porous media, heat transfer in plasmas, and population dynamics. It is known for its nonlinear diffusion and finite propagation speed. In this paper, we study numerical solutions of the one-dimensional direct and inverse PME using Physics-Informed Neural Networks (PINNs), and compare them with classical numerical methods and available analytical and manufactured solutions. While PINNs provide a flexible framework for solving both forward and inverse problems, we show that the standard inverse formulation suffers from a strong sensitivity to the initial guess, leading to only local convergence. To address this issue, we propose a novel two-stage PINN training framework for the inverse problem, which significantly improves convergence stability and allows reliable recovery of the unknown parameter even for poor initial guesses. Overall, the proposed approach demonstrates that PINNs are a flexible and accurate alternative to classical methods for the 1D PME, and the introduced two-stage training strategy substantially improves their robustness in inverse problems, providing a solid basis for extensions to more complex geometries and higher-dimensional cases.
\end{abstract}

\vspace{0.5em}
\newpage
\hrule

\tableofcontents

\vspace{0.75cm}

\hrule

\section{Introduction}

Inverse problems governed by nonlinear partial differential equations (PDEs) arise in many scientific and engineering applications ~\cite{Moufawad, gholami2016inverse, Caro}, where unknown parameters or coefficients must be inferred from observed data. In this paper, we solve the inverse Porous Medium Equation (PME) Eq.~\ref{eq: PME} problem using Physics-Informed Neural Networks (PINNs), which is a mesh-free framework for solving such problems through the incorporation of governing physical laws into the learning. 
PINNs have several important advantages compared to classical numerical methods. They are mesh-free, which makes them easier and more flexible to extend to higher-dimensional problems, since classical approaches usually require new and more complex discretization schemes and more advanced and expensive numerical methods when moving from 1D to higher dimensions. In addition, inverse problems in classical frameworks are generally more expensive and sometimes ill-conditioned than direct problems, while PINNs can handle inverse problems in a more efficient way within the same framework. Finally, PINNs can also work with sparse data in inverse problems, without the need for full-field or dense measurements. These advantages, along with the growing development of scientific machine learning methods, were the main motivation behind using PINNs as the main approach in this paper. 

An elementary attempt to solve the 1D PME using PINNs only for $m = 3$ was done in \cite{NSG}. In this paper, we aim to improve the results produced by \cite{NSG} by improving the PINN model and focusing on the inverse problem. We also generalize the solution for any $m$, and we manufactured three different solutions to further verify the work and test its robustness. For the inverse problem, we propose a novel two-stage learning approach that achieves convergence even when the initial guess is far from the true solution, which solves the problem of local convergence when only standard PINN techniques were involved. Throughout, we compare the performance of PINNs to traditional numerical solvers in terms of accuracy, stability, and convergence, with the broader aim of understanding the strengths and limitations of PINNs in solving nonlinear diffusion equations like the PME. These classical numerical methods are used in this work as a validation benchmark rather than a competing approach, since the goal is not to outperform them in all cases, but to develop a robust and reliable PINN framework that remains stable and accurate across different settings.

To elaborate, the motivation behind writing this paper and doing this work is to build strong foundations where the 1D PME PINN performs very well, in order to prepare for the main topic and objective: the 2D and higher-dimensional PME in more complex geometries, where classical methods normally struggle, and to explore the advantages of PINNs. The 1D conditions presented in this paper represent the setting where classical methods usually excel. Our aim is to develop a robust PINN that can be as accurate as feasible compared to the classical method in these basic geometries, in order to ensure robustness before moving to higher and more complex geometries, where we expect PINNs to become more advantageous. In particular, establishing a reliable inverse PINN framework in 1D is a necessary step before moving to higher-dimensional cases, because unlike classical numerical methods, which usually need new discretization techniques, new mesh generation, and sometimes even different solvers when going to higher dimensions, PINNs keep the same main structure. Only small changes are needed, mainly in the input dimension and the way collocation points are sampled. This is why the paper is divided into two sections: the first section, Section ~\ref{sec: sec1}, focuses on the direct PME problem, while the second, Section ~\ref{sec:inverse_pme1} which represents the main contribution of this work, focuses on the inverse problem. In brief, the main objectives of this work are:
\begin{enumerate}
    \item To develop a robust PINN framework for the inverse PME problem in 1D.
    
    \item To establish a strong and reliable foundation for extending the methodology to higher-dimensional problems.
\end{enumerate}

Summarizing the main results of this paper, for both the direct and inverse problems, the PINN framework was found to be accurate, with most errors below 5\%. It is also more flexible than the classical methods, although it comes at a higher computational cost.  Overall, PINNs prove to be a viable and flexible alternative to classical methods for the 1D PME, with clear room for further improvement, and give a solid foundation for extending the framework to higher-dimensional geometries.

The main contributions of this work are summarized as follows:
\begin{itemize}
    \item We develop an improved PINN framework for the 1D Porous Medium Equation covering both direct and inverse problems.

    \item We investigate the effect of initialization, training strategy, stability techniques, and loss scheduling on inverse PINN performance (Section ~\ref{sec: experiments}).
    
    \item We identify the local convergence limitation of standard inverse PINNs for the PME, where accurate parameter recovery strongly depends on the initial guess (Section ~\ref{sec: standard}).
    
    \item We propose a novel two-stage training PINN framework for the inverse PME problem that improves convergence robustness and stabilizes parameter estimation across different initial guesses (Section ~\ref{sec:proposed}).

    \item We compare the proposed framework against classical numerical inverse solvers across several manufactured solutions and values of \(m\).

    \item We establish a robust and well-behaved inverse PINN framework for the 1D PME, which provides a solid foundation for extending the methodology to higher-dimensional and more complex geometries.

\end{itemize}

The proposed two-stage learning strategy is not restricted to the PME and can be naturally extended to other PDE-based inverse problems, where the number of stages and their structure can be adapted depending on the specific problem, the data availability, and the convergence behavior, potentially leading to more general multi-stage learning frameworks tailored to each application.\vspace{2mm}

To better situate this work, we now briefly review relevant literature on machine learning and Physics-Informed Neural Networks for solving PDEs and inverse problems.
Machine learning (ML) has become increasingly prominent with the data revolution in the digital age, with huge success in applications such as image classification \cite{krizhevsky2012imagenet} and natural language processing \cite{lecun2015deep}. These successes are usually driven by large, high-quality data sets. However, a lot of applied science is in the small-data regime with large uncertainty. In such cases, ML-based methods, e.g., PINNs, are particularly valuable \cite{raissi2017physics1}.\vspace{2mm}

\noindent Neural networks in particular possess several properties that make them important for scientific computing. They are universal function approximators \cite{hornik1991universal} and benefit from automatic differentiation, which makes calculating the derivatives that occur in differential equation-solving more straightforward \cite{baydin2018automatic}. The use of neural networks to approximate ordinary and partial differential equations  started in the 1990s \cite{meade1994numerical, vanmilligen1995neural, lagaris1998artificial} and has continued as an active area of research with deep learning techniques. Early approaches blended neural networks with traditional methods, e.g., weighted residuals, constrained backpropagation, and Galerkin methods, for the solution of initial and boundary value problems \cite{lagaris2000irregular}.
These developments naturally led to physics-informed formulations. More recently, Raissi et al. \cite{raissi2017physics1, raissi2019physics} proposed Physics-Informed Neural Networks, which directly incorporate physical laws into the training procedure.  It is a mesh-free approach, avoids interpolation errors, and offers a unified treatment of both forward and inverse problems. PINNs have been applied to diverse areas such as heat transfer \cite{zobeiry2021heat, mishra2021radiative} and manufacturing \cite{zhu2021metal}. Another advantage is the availability of several Python libraries \cite{paszke2019pytorch, lu2021deepxde, chen2020neurodiffeq} that facilitate the implementation of PINNs. Most notably, PINNs provide a unified framework for both direct (forward) and inverse problems. These advantages explain the rapid adoption and exploration of PINNs in scientific applications.\vspace{2mm}

\noindent PINNs do come with some issues, though. Their loss landscapes may be very non-convex, which can induce optimization difficulties and decreased accuracy. More sophisticated physics typically require deeper networks, which increase the computation cost, memory requirements, and training time. Also, PINNs are confronted with PDEs with parameter discontinuities or domain heterogeneities, e.g., conjugate heat transfer in heterogeneous media.
To address these limitations, several enhanced PINN variants have been developed. Domain decomposition techniques have been proposed to overcome this difficulty \cite{shukla2021parallel, moseley2023finite}. These divide the computational domain into subdomains, for which local PINNs are separately defined and coupled together with interface conditions. Adaptive activation functions \cite{jagtap2020locally} and even variational principles have been examined to improve accuracy.
In order to further improve performance, Transfer Physics-Informed Neural Networks (TPINNs) were introduced. Unlike training separate neural networks for each of the subdomains, TPINNs transfer layers between subdomains but allow certain layers to learn to adjust to local physics. Besides increasing efficiency, it also shows some effectiveness in solving stiff differential equations \cite{seiler2025stifftransferlearningphysicsinformed, wang2025transferlearningphysicsinformedneural}. Beyond these improvements, several other neural network-based approaches have been proposed to handle complex scientific computing problems where conventional numerical methods struggle, especially in highly nonlinear or high-dimensional systems. Examples include geometry-aware PINNs for more complex domains \cite{sukumar2022exact}, specific activation functions for high-dimensional nonlinear wave equations \cite{zhou2023deep}, and convolutional PINNs for irregular geometries \cite{gao2021phygeonet}.The demand for inverse PINNs is especially urgent in porous media simulations, as highlighted in recent studies such as \cite{berardi2025inverse}, where Physics-Informed Neural Networks are used for transport models in porous materials. In this context, correct parameter estimation, such as permeability and dispersivity, is pivotal for accurately modeling large-scale dynamics. Direct measurement of these parameters is typically expensive, time-consuming, and spatially limited. PINNs address this issue by enabling parameter estimation from sparse and noisy measurements of state variables of PDE systems. This makes them particularly valuable not only in porous media applications, but also in agronomy, soil sciences, and hydrology \cite{berardi2023optimal, celia1990mass, marinoschi2006functional}, as well as in materials science and biomedical engineering \cite{wein2019topology}.\vspace{2mm}

\noindent In summary, PINNs were effective forward and inverse PDE problem solvers, overcoming many of the limitations of conventional numerical methods. In particular, inverse problems are often ill-posed, meaning they can be non-unique or unstable, and are typically more computationally complex than forward problems. While recent studies on PINNs provided improved accuracy and efficiency, issues are still open. This requires further investigation into the application and extension of PINN approaches to problems such as the Porous Medium Equation (PME), where both stable forward solutions and accurate parameter estimation are paramount. \vspace{2mm}

\noindent The one-dimensional PME has been studied using classical numerical approaches, including finite difference schemes with various sweep and discretization strategies, as well as iterative solvers such as SOR and Newton-based methods \cite{chew2021half, chew2020quartersweep, chew2016newton}, and variational approaches \cite{duan2018pme}. As for the inverse problem, see \cite{blessing2025pme, carstea2021inversepme}. 
However, despite this rich body of work, these approaches do not incorporate machine learning data-driven learning frameworks. This highlights the potential interest in exploring the effectiveness of PINNs for both direct and inverse formulations of the PME. 

\section{The Direct Problem of the 1D PME}
\label{sec: sec1}
The Porous Medium Equation (PME) is a nonlinear degenerate parabolic partial differential equation of the form
\begin{equation}
\begin{cases}
\partial_t u = \Delta(u^m) \quad & m > 1 \\
u(t,a) = 0,\quad u(t,b) = 0 & t > t_0, \\
u(t_0,x) = u_0(x) & x \in (a,b),
\label{eq: PME}
\end{cases}
\end{equation}
where \( u = u(x,t) \geq 0 \) denotes the density of a diffusing quantity, and \( m \) is the nonlinearity (or polytropic) exponent.

\noindent The equation is degenerate due to the fact that the diffusion coefficient \( D(u) = mu^{m-1} \) becomes zero at \( u = 0 \) if the equation is equivalently written in the form \( \partial_t u = \partial_x (m u^{m-1} \partial_xu) \).
The PME unifies a wide range of nonlinear diffusion phenomena occurring in physics, biology, and engineering, including gas transmission in porous media, population growth, and heat transfer in nonlinear materials. The PME can be considered as a nonlinear version of the classical heat equation, which is regained formally in the case \( m = 1 \).

\noindent The mathematical theory of the PME developed over the years, beginning with foundational work of Oleinik et al. in the late 1950s, and finally reaching maturity in the 1980s, partly due to developments in functional analysis and nonlinear PDE theory~\cite{vazquez2007porous}.

\noindent Finite speed of propagation is also one of the typical features of the PME, meaning that disturbances propagate at finite speed, unlike in the classical heat equation, where infinite speed of propagation occurs. Indeed, if the initial data \( u_0 \) is compactly supported, then so is the solution \( u(t) \) for all \( t > 0 \), with a sharp free boundary separating the region where \( u > 0 \) from the vacuum. That results in the formation of moving interfaces and makes the maximum principle inapplicable. A fundamental exact solution of the PME is the Barenblatt solution. It describes the evolution of a compactly supported density from a Dirac-type initial condition and serves as a paradigm to study spreading rates, self-similarity, and asymptotic behavior. We will utilize such a solution to verify our numerical and PINNs methods. For a complete treatment of the PME and other related nonlinear diffusion equations, we would like to refer the reader to Vázquez's monograph~\cite{vazquez2007porous}, which will be the principal theoretical foundation of this work.

In this section, we solve the one-dimensional direct problem of the PME given in Eq.~\ref{eq: PME} written in the following form: 
\begin{equation}
    \begin{cases}
\partial_t u(t,x) = \partial_x\left( \beta u^{\beta - 1} \partial_x u(t,x) \right), & x \in [-1,1],\; t \in [0,1], \\
u(t, -1) = u(t, 1) = 0, & t \in [0,1], \\
u(0, x) = u_0(x), & x \in [-1,1],
\end{cases}
\label{eq: PME2}
\end{equation}
We solve this problem using a classical numerical method (Section~\ref{sec: classical}) and a PINN approach (Section~\ref{sec: direct pinn}). We then introduce four different reference solutions of the PME (Section~\ref{sec: ref sol}), and finally, we present the results for the 1D direct problem in Section~\ref{sec: test direct}.

\subsection{Methodologies}
The inverse framework developed in this paper builds upon the direct PME solvers introduced in our earlier work \cite{NSG}. For completeness, we briefly summarize the classical and PINN formulations used throughout this work. Additional implementation details and algorithmic descriptions can be found in \cite{NSG}.

\subsubsection{Classical Numerical Scheme}
\label{sec: classical}

As a reference solution, we employ the finite difference formulation introduced in \cite{NSG}, based on a backward Euler discretization in time and central differences in space. For each interior grid point, the resulting nonlinear residual is

\[
F_i(u^{n+1}) =
u_i^{n+1} - u_i^n
-
\frac{m \Delta t}{\Delta x^2}
\left[
\left(\frac{u_{i+1}^{n+1}+u_i^{n+1}}{2}\right)^{m-1}
( u_{i+1}^{n+1}-u_i^{n+1} )
-
\left(\frac{u_i^{n+1}+u_{i-1}^{n+1}}{2}\right)^{m-1}
( u_i^{n+1}-u_{i-1}^{n+1} )
\right].
\]
with the vector of unknowns at time step \( n+1 \) as
\[
u^{n+1} = (u_1^{n+1},\dots,u_{N-1}^{n+1})^\top \in \mathbb{R}^{m},
\quad m = N-1,
\]
and the residual vector
\[
F(u^{n+1}) = (F_1,\dots,F_{N-1})^\top.
\]
The nonlinear system is solved at each time step using Newton's method. The update rule is

\[
J (u^{(k)}) \, \delta u^{(k)} = -F(u^{(k)}), \quad 
u^{(k+1)} = u^{(k)} + \delta u^{(k)},
\]
where the Jacobian matrix is given by

\[
\bigl[J\bigl(u^{(k)}\bigr)\bigr]_{ij}
= \frac{\partial F_i}{\partial u_j}\Bigg|_{u=u^{(k)}},
\qquad i,j=1,\dots,s.
\] 
and is approximated using finite differences. 

\noindent Table ~\ref{tab:parameters} summarizes the parameters used. The pseudocode of the classical PME direct problem, using the Barenblatt solution~\ref{sec: barenblatt} for the initial and boundary conditions, is summarized in Algorithm 11 of \cite{NSG}. 

\begin{table}[H]
\centering
\begin{tabular}{|l|l|}
\hline
\textbf{Parameter} & \textbf{Value} \\ \hline
Spatial domain & \(x\in[-1,1]\) with \(N=100\) intervals (\(N+1\) grid points) \\ \hline
Time domain & \(t\in[0,1]\) with \(\Delta t=0.01\) (\(100\) time steps) \\ \hline
Tolerance (Newton) & \(10^{-6}\) \\ \hline
Maximum Newton iterations & \(20\) per time step \\ \hline
Jacobian approximation & Finite differences with perturbation \(h=10^{-6}\) \\ \hline
\end{tabular}
\caption{Parameters for Direct Classical Method}
\label{tab: 8}
\label{tab:parameters}
\end{table}

\subsubsection{PINN Approach} \label{sec: direct pinn}

Physics-Informed Neural Networks (PINNs) are built upon standard feedforward neural networks, which act as universal function approximators \cite{hornik1991universal} by learning mappings between inputs and outputs through trainable weights and biases. Their training relies on automatic differentiation together with optimization algorithms such as ADAM and L-BFGS \cite{liu1989limited, kingma2014adam}, to minimize the loss function and accurately approximate PDE solutions. What distinguishes PINNs from traditional neural networks is the fact that they include physical laws, expressed as differential equations, directly in the loss function. To explain, PINNs use the residuals of the governing PDE as part of the loss. This means the network is penalized when it mispredicts known data and also when its predictions violate the underlying physics.

\noindent For instance, suppose we have a PDE of the form:
\[
\mathcal{N}[u](\mathbf{x}, t;c) = 0,
\]
where \( \mathcal{N} \) is a differential operator acting on the unknown function \( u(\mathbf{x}, t) \) and c is an underlying parameter of the PDE. A PINN seeks to approximate \( u \) with a neural network \( u_\theta(\mathbf{x}, t) \), parameterized by weights \( \theta \). During training,

\begin{enumerate}
    \item In solving the direct problem, the total loss is composed of the PDE, initial value, and boundary conditions residuals. Also, it may include any additional physical constraints. 
    \[
\mathcal{L}_{direct}(\theta) = \mathcal{R}_{\text{PDE}}(\theta) + \mathcal{R}_{\text{initial value}}(\theta) + \mathcal{R}_{\text{boundaries}}(\theta) + \lambda \cdot \mathcal{R}_{\text{constraint}}(\theta) .
\]
where $\lambda$ is mathematically speaking the Lagrange multiplier, but in our computational context, we refer to it as a hyperparameter. It is worth noting that any other residual in the total loss can also be multiplied by another hyperparameter. 

\item In solving the inverse problem, i.e. recovering c, the total loss is composed of the direct problem loss, along with a data loss. In this way, the network will be penalized for mispredicting both the physics and the data.
\[
\mathcal{L}_{inverse}(\theta) = \mathcal{L}_{\text{direct}}(\theta) + \mathcal{R}_{\text{data}}(\theta) 
\]
\end{enumerate}
In this section, we give the details of the PINN constructed to solve the direct PME problem that gave the results in Section ~\ref{sec: test direct}.

\begin{figure}[H]
    \centering
    \includegraphics[width=1\textwidth]{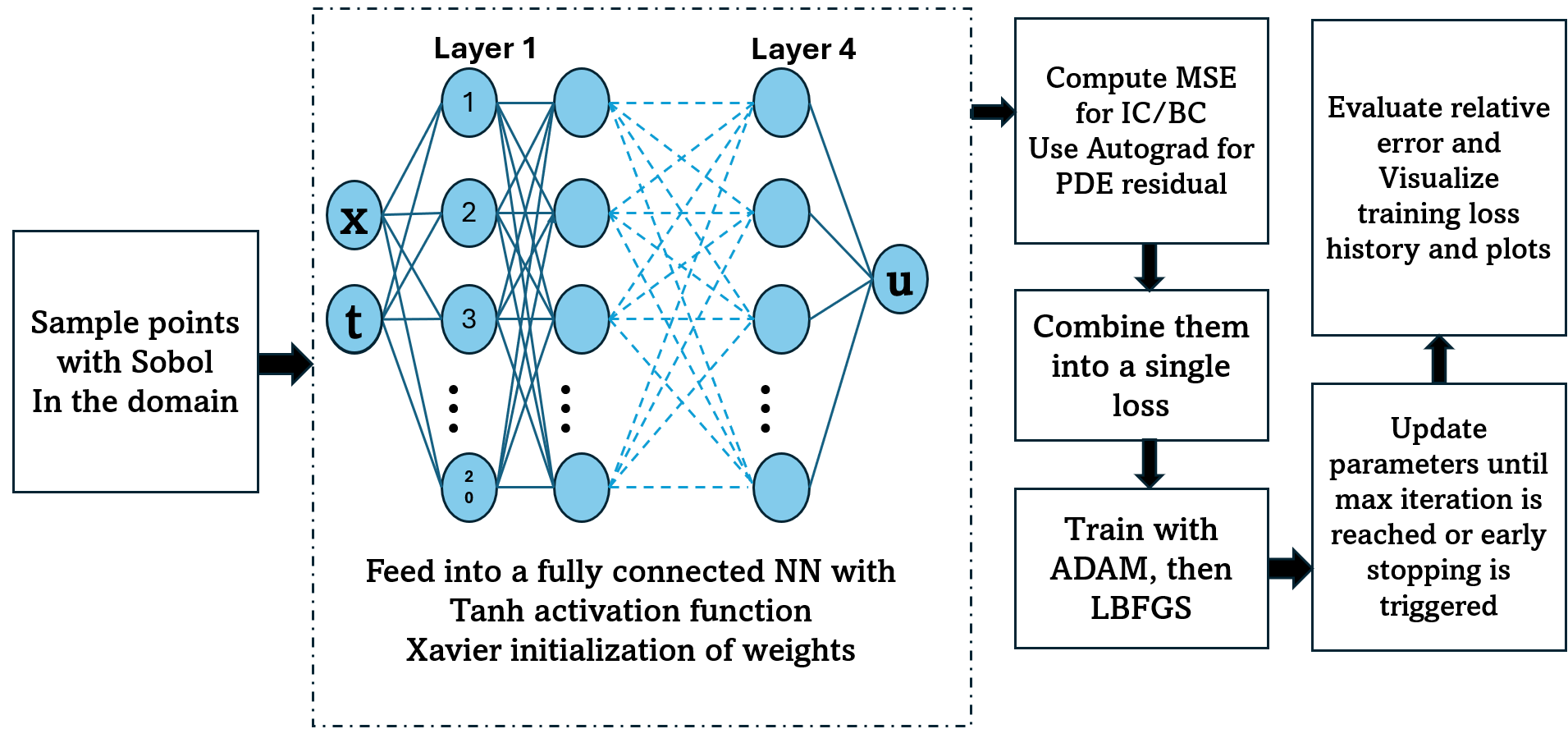}\vspace{-2mm}
\caption{Schematic of the direct PINN}
    \label{direct schematic}
\end{figure}

The direct PINN framework follows the formulation presented in \cite{NSG} and is summarized in the following section. The direct problem serves as the foundation for the inverse methodology developed later in this paper. A schematic representation of the direct PINN is shown in Figure~\ref{direct schematic}.\vspace{-1mm}

\paragraph*{Neural Network Architecture} 

The solution \(u(t,x)\) is approximated using a fully connected neural network \(u_\theta(t,x)\). The network architecture used throughout this work is summarized in Table~\ref{tab: arch}.

\begin{table}[H]
\centering
\begin{tabular}{|l|l|}
\hline
\textbf{Component} & \textbf{Specification} \\ \hline
Input layer & Size 2 (for \(t\) and \(x\)) \\ \hline
Hidden layers & 4 layers, each with 20 neurons \\ \hline
Activation function & Hyperbolic tangent \( \tanh \) \\ \hline
Output layer & Single neuron returning \(u_\theta(t,x)\) \\ \hline
\end{tabular}
\caption{Architecture of the neural network.}
\label{tab: arch}
\end{table}

\paragraph*{weights} All weights of the neural network are initialized with the \textit{Xavier initialization} procedure that works best with \( \tanh \) activation function. This process initializes the weights by sampling random values from a uniform distribution over the range:
\[
W \sim \mathcal{U}\left( -\sqrt{\frac{6}{n_{\text{in}} + n_{\text{out}}}},\ \sqrt{\frac{6}{n_{\text{in}} + n_{\text{out}}}} \right),
\]
\noindent where \( n_{\text{in}} \) and \( n_{\text{out}} \) are the number of input and output units of a layer, respectively \cite{glorot2010understanding}. Such initialization helps maintain stable gradients during training, especially when activation functions like \(\tanh \) are employed since it maintains the signal variance across layers constant.

\paragraph*{Loss Function and Optimization}

The training objective is to minimize a composite loss \( \mathcal{L} \) that enforces agreement with the initial and boundary conditions as well as the underlying PDE. Let \( \mathcal{D}_t \), \( \mathcal{D}_b \), and \( \mathcal{D}_\Omega \) denote sampled points on the initial boundary \( t=0 \), spatial boundaries \( x = \pm 1 \), and interior domain \( \Omega \), respectively.\vspace{2mm}

\noindent The total loss is:\vspace{-4mm}
\begin{equation}
\label{eq: loss direct}
\mathcal{L}_\theta = \log_{10} \left( \lambda_u \cdot (\mathcal{L}_b +\mathcal{L}_t )+ \mathcal{L}_{\text{PDE}}\right)
\end{equation}
with \vspace{-5mm}
\begin{align}
\mathcal{L}_b &= \mathbb{E}_{(t,x)\in\mathcal{D}_b} \left[ \left( u_\theta(t,x) - u_{\text{BC}}(t,x) \right)^2 \right], \\
\mathcal{L}_t &= \mathbb{E}_{(t,x)\in\mathcal{D}_t} \left[ \left( u_\theta(t,x) - u_0(x) \right)^2 \right], \\
\mathcal{L}_{\text{PDE}} &= \mathbb{E}_{(t,x)\in\mathcal{D}_\Omega} \left[ \left( \partial_t u_\theta - \partial_x \left( mu_\theta^{m-1} \partial_x u_\theta \right) \right)^2 \right].\vspace{-9mm}
\end{align}
\noindent where \(\mathbb{E}[\cdot]\) denotes the empirical mean over the sampled training points, i.e.,
$\displaystyle
\mathbb{E}[f(x)] = \frac{1}{N}\sum_{i=1}^{N} f(x_i),
$
which corresponds to the Mean Squared Error (MSE) used in the loss computation.\\ Gradients are computed using automatic differentiation.
\paragraph*{Training Procedure} \label{sec: training}

The neural network is trained using the following optimization approach:

\begin{enumerate}
    \item ADAM Optimizer for up to 2000 epochs with a learning rate $ = 0.001$.
    \item followed by L-BFGS for up to 500 epochs with a learning rate $ = 0.05$.
    \item Early Stopping: The training halts if the validation loss does not improve by $10^{-6}$ over a predefined patience window equals to 200 epochs in our case.
\end{enumerate}

\noindent The optimizer minimizes the logarithmic loss for stability and to better visualize loss curves, and training points are generated using a Sobol sequence, which ensures a quasi-random coverage of the space-time domain.The network was trained with $n_{\text{int}} = 256$ interior points, $n_{\text{sb}} = 64$ points on each spatial boundary, and $n_{\text{tb}} = 64$ points for the initial condition with a total of 448 points. 
\paragraph*{Evaluation and Accuracy}
\noindent The solution \( u_\theta(t,x) \) is evaluated against the analytical form. We compute the relative \( L^2 \)-error over the same grid points used by the classical Newton method which will be put in comparison with our PINN.

\begin{equation*}
\epsilon_{rel} = \frac{ \left\| u_\theta - u \right\|_{L^2} }{ \left\| u \right\|_{L^2} } = \sqrt{ \frac{ \sum_{i=1}^N \left( u_\theta(t_i, x_i) - u(t_i, x_i) \right)^2 }{ \sum_{i=1}^N u(t_i, x_i)^2 } }.\vspace{-5mm}
\end{equation*}

\subsection{Reference Solutions}
\label{sec: ref sol}
In this section, we present the four different PME solutions on which we will test our methods.
\subsubsection{Barenblatt Solution}
\label{sec: barenblatt}

The Barenblatt solution, thoroughly explained in \cite{vazquez2007porous}, is a well-known self-similar solution to the PME and serves as a fundamental benchmark for validating our numerical methods. It represents the evolution of an initial mass concentrated at a point and is particularly useful due to its explicit analytical form.

\noindent The Barenblatt solution in one spatial dimension is given by:
\[
u(t,x) = t^{-\alpha} F\left(x t^{- \alpha} \right),
\]
where the similarity profile \( F \) is defined as:
\[
F(\eta) = \left(C - \kappa \eta^2\right)_{+}^{\frac{1}{m-1}},
\quad \text{with} \quad \eta = x t^{-\alpha}.
\]

\noindent Here, \( m > 1 \) is the nonlinearity parameter of the PME, and the constants are:
\[
\alpha = \frac{1}{m+1}, \quad \kappa = \frac{(m - 1)\alpha}{2m} = \frac{m - 1}{2m(m+1)}.
\]

\noindent The parameter \( C \) is a constant that depends on the total mass of the solution and determines the support of \( u(t,x) \). The notation \( (\cdot)_+ \) indicates that the expression is set to zero whenever the argument is negative to ensure compact support. This solution will not involve the homogeneous Dirichlet boundary conditions.

\subsubsection{Manufactured Solution}

To further test the accuracy of our solver, we use two manufactured solutions. These are artificial solutions that we plug into the equation to compute the source term \( f(x,t) \), so that the equation is exactly satisfied. This helps us evaluate how well the solver can recover known solutions, even when no closed-form analytical solution exists for the original PME.

\subsubsubsection*{Polynomial Solution}
\label{sec: poly sol}

We first consider a smooth polynomial function in both time and space:
\[
u(x,t) = (1 - t)^2 (1 - x^2)^2.
\]
By plugging this into the equation \( \partial_t u = \partial_x \left({m} u^{m-1} \partial_x u \right) + f(x,t) \), we obtain the forcing term:
\[
f(x,t) = -2(1 - t)(1 - x^2)^2 
+ 4m (1 - t)^{2m} \left[
(1 - x^2)^{2m - 1} - 2(2m - 1)x^2 (1 - x^2)^{2m - 2}
\right].
\]
This solution is zero at the boundaries, so it satisfies the same boundary conditions as our test setup.

\subsubsubsection*{Damped Harmonic Solution}
\label{sec: harmonic}
We also test a harmonic function that decays over time:
\[
u(x,t) = e^{-t} \cos\left(\frac{\pi}{3}x\right),
\]
with the corresponding source term:
\[
f(x,t) = -e^{-t} \cos\left(\frac{\pi}{3}x\right) 
+ m \left(\frac{\pi}{3}\right)^2 e^{-tm} \left[
-(m - 1) \cos^{m - 2}\left(\frac{\pi}{3}x\right) \sin^2\left(\frac{\pi}{3}x\right)
+ \cos^m\left(\frac{\pi}{3}x\right)
\right].
\]
 \noindent This solution will not involve the homogeneous Dirichlet boundary conditions.

\subsubsubsection*{m-Dependent manufactured solution} \label{sec: m-dep}
As one can notice, all previously manufactured solutions did not depend on 
m. Here, we suggest an additional 
m-dependent solution: 
\[
u(x,t) = (\cos{\frac{\pi x}{2}})^{2m} (1+t)^m
\]
with its corresponding source term: 
\[
f(x, t) = m (1 + t)^{m - 1} \cos^{2m}\left( \frac{\pi x}{2} \right)
- \frac{\pi^2}{2} m^2 (1 + t)^{m^2} \left[
(2m^2 - 1) \cos^{2m^2 - 2}\left( \frac{\pi x}{2} \right) \sin^2\left( \frac{\pi x}{2} \right)
- \cos^{2m^2}\left( \frac{\pi x}{2} \right)
\right]
\]

\noindent All these manufactured solutions are used to validate our classical and PINN-based solvers by directly comparing the computed results with the exact ones, as well as to compare the results with those of the classical solver.

\subsection{Direct PINN Testing without Data} \label{sec: test direct}

\begin{figure}[H]
    \centering
    
    \begin{subfigure}[b]{0.32\textwidth}
        \centering
        \includegraphics[width=\textwidth]{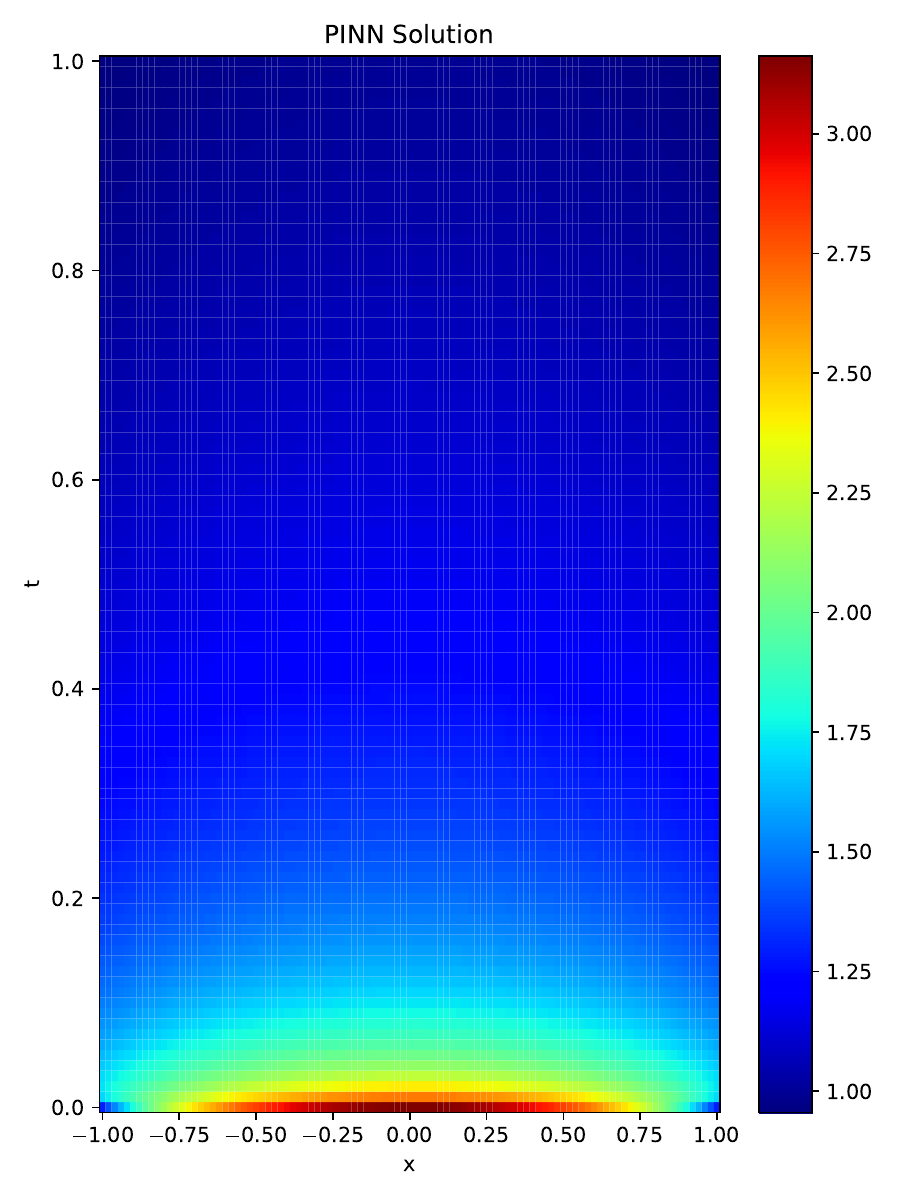}
        \caption{PINN solution}
    \end{subfigure}
    \hfill
    \begin{subfigure}[b]{0.32\textwidth}
        \centering
        \includegraphics[width=\textwidth]{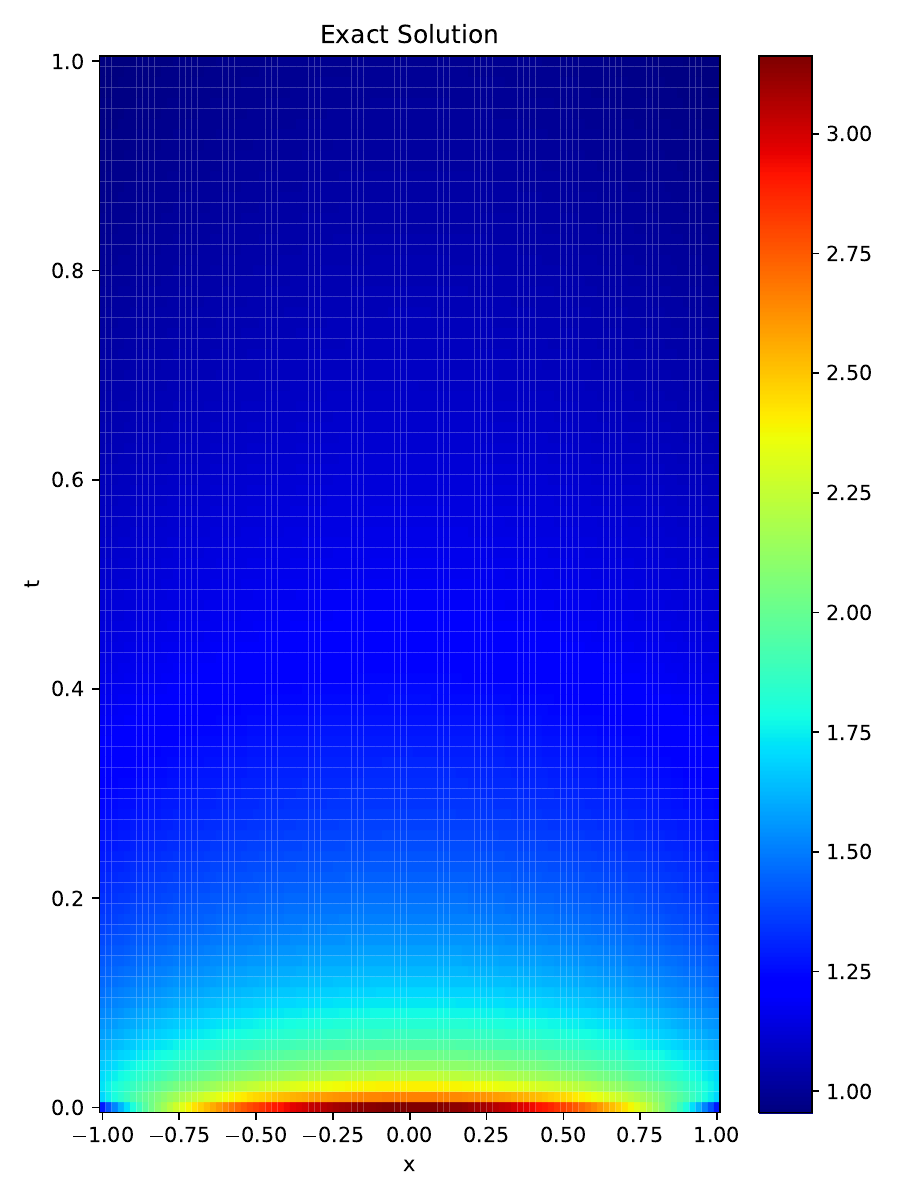}
        \caption{Exact solution}
    \end{subfigure}
    \hfill
    \begin{subfigure}[b]{0.32\textwidth}
        \centering
        \includegraphics[width=\textwidth]{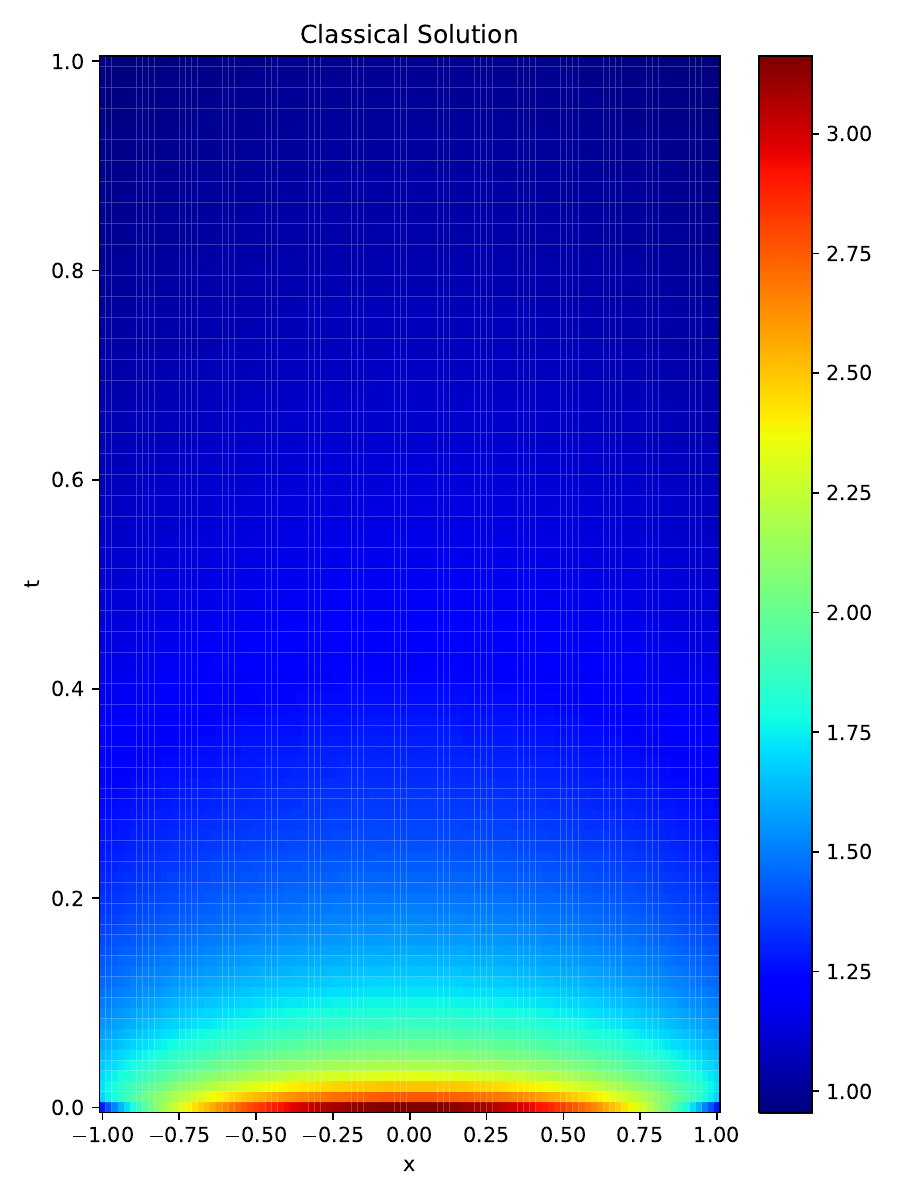}
        \caption{Classical solution}
    \end{subfigure}
    
    \caption{visual comparison of Barenblatt numerical solutions for $m = 3 $}
    \label{fig:solutions_baren}\vspace{-2mm}
\end{figure}

\begin{figure}[H]
    \centering
    
    \begin{subfigure}[b]{0.32\textwidth}
        \centering
        \includegraphics[width=\textwidth]{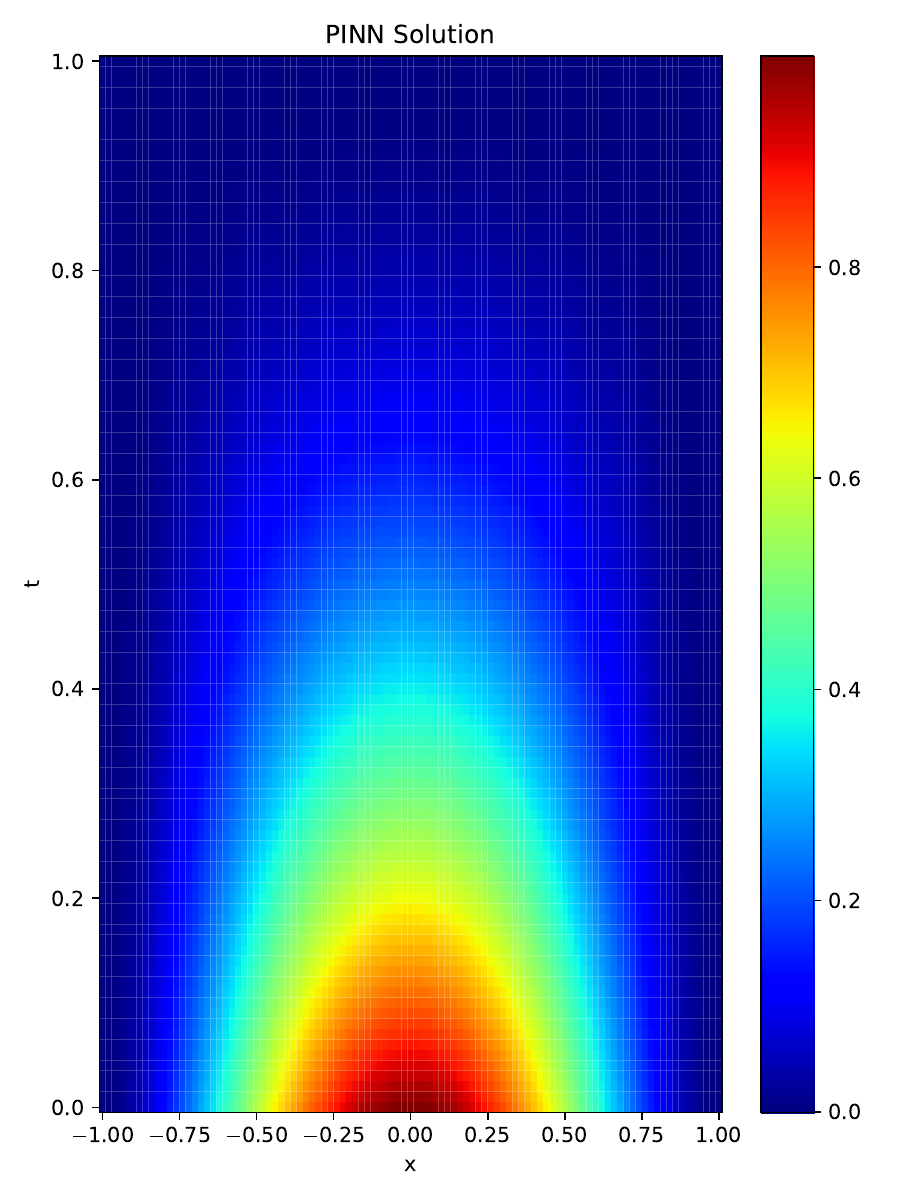}
        \caption{PINN solution}
    \end{subfigure}
    \hfill
    \begin{subfigure}[b]{0.32\textwidth}
        \centering
        \includegraphics[width=\textwidth]{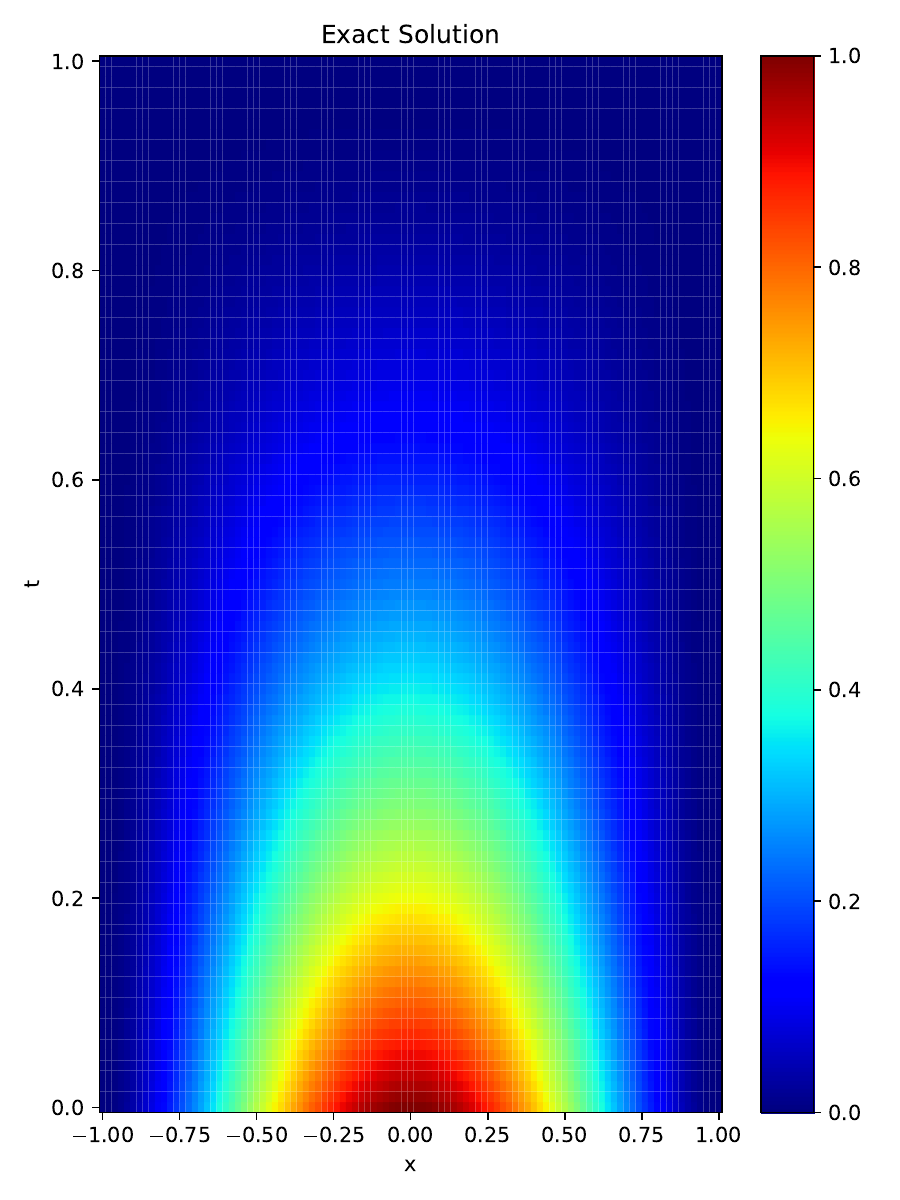}
        \caption{Exact solution}
    \end{subfigure}
    \hfill
    \begin{subfigure}[b]{0.32\textwidth}
        \centering
        \includegraphics[width=\textwidth]{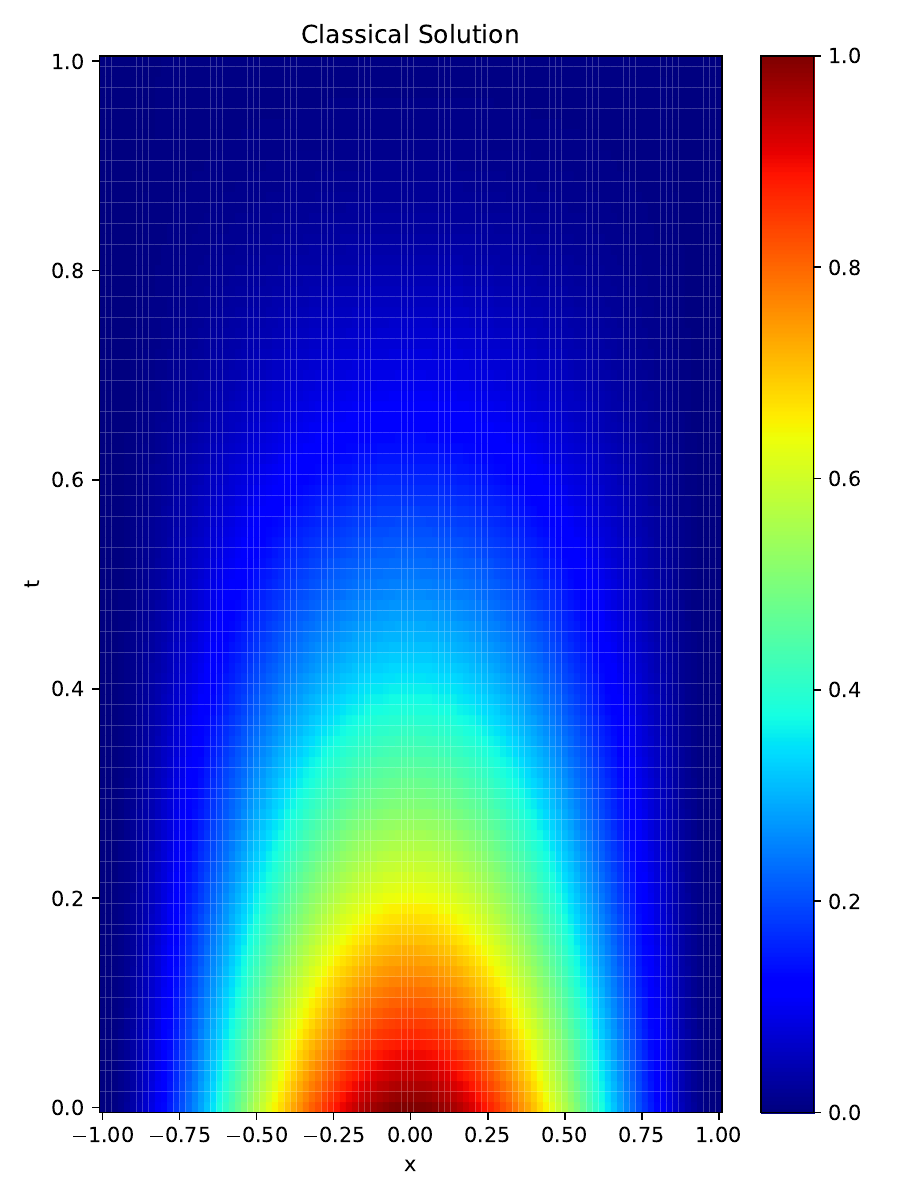}
        \caption{Classical solution}
    \end{subfigure}
    
    \caption{visual comparison of the Polynomial numerical solutions }
    \label{fig:solutions_poly}
\end{figure}

In this section, we validate our classical numerical solver (~\ref{sec: classical}) and the PINN approach (~\ref{sec: direct pinn}) described against our four reference solutions provided in ~\ref{sec: ref sol}. The PINN testing in this section follows the standard approach presented in Section~\ref{sec: direct pinn}, where no data are incorporated into the loss function. A variation of this standard testing procedure, in which measurement data are added as an additional term in the loss function, is considered later in Section~\ref{direct noise}. First, we test the performance for different values of the polytropic exponent \( m \in \{2, 3, 4, 5\} \). Then, we compute the relative \( L^2 \)-error between the approximated solutions and the exact solution over the 10201 points of the classical mesh in Table ~\ref{tab:parameters}:

\begin{equation*}
\epsilon_{\text{rel}} = \frac{ \| u_{\text{approx}} - u_{\text{exact}} \|_{L^2} }{ \| u_{\text{exact}} \|_{L^2} }.
\end{equation*} and provide a table of the errors for each solution (Tables \ref{tab:barenblatt_errors}, \ref{tab:polynomial_error}, \ref{tab:harmonic_error}, \ref{tab: m-dep error}). Other tables are also provided for each solution to compare the time (Tables \ref{tab:barenblatt_time}, \ref{tab:poly_time} , \ref{tab:har_time}, \ref{tab:mdep_time}) and the epochs (Tables \ref{tab:barenblatt_epochs}, \ref{tab:poly_epochs} , \ref{tab:har_epochs}, \ref{tab:mdep_epochs}) used by each method in each case. The reported epochs are approximate values estimated visually from the loss curves (Appendix ~\ref{app:A}), reflecting the general convergence behavior of each method. Rather than reporting the final training epoch, we report the approximate epoch at which the loss first reached its plateau (i.e., the last meaningful improvement). Consequently, epochs corresponding to a plateau with no significant reduction in the loss, including any patience period preceding early stopping, were not counted. 
Also, only a visual comparison for the case of $m=3$ is presented in this section (Figures \ref{fig:solutions_baren}, \ref{fig:solutions_poly}, \ref{fig:solutions_harm}, \ref{fig:solutions_mdep}); for the visuals of the other cases, see Appendix ~\ref{app:A}. For this section, all codes were run in Python using Google Colab with an NVIDIA Tesla T4 GPU \cite{googlecolab}. Lastly, for each type of solution, we provide some observations and comments when needed.

\subsubsection*{Barenblatt Solution Results ~\ref{sec: barenblatt}}
While both methods are scientifically accurate (relative error below 5\% in all cases, as shown in Table \ref{tab:barenblatt_errors}), the PINN consistently outperformed the classical numerical method in terms of accuracy. On the other hand, the accuracy improved as the nonlinearity parameter \( m \) increased for both methods. This is expected, as larger \( m \) values lead to more localized Barenblatt profiles, which are easier to approximate numerically and by neural networks.
Note that for \( m = 2 \), both methods had their largest errors. 

\begin{table}[H]
\centering
\caption{Relative \( L^2 \)-errors for the Barenblatt Solution}
\begin{tabular}{|c|c|c|c|}
\hline
\( m \) & PINN vs Exact & Classical vs Exact & PINN vs Classial \\
\hline
2 & 1.774 \( \times 10^{-2}\) & 2.475 \( \times 10^{-2}\) &  2.591 \( \times 10^{-2}\) \\
3 & 5.378 \( \times 10^{-4}\) & 1.557  \( \times 10^{-2}\)   & 1.315 \( \times 10^{-2}\) \\
4 & 3.394 \( \times 10^{-4}\) & 7.232 \( \times 10^{-3}\)  & 6.404 \( \times 10^{-3}\) \\
5 & 1.575 \( \times 10^{-4}\)  & 3.652 \( \times 10^{-3}\) & 3.322 \( \times 10^{-3}\) \\
\hline
\end{tabular}
\label{tab:barenblatt_errors}\vspace{-3mm}
\end{table}
\begin{table}[H]
\centering

\begin{minipage}{0.48\textwidth}
\centering
\caption{Training epochs for the Barenblatt solution}
\begin{tabular}{|c|c|c|}
\hline
$m$ & Adam epochs & LBFGS epochs \\ \hline
2 & 2000 & 500 \\ \hline
3 & 2000 & 500 \\ \hline
4 & 2000 & 500 \\ \hline
5 & 2000 & 500 \\ \hline
\end{tabular}
\label{tab:barenblatt_epochs}
\end{minipage}
\hfill
\begin{minipage}{0.48\textwidth}
\centering
\caption{time (sec) for the Barenblatt solution}
\begin{tabular}{|c|c|c|}
\hline
$m$ & PINN training time & Classical time \\ \hline
2 & 207.37 & 4.62 \\ \hline
3 & 194.96 & 5.31 \\ \hline
4 & 154.08 & 4.95 \\ \hline
5 & 159.39 & 5.28 \\ \hline
\end{tabular}
\label{tab:barenblatt_time}
\end{minipage}\vspace{-6mm}
\end{table}
To better illustrate the evolution of the diffusion process, a zoomed-in view of the solution is presented over a reduced time interval.\vspace{-3mm}
\begin{figure}[H]
    \centering
    \includegraphics[width=0.78\linewidth]{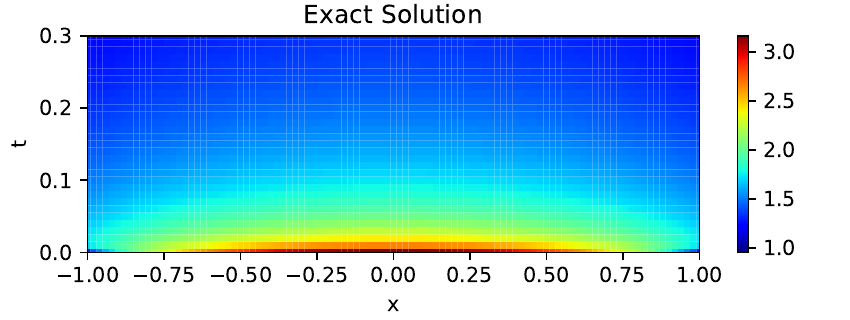}\vspace{-5mm}
    \caption{Zoomed-in view of the exact Barenblatt solution for $m = 3 $ for t $\in [0,0.3]$}\label{fig:bar_zoom}
\end{figure}

\subsubsection*{Polynomial Solution Results ~\ref{sec: poly sol}}
\label{sec: poly res}
 Although the exact solution of the manufactured polynomial case does not depend on \( m \), we still test the PINN for different values of \( m \in \{2, 3, 4, 5\} \). The goal is not to compare the solution shapes, since they are visually the same, but to observe how the PINN handles increasing nonlinearity as \( m \) increases.
 
The visuals are all the same over different values of m, as expected, since the solution is independent of m. However, from the training loss graphs shown in Appendix ~\ref{app:A}, we can observe that as \( m \) increases, the number of epochs needed to reach the minimum training loss also increases. This suggests that for higher values of \( m \), the solution becomes more complex and requires more epochs for the model to converge. Despite this, the convergence behavior remains consistent, with the loss gradually decreasing and eventually stabilizing. This indicates that the PINN is able to learn the solution for all values of \( m \), but the learning process takes longer as the nonlinearity strengthens as verified from the Table ~\ref{tab:poly_epochs} and ~\ref{tab:poly_time}.

 Table ~\ref{tab:polynomial_error} shows that the relative \( L^2 \)-errors of the PINN and exact solutions are low and accurate for all values of \( m \). While the error for \( m = 2 \) is slightly higher, all cases are considered accurate, with the PINN performing well and better than the classical method across all other values of \( m \), except for a negligibly higher error in the case of \( m = 2 \).

\begin{table}[H]
\centering
\caption{Relative \( L^2 \)-errors for the Polynomial solution}
\label{tab:polynomial_error}
\begin{tabular}{|c|c|c|c|}
\hline
\( m \) & PINN vs Exact & Classical vs Exact & PINN vs Classical \\
\hline
2 & 1.548 \( \times 10^{-2} \) & 1.012 \( \times 10^{-2}\) &  2.255 \( \times 10^{-2}\) \\
3 & 8.050 \( \times 10^{-4} \) & 1.134 \( \times 10^{-2}\)   & 1.127 \( \times 10^{-2}\) \\
4 & 1.255 \( \times 10^{-3} \) & 1.186 \( \times 10^{-2}\)  & 1.157 \( \times 10^{-2}\) \\
5 & 1.740 \( \times 10^{-3} \)  & 1.215\( \times 10^{-2}\) & 1.180 \( \times 10^{-2}\) \\
\hline
\end{tabular}

\end{table}

\begin{table}[H]
\centering

\begin{minipage}{0.48\textwidth}
\centering
\caption{Training epochs for the Polynomial solution}
\begin{tabular}{|c|c|c|}
\hline
$m$ & Adam epochs & LBFGS epochs \\ \hline
2 & 2000 & 1 \\ \hline
3 & 2000 & 310 \\ \hline
4 & 2000 & 410 \\ \hline
5 & 2000 & 460 \\ \hline
\end{tabular}
\label{tab:poly_epochs}
\end{minipage}
\hfill
\begin{minipage}{0.48\textwidth}
\centering
\caption{Time (sec) for the Polynomial solution}
\begin{tabular}{|c|c|c|}
\hline
$m$ & PINN training time & Classical time \\ \hline
        2 & 38.65 & 5.32 \\ \hline
        3 & 157.23 & 5.82 \\ \hline
        4 & 185.16 & 5.59 \\  \hline
        5 & 194.58 & 4.33 \\   \hline

\end{tabular}
\label{tab:poly_time}
\end{minipage}

\end{table}

\subsubsection*{Damped Harmonic solution Results~\ref{sec: harmonic}}
The damped harmonic solution does not depend on \( m \) either. Therefore, we perform the same testing as in the polynomial solution case (Section~\ref{sec: poly res}) for the same reason. \vspace{-2mm}
\begin{figure}[H]
    \centering
    
    \begin{subfigure}[b]{0.32\textwidth}
        \centering
        \includegraphics[width=\textwidth]{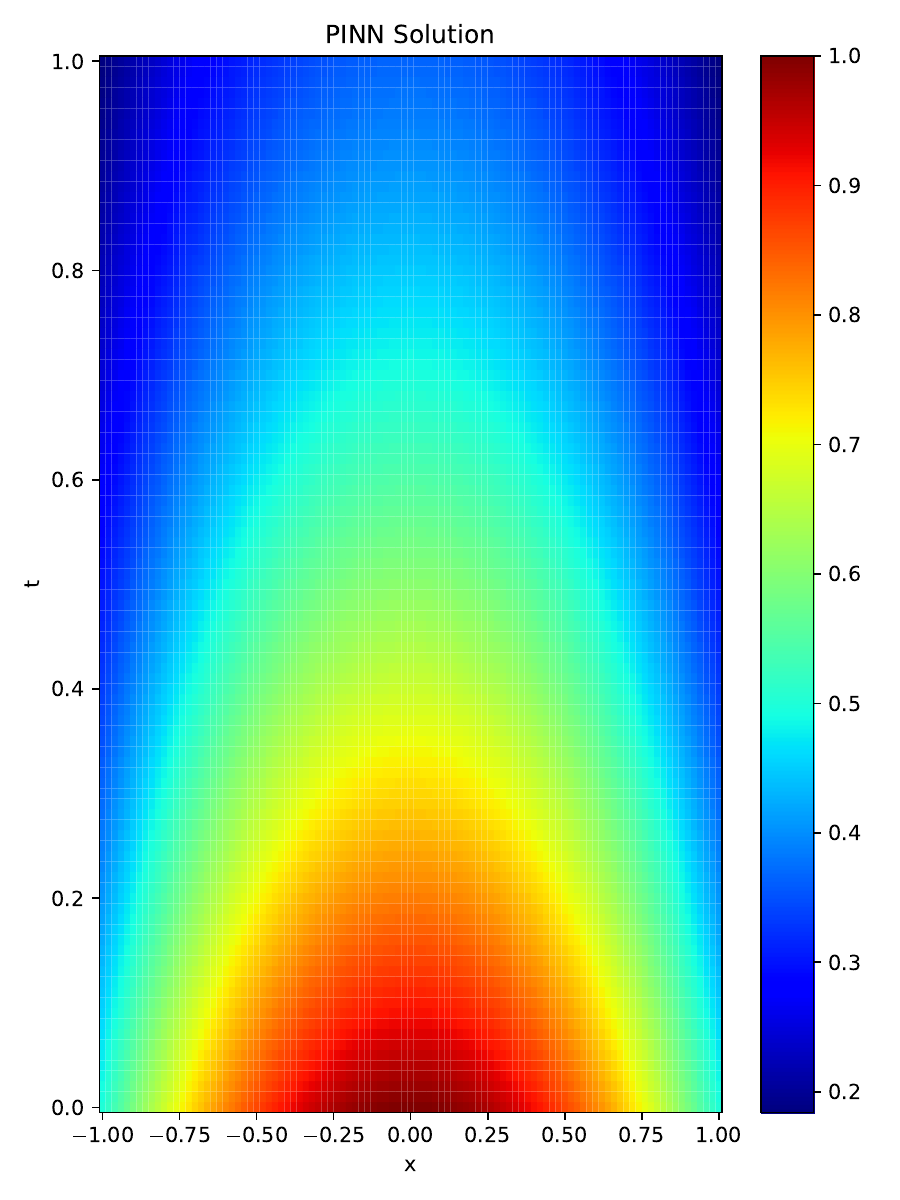}
        \caption{PINN solution}
    \end{subfigure}
    \hfill
    \begin{subfigure}[b]{0.32\textwidth}
        \centering
        \includegraphics[width=\textwidth]{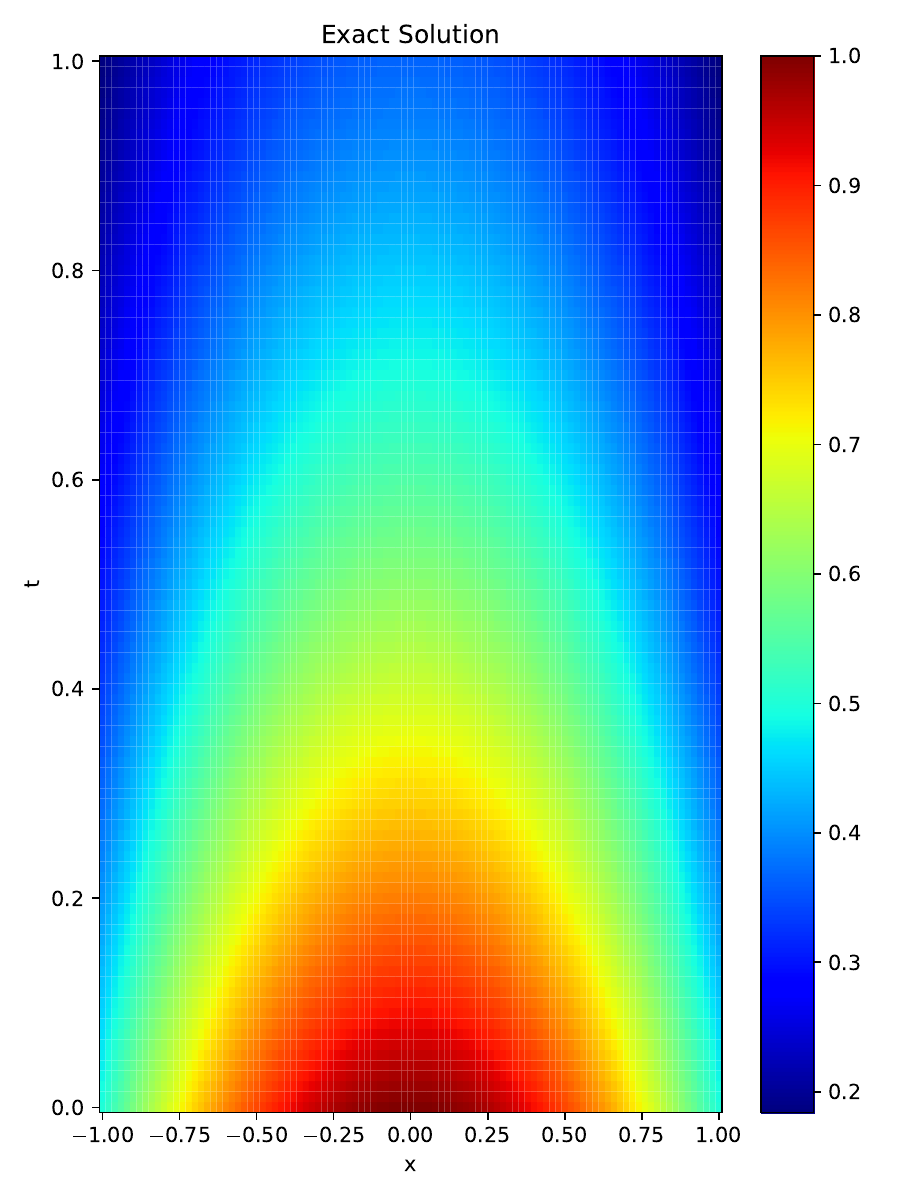}
        \caption{Exact solution}
    \end{subfigure}
    \hfill
    \begin{subfigure}[b]{0.32\textwidth}
        \centering
        \includegraphics[width=\textwidth]{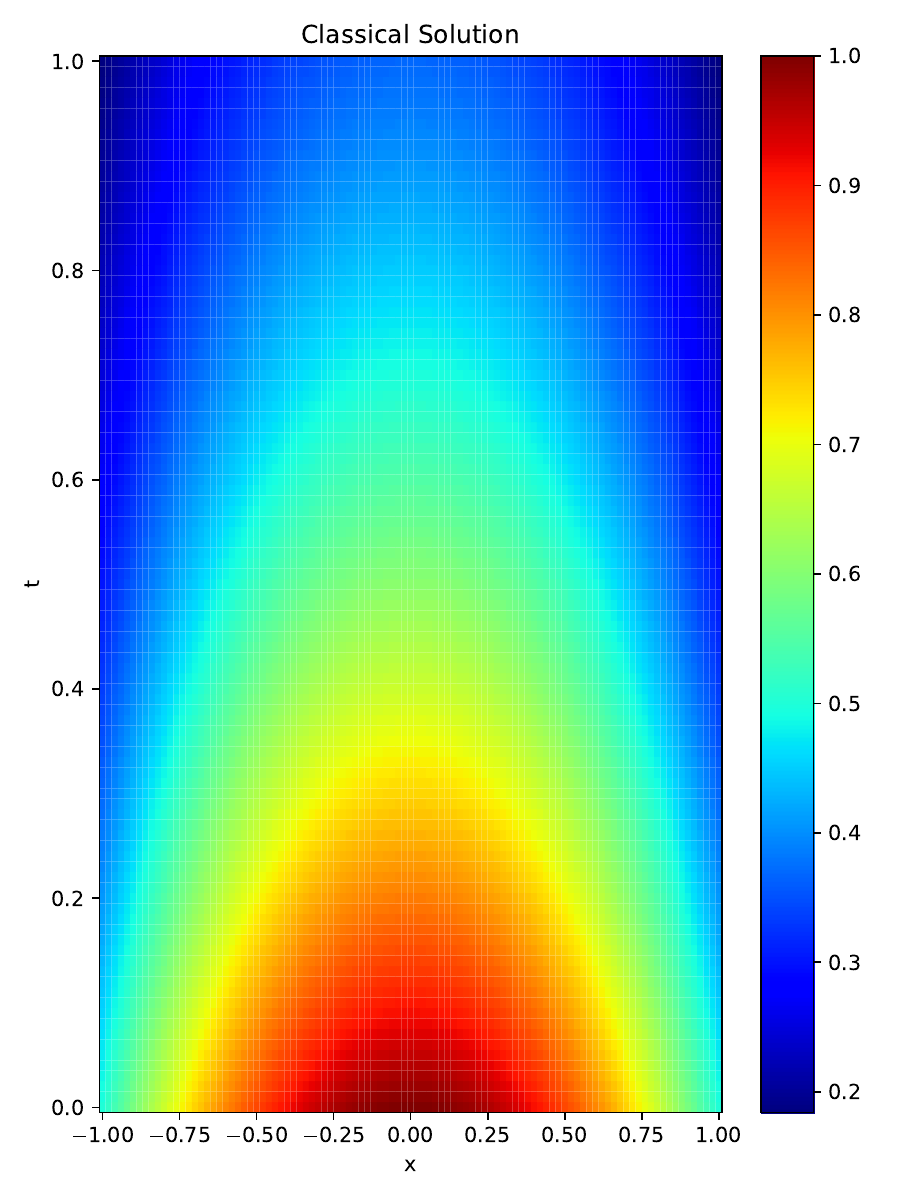}
        \caption{Classical solution}
    \end{subfigure}   
    \caption{visual comparison of the Harmonic numerical solutions }
\label{fig:solutions_harm}\vspace{-2mm}
\end{figure}
\noindent Table ~\ref{tab:harmonic_error} shows that the PINN outperforms the classical method in all cases by two orders of magnitudes. 
\begin{table}[H]
\centering
\caption{Relative \( L^2 \)-errors for the Damped Harmonic solution}
\label{tab:harmonic_error}\vspace{-1mm}
\begin{tabular}{|c|c|c|c|}
\hline
\( m \) & PINN vs Exact & Classical vs Exact & PINN vs Classical \\
\hline
2 & 2.081\( \times 10^{-5} \) & 1.292 \( \times 10^{-3}\) &  1.287 \( \times 10^{-3}\) \\
3 & 2.096 \( \times 10^{-5} \) & 1.449  \( \times 10^{-3}\)   & 1.451 \( \times 10^{-3}\) \\
4 & 3.656\( \times 10^{-5} \) & 1.723 \( \times 10^{-3}\)  & 1.722 \( \times 10^{-3}\) \\
5 & 5.156\( \times 10^{-5} \)  & 2.049\( \times 10^{-3}\) & 1.566 \( \times 10^{-3}\) \\
\hline
\end{tabular}
\end{table}\vspace{-2mm}
\begin{table}[H]
\centering
\begin{minipage}{0.48\textwidth}
\centering
\caption{Training epochs for the Harmonic solution}\vspace{-1mm}
\begin{tabular}{|c|c|c|}
\hline
$m$ & Adam epochs & LBFGS epochs \\ \hline
2 & 2000 & 210 \\ \hline
3 & 2000 & 140 \\ \hline
4 & 2000 & 65 \\ \hline
5 & 2000 & 130 \\ \hline
\end{tabular}
\label{tab:har_epochs}
\end{minipage}
\hfill
\begin{minipage}{0.48\textwidth}
\centering
\caption{Time (sec) for the Harmonic solution}\vspace{-1mm}
\begin{tabular}{|c|c|c|}
\hline
$m$ & PINN training time & Classical time \\ \hline
        2 & 149.95 & 5.38 \\ \hline
        3 & 112.42 & 6.33 \\ \hline
        4 & 84.34 & 5.29\\  \hline
        5 & 108.06 & 6.21\\  \hline
\end{tabular}
\label{tab:har_time}
\end{minipage}

\end{table}

\subsubsection*{m-dependent solution Results ~\ref{sec: m-dep}}
This solution was harder to approximate using both the classical and the PINN methods. Classically, the higher the m is, the more time it takes to give a solution (Table \ref{tab:mdep_time}). Meanwhile, the early stopping was triggered in all cases in the PINN case. Moreover, the error increases as m increases in both methods. Table \ref{tab: m-dep error} shows that for \( m = 2 \), the smallest tested value of \( m \), both the PINN and the classical method have errors of the same order of magnitude, with the PINN being slightly more accurate. As \( m \) increases, the errors of both methods also increase. However, the classical method remains more accurate, while the PINN error grows more rapidly. For values of \( m \) greater than 2, the classical method gives errors in an approximate range of 0.8\% to 2\%, whereas the PINN errors range from about 2\% to 15\%.\vspace{-2mm}
\begin{figure}[H]
    \centering
\begin{subfigure}[b]{0.32\textwidth}
        \centering
        \includegraphics[width=\textwidth]{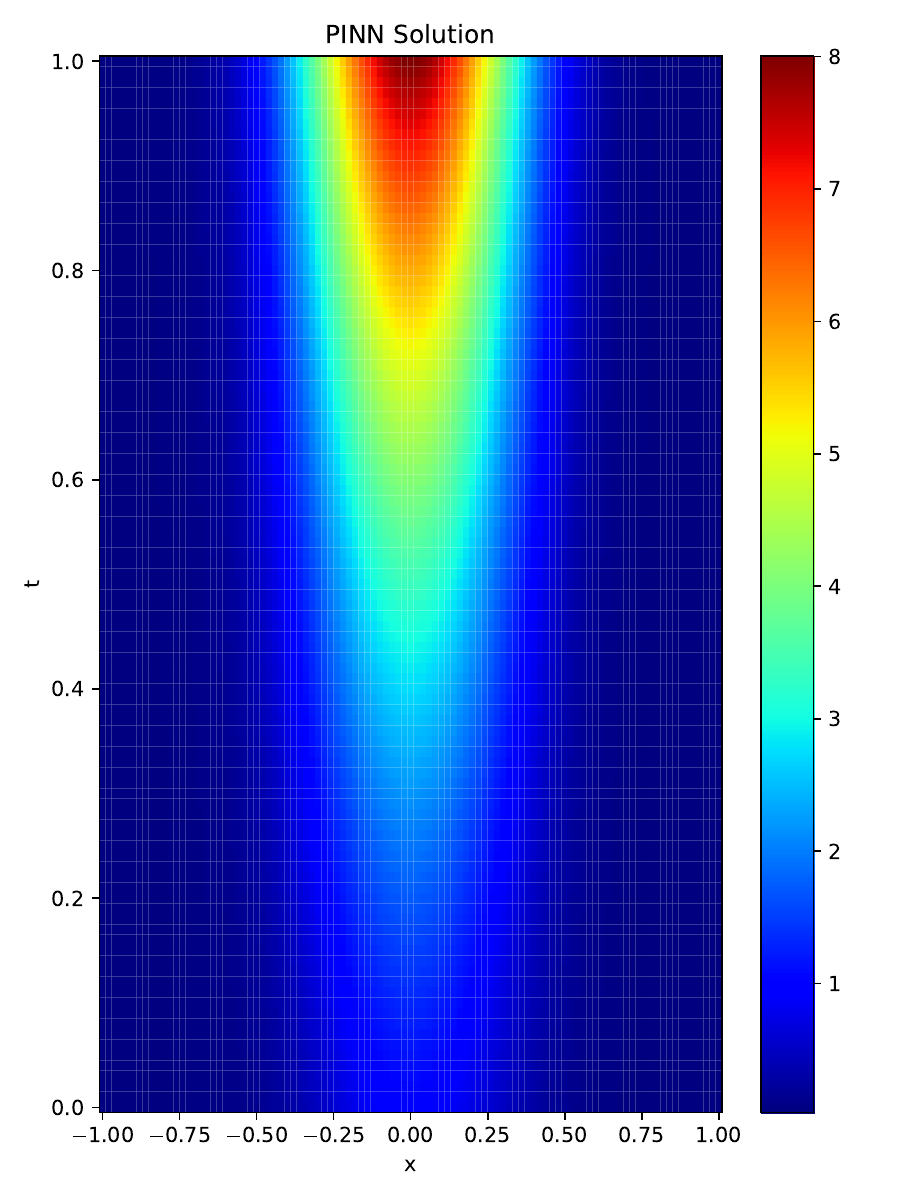}
        \caption{PINN solution}
    \end{subfigure}
    \hfill
    \begin{subfigure}[b]{0.32\textwidth}
        \centering
        \includegraphics[width=\textwidth]{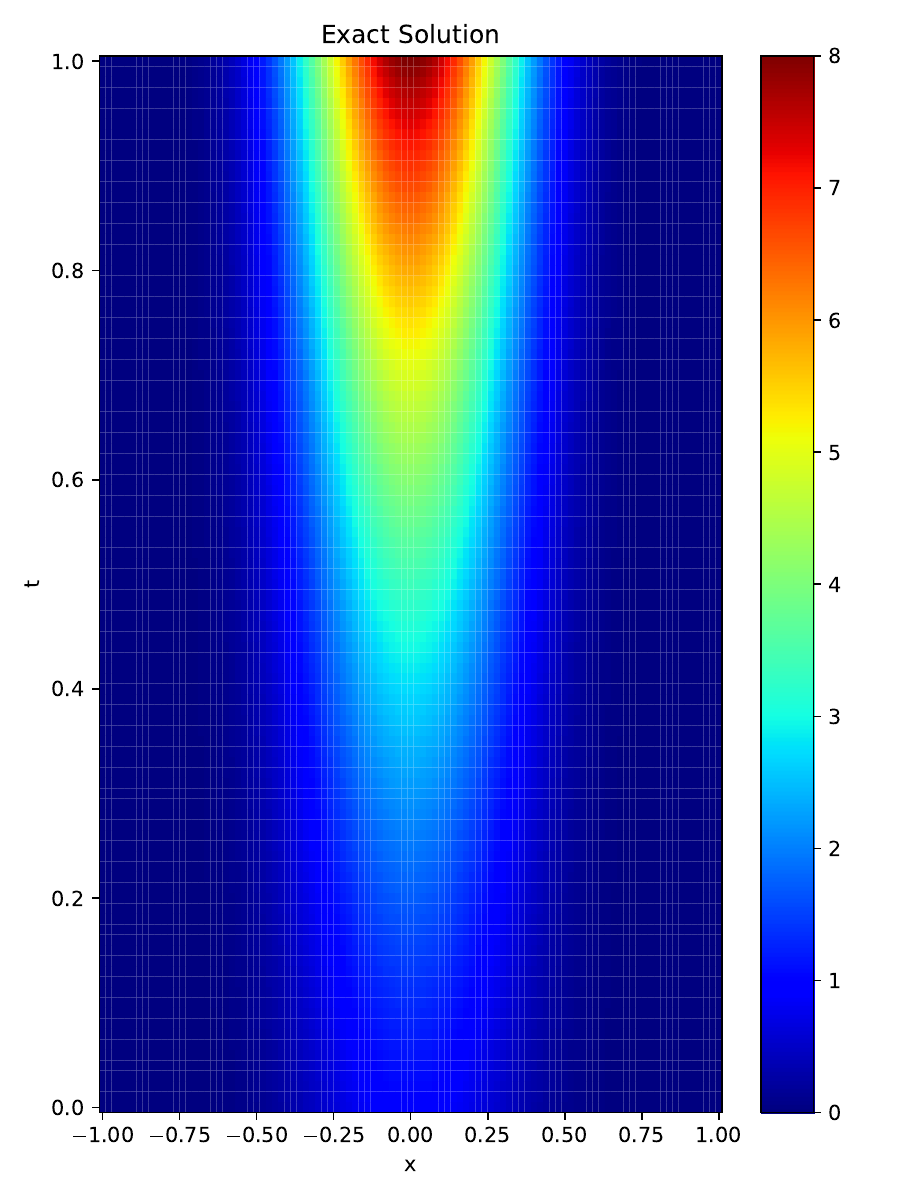}
        \caption{Exact solution}
    \end{subfigure}
    \hfill
    \begin{subfigure}[b]{0.32\textwidth}
        \centering
        \includegraphics[width=\textwidth]{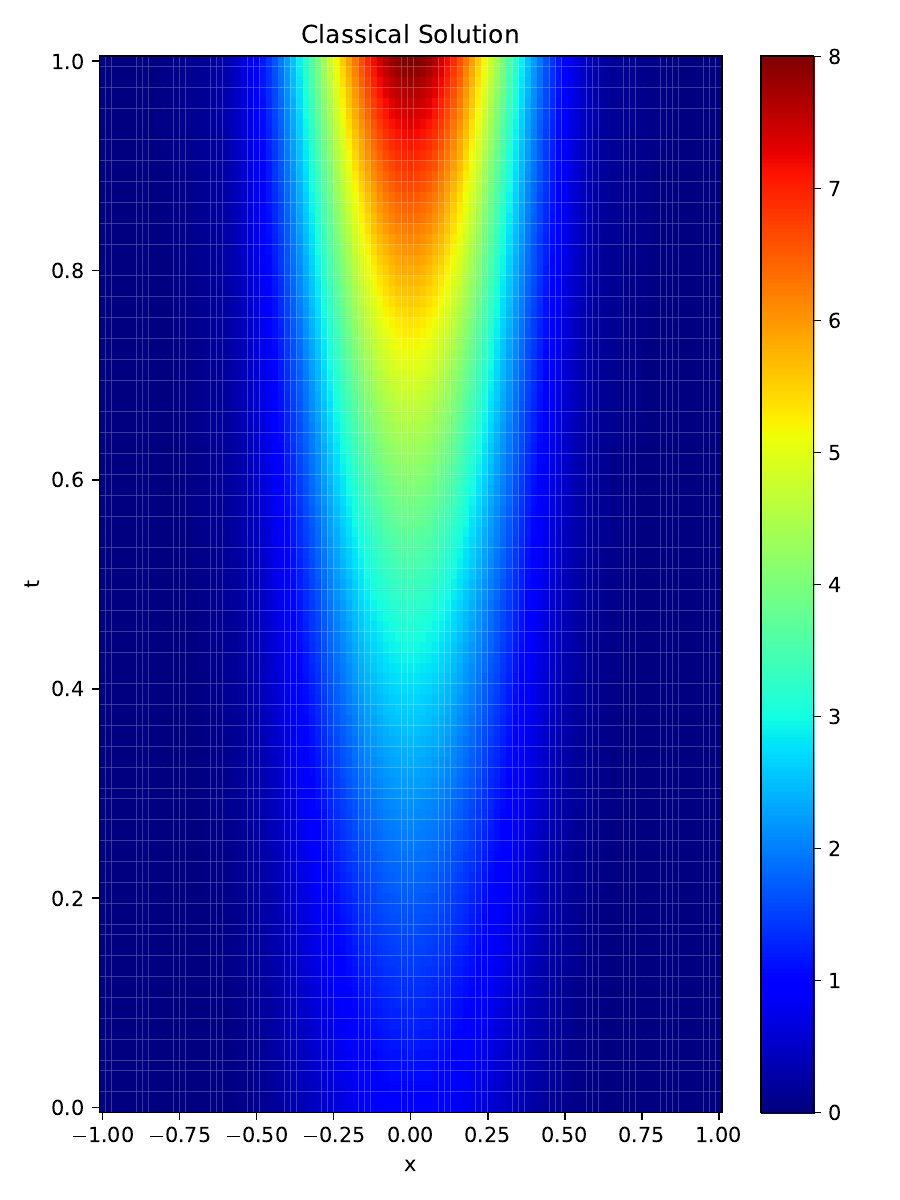}
        \caption{Classical solution}
    \end{subfigure}
    \caption{visual comparison of the m-dependent numerical solutions }
\label{fig:solutions_mdep}\vspace{-4mm}
\end{figure}
\noindent Note that we can always keep improving the PINN to our problem by tuning more of the hyperparameters, changing architectures, or doing some computational tricks that can reduce the effect of nonlinearity. However, for the tests and comparisons to be adequate and fair across all the test cases and all the different solutions we have, we decided to fix our model to what is detailed in section ~\ref{sec: direct pinn} and perform the tests. \vspace{-4mm}
\begin{table}[H]
\centering
\caption{Relative \( L^2 \)-errors for the m-dependent solution}
\label{tab: m-dep error}\vspace{-2mm}
\begin{tabular}{|c|c|c|c|}
\hline
\( m \) & PINN vs Exact & Classical vs Exact & PINN vs Classical \\
\hline
2 & 1.586 \( \times 10^{-3} \) & 2.715 \( \times 10^{-3}\) &  1.636\( \times 10^{-3}\) \\
3 & 1.651 \( \times 10^{-2} \) & 8.514 \( \times 10^{-3}\)   & 1.311\( \times 10^{-2}\) \\
4 & 7.679 \( \times 10^{-2} \) & 1.437 \( \times 10^{-2}\)  & 7.074  \( \times 10^{-2}\) \\ 
5 & 1.505\( \times 10^{-1} \)  & 1.802 \( \times 10^{-2} \)
& 1.420\( \times 10^{-1} \)
\\
\hline
\end{tabular}
\end{table}\vspace{-4mm}
\begin{table}[H]
\centering
\begin{minipage}{0.48\textwidth}
\centering
\caption{Training epochs for the m-dependent sol.}\vspace{-2mm}
\begin{tabular}{|c|c|c|}
\hline
$m$ & Adam epochs & LBFGS epochs \\ \hline
2 & 2000 & 200 \\ \hline
3 & 5 & 150 \\ \hline
4 & 5 & 120 \\ \hline
5 & 5 & 90 \\ \hline
\end{tabular}
\label{tab:mdep_epochs}
\end{minipage}
\hfill
\begin{minipage}{0.48\textwidth}
\centering
\caption{Time (sec) for the m-dependent solution}\vspace{-2mm}
\begin{tabular}{|c|c|c|}
\hline
$m$ & PINN training time & Classical time \\ \hline
        2 & 62.89 & 5.66 \\ \hline
        3 & 56.53 & 10.09 \\ \hline
        4 & 43.45 & 11.51 \\ \hline
        5 & 36.141 & 24.33 \\ \hline
\end{tabular}
\label{tab:mdep_time}
\end{minipage}
\end{table}

\subsection{Direct PINN Testing with Data}
\label{direct noise}

A natural extension of the direct PINN is to include additional measurement data in the training process rather than minimizing the loss only over the PDE residual together with the initial and boundary conditions. This setting corresponds to the practical scenario in which the governing PDE is known, but a set of observations of the solution is also available. In such a case, the PINN can exploit both the physical constraints imposed by the PDE and the available measurements, potentially improving the quality of the approximation or the time complexity. This experiment is performed only for the polynomial solution. The polynomial test case was selected because, for different values of $m$, the relative errors obtained by the standard direct PINN range from approximately $10^{-2}$ to $10^{-4}$. This wide range is representative of the error magnitudes observed for the other manufactured solutions and therefore provides a suitable benchmark for assessing the influence of incorporating measurement data.

The only modification with respect to the loss function defined in~\eqref{eq: loss direct} is the addition of a measurement loss term. The total loss becomes
\[
L=\log_{10}\!\left(\lambda_u L_u + L_{\mathrm{int}} + \lambda_{\mathrm{meas}}L_{\mathrm{meas}}\right),
\]
where $\lambda_{\mathrm{meas}}$ is a hyperparameter chosen to be $10$, and
\[
L_{\mathrm{meas}}
=
\frac{1}{N_{\mathrm{meas}}}
\sum_{i=1}^{N_{\mathrm{meas}}}
\left(
u_\theta(t_i,x_i)-u_i^{\mathrm{data}}
\right)^2.
\]
We choose $\lambda_{\mathrm{meas}}=10$, equal to the weight assigned to the boundary and initial condition loss ($\lambda_u=10$), since both terms enforce agreement with known pointwise values of the solution, whereas the interior loss enforces the governing differential equation. The measurement dataset consists of $N_{\mathrm{meas}}=225$ exact solution values sampled on a uniform $15\times15$ grid over the computational domain $(x,t)\in[-1,1]\times[0,1]$. 

The results are summarized in Table~\ref{tab: direct poly w meas} and compared with the standard direct PINN results reported in Tables~\ref{tab:poly_time} and~\ref{tab:polynomial_error}. The inclusion of the measurement loss does not lead to a significant improvement in the approximation accuracy. The relative \(L^2\)-errors remain of the same order of magnitude as those obtained with the standard PINN. For \(m=4\) and \(m=5\), the errors slightly decrease from \(1.255\times10^{-3}\) to \(1.008\times10^{-3}\) and from \(1.740\times10^{-3}\) to \(1.622\times10^{-3}\), respectively. However, for \(m=2\) and \(m=3\), the addition of measurement data slightly increases the error from \(1.548\times10^{-2}\) to \(1.716\times10^{-2}\) and from \(8.050\times10^{-4}\) to \(8.214\times10^{-4}\), respectively. Also, the inclusion of measurement data introduces an additional computational cost, as the network must minimize an extra loss term.

\begin{table}[h!]
\centering

\begin{tabular}{|l|l|l|}
\hline
m & PINN relative error                                    & PINN training time \\ \hline
2 & 1.716  $ \times 10^{-2}$ & 60.74              \\ \hline
3 & 8.214 $ \times 10^{-4}$ & 207.00             \\ \hline
4 & 1.008 $ \times 10^{-3}$ & 193.70             \\ \hline
5 & 1.622 $ \times 10^{-3}$ & 191.69             \\ \hline
\end{tabular}
\caption{Relative \( L^2 \)-errors and training time (second) for the direct PINN applied to the polynomial solution using exact data}
\label{tab: direct poly w meas}
\end{table}
\noindent Therefore, for this test case, where the standard direct PINN already provides accurate solutions, adding additional measurement data is not necessary, as it does not provide a meaningful improvement in accuracy while increasing the training time. \\
\noindent It was also worth performing this test on the m-dependent solution, where the standard direct PINN produced larger errors, exceeding \(5\%\) for \(m=4\) and \(m=5\). The aim is to check whether adding measurement data helps in the cases where the standard PINN struggles to converge. The results are reported in the following Table~\ref{tab: mdep w meas}: 

\begin{table}[H]
\centering
\caption{Relative \(L^2\)-errors and training time (second) for the direct PINN applied to the m-dependent solution using exact data}
\label{tab: mdep w meas}
\begin{tabular}{|l|l|l|}
\hline
m & PINN relative error & PINN training time \\ \hline
2 & 1.109 \(\times 10^{-3}\) & 69.66 \\ \hline
3 & 1.283 \(\times 10^{-2}\) & 71.18 \\ \hline
4 & 6.747 \(\times 10^{-2}\) & 71.32 \\ \hline
5 & 2.284 \(\times 10^{-1}\) & 46.68 \\ \hline
\end{tabular}
\end{table}

\noindent The addition of measurement data does not resolve the difficulty observed for the m-dependent solution. For \(m=2\), \(m=3\), and \(m=4\), the relative error decreases slightly, but it remains of the same order of magnitude as before, and the errors for \(m=4\) and \(m=5\) stay well above \(5\%\). For \(m=5\), the error even increases, from \(1.505\times10^{-1}\) to \(2.284\times10^{-1}\). At the same time, the training time increases in all four cases. These observations are consistent with the conclusion drawn from the polynomial test: adding measurement data does not lead to a meaningful improvement in accuracy and does not rescue the cases where the standard PINN fails, while it increases the training time.\vspace{2mm}

\noindent In summary, for the direct problem, the PINN achieves accurate solutions close to the classical methods, where the PINN performs better in some cases and the classical method performs better in others, except for the \(m\)-dependent solution case, where the difference between the PINN and the classical method is relatively larger than in the other solution cases. Also, in terms of time, the classical method is always faster. However, one should note that the time data presented in the tables for the PINN include the time taken to trigger early stopping, which is equal to the time required to run 200 epochs without notable improvement. Nevertheless, this would still be negligible, and the classical method remains faster. 

\section{The 2-stage Learning Inverse PINN } \label{sec:inverse_pme1}

Since the direct PINN approach provided accurate and stable results, we now extend it to the inverse problem, highlighting the flexibility of PINNs, where only minor modifications are needed to move from solving the direct problem to solving the inverse one.
The inverse problem of the PME consists of determining the polytropic parameter $m$ from given measurements. In our study, these measurements are generated from exact reference solutions. 

\noindent We begin by applying a classical optimization-based method using MATLAB \texttt{fmincon}  in Section~\ref{sec:classical}. Next, we present a PINN to solve the same inverse problem in Section~\ref{sec:pinnInv}, where we propose our novel inverse PINN model ~\ref{sec:proposed}. Finally, we compare and analyze the performance of these methods through numerical testing in Section~\ref{sec:testing}.

\subsection{Classical Approach:  \texttt{fmincon}}\label{sec:classical}

MATLAB's  \texttt{fmincon} is an optimization solver that is designed to find the constrained minimum of a multi-variable scalar function from an initial guess. Such a procedure, generally called nonlinear programming or constrained nonlinear optimization, is most appropriate whenever parameters must be estimated within physically reasonable constraints. In this work, we use it to estimate the parameter \( m \) of the PME by minimizing the difference between synthetic reference data and our numerical solution, under the constraint \( 1 \le m \le 10 \). 

\noindent The implementation can be summarized as follows:

\begin{enumerate}
    \item We generate synthetic data from solutions ~\ref{sec: ref sol} with a known \( m \).
    \item Define an objective function that computes the error between the numerical and exact solutions 
    \[
J(m) = \sum_{i,k} \big[ u_{\text{num}}(x_i,t_k;m) - u_{\text{exact}}(x_i,t_k) \big]^2,
\]
    \item For a given \( m \), solve the nonlinear PDE system using Newton’s method at each time step to obtain \( u_{\text{num}}(x,t;m) \).
    \item Use MATLAB’s \texttt{fmincon} to iteratively update \( m \) and minimize the objective \( J(m) \), subject to the bounds \( 1 \le m \le 10 \).
\end{enumerate}

\noindent A pseudo-code of the corresponding inverse problem algorithm using the Barenblatt solution is found in Algorithm 8 in \cite{NSG}. 

\subsection{PINN Approach}\label{sec:pinnInv}
Physics-informed neural networks are particularly flexible because once a direct PINN is implemented, only minor modifications are needed to adapt it to an inverse problem. Instead of restating the full direct PINN formulation, in this section we focus only on the adjustments required to estimate unknown parameters from data. Below, we highlight the differences in the loss function, network outputs, and training procedure that distinguish the inverse PINN from the direct case.

\subsubsection*{Standard Inverse PINN Model} \label{sec: standard}

\begin{figure}[H]
\begin{subfigure}{0.5\textwidth}
\includegraphics[width= 1\linewidth, height=6cm]{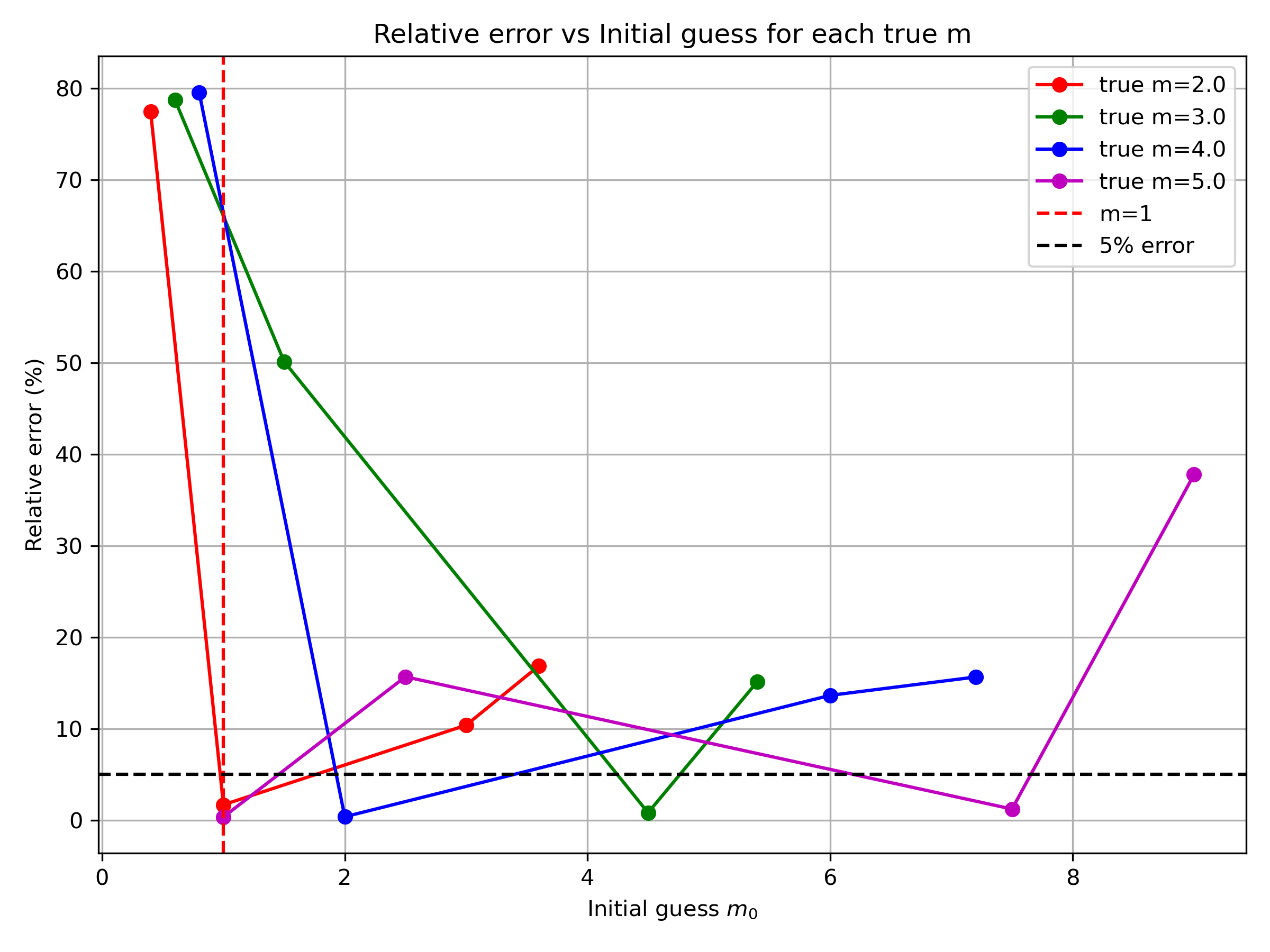} 
\caption{Standard Inverse PINN}
\label{fig:subim1}
\end{subfigure}
\begin{subfigure}{0.5\textwidth}
\includegraphics[width=1\linewidth, height=6cm]{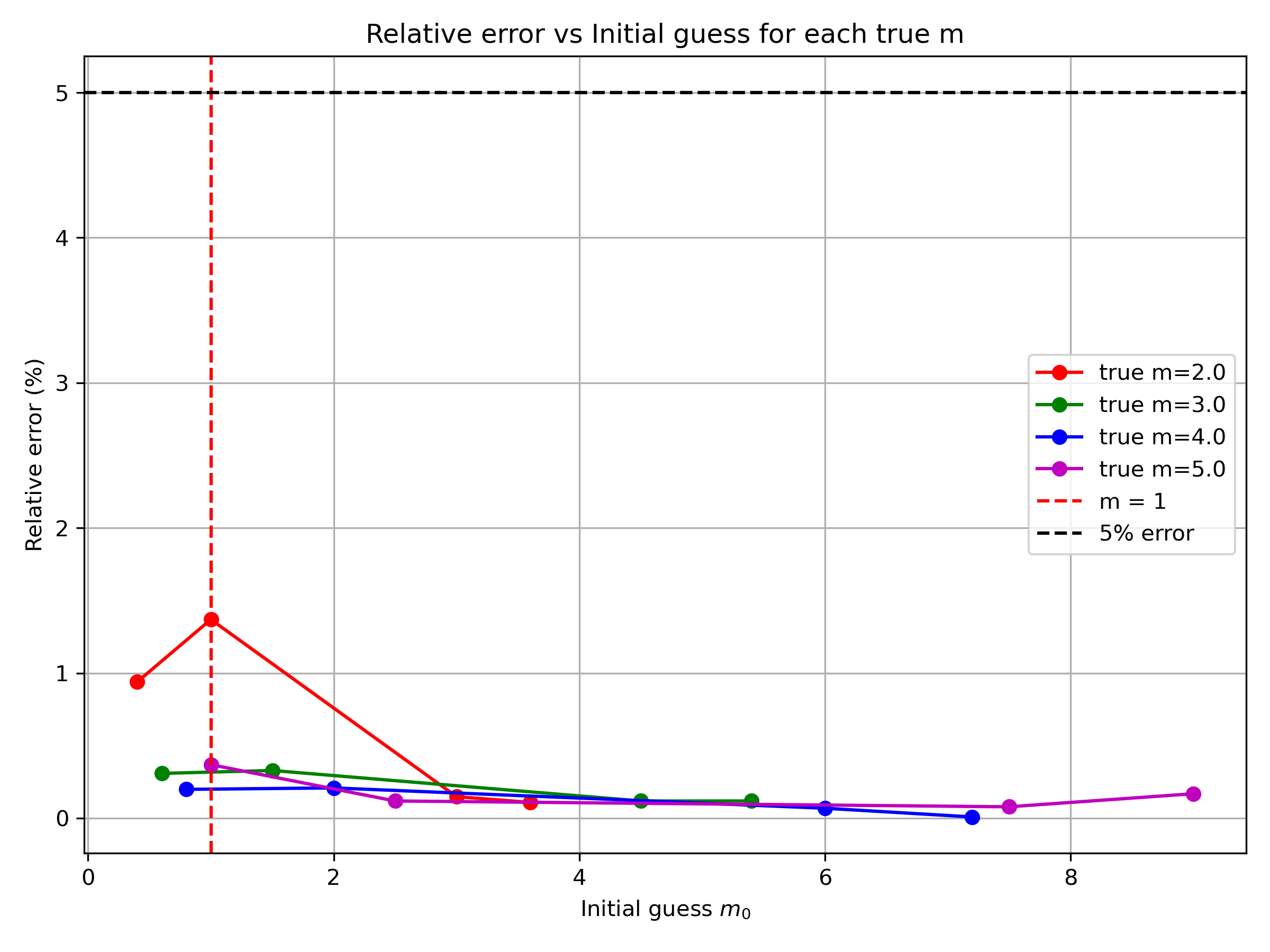}
\caption{ 2-Stage Training PINN }
\label{fig:subim2}
\end{subfigure}

\caption{Relative error of the estimated parameter \(m\) with respect to different initial guesses. Each color corresponds to a different true value of \(m\). The x-axis represents the initial guess \(m_0\), taken around the true value ($\pm 50\%$ and $80\%$ perturbations). The y-axis shows the relative error (\%). The left plot shows the standard inverse PINN, where the solution is accurate only when the initial guess is close to the true value, while the right plot shows the results after applying the 2-stage training PINN, where the method achieves significantly lower and more stable errors across different initial guesses.}
\label{fig: 2stage}
\end{figure}

\begin{enumerate}
\label{standard}
    \item \textbf{Neural Network Architecture:} Same architecture and initialization, but we treat $m$ as trainable real scalar parameter.
    \item \textbf{Loss function and optimization:} In addition to the previous loss elements, we add the measurement loss 
    \begin{equation}
    \mathcal{L}_\theta = \log_{10} \left( 
    \lambda_u \cdot (\mathcal{L}_b + \mathcal{L}_t) +
    \mathcal{L}_{\text{PDE}} +
    \lambda_m \cdot \mathcal{L}_{\text{meas}}
    \right),
    \end{equation}
     with \( \mathcal{L}_{\text{meas}} \)the data mismatch term:
    \[
    \mathcal{L}_{\text{meas}} = \mathbb{E}_{(t,x)\in\mathcal{D}_{\text{meas}}} \left[ \left( u_\theta(t,x) - u_{\text{meas}}(t,x) \right)^2 \right],
    \]
    where \( \mathcal{D}_{\text{meas}} \) denotes the set of measurement points where we evaluate the exact solution. These measurement points are generated using the corresponding exact solution.
    \item \textbf{Training Procedure} Same as in the direct PINN ~\ref{sec: training}    
\end{enumerate}
After these standard modifications for the inverse PINN, the neural network was able to converge only locally, where it provided a good approximation of \(m\) only when the initial guess was close to the true value, and large approximation errors otherwise. Therefore, aiming for a global convergence, we started experimenting with several aspects, including different initialization, architecture, optimization techniques, and other details. Among all tested approaches, the most significant improvement came from the introduction of a novel two-stage training PINN framework, which substantially stabilized the inverse optimization process and improved convergence across different initial guesses.
Figure ~\ref{fig: 2stage} shows the results of the inverse problem for the Harmonic solution before and after the two-stage Training Model. For a comparison of the other solutions, see Appendix ~\ref{app:C}. The standard approach is the one in which the only changes over the direct PINN model are the three enumerated points in ~\ref{standard}, while the two-stage learning model will be detailed in the following section.

The comparison in Figure ~\ref{fig: 2stage} highlights the main contribution of this work. While the standard inverse PINN behaves as a locally convergent method that strongly, the proposed two-stage framework significantly improves robustness and allows accurate parameter recovery across a much wider range of initial guesses. This demonstrates that the training schedule itself plays a critical role in stabilizing inverse PINN optimization for nonlinear PDEs such as the PME.

\subsection{Inverse PINN Design Experiments} \label{sec: experiments}

In the following parts of this section, for each PINN component, we explain the different trials we conducted and highlight what worked best and was used in the final proposed two-stage inverse PINN framework ~\ref{sec: experiments} that produced the results in section ~\ref{sec:testing}. The following modifications collectively led to the final inverse PINN model ~\ref{sec:proposed} proposed in this work.

\subsubsection*{Neural Network Architecture}
For the inverse PINN, we first started experimenting by increasing and decreasing the length (i.e., number of neurons per layer) and depth (i.e., the number of layers of the neural network) of our neural network, but no clear improvements were shown. Then, as suggested and discussed in the literature \cite{chen2025rethinkingshapeconventionmlp}, we tried the narrow-wide-narrow and wide-narrow-wide architectures, but again the initial configuration used in the direct problem ~\ref{tab: arch} gave the best results, so we kept it. 

\subsubsection*{Weights}

We experimented with several weight initializations to improve training dynamics and reduce the error in estimating $m$:
\begin{itemize}
    \item Xavier normal: Xavier Glorot normal initialization
    \item Xavier uniform: Xavier Glorot uniform initialization
    \item kaiming\_normal: He normal initialization
    \item orthogonal: Orthogonal initialization
    \item scaled\_normal: Custom scaled normal initialization
\end{itemize}

\noindent Initially, in the direct problem,  we used the Xavier uniform scheme, but it was shown that the Xavier normal initialization consistently produced better convergence and more accurate parameter recovery in the inverse setting. Even though both Xavier uniform and normal initializations have the same variance and aim to prevent vanishing or exploding gradients, the normal distribution, due to its unbounded support and variance characteristics \cite{glorot2010understanding}, proved more effective for our PINN training.

\noindent Mathematically, the Glorot initialization selects the variance of each layer’s weights such that forward activations and backward gradients remain stable across layers. For a layer with n-in and n-out number of neurons, the target variance is given by\vspace{-2mm}
\[
\mathrm{Var}(w) = \frac{2}{\text{n\_in} + \text{n\_out}}.
\]
This variance can be realized either through a uniform distribution,
\[
w \sim U\left(-g\sqrt{\frac{6}{\text{n\_in} + \text{n\_out}}},\;
g\sqrt{\frac{6}{\text{n\_in} + \text{n\_out}}}\right),
\]
or through a normal distribution,\vspace{-4mm}
\[
w \sim \mathcal{N}\left(0,\; g^2 \frac{2}{\text{n\_in} + \text{n\_out}}\right),
\]
with $g=1$ in our case. The Xavier normal initialization gave us better results in terms of errors.
\subsubsection*{Stability}

We also did the following modifications to enhance stability:
\begin{itemize}
    \item Since $m$ is known to be positive, we parameterize it as
\[
m = \exp(\widehat{m}),
\]
This formulation guarantees $m>0$ by construction and avoids the need for manual clipping.

\item  we applied L2 gradient clipping to prevent instability caused by large gradients during backpropagation. At each training step, the total gradient norm $\|\mathbf{g}\|$ over all trainable parameters is computed. If
\[
\|\mathbf{g}\| > \text{max\_norm},
\]
the gradients are rescaled according to
\[
\mathbf{g} \leftarrow \mathbf{g}\cdot\frac{\text{max\_norm}}{\|\mathbf{g}\|},
\]
with $\text{max\_norm}=1$ in our implementation. This limits excessively large parameter updates specifically encountered when $m $ is big,  and prevents the occurrence of NaN values during training.

\end{itemize}

\subsubsection{The Proposed Two-Stage Training Inverse PINN Model}
\label{sec:proposed}
We increased the number of training points, so the model used 2048 interior points \((n_{\text{int}})\), 256 spatial boundary points per side \((n_{\text{sb}})\), 256 temporal boundary points \((n_{\text{tb}})\), and 225 measurement points \((n_{\text{meas}} = 15 \times 15 )\) uniformly spread over the domain.
\noindent We initially tried to use the same training procedure as in the direct problem (jointly minimizing PDE, initial/boundary, and measurement losses with fixed hyperparameters (i.e., loss weights)). However, we observed that the recovered $m$ was highly sensitive to the manual choice of loss weights, and we could not obtain a reliable convergence behavior across different initial guesses $m_0$ (i.e., the model did not consistently converge to the true $m$ for different initial guesses). A natural idea was then to make the loss weights learnable parameters and optimize them during training, but this behaved poorly in practice: it led to unstable training and catastrophic performance even for the direct problem. This suggested that, in our setting, adding more trainable degrees of freedom makes optimization harder rather than easier. 

Motivated by these observations, we propose the two-stage training PINN, which is specifically designed to address these problems. The main idea behind the proposed two-stage training PINN is to separate the learning of the solution structure from the learning of the physical parameter \(m\). Instead of optimizing both \(u_\theta\) and \(m\) simultaneously from the beginning, the network is first guided toward learning a physically meaningful solution profile before updating the unknown parameter.
\paragraph*{Two-stage strategy}
Since we already validated in the direct problem that the network can approximate accurate solutions when the PDE parameter is known, we adopted the following training schedule:
\begin{enumerate}
    \item \textbf{Stage 1 (Warm-up, 20\% of epochs):} $m$ is frozen at the initial guess $m_0$. The network $u_\theta(t,x)$ is trained to satisfy the PDE and initial/boundary conditions as if it is solving the direct problem with m = $m_0$ the initial guess.
    \item \textbf{Stage 2 (Main training):} $m$ is unfrozen and we jointly optimize $(\theta,m)$.
\end{enumerate}
Intuitively, Stage 1 helps the model learn a reasonable solution structure before it starts adjusting the physical parameter. This reduces the tendency of the optimizer to use $m$ to compensate for a poorly learned $u_\theta$ through incorrect updates of \(m\), which was one of the main causes of unstable convergence in the standard inverse PINN formulation. Also, in this inverse training, we only used ADAM optimizer with 10,000 total epochs. 
\paragraph*{Curriculum via hyperparameter annealing for measurements}
Even with two-stage training, the choice of $\lambda_m$ remains important. Instead of selecting a single fixed value, we used a simple curriculum learning (sometimes referred to as loss-weight annealing): during the warm-up, $\lambda_m$ is kept small at 10, and after warm-up it is increased gradually to its final value at 100. Concretely, letting $E$ be the total number of epochs and $E_w$ the warm-up epochs, we use:
\[
\lambda_m(e) =
\begin{cases}
\lambda_{m}^{\text{init}}, & e < E_w,\\
\lambda_{m}^{\text{init}} + \dfrac{e-E_w}{E-E_w}\Big(\lambda_{m}^{\text{final}}-\lambda_{m}^{\text{init}}\Big), & e \ge E_w,
\end{cases}
\]

\begin{figure}[H]
    \centering
    \includegraphics[width=0.8\linewidth]{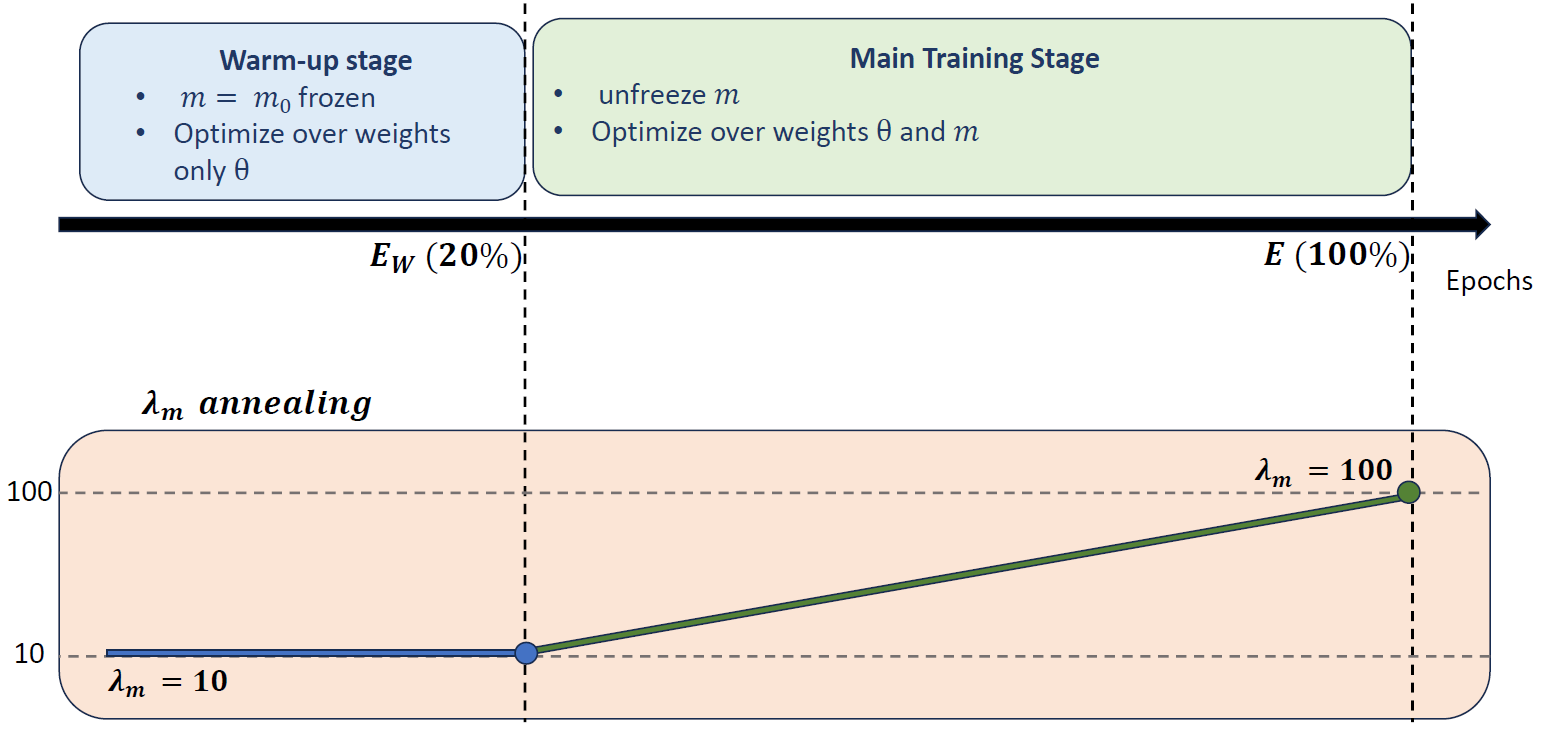}
    \caption{Training Diagram}
    \label{fig:placeholder}
\end{figure}

The pseudo-code of the full procedure is summarized in Algorithm~\ref{alg:two_stage}.

\begin{algorithm}[H]
\caption{Two-Stage Training Pseudo-code for the Inverse PINN}
\label{alg:two_stage}
\begin{algorithmic}[1]
\Require Initial guess $m_0$; total epochs $E$; warm-up fraction $f_w$; $\lambda_u$; $\lambda_m^{\text{init}}, \lambda_m^{\text{final}}$
\State $E_w \leftarrow \lfloor f_w \cdot E \rfloor$

\Statex
\Statex \textbf{Stage 1:  Warm-up} ($e = 0, \ldots, E_w - 1$)
\State Freeze $m$ at $m_0$ \quad ($\nabla_m \leftarrow 0$)
\For{$e = 0$ \textbf{to} $E_w - 1$}
    \State Update $\theta$ only by minimizing $\log_{10}\!\left(\lambda_u(\mathcal{L}_b + \mathcal{L}_t) + \mathcal{L}_{\text{PDE}} + \lambda_m^{\text{init}}\,\mathcal{L}_{\text{meas}}\right)$
\EndFor

\Statex
\Statex \textbf{Stage 2: Main Training} ($e = E_w, \ldots, E - 1$)
\State Unfreeze $m$
\For{$e = E_w$ \textbf{to} $E - 1$}
    \State $\lambda_m(e) \leftarrow \lambda_m^{\text{init}} + \dfrac{e - E_w}{E - E_w}\left(\lambda_m^{\text{final}} - \lambda_m^{\text{init}}\right)$
    \State Update $(\theta,\, m)$ jointly by minimizing $\log_{10}\!\left(\lambda_u(\mathcal{L}_b + \mathcal{L}_t) + \mathcal{L}_{\text{PDE}} + \lambda_m(e)\,\mathcal{L}_{\text{meas}}\right)$
\EndFor

\Return Trained $\theta$, estimated $m$
\end{algorithmic}
\end{algorithm}

\paragraph*{Relation to prior work.}
The idea of splitting training into stages is consistent with existing PINN practices. For example, a two-stage PINN has been used in the direct setting where a second stage adds extra physical structure (e.g., conserved quantities) to refine the solution \cite{lin2022twoStageConserved}. In a related direction, a two-stage learning strategy has been proposed where a pre-training phase helps boundary losses converge, followed by a main phase that balances the full objective \cite{shi2023atpinn}. Our warm-up (freeze-unfreeze) schedule follows the same general motivation: stabilize training by learning easier constraints first.

\noindent Our gradual increase of $\lambda_m$ is also aligned with curriculum learning ideas. Curriculum learning was introduced as training models by presenting easier samples/constraints first and progressively increasing difficulty \cite{bengio2009curriculum}. Recently, curriculum strategies have been adapted to PINNs and shown to improve convergence and solution quality \cite{duffy2024dynamicCurriculumPINN}. Similar curriculum training has also been used in time-dependent physics problems by splitting training along temporal intervals and training progressively over time \cite{bekele2024pinnCurriculumPoroelastic}. In our inverse PME setting, we apply this principle at the loss level (annealing $\lambda_m$) to smoothly introduce measurement data while maintaining stable PDE/BC learning in our inverse problem.

\subsection{Testings and results} \label{sec:testing}
We test our classical and PINN approaches on all the different reference solutions for \(m \in \{2, 3, 4, 5\}\). 
To study the sensitivity of each method to the initial condition, each value of \(m\) is tested using four different initial conditions corresponding to \(\pm 50\%\) and \(\pm 80\%\) of the value. 
To better analyze the results, they are represented by plots showing the relative error versus the initial guess for every value of \(m\) and for all four reference solutions. 
For this section, the PINN codes were run in Python on Google Colab with an NVIDIA Tesla T4 GPU \cite{googlecolab}, while the classical codes were run locally in MATLAB R2024a (Version 24.1) using MATLAB's \texttt{fmincon}. When small adjustments are made in either our classical or PINN approaches to improve results or overcome failures and errors, they are noted in their corresponding sections. 
If no comments are present, this means that the results were obtained using the following configurations for the classical and PINN methods, along with the configuration of the  \texttt{fmincon}, in tables ~\ref{tab: 8}, ~\ref{tab: 7}, and ~\ref{tab:fmincon_config} respectively.

\begin{table}[h]
\centering
\caption{ \texttt{fmincon} Configuration}
\label{tab:fmincon_config}
\begin{tabular}{|l|l|}
\hline
\textbf{Parameter} & \textbf{Value} \\
\hline
Optimization algorithm & Interior-point (MATLAB default) \\
\hline
Maximum iterations & 50 \\
\hline
Optimality tolerance & $10^{-6}$ \\
\hline
Lower bound & $m = 1$ \\
\hline
Upper bound & $m = 10$ \\
\hline
\end{tabular}
\end{table}

\begin{figure}[H]
    \centering
\includegraphics[width=1\linewidth]{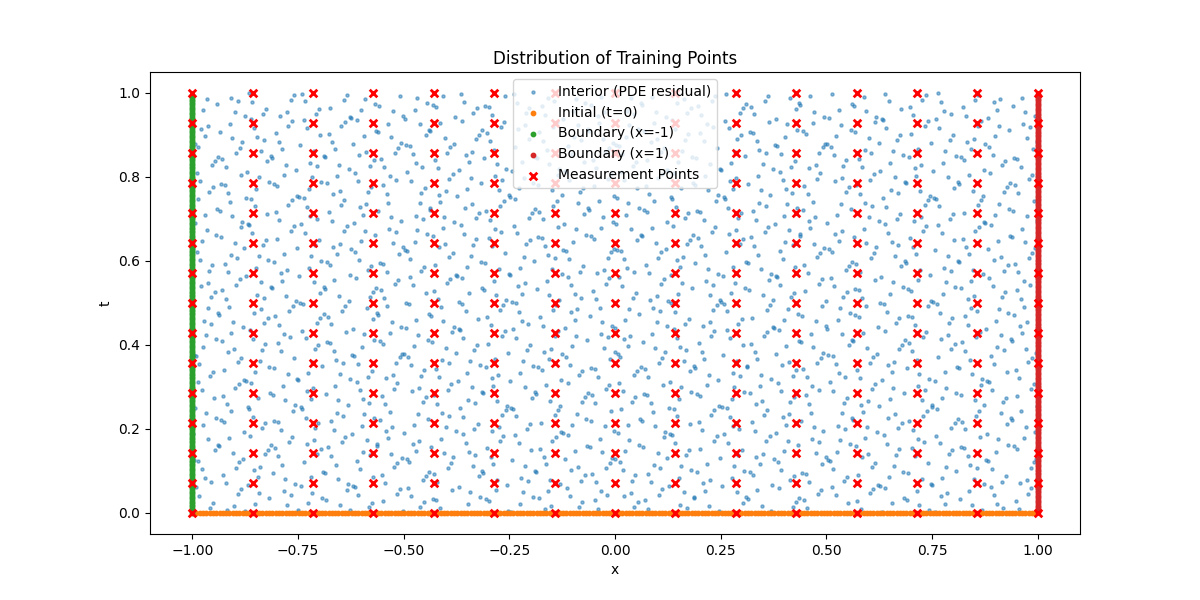}\vspace{-2mm}
\caption{Distribution of Training Points}
\label{fig:trainingpnt}
\end{figure}

\noindent To ensure a fair comparison between the classical inverse solver and the PINN approach, we used the same uniform distribution of measurements for the classical method (Figure \ref{fig:trainingpnt}). 
Upon testing, some classical objective functions were undefined for some initial guesses, others produced numerical instabilities and sigularities and this resulted in a failed run. Therefore, we wrapped the call inside a try–catch block so that in case something went wrong during the run, the run was tagged as 'failed' in the results file so we don't crash the whole script. This is translated into non-existent points in the following plots. On the contrary, the PINN does not implement the PDE step-by-step but learns a function and adjusts the parameter; therefore, this allows it to handle a very wide range of initial values with better stability and without crashing.
In addition to the tests using data generated from the exact reference solutions, the robustness of the inverse PINN is further investigated by testing it with noisy data, where the measurement data is generated from the direct PINN solution and then perturbed with different levels of noise. This analysis is presented in Section~\ref{sec: inverse without noise} and Section~\ref{sec: inverse with noise} respectively.
\begin{table}[htbp]
\centering
\caption{Configurations and Parameters for Inverse PINN Method}
\label{tab: 7}
\resizebox{\textwidth}{!}{%
\begin{tabular}{|l|l|l|}
\hline
\textbf{Category} & \textbf{Parameter} & \textbf{Value/Description} \\
\hline
Training Points ( see Fig ~\ref{fig:trainingpnt} )& n\_int (interior points) & 2048 \\
& n\_sb (spatial boundary points per side) & 256 \\
& n\_tb (temporal boundary points) &  256\\
& n\_meas (measurement points) & 15x15 = 225\\
\hline
Domain & Time domain & [0, 1] \\
& Space domain & [-1, 1] \\
\hline
Loss Weights & lambda\_u (boundary/initial loss weight) & 10 \\
& lambda\_m (measurement loss weight) & 100 \\
\hline
Neural Network (approximate\_solution) & Input dimension & 2 (t, x) \\
& Output dimension & 1 (u) \\
& Hidden layers & 4  \\
& Neurons per layer & 20 \\
& Activation function & nn.Tanh() \\
& Initialization & Xavier normal , bias=0 \\
& Retrain seed & 42 \\
\hline
Sampling& Generator & torch.quasirandom.SobolEngine(dimension=2) \\
\hline
Optimizer - ADAM & Learning rate (network) & 5e-4 \\
& Learning rate (m parameter) & 1e-3 \\
\hline
Training Config & num\_epochs & 10000 \\
& early\_stopping\_patience & 500 \\
& early\_stopping\_min\_delta & 1e-7 \\
\hline
\end{tabular}%
}\vspace{-6mm}
\end{table}

\subsubsection{Inverse PINN testing without Noisy Data}
\label{sec: inverse without noise}
In this section, the inverse PINN is first tested using data generated directly from the exact solution, without adding any additional perturbation or considering the direct PINN solution.

\begin{figure}[H]
\begin{subfigure}{0.5\textwidth}

\includegraphics[width= 1\linewidth, height=6cm]{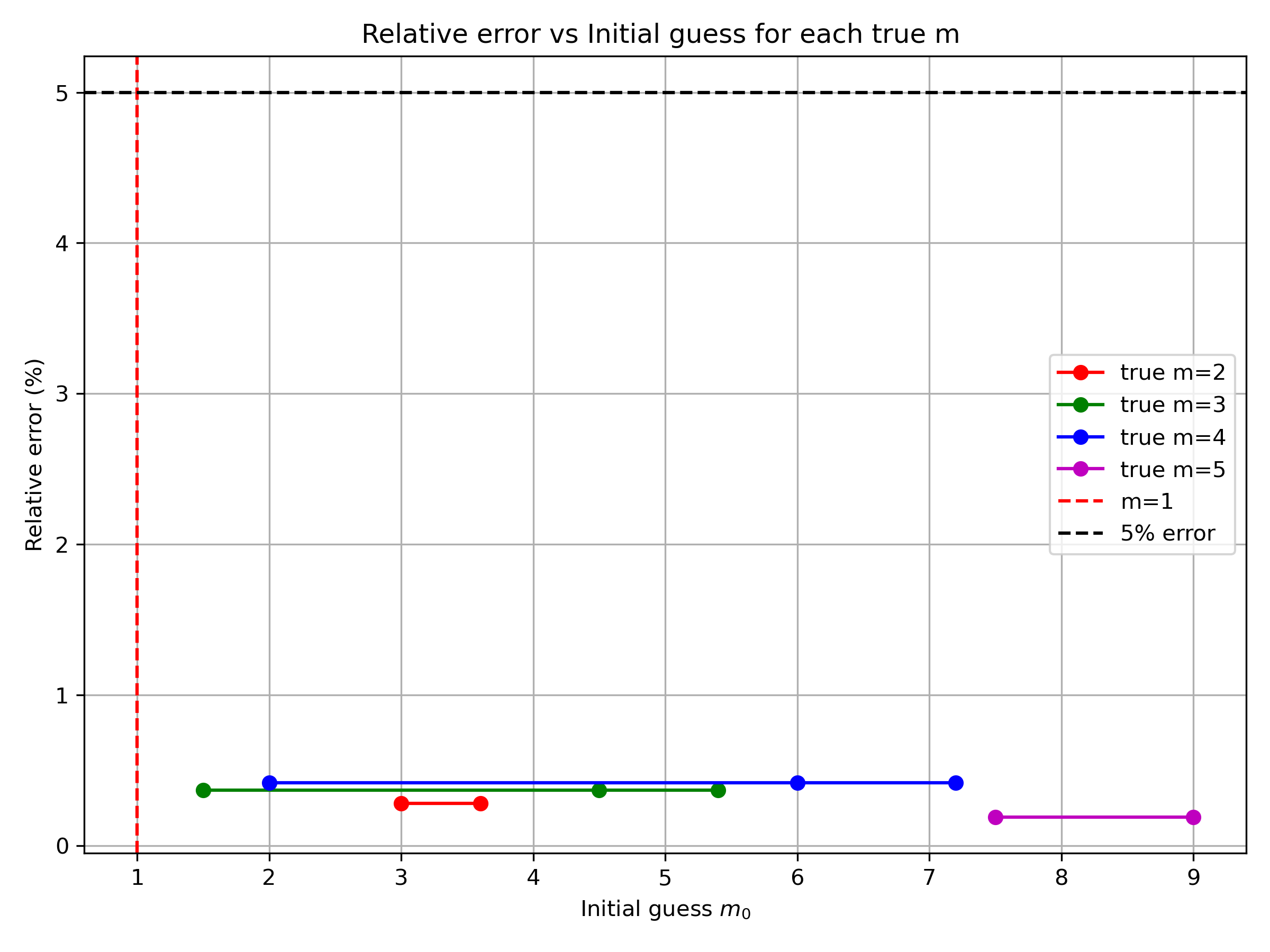} 
\caption{ \texttt{fmincon} results}
\label{fig:subim1}
\end{subfigure}
\begin{subfigure}{0.5\textwidth}
\includegraphics[width=1\linewidth, height=6cm]{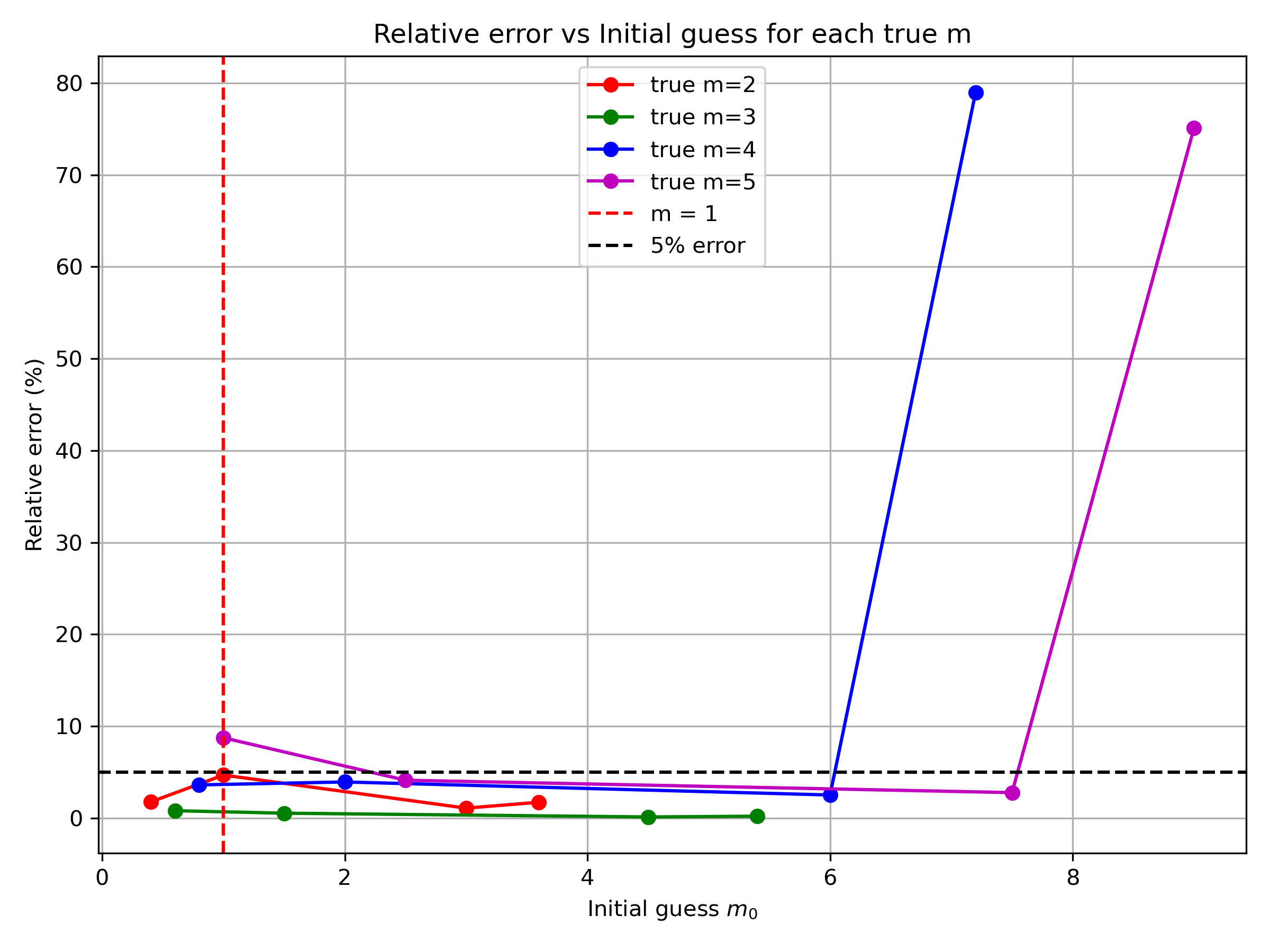}
\caption{ PINN results }
\label{fig:subim2}
\end{subfigure}

\caption{Inverse Barenblatt solution results plots}
\label{fig:image2}
\label{fig: inv bar}
\end{figure}

\paragraph*{Barenblatt Solution} The testing results for the Barenblatt solution are presented and visualized in the following Table~\ref{tab: inv bar} and Figure~\ref{fig: inv bar}:

\begin{table}[H]
\begin{tabular}{|l|l|l|l|l|l|}
\hline
true m & initial guess &  \texttt{fmincon} relative error & PINN relative error &  \texttt{fmincon}time & PINN training time \\ \hline

\multirow{4}{*}{2} 
 & 0.4           & failed                   & 1.74\%              & 2.30           & 228.70             \\ \cline{2-6}
 & 1             & failed                   & 4.67\%              & 0.70           & 238.43             \\ \cline{2-6}
 & 3             & 0.28\%                   & 1.07\%              & 14.50          & 229.50  \\           \cline{2-6}
 & 3.6           & 0.28\%                   & 1.71\%              & 19.88          & 124.91             \\ \hline

 \multirow{4}{*}{3}
& 0.6           & failed                   & 0.79\%              & 0.76           & 230.82             \\ \cline{2-6}
& 1.5           & 0.37\%                   & 0.52\%              & 12.52          & 238.57             \\ \cline{2-6}
& 4.5           & 0.37\%                   & 0.11\%              & 11.16          & 216.04             \\ \cline{2-6}
& 5.4           & 0.37\%                   & 0.19\%              & 10.21          & 212.74             \\ \hline

 \multirow{4}{*}{4}
 & 0.8           & failed                   & 3.60\%              & 0.72           & 215.14             \\ \cline{2-6}
 & 2             & 0.41\%                   & 3.92\%              & 8.41           & 216.10             \\ \cline{2-6}
 & 6             & 0.41\%                   & 2.50\%              & 8.16           & 206.24             \\ \cline{2-6}
 & 7.2           & 0.41\%                   & 79.00\%             & 7.00           & 58.49              \\ \hline

 \multirow{4}{*}{5} 
& 1             & failed                   & 8.73\%              & 0.68           & 162.13             \\ \cline{2-6}
& 2.5           & failed                   & 4.11\%              & 0.56           & 196.47             \\ \cline{2-6}
& 7.5           & 0.19\%                   & 2.76\%              & 6.92           & 58.37              \\ \cline{2-6}
 & 9             & 0.19\%                   & 75.11\%             & 7.59           & 78.01              \\ \hline
\end{tabular}
\caption{Relative error and computational time (seconds) comparison between the classical method and the PINN approach for the Barenblatt solution with different values of \(m\) and initial guesses.}
\label{tab: inv bar}
\end{table}

As explained before, many classical cases failed, while the PINN was more stable and produced results for all cases. However, when the classical method did run successfully, it mainly gave more accurate results (smaller relative error). Still, the errors were generally considered accurate, being below 5\% for all successful classical and PINN cases, with a few exceptions for the PINN.
 
\noindent After seeing these results, we noticed that the PINN relative error increased with \(m\) , and we wanted to understand why the errors for the last cases of \(m = 4\) and \(m = 5\), with their corresponding largest initial guesses, gave such off results. So we repeated the runs, but this time we plotted the evolution of both the training loss and the approximation of \(m\) versus the epochs.

\noindent We noticed that the loss in all cases of \(m\) followed a similar trend, where we have a sudden decrease, then a plateau, and then a new decrease, as shown in the example case of \(m = 3\) with the initial guess 4.5.

\begin{figure}[H]
\centering
\includegraphics[width=1\linewidth]{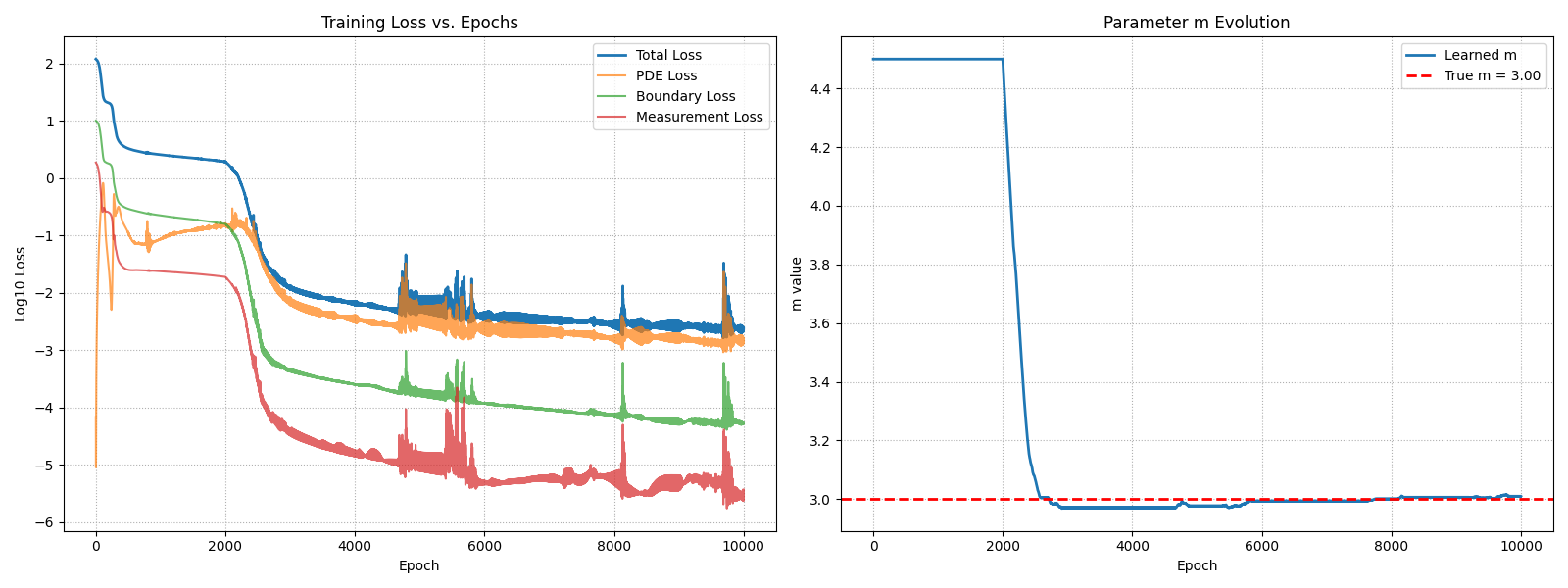}
\caption{Training loss and m evolution for case m=3 with initial guess =4.5}
\end{figure}

\noindent However, we noticed that in the case of m = 4 with an initial guess of 7.5, the training stopped in the plateau region at 4000 epochs due to the early stopping function. We expected it to follow the same loss trend as the other cases. Therefore, when we removed early stopping, the error in this case decreased from 79\% to 1.38\%.

\begin{figure}[H]
\begin{subfigure}{0.5\textwidth}
\includegraphics[width=1\linewidth, height=6cm]{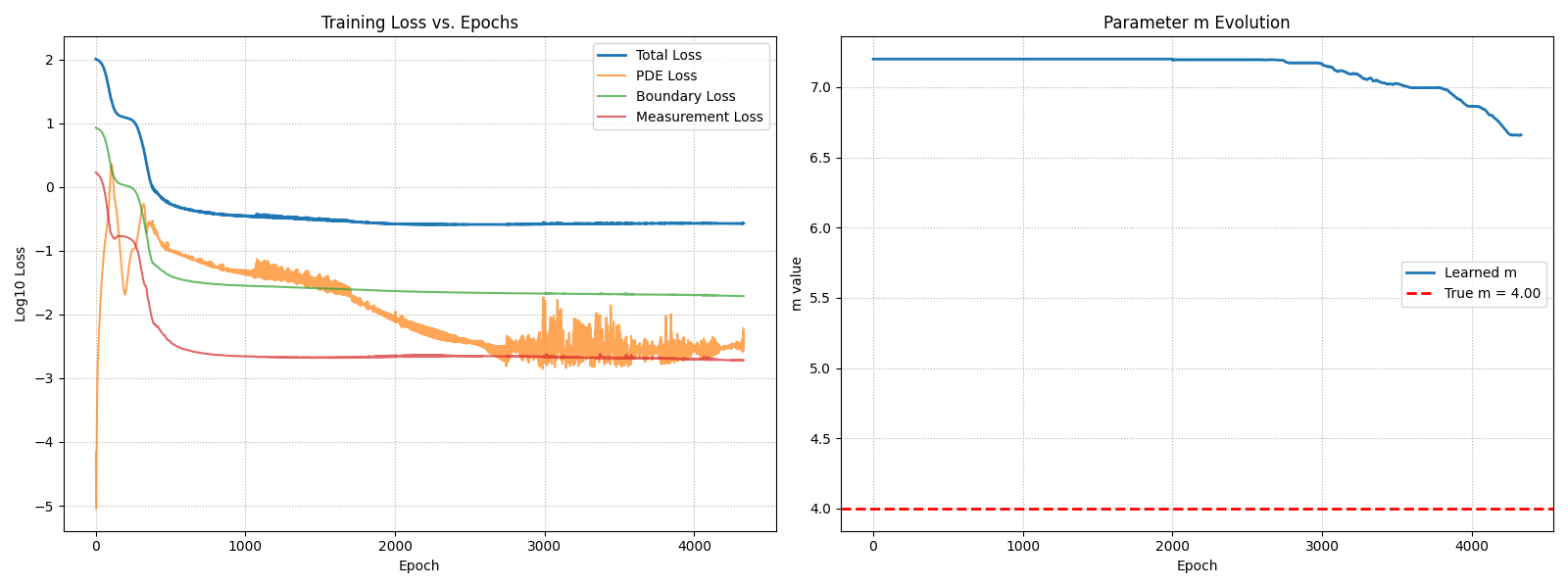}
\caption{with early stopping}
\label{fig:subim1}
\end{subfigure}
\begin{subfigure}{0.5\textwidth}
\includegraphics[width=1\linewidth, height=6cm]{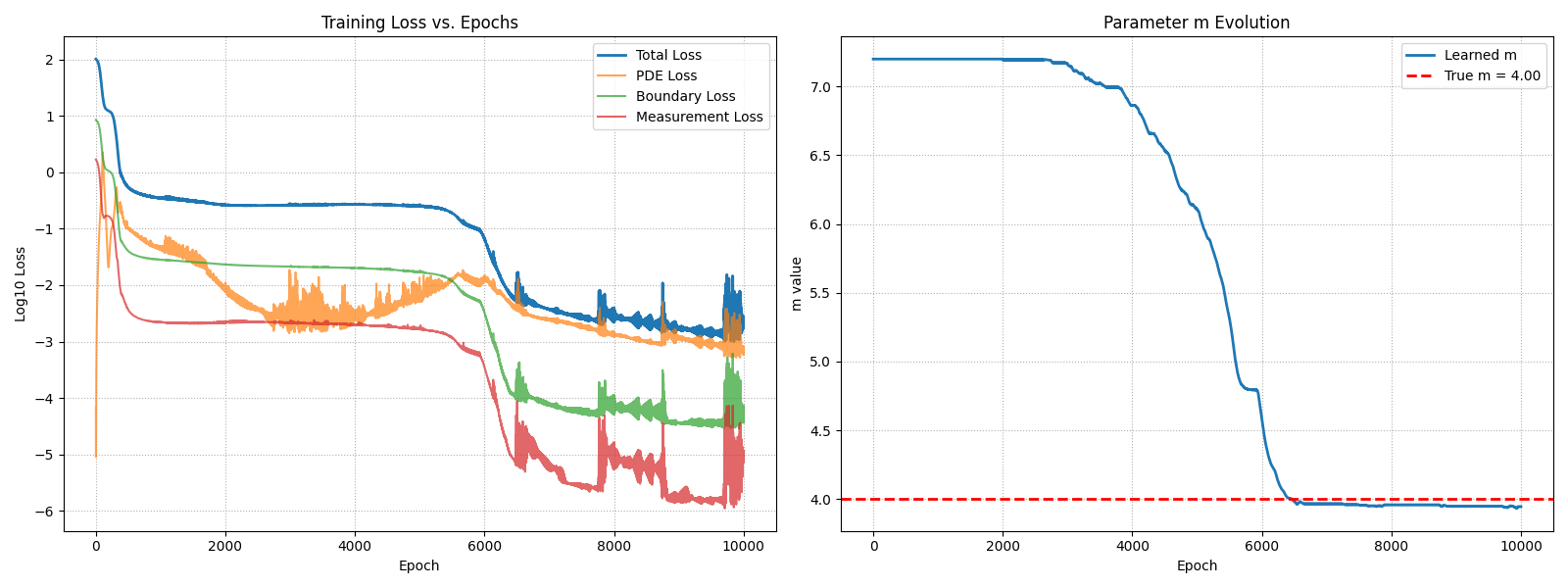}
\caption{without early stopping}
\label{fig:subim2}
\end{subfigure}
\caption{Training loss and m evolution for case m=4 with initial guess =7.5, before and after removing the early stopping function}
\end{figure}

\paragraph*{m-independent solutions}
A first try to run the classical solver ended up with numerical errors and singularities and failed for all cases. Some updated values of \( u \) may become zero or negative, which causes instability when computing terms such as \(u^{m-1}\), this is because the bounds  \texttt{fmincon} ensures that the parameter m stays between 1.1 and 10, but it does not guarantee that the numerical solution \(u\) remains positive during Newton iterations. Therefore, to avoid this issue, we imposed a small lower bound  using \(u=\max(u,10^{-12})\), which keeps the solution positive and improves stability. Everything worked perfectly after this change as shown in Tables \ref{tab:inv_poly} and \ref{tab:inv_har}.

\begin{table}[H]
\begin{tabular}{|l|l|l|l|l|l|}
\hline
true m & initial guess &  \texttt{fmincon} relative error & PINN relative error &  \texttt{fmincon} time & PINN training time \\ \hline

\multirow{4}{*}{2}
& 0.4           & 1.10\%                   & 0.51\%              & 30.05          & 360.64             \\ \cline{2-6}
 & 1             & 1.10\%                   & 0.28\%              & 29.46          & 356.92             \\ \cline{2-6}
 & 3             & 1.10\%                   & 0.15\%              & 16.38          & 298.42             \\ \cline{2-6}
& 3.6           & 1.10\%                   & 0.22\%              & 34.53          & 354.22             \\ \hline

\multirow{4}{*}{3}
& 0.6           & 0.07\%                   & 0.14\%              & 33.93          & 360.07             \\ \cline{2-6}
& 1.5           & 0.07\%                   & 0.20\%              & 24.96          & 362.01             \\ \cline{2-6}
& 4.5           & 0.07\%                   & 0.13\%              & 41.15          & 363.28             \\ \cline{2-6}
 & 5.4           & 0.07\%                   & 0.11\%              & 52.86          & 363.28             \\ \hline

 \multirow{4}{*}{4}
 & 0.8           & 0.22\%                   & 0.54\%              & 25.95          & 358.70             \\ \cline{2-6}
 & 2             & 0.22\%                   & 0.25\%              & 38.87          & 349.53             \\ \cline{2-6}
& 6             & 0.22\%                   & 0.45\%              & 39.00          & 360.74             \\ \cline{2-6}
& 7.2           & 0.22\%                   & 0.26\%              & 35.97          & 359.90             \\ \hline

\multirow{4}{*}{5}
 & 1             & 0.14\%                   & 76.74\%             & 35.50          & 143.21             \\ \cline{2-6}
& 2.5           & 0.14\%                   & 0.11\%              & 63.15          & 298.59             \\ \cline{2-6}
 & 7.5           & 0.14\%                   & 0.45\%              & 35.70          & 356.80             \\ \cline{2-6}
& 9             & 0.14\%                   & 0.31\%              & 46.35          & 356.47             \\ \hline
\end{tabular}
\caption{Relative error and computational time (seconds) comparison between the classical method and the PINN approach for the polynomial solution with different values of \(m\) and initial guesses.}\label{tab:inv_poly}
\end{table}

\begin{figure}[H]
\begin{subfigure}{0.5\textwidth}
\includegraphics[width= 1\linewidth, height=6cm]{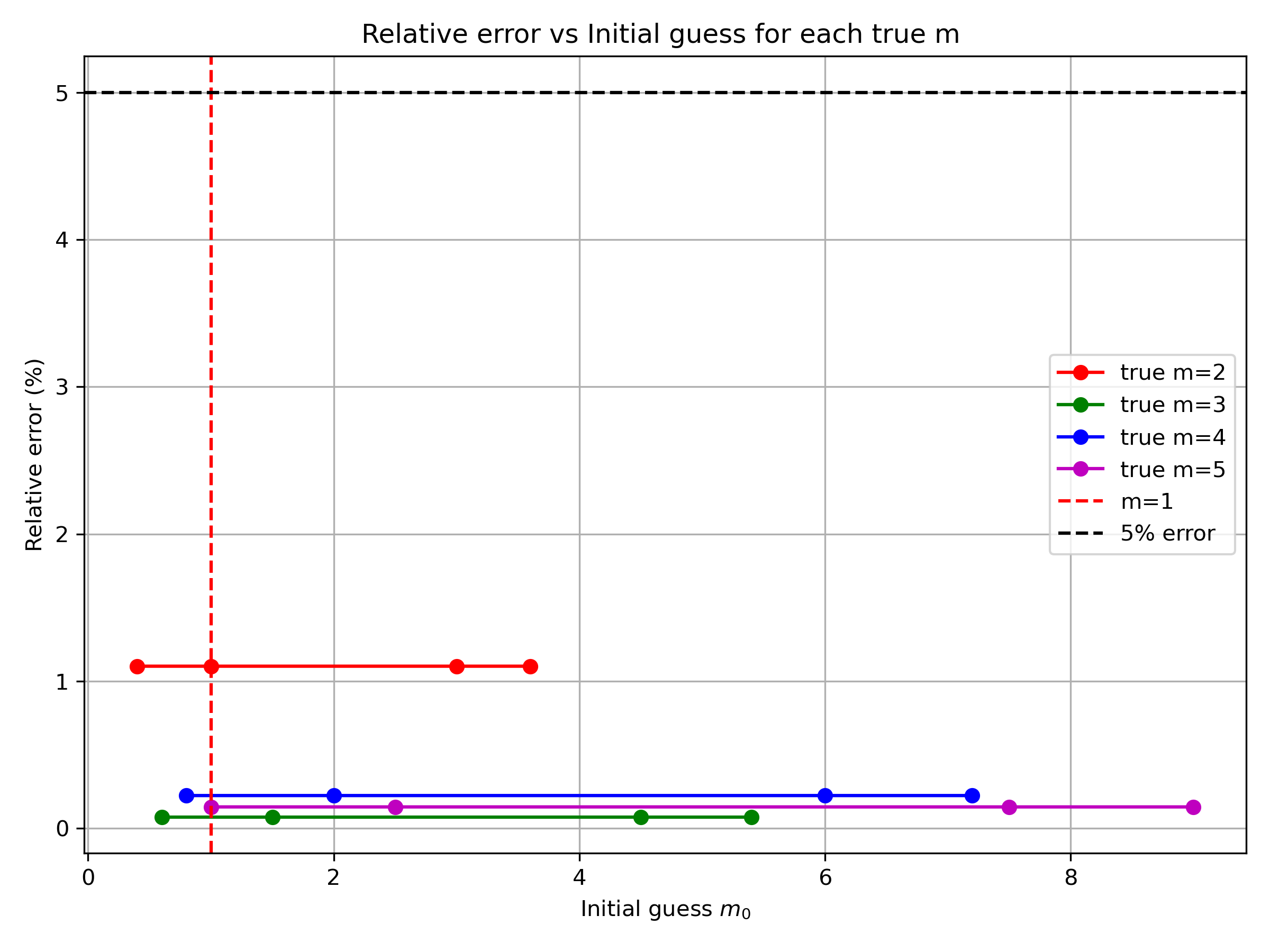} 
\caption{ \texttt{fmincon} results}
\label{fig:subim1}
\end{subfigure}
\begin{subfigure}{0.5\textwidth}
\includegraphics[width=1\linewidth, height=6cm]{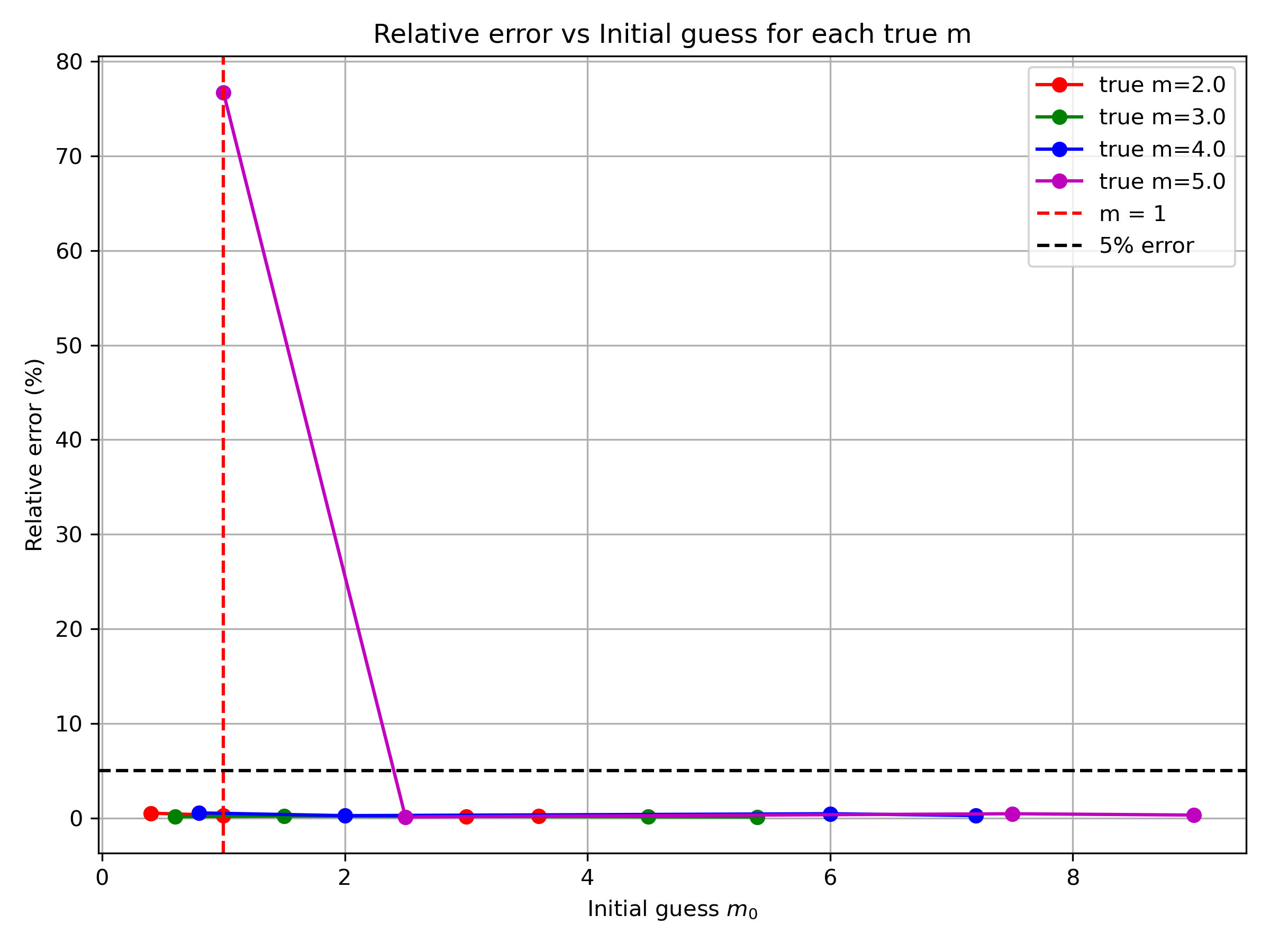}
\caption{ PINN results }
\label{fig:subim2}
\end{subfigure}

\caption{Inverse Polynomial solution plots}
\label{fig:image2}
\end{figure}

\begin{table}[H]
\begin{tabular}{|l|l|l|l|l|l|}
\hline
true m & initial guess &  \texttt{fmincon} relative error & PINN relative error &  \texttt{fmincon} time & PINN training time \\ \hline

\multirow{4}{*}{2} 
& 0.4           & 0.85\%                    & 0.94\%              & 5.38               & 338.86                 \\ \cline{2-6}
  & 1             & 0.85\%                    & 1.37\%              & 7.32               & 338.78                 \\ \cline{2-6}
 & 3             & 0.85\%                    & 0.15\%              & 12.32              & 250.15                 \\ \cline{2-6}
  & 3.6           & 0.85\%                    & 0.11\%              & 14.43              & 341.01                 \\ \hline

\multirow{4}{*}{3} 
 & 0.6           & 0.70\%                    & 0.31\%              & 9.31               & 352.13                 \\ \cline{2-6}
 & 1.5           & 0.70\%                    & 0.33\%              & 7.77               & 349.41                 \\ \cline{2-6}
 & 4.5           & 0.70\%                    & 0.12\%              & 7.46               & 349.22                 \\ \cline{2-6}
 & 5.4           & 0.70\%                    & 0.12\%              & 9.19               & 309.47                 \\ \hline

\multirow{4}{*}{4} 
& 0.8           & 0.71\%                    & 0.20\%              & 7.59               & 343.31                 \\ \cline{2-6}
 & 2             & 0.71\%                    & 0.21\%              & 7.17               & 342.84                 \\ \cline{2-6}
 & 6             & 0.71\%                    & 0.07\%              & 7.63               & 284.96                 \\ \cline{2-6}
 & 7.2           & 0.71\%                    & 0.01\%              & 7.61               & 343.86                 \\ \hline

\multirow{4}{*}{5} 
& 1             & 0.80\%                    & 0.37\%              & 8.58               & 344.33                 \\ \cline{2-6}
 & 2.5           & 0.80\%                    & 0.12\%              & 7.57               & 345.70                 \\ \cline{2-6}
& 7.5           & 0.80\%                    & 0.08\%              & 7.58               & 343.78                 \\ \cline{2-6}
 & 9             & failed                    & 0.17\%              & 1.57               & 344.36                 \\ \hline
\end{tabular}
\caption{Relative error and computational time (seconds) comparison between the classical method and the inverse PINN approach for the harmonic solution with different values of \(m\) and initial guesses.}\label{tab:inv_har}
\end{table}

\begin{figure}[H]
\begin{subfigure}{0.5\textwidth}
\includegraphics[width= 1\linewidth, height=6cm]{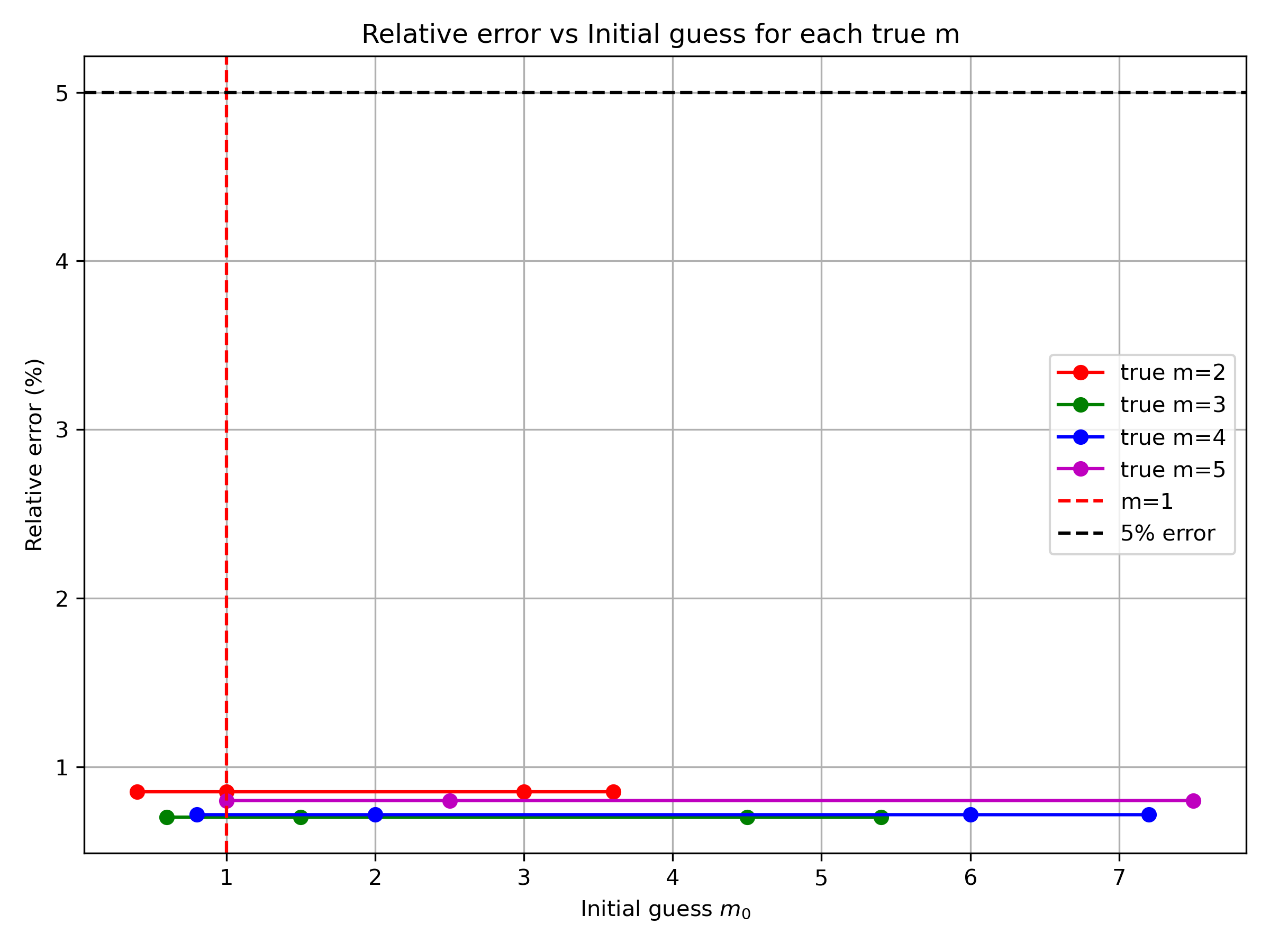} 
\caption{ \texttt{fmincon} results}
\label{fig:subim1}
\end{subfigure}
\begin{subfigure}{0.5\textwidth}
\includegraphics[width=1\linewidth, height=6cm]{PINN_inverse_plots/HPINN.png}
\caption{ PINN results }
\label{fig:subim2}
\end{subfigure}

\caption{Inverse Harmonic solution plots}
\label{fig:image2}
\end{figure}

\paragraph*{m-dependent solution}
The positivity of \(u\) in the classical solver was also forced in this case. 

\noindent As shown in \ref{tab: inverse mdep}, the PINN model clearly performs differently depending on the nonlinearity parameter \(m\). For \(m = 2\), almost all initial guesses give relative errors approximately within the 5\% threshold, except for the case where the initial guess was 50\% above the true value. As m increases, the errors grow significantly, reaching between 20\% and 80\% across the tested initial guesses. For \(m = 5\) across all initial guesses, and for \(m = 3\) and 4 with initial guesses that are $+$50\% and $+$80\% away from the true value, the relative errors are equal to the initial guess offset itself. This means the PINN did not learn anything and stayed stuck at its initial guess. This confirms that the model struggles as nonlinearity increases, learning best at \(m = 2\) and worst at \(m = 5\).

\begin{table}[H]
\begin{tabular}{|l|l|l|l|l|l|}
\hline
true m & initial guess &  \texttt{fmincon} relative error & PINN relative error &  \texttt{fmincon} time & PINN training time \\ \hline

\multirow{4}{*}{2}
 & 0.4            & 0.08\%                     & 6.79\%                & 107.29             & 379.69                  \\ \cline{2-6}
& 1              & 0.08\%                     & 2.83\%                & 110.53             & 382.94                  \\ \cline{2-6}
 & 3              & 0.08\%                     & 49.85\%               & 71.56              & 266.39                  \\ \cline{2-6}
& 3.6            & 0.08\%                     & 0.15\%                & 55.11              & 366.76                  \\ \hline

\multirow{4}{*}{3}
 & 0.6            & failed                     & 9.49\%                & 26.84              & 370.63                  \\ \cline{2-6}
 & 1.5            & 0.27\%                     & 29.66\%               & 1248.56            & 376.23                  \\ \cline{2-6}
 & 4.5            & failed                     & 51.31\%               & 0.69               & 369.48                  \\ \cline{2-6}
 & 5.4            & failed                     & 81.89\%               & 12.34              & 369.26                  \\ \hline

\multirow{4}{*}{4}
 & 0.8            & failed                     & 22.32\%               & 39.68              & 377.16                  \\ \cline{2-6}
 & 2              & failed                     & 19.84\%               & 27.14              & 370.12                  \\ \cline{2-6}
 & 6              & failed                     & 50.46\%               & 0.28               & 364.21                  \\ \cline{2-6}
 & 7.2            & failed                     & 79.71\%               & 0.29               & 363.13                  \\ \hline
\multirow{4}{*}{5}
 & 1              & failed                     & 80.00\%               & 9.76               & 179.46                  \\ \cline{2-6}
& 2.5            & failed                     & 50.00\%               & 0.95               & 179.96                  \\ \cline{2-6}
& 7.5            & failed                     & 50.00\%               & 4.73               & 177.85                  \\ \cline{2-6}
 & 9              & failed                     & 80.00\%               & 0.11               & 178.40                  \\ \hline
\end{tabular}
\caption{Relative error and computational time (seconds) comparison between the classical method and the PINN approach for the m\_dependent solution with different values of \(m\) and initial guesses.}
\label{tab: inverse mdep}
\end{table}

\begin{figure}[H]
\begin{subfigure}{0.5\textwidth}
\includegraphics[width= 1\linewidth, height=6cm]{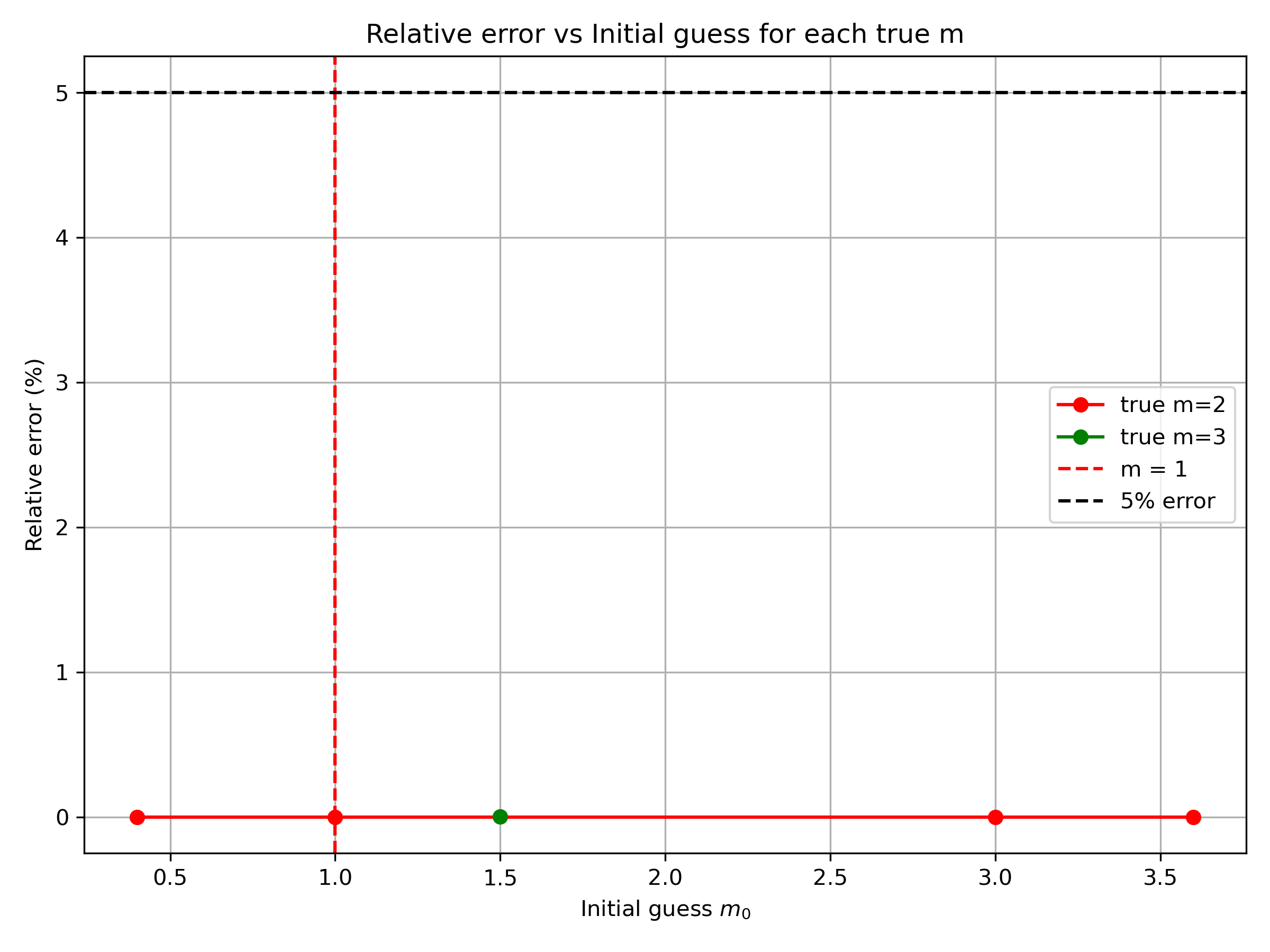} 
\caption{ \texttt{fmincon} results}
\label{fig:subim1}
\end{subfigure}
\begin{subfigure}{0.5\textwidth}
\includegraphics[width=1\linewidth, height=6cm]{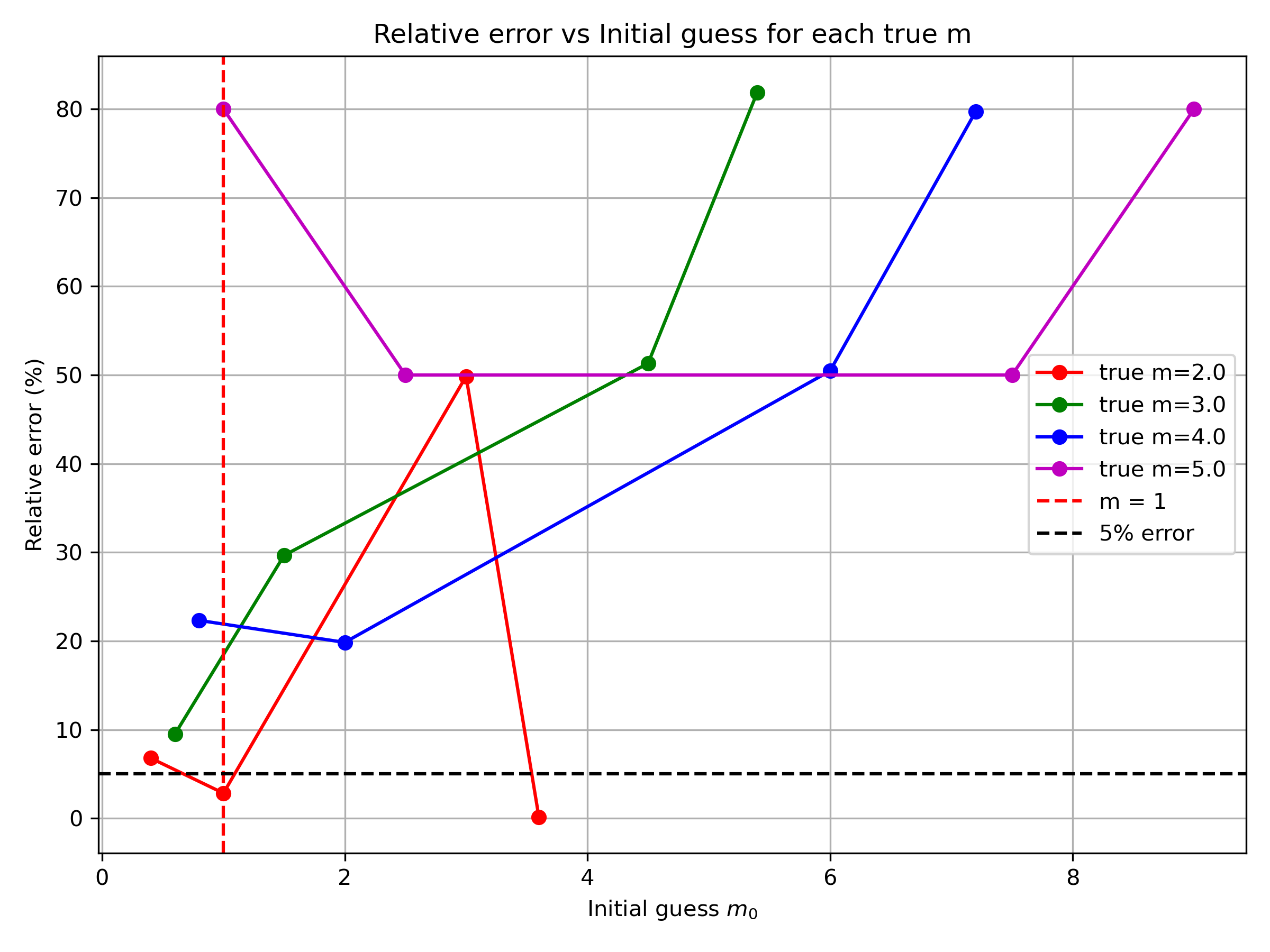}
\caption{ PINN results }
\label{fig:subim2}
\end{subfigure}

\caption{Inverse m-dependent solution plots}
\label{fig:image2}
\end{figure}

\begin{figure}[H]
\begin{subfigure}{0.5\textwidth}
\includegraphics[width=1\linewidth, height=6cm]{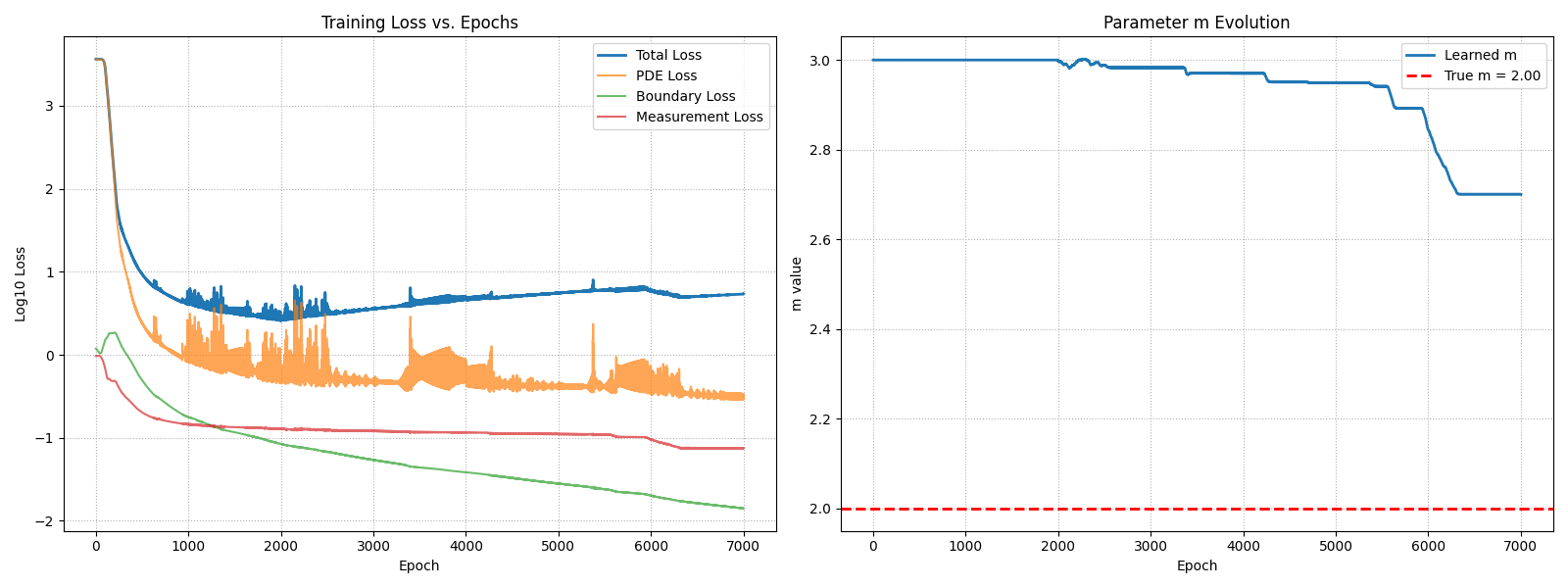}
\caption{with early stopping}
\label{fig:subim1}
\end{subfigure}
\begin{subfigure}{0.5\textwidth}
\includegraphics[width=1\linewidth, height=6cm]{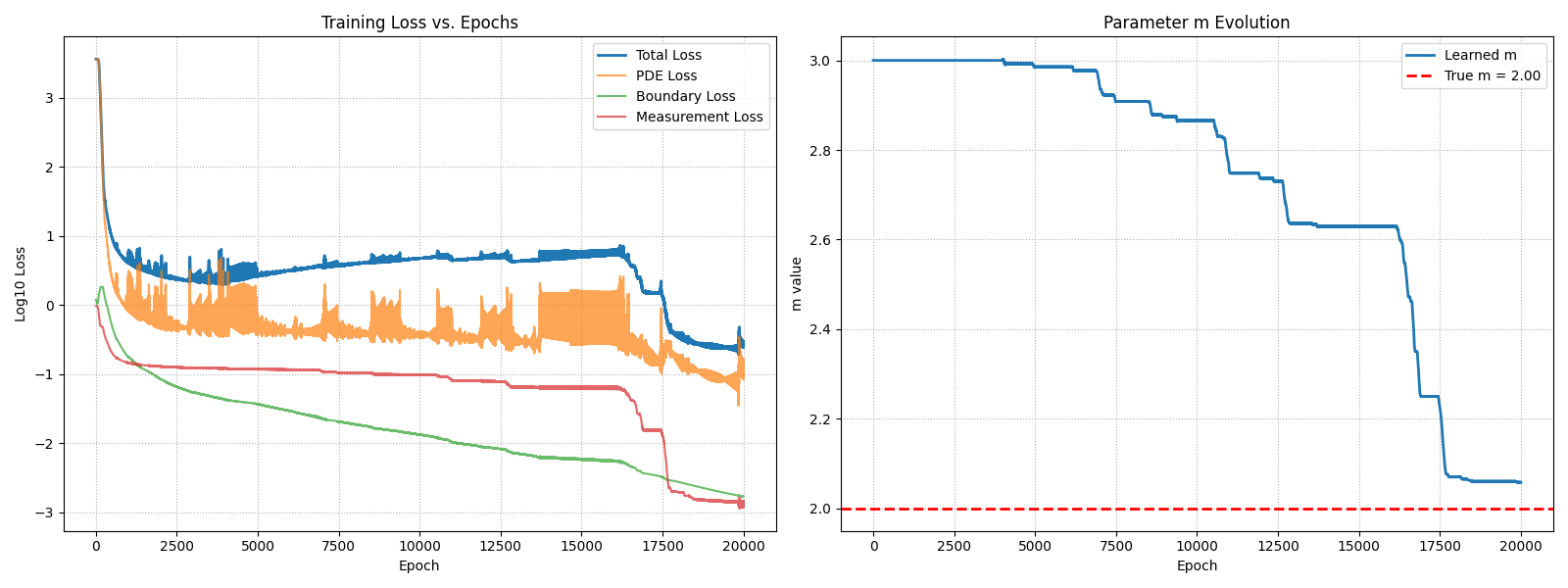}
\caption{without early stopping}
\label{fig:1}
\end{subfigure}
\caption{Training loss and m evolution for case m=2 with initial guess =3, after increasing the epochs and removing the early stopping, for the m-dependent solution}
\label{fig:m2 mdep}
\end{figure}

Looking at the training loss evolution shown in Appendix ~\ref{appendix B}, it can be seen that in many cases either early stopping was triggered or the epochs ran out right when the loss was starting to go down. To check this, the number of epochs was increased from 10,000 to 20,000 for the case of \(m = 2\) with an initial guess 50\% above the true value, and the PINN kept learning, bringing the error down to 2.9\% as shown in figure ~\ref{fig:m2 mdep}.

The loss trends for \(m = 3\) and 4 with initial guesses of -50\% and -80\% suggest that more epochs could give better results in those cases as well. However, for cases where the loss does not improve at all and stays flat, other changes are needed, such as rescaling the problem, using the L-BFGS optimizer, or trying other improved training techniques.
On the other hand, the classical method was able to accurately estimate $m$ only for \(m = 2\) with all initial guesses and \(m = 3\) with an initial guess of 1.5 (-50\%), while failing to produce results for all other cases due to numerical instabilities and matrix singularities.
\\
\noindent Overall, for the inverse problem, the PINN model proposed in this work achieved results comparable to the classical method in terms of accuracy across all tested cases. Even in the m-dependent solution case, which was the most challenging for both approaches, the PINN was still able to produce results where the classical method completely failed due to numerical instabilities. While not all PINN results were accurate, the model still has room for improvement, as one can explore different optimization strategies, architectural changes, and rescaling tricks to push the results further. The classical method, on the other hand, offers no such flexibility.

\subsubsection*{Time comparison }
\begin{figure}[H]
    \centering
\includegraphics[width=0.9\textwidth, height = 5cm]{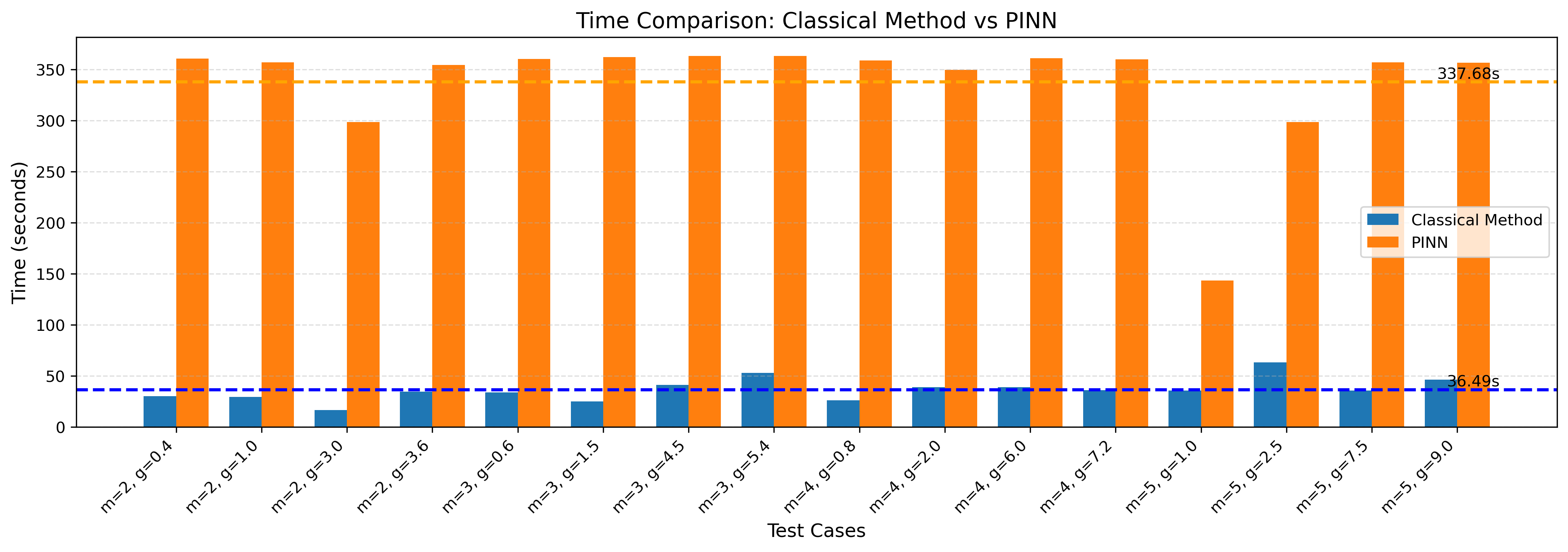}\vspace{-4mm}
     \caption*{Polynomial Solution}
    \vspace{0.3cm}
\includegraphics[width=0.9\textwidth, height = 5cm]{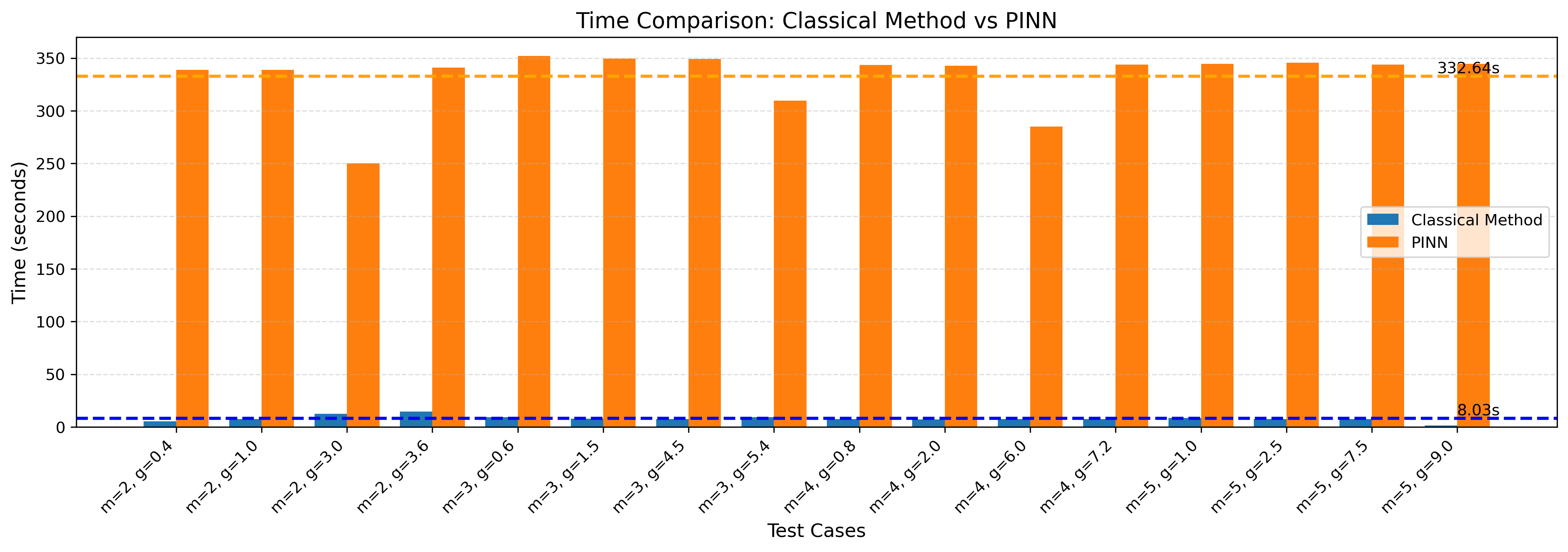}\vspace{-4mm}
     \caption*{Harmonic Solution}  \vspace{0.3cm}  
\includegraphics[width=0.9\textwidth, height = 5cm]{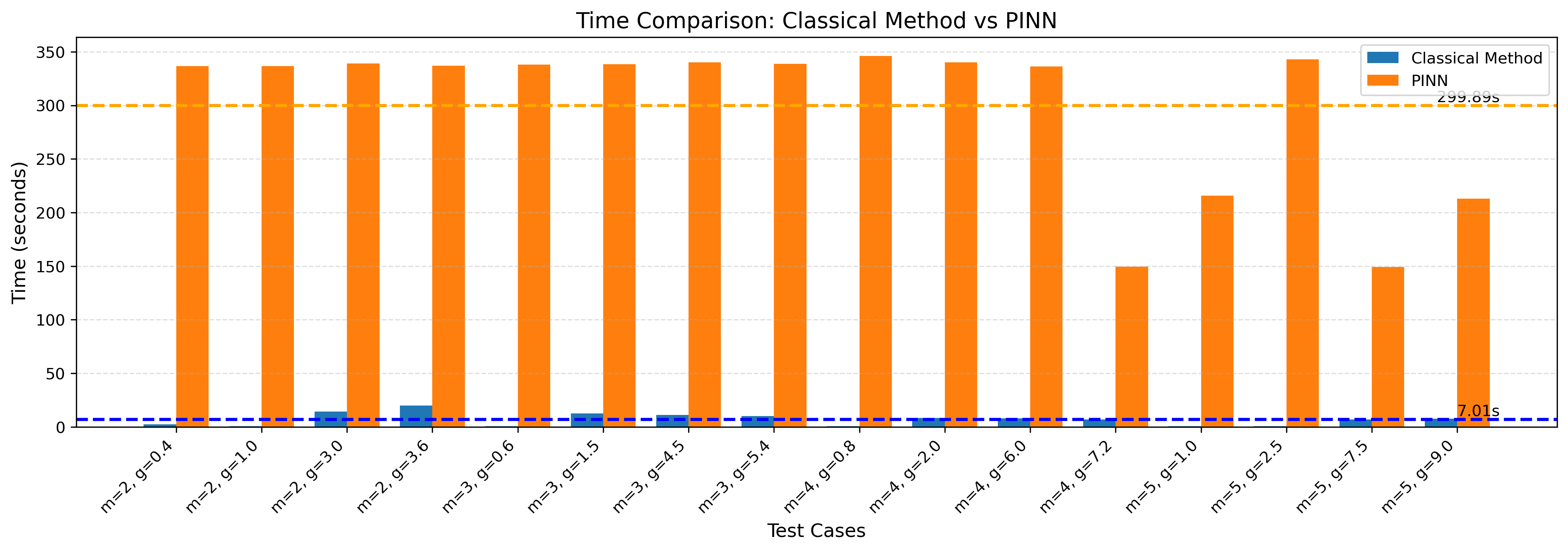}\vspace{-4mm}
    \caption*{Barenblatt Solution}    \vspace{0.3cm}
\includegraphics[width=0.9\textwidth, height = 5cm]{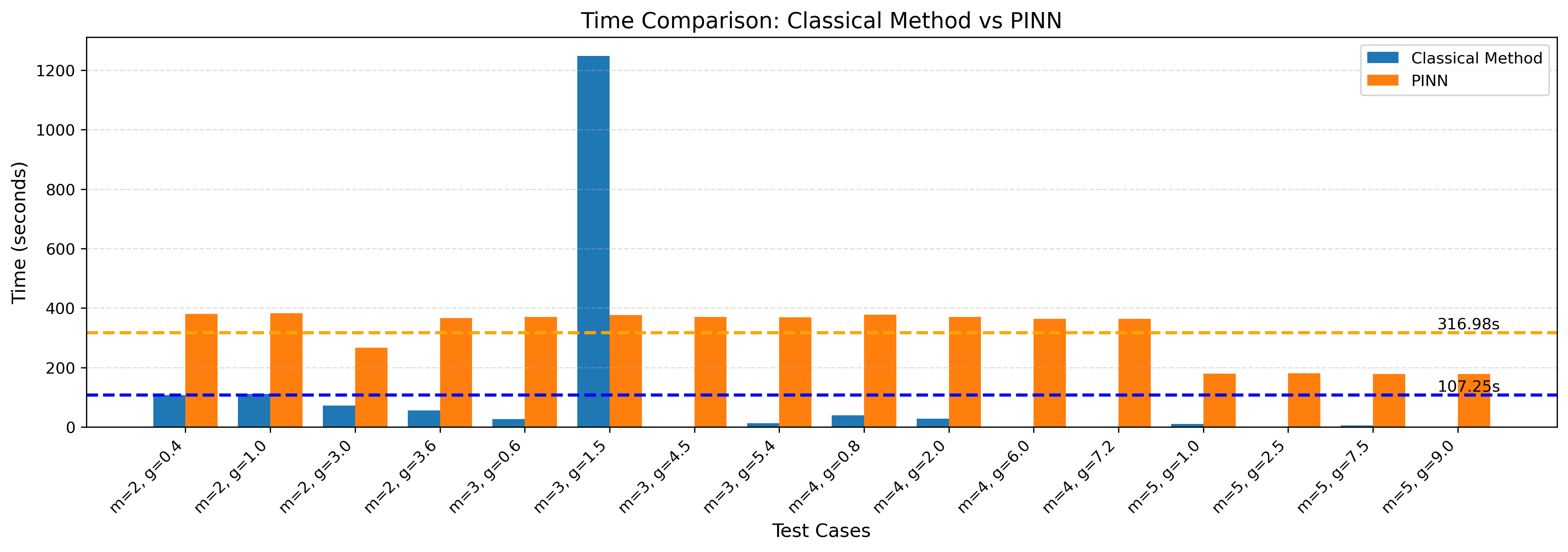}\vspace{-4mm}
         \caption*{m-depenednt}
\caption{Computational time comparison of PINN and the classical method across m-dependent and m-independent solutions. Dashed lines show average times.}
    \label{fig: time mind}
\end{figure}

In this section, the computational times of the PINN approach and the classical method are compared for different test cases. The following plots (Figure ~\ref{fig: time mind}) illustrate the execution time of each method and highlight their average performance.\vspace{2mm}

\noindent Overall, the \texttt{fmincon} method is consistently faster than the PINN approach in all cases. Its average computational time ranges from about 7 s to 36 s depending on the solution, except for the m-dependent solution where it rises to around 107 s, while the PINN method requires roughly 300 s to 340 s on average.

At first glance, the substantial difference in computational time between the two approaches was unexpected. Indeed, the classical inverse solver relies on \texttt{fmincon}, which repeatedly solves the direct problem during the optimization process. Based on the direct problem experiments presented in Section~\ref{sec: test direct}, solving a single direct problem using the classical Newton method required between approximately \(4\) s and \(24\) s in the Python implementation, depending on the considered test case. For instance, for the harmonic solution, a single direct solve required on average about \(6\) s.

Since the inverse optimization typically converges in approximately eight iterations and each iteration requires at least one direct solve, one would expect the overall computational time of the classical inverse method to be at least \(48\) s for the harmonic solution. However, the measured execution time in MATLAB was only about \(7\) s, indicating a significant discrepancy.

To investigate this observation, the direct classical solver was implemented and executed independently in MATLAB using the same numerical method and discretization parameters as in the Python implementation. The comparison revealed that MATLAB executes the direct solver considerably faster. For example, solving the harmonic direct problem for all tested values of \(m\) required \(26.34\) s in Python, whereas the corresponding MATLAB implementation required only \(1.21\) s. This comparison indicates that the large difference observed in the inverse problem timings is primarily due to the execution speed of the two programming environments rather than the optimization procedure itself. In particular, the Python implementation contains several nested loops that are executed much less efficiently than the corresponding MATLAB implementation, leading to a significantly higher computational cost despite both implementations employing the same numerical algorithm. Table ~\ref{tab:classical_times} summarizes the computation times (in seconds) of the same classical method presented in ~\ref{sec: classical} for the four manufactured solutions using MATLAB instead of Python for the direct problem. 
\begin{table}[ht]
\centering
\caption{Computation time (seconds) for solving the direct problem using MATLAB for different manufactured solutions.}\vspace{-2mm}
\label{tab:classical_times}
\renewcommand{\arraystretch}{1.2}
\begin{tabular}{|c|c|c|c|c|}
\hline
$m$ & Barenblatt & Polynomial & Harmonic & m-dependent \\
\hline
2 & 0.078 & 0.087 & 0.182 & 0.134 \\
3 & 0.040 & 0.053 & 0.033 & 0.054 \\
4 & 0.388 & 0.529 & 0.496 & 0.938 \\
5 & 0.378 & 0.524 & 0.491 & 1.133 \\
\hline
\end{tabular}
\end{table}
\subsubsection{Inverse PINN testing with Noisy Data}
\label{sec: inverse with noise}
Another aspect that should be investigated is the performance of the inverse PINN when the available solution data is not obtained from the exact solution, but rather from a perturbed dataset. In this section, the inverse PINN is tested under three different data scenarios. First, the measurement dataset is generated from the direct PINN solution, where the approximation error of the direct PINN introduces a deviation from the exact solution. Second, an additional \(1\%\) noise is added to the direct PINN-generated dataset, and finally, the same test is performed with \(5\%\) noise. These experiments aim to evaluate the robustness of the inverse PINN when the available data becomes increasingly different from the exact solution. As in Section~\ref{direct noise}, this testing is performed only for the polynomial solution due to the range of relative errors obtained for different values of \(m\), which provides a suitable benchmark for testing the effect of data perturbations.

The noise is added to the measurement dataset as follows. Let \(u^{\text{PINN}}_i\) denote the value of the direct PINN solution at the measurement point \((t_i, x_i)\), for \(i = 1, \dots, N_{\text{meas}}\), where \(N_{\text{meas}} = 225\) corresponds to the \(15 \times 15\) grid described in Section~\ref{direct noise}. The perturbed dataset is obtained by adding an independent Gaussian perturbation to each measurement,\vspace{-1mm}
\begin{equation}
    u^{\text{obs}}_i = u^{\text{PINN}}_i + \varepsilon_i, \qquad
    \varepsilon_i \sim \mathcal{N}(0, \sigma^2), \qquad i = 1, \dots, N_{\text{meas}},\vspace{-1mm}
\end{equation}
where the standard deviation \(\sigma\) is chosen relative to the magnitude of the data itself,
\begin{equation}
    \sigma = \delta \sqrt{\frac{1}{N_{\text{meas}}} \sum_{i=1}^{N_{\text{meas}}} \left( u^{\text{PINN}}_i \right)^2 },
\end{equation}
and \(\delta\) is the prescribed noise level, taken here as \(\delta = 0.01\) and \(\delta = 0.05\), corresponding to \(1\%\) and \(5\%\) noise respectively.

Two remarks are in order regarding this choice. First, the perturbation is additive rather than multiplicative. Since the polynomial solution \(u(t,x) = (1-t)^2(1-x^2)^2\) vanishes at \(x = \pm 1\) and at \(t = 1\), a perturbation proportional to the local value of the solution would leave these measurements unperturbed, which does not reflect the behaviour of a physical measurement device. Scaling \(\sigma\) by the root mean square of the data instead assigns the same absolute perturbation to every measurement point, independently of the local magnitude of the solution. Second, the noise is applied only to the measurement dataset; the initial and boundary conditions are kept exact, as these are assumed to be known.

The random perturbations are generated with a fixed seed, so that the same realisation of the noise is used across all runs and the results reported below are reproducible. The remaining components of the inverse problem, namely the network architecture, the optimizer, the loss weights and the number of training epochs, are left unchanged from Section~\ref{sec: standard}, so that any variation in the recovered value of \(m\) can be attributed to the perturbation of the data alone.The parameters used for the data perturbation are summarized in Table~\ref{tab:noise_config}.

\begin{table}[H]
\centering
\caption{Parameters of the data perturbation used in the inverse PINN experiments.}
\label{tab:noise_config}
\begin{tabular}{|l|l|}
\hline
\textbf{Parameter} & \textbf{Value} \\
\hline
Noise model                        & Additive Gaussian, \(\varepsilon_i \sim \mathcal{N}(0,\sigma^2)\) \\
Noise levels \(\delta\)            & \(0\%,\ 1\%,\ 5\%\) \\
Standard deviation \(\sigma\)      & \(\delta \cdot \text{RMS}\left(u^{\text{PINN}}\right)\) \\
Perturbed quantity                 & Measurement dataset only (from direct PINN) \\
Number of perturbed points         & \(N_{\text{meas}} = 225\) \\
Random seed                        & 42 \\ \hline
\end{tabular}
\end{table}
\noindent The results obtained for the three data scenarios are reported in Tables ~\ref{tab:zero per}, \ref{tab:one per}, and \ref{tab:five pct}
\begin{table}[H]
\centering
\caption{PINN inverse results for the polynomial solution using measurement data generated from the direct PINN with 0\% additional noise}\label{tab:zero per}
\begin{tabular}{|l|l|l|l|l|l|}
\hline
true m & initial guess     & recovered m       & PINN relative error & PINN training time s\\ \hline
\multirow{4}{*}{2}
 & 0.4                & 2.00 & 0.22\%  & 338.66 \\ \cline{2-5}
& 1.0                & 2.00  & 0.01\%  & 285.76 \\ \cline{2-5}
& 3.0                & 2.00  & 0.30\%  & 224.85  \\ \cline{2-5}
& 3.6                & 2.00 & 0.35\%  & 244.06  \\ \hline
\multirow{4}{*}{3}

 & 0.6 & 3.00 & 0.28\%   & 341.31 \\ \cline{2-5}
 & 1.5                & 3.00  & 0.19\%    & 337.40 \\ \cline{2-5}
& 4.5                & 3.00 & 0.16\%   & 310.95 \\ \cline{2-5}
 & 5.4                & 3.00 & 0.09\%  & 310.50  \\ \hline

 \multirow{4}{*}{4}
 & 0.8                & 4.02 & 0.56\%   & 343.50\\ \cline{2-5}
& 2.0                & 4.02  & 0.53\%   & 342.98 \\ \cline{2-5}
 & 6.0                & 4.02 & 0.53\%    & 346.10 \\ \cline{2-5}
 & 7.2                & 4.01 & 0.37\%   & 344.70 \\ \hline
 
 \multirow{4}{*}{5}
& 1.0                & 1.14 & 77.17\%  & 136.88 \\ \cline{2-5}
 & 2.5                & 4.99  & 0.02\%  & 340.89 \\ \cline{2-5}
 & 7.5                & 5.02  & 0.51\%   & 343.42   \\ \cline{2-5}
  & 9.0                & 5.02  & 0.42\%  & 343.97 \\ \hline
\end{tabular}
\end{table}

\begin{table}[H]
\centering
\caption{PINN inverse results for the polynomial solution using measurement data generated from the direct PINN with 1\% additional noise}\label{tab:one per}
\begin{tabular}{|l|l|l|l|l|}
\hline
true m & initial guess & recovered m & PINN relative error pct & PINN training time s \\
\hline
\multirow{4}{*}{2}
& 0.4            & 2.00    & 0.04\%               & 355.85          \\ \cline{2-5}
 & 1              & 1.99     & 0.22\%               & 254.72         \\
 \cline{2-5}
 & 3              & 1.99     & 0.02\%               & 201.63          \\ \cline{2-5}
  & 3.6            & 2.00     & 0.40\%               & 223.58          \\ \hline
  \multirow{4}{*}{3}
 & 0.6            & 3.00     & 0.04\%               & 359.60         \\ \cline{2-5}
 & 1.5            & 3.00     & 0.07\%               & 360.32          \\ \cline{2-5}
  & 4.5            & 3.00   & 0.12\%               & 227.98          \\ \cline{2-5}
 & 5.4            & 3.00     & 0.02\%               & 233.90          \\ \hline
 \multirow{4}{*}{4}
 & 0.8            & 4.01    & 0.41\%               & 356.59         \\ \cline{2-5}
 & 2              & 4.01    & 0.30\%               & 329.19        \\ \cline{2-5}
 & 6              & 4.01    & 0.40\%               & 246.92         \\ \cline{2-5}
& 7.2            & 4.00    & 0.24\%               & 326.11         \\ \hline
\multirow{4}{*}{5}
 & 1              & 1.15    & 76.96\%              & 141.84        \\ \cline{2-5}
 & 2.5            & 4.99     & 0.13\%               & 354.21         \\ \cline{2-5}
& 7.5            & 5.01    & 0.39\%               & 324.37        \\ \cline{2-5}
 & 9              & 5.01     & 0.34\%               & 324.97  \\ \hline
\end{tabular}
\end{table}

\begin{table}[H]
\centering
\caption{PINN inverse results for the polynomial solution using measurement data generated from the direct PINN with 5\% additional noise}\label{tab:five pct}.
\begin{tabular}{|l|l|l|l|l|}
\hline
true m & initial guess & recovered m & PINN relative error & PINN training time s \\ \hline
\multirow{4}{*}{2}
 & 0.4            & 1.98   & 0.58\%               & 199.68         \\ \cline{2-5}
 & 1              & 1.99     & 0.38\%               & 178.59          \\ \cline{2-5}
 & 3              & 1.99     & 0.39\%               & 165.93         \\ \cline{2-5}
 & 3.6            & 1.99   & 0.16\%               & 176.49          \\ \hline
\multirow{4}{*}{3}
& 0.6            & 2.96     & 1.27\%               & 210.94          \\ \cline{2-5}
& 1.5            & 2.96     & 1.20\%               & 197.06         \\ \cline{2-5}
 & 4.5            & 2.99     & 0.27\%               & 172.51        \\ \cline{2-5}
 & 5.4            & 2.99     & 0.25\%               & 176.41          \\ \hline
\multirow{4}{*}{4}
  & 0.8            & 3.98    & 0.49\%               & 222.55        \\ \cline{2-5}
  & 2              & 3.99    & 0.12\%               & 191.97        \\ \cline{2-5}
 & 6              & 4.00     & 0.14\%               & 171.97         \\ \cline{2-5}
  & 7.2            & 4.00    & 0.12\%               & 177.12       \\ \hline
\multirow{4}{*}{5}
 & 1              & 1.17     & 76.48\%              & 142.86           \\ \cline{2-5}
 & 2.5            & 4.98     & 0.37\%               & 205.75        \\ \cline{2-5}
 & 7.5            & 5.03    & 0.61\%               & 164.99         \\ \cline{2-5}
 & 9              & 5.02     & 0.56\%               & 173.50         \\ \hline
\end{tabular}
\end{table}

Overall, across all three data scenarios, the recovered value of m remains accurate, with the relative error not exceeding 1.5\% in every case. The only exception is the case where the true m is 5 and the initial guess is 1, for which the error stays between 76\% and 77\% across all perturbation levels. This case, however, already fails in the noise-free scenario, where the error is 76.74\%; and since the noisy datasets are obtained from the same direct PINN solution, this behaviour is carried over to the 1\% and 5\% cases as well. Apart from this case, the recovered m remains stable as the noise level increases, which indicates that the inverse PINN is indeed robust to perturbations in the measurement data.

\section{Conclusion}

In this paper, we developed a PINN framework for solving the inverse one-dimensional Porous Medium Equation. The main objective was to establish a robust inverse PINN methodology and provide a strong foundation for future extensions to higher-dimensional PME problem. To achieve this, we investigated the factors affecting inverse PINN performance and proposed a novel two-stage learning framework to improve convergence robustness and parameter recovery.

The results showed that the proposed framework successfully solves both the direct and inverse PME problems and produces accurate solutions across different manufactured solutions and values of \(m \) . In particular, the proposed two-stage strategy improves the robustness of the inverse problem and reduces the dependence on the initial parameter guess.

Future work will focus on extending the framework to two-dimensional and higher-dimensional PME problems, including Barenblatt solutions and solutions generated using finite difference or other numerical methods for validation, comparison, and data collection. At this point, one clear advantage of the proposed PINN framework is that these high-dimensional extensions require only minor modifications to the existing code, mainly through adding extra spatial coordinates to the network inputs, updating the PDE residual, and sampling collocation points in the new domains. Overall, the developed framework provides a solid foundation for future studies of more complex PME problems, where the advantages of PINNs are expected to become more significant.

\newpage
\appendix
\section{Visual Comparisons of Direct Solutions}
\label{app:A}

For each Solution and \( m \), we visualize the heatmap of the PINN, exact, and classical solution, along with the evolution of the PINN training loss.  \vspace{-2mm}
\subsection{Barenblatt solution}
\vspace{1mm}
\textbf{m = 2}\vspace{-2mm}
\begin{figure}[H]
    \centering
    \begin{subfigure}[b]{0.32\textwidth}
        \centering
        \includegraphics[width=\textwidth, height = 6.5 cm]{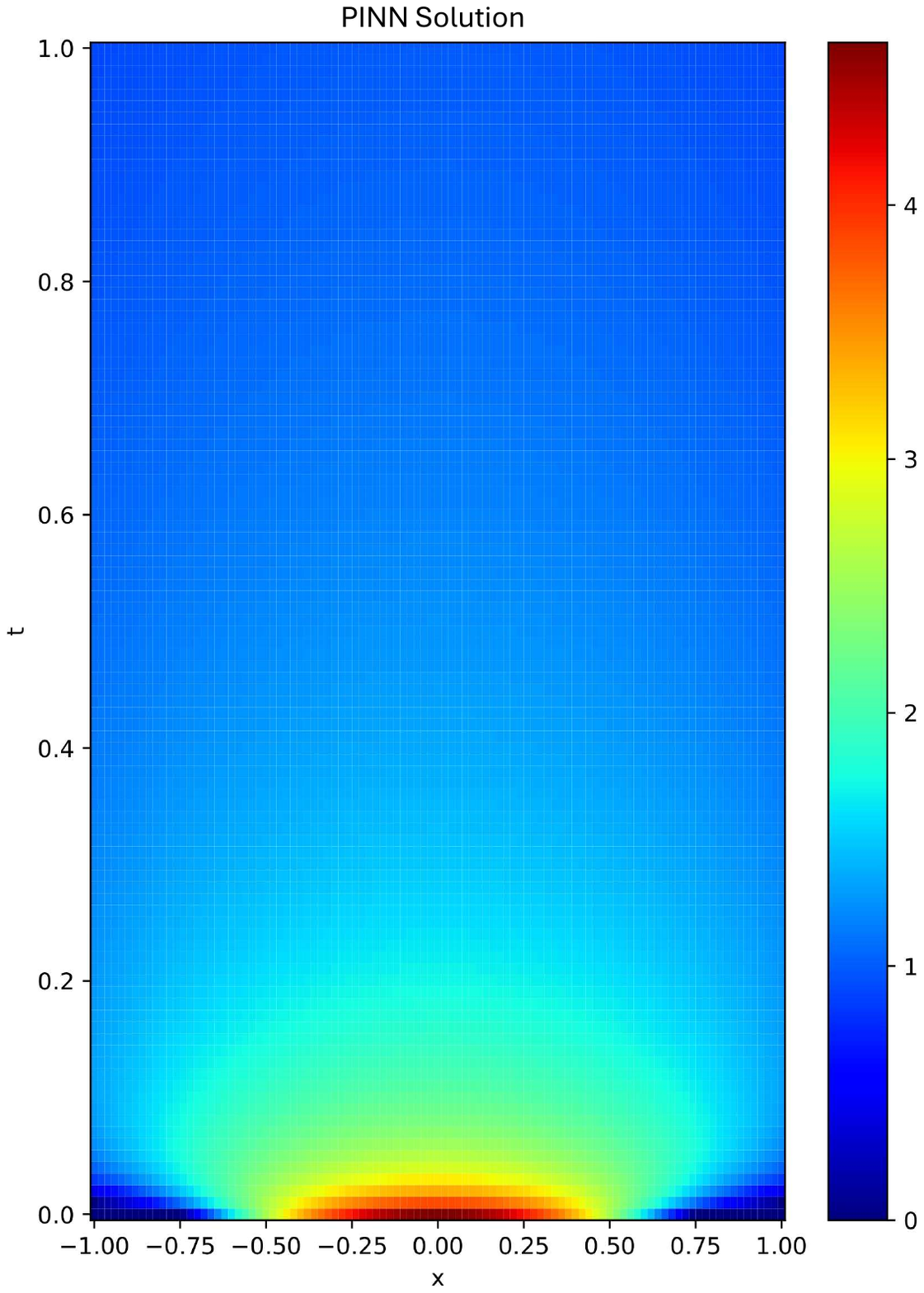}
        \caption{PINN solution}
    \end{subfigure}
    \hfill
    \begin{subfigure}[b]{0.32\textwidth}
        \centering
        \includegraphics[width=\textwidth, height = 6.5 cm]{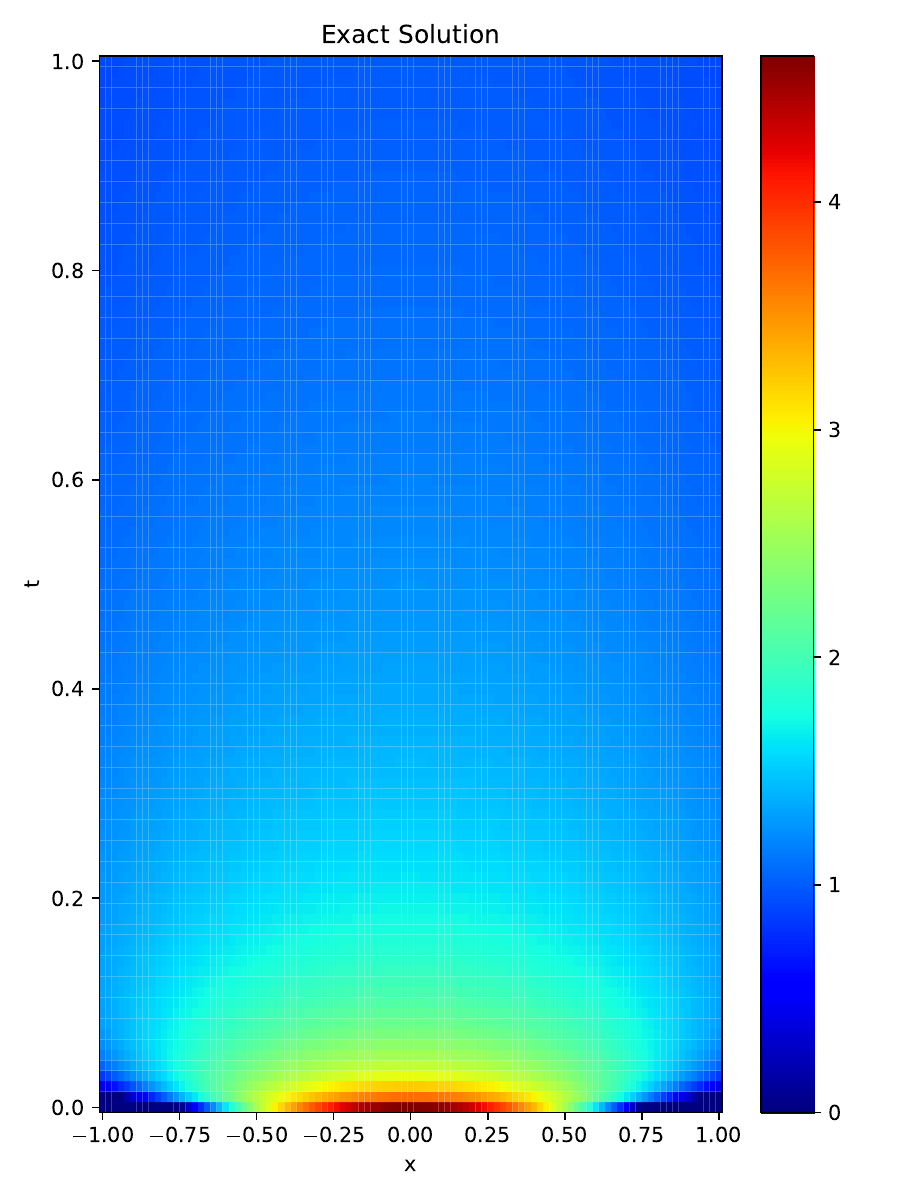}
        \caption{Exact solution}
    \end{subfigure}
    \hfill
    \begin{subfigure}[b]{0.32\textwidth}
        \centering
        \includegraphics[width=\textwidth, height = 6.5 cm]{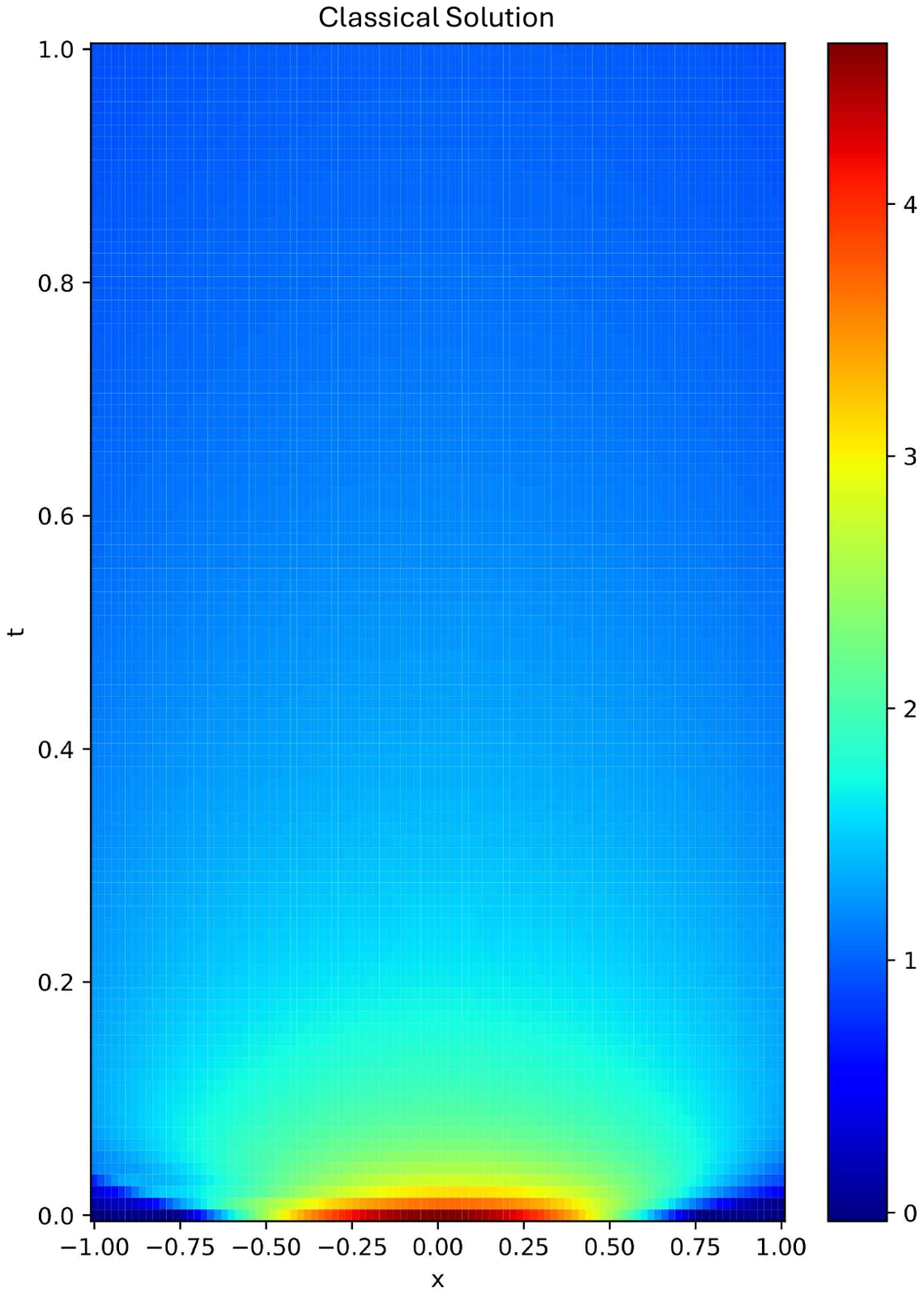}
        \caption{Classical solution}
    \end{subfigure}\vspace{-4mm}
\end{figure}
 \begin{figure}[H]
    \centering
\includegraphics[width=0.8\linewidth]{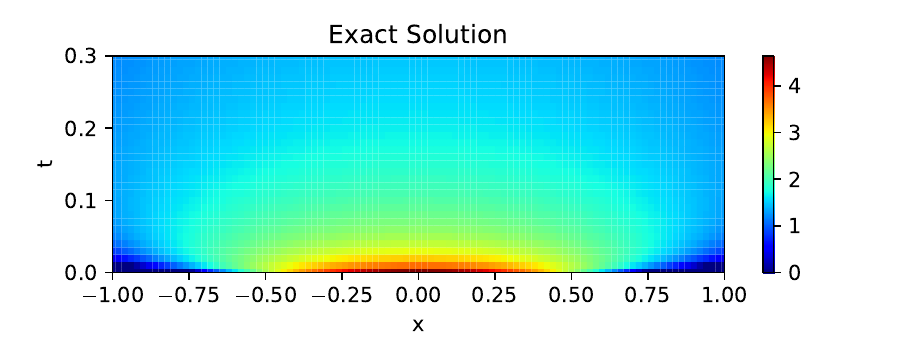}\vspace{-6mm}
    \caption{Zoomed-in view of the exact Barenblatt solution for $m = 2 $ for t $\in [0,0.3]$}
\end{figure}
\begin{figure}[H]
\begin{subfigure}{0.50\textwidth}
\includegraphics[width=0.9\linewidth, height=5cm]{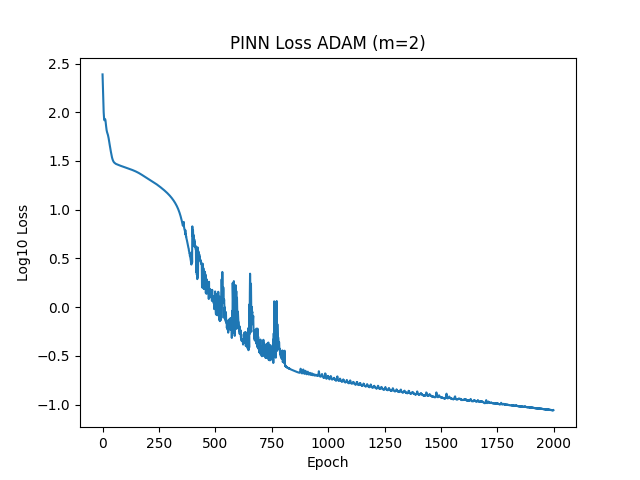} 
\caption{PINN ADAM loss}
\label{fig:subim1}
\end{subfigure}
\begin{subfigure}{0.50\textwidth}
\includegraphics[width=0.9\linewidth, height=5cm]{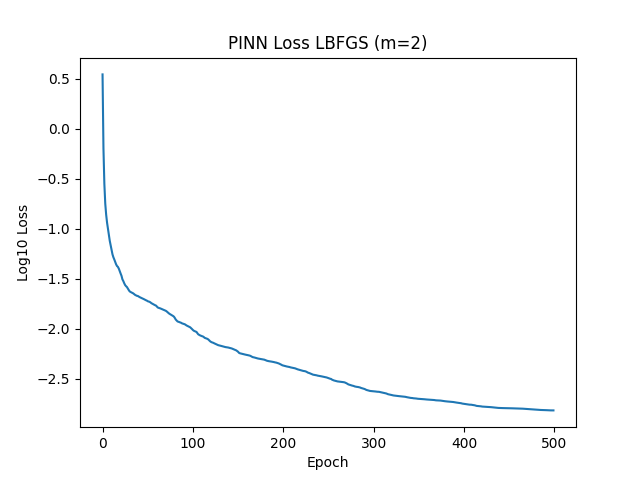}
\caption{PINN LBFGS loss }
\label{fig:subim2}
\end{subfigure}
\end{figure}
%\hrule
%\vspace{3mm}
\newpage
\textbf{m = 3}

\begin{figure}[H]
    \centering
    
    \begin{subfigure}[b]{0.32\textwidth}
        \centering
        \includegraphics[width=\textwidth, height = 6.5 cm]{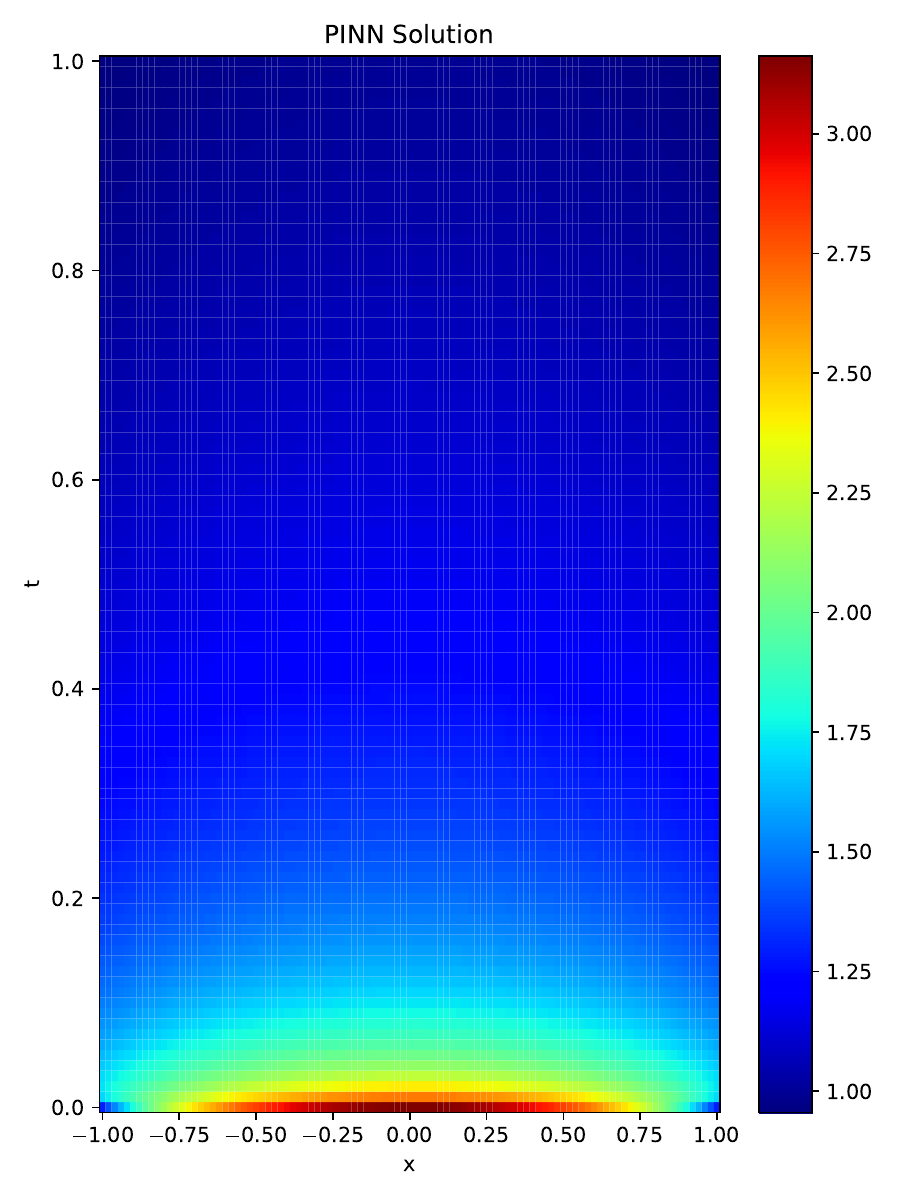}
        \caption{PINN solution}
    \end{subfigure}
    \hfill
    \begin{subfigure}[b]{0.32\textwidth}
        \centering
        \includegraphics[width=\textwidth, height = 6.5 cm]{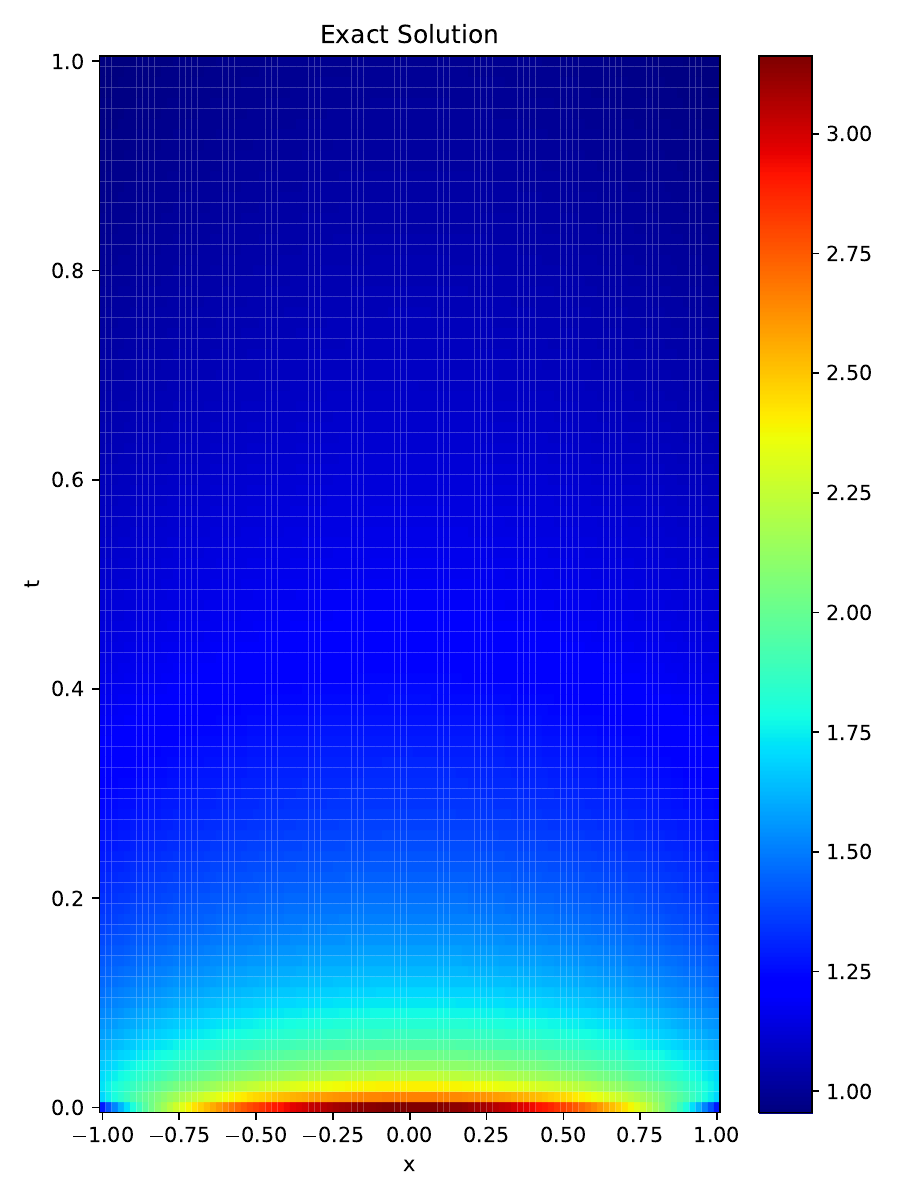}
        \caption{Exact solution}
    \end{subfigure}
    \hfill
    \begin{subfigure}[b]{0.32\textwidth}
        \centering
        \includegraphics[width=\textwidth, height = 6.5 cm]{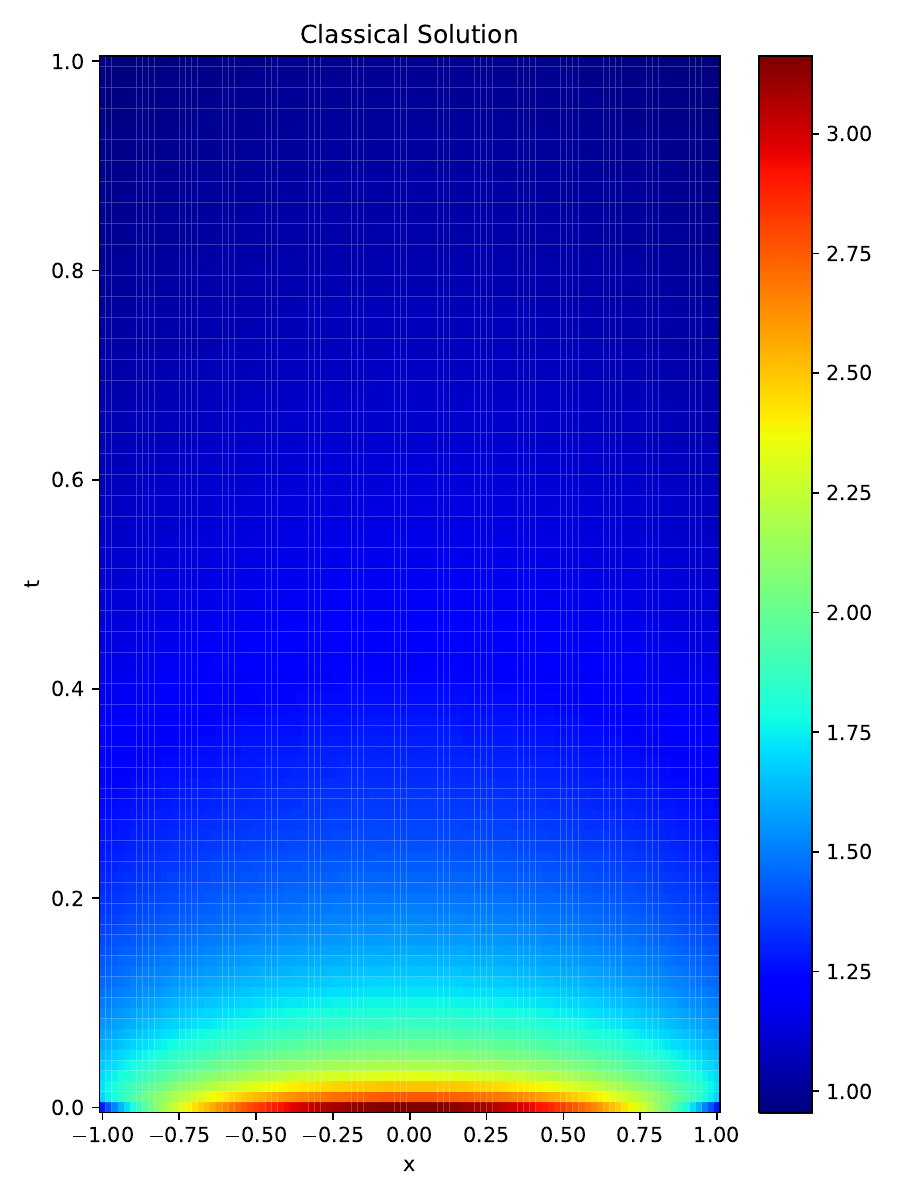}
        \caption{Classical solution}
    \end{subfigure}
\end{figure}
\begin{figure}[H]
    \centering
    \includegraphics[width=0.8\linewidth]{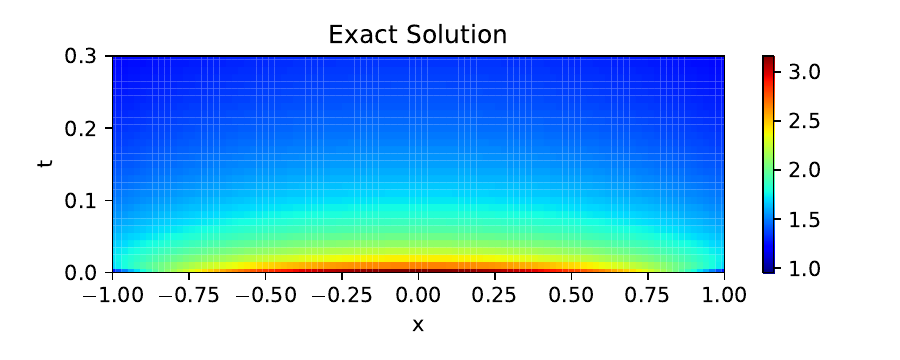}
    \caption{Zoomed-in view of the exact Barenblatt solution for $m = 3 $ for t $\in [0,0.3]$}
\end{figure}

\begin{figure}[H]

\begin{subfigure}{0.50\textwidth}
\includegraphics[width=0.9\linewidth, height=5cm]{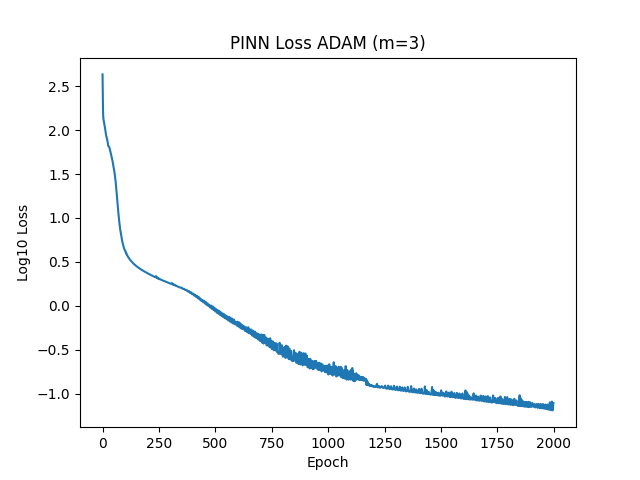} 
\caption{PINN ADAM loss}
\label{fig:subim1}
\end{subfigure}
\begin{subfigure}{0.50\textwidth}
\includegraphics[width=0.9\linewidth, height=5cm]{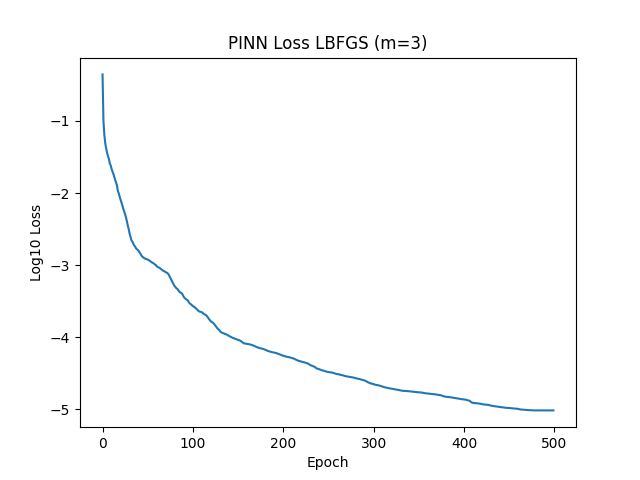}
\caption{PINN LBFGS loss }
\label{fig:subim2}
\end{subfigure}

\end{figure}

\newpage
\textbf{m = 4}
\begin{figure}[H]
    \centering
    
    \begin{subfigure}[b]{0.32\textwidth}
        \centering
        \includegraphics[width=\textwidth, height =6.5 cm]{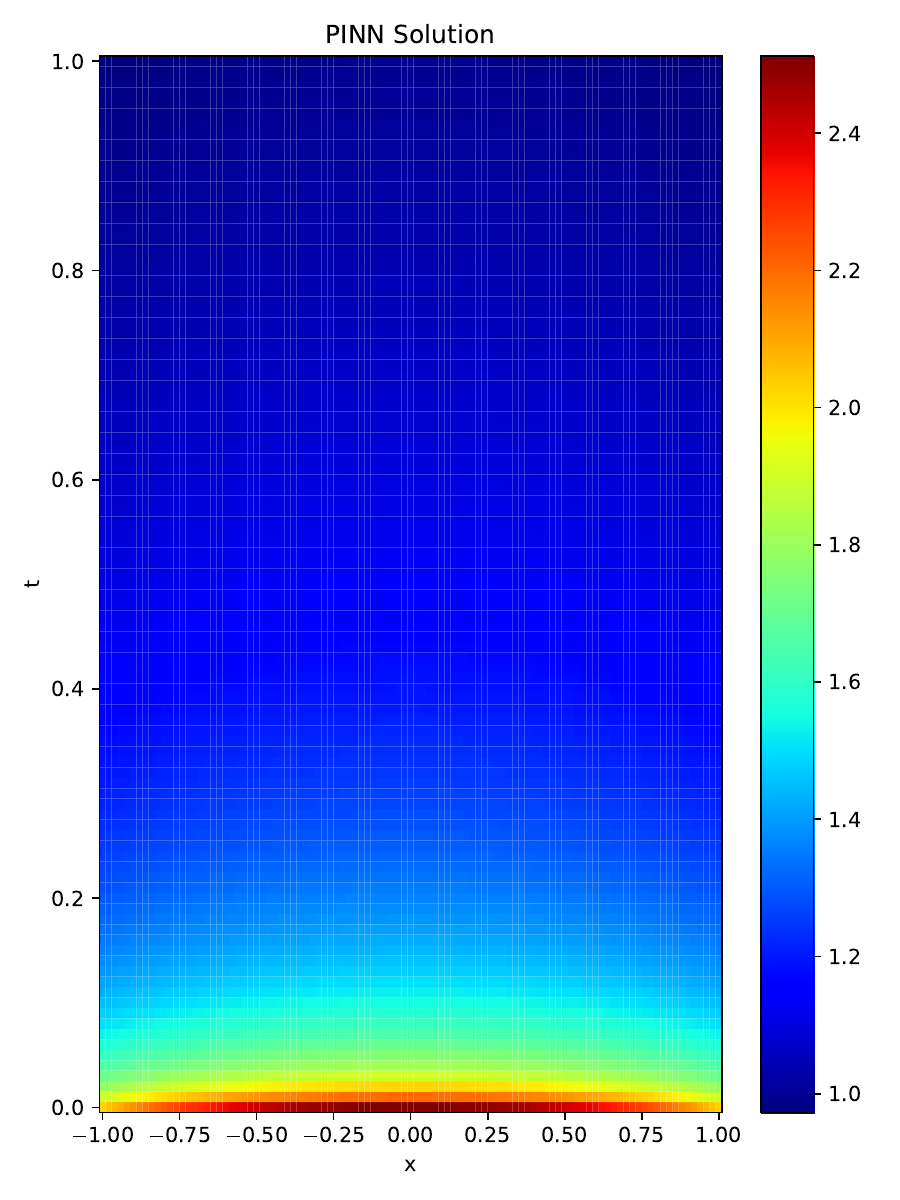}
        \caption{PINN solution}
    \end{subfigure}
    \hfill
    \begin{subfigure}[b]{0.32\textwidth}
        \centering
        \includegraphics[width=\textwidth, height = 6.5 cm]{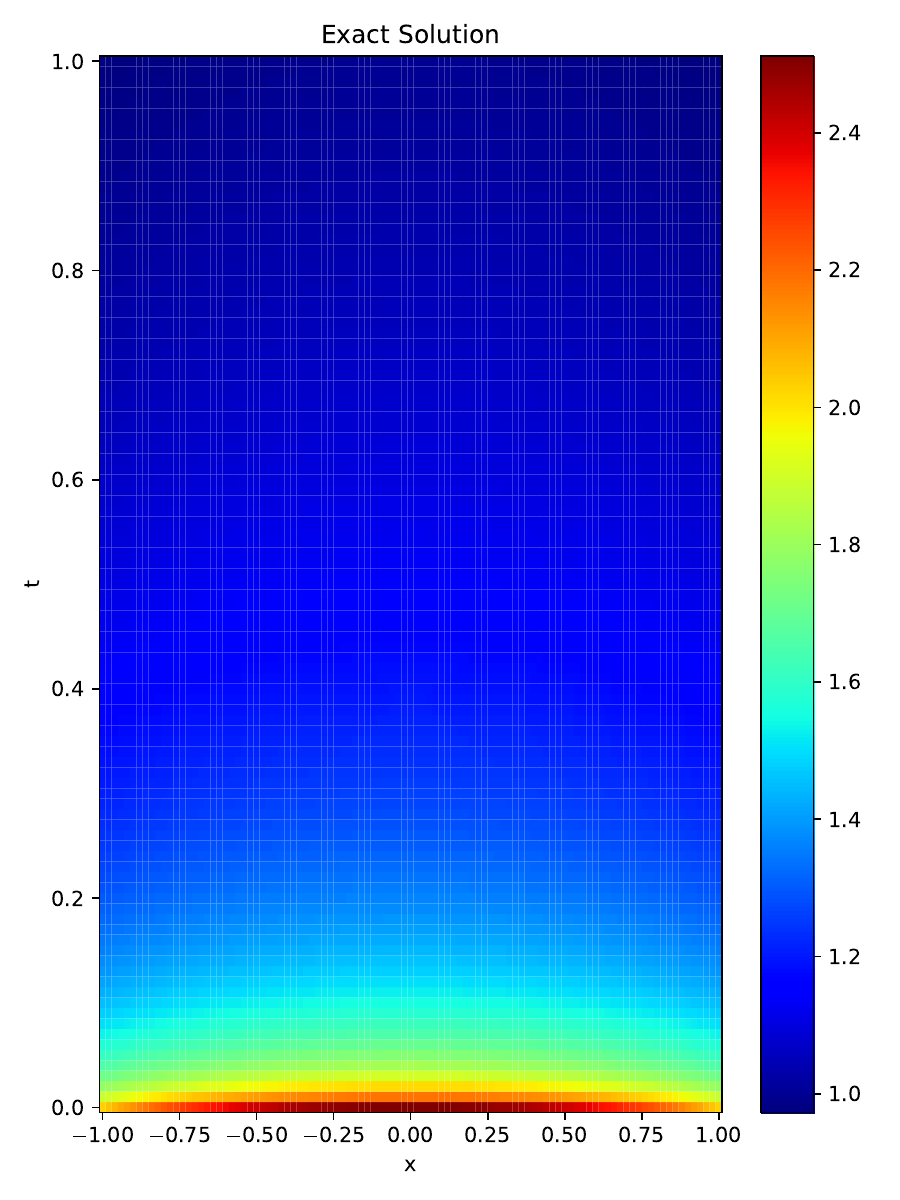}
        \caption{Exact solution}
    \end{subfigure}
    \hfill
    \begin{subfigure}[b]{0.32\textwidth}
        \centering
        \includegraphics[width=\textwidth, height = 6.5 cm ]{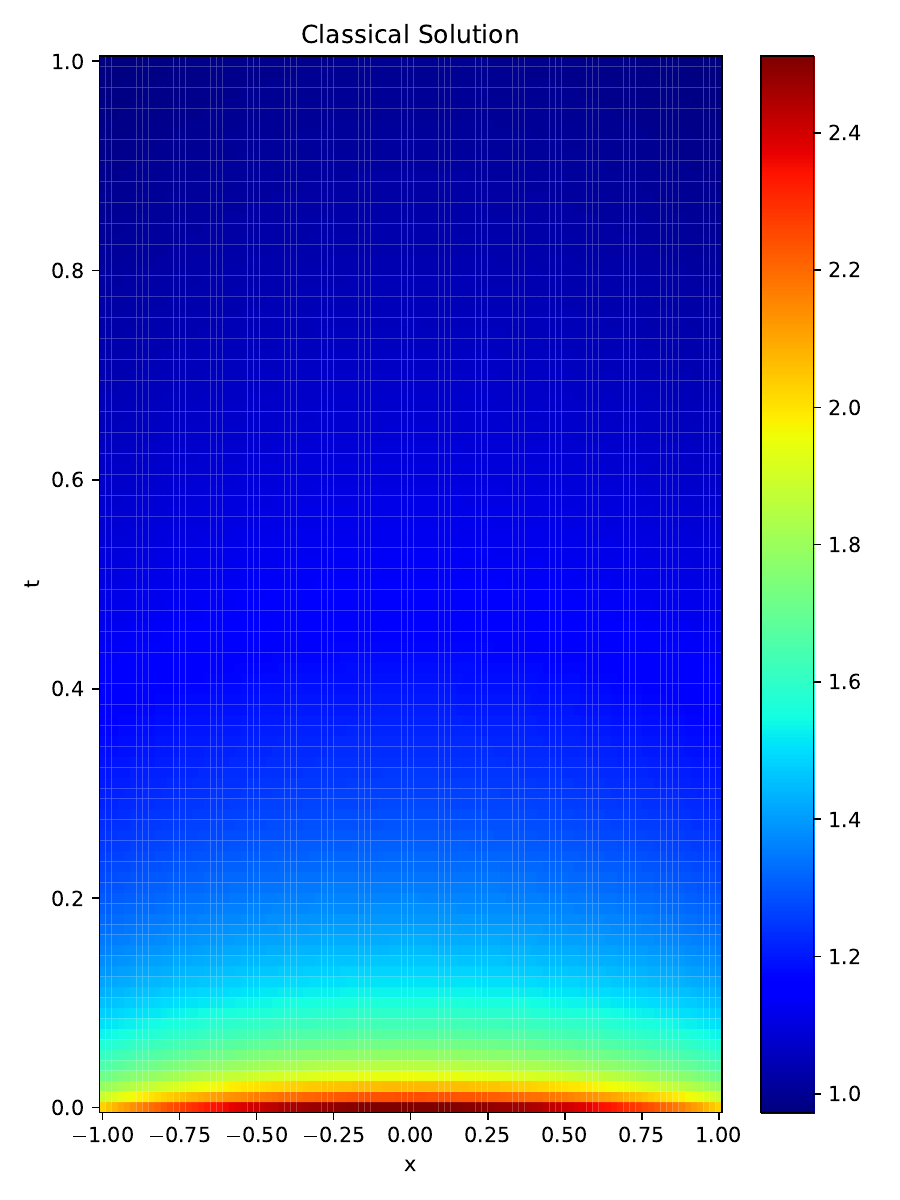}
        \caption{Classical solution}
    \end{subfigure}
\end{figure}

\begin{figure}[H]
    \centering
    \includegraphics[width=0.8\linewidth]{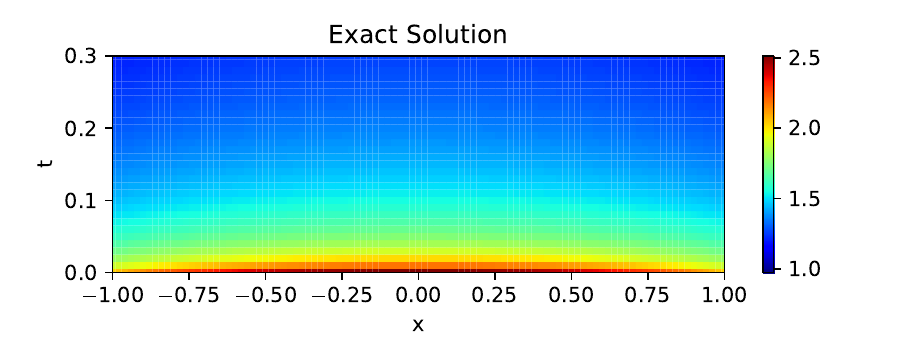}
    \caption{Zoomed-in view of the exact Barenblatt solution for $m = 4 $ for t $\in [0,0.3]$}
\end{figure}

\begin{figure}[H]

\begin{subfigure}{0.50\textwidth}
\includegraphics[width=0.9\linewidth, height=5cm]{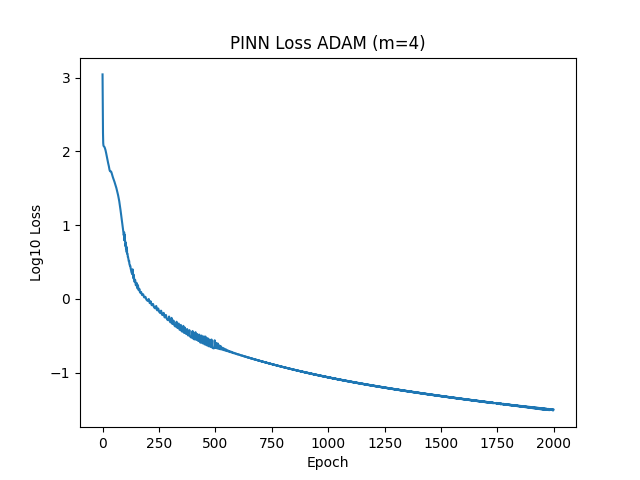} 
\caption{PINN ADAM loss}
\label{fig:subim1}
\end{subfigure}
\begin{subfigure}{0.50\textwidth}
\includegraphics[width=0.9\linewidth, height=5cm]{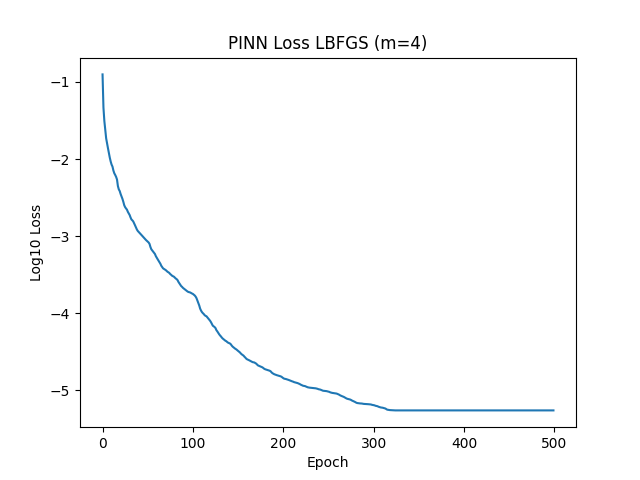}
\caption{PINN LBFGS loss }
\label{fig:subim2}
\end{subfigure}

\end{figure}

\newpage

\textbf{m = 5}

\begin{figure}[H]
    \centering
    
    \begin{subfigure}[b]{0.32\textwidth}
        \centering
        \includegraphics[width=\textwidth, height =6.5 cm]{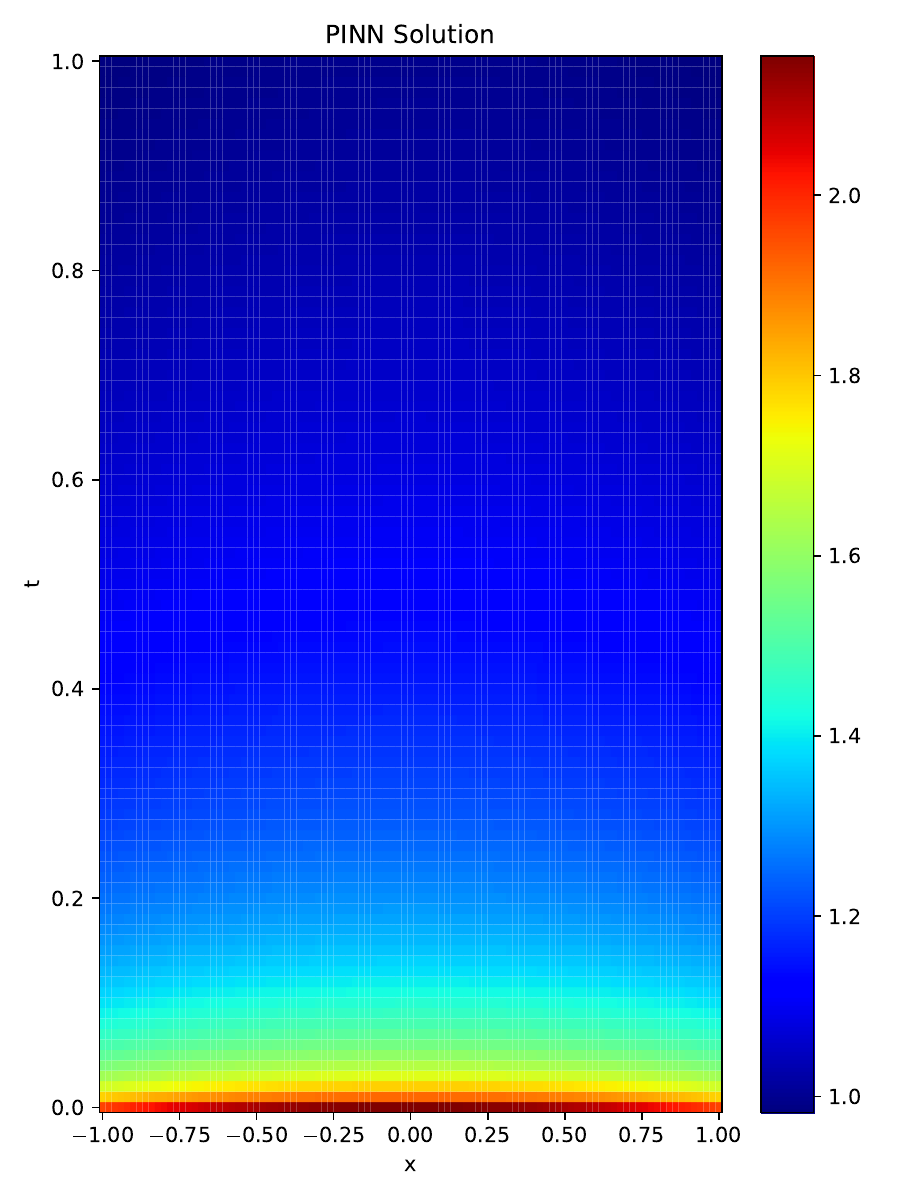}
        \caption{PINN solution}
    \end{subfigure}
    \hfill
    \begin{subfigure}[b]{0.32\textwidth}
        \centering
        \includegraphics[width=\textwidth, height = 6.5 cm]{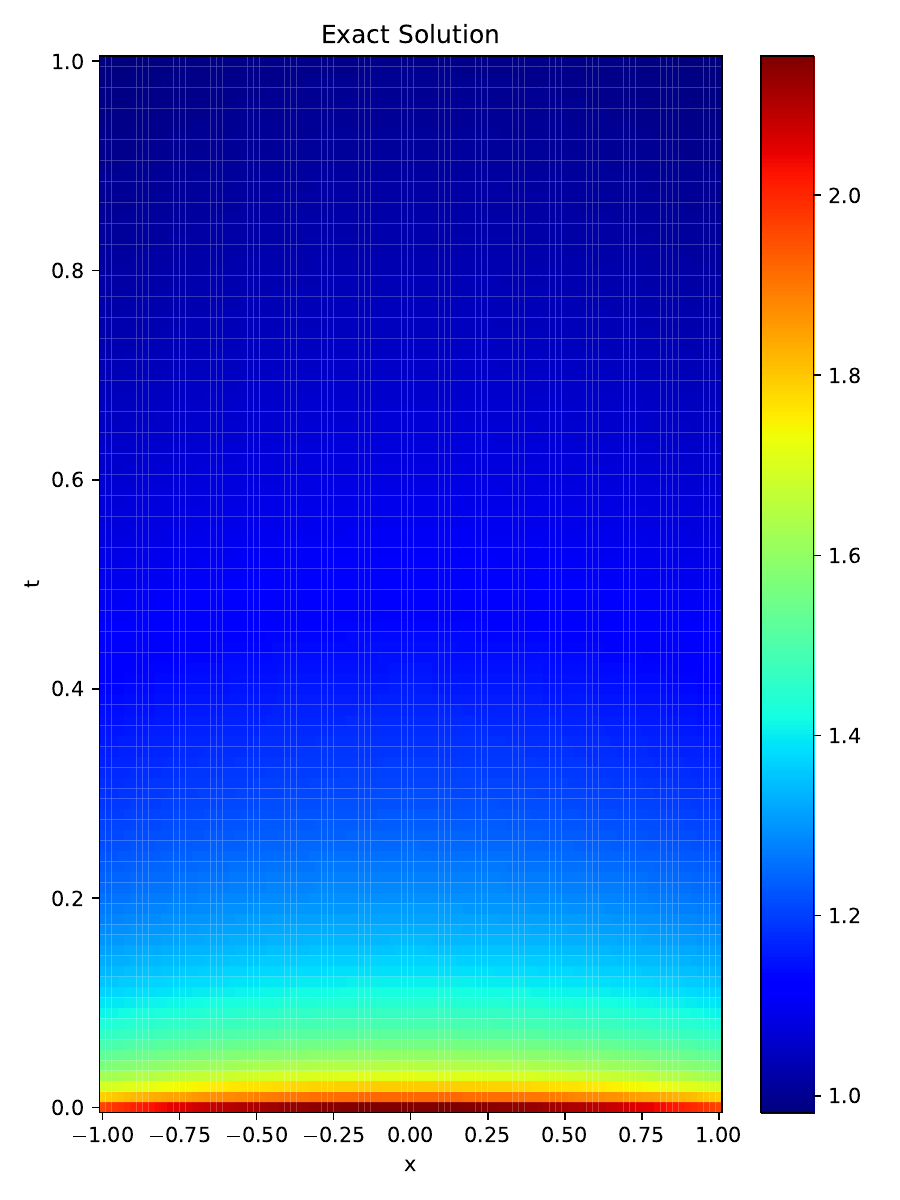}
        \caption{Exact solution}
    \end{subfigure}
    \hfill
    \begin{subfigure}[b]{0.32\textwidth}
        \centering
        \includegraphics[width=\textwidth, height = 6.5 cm ]{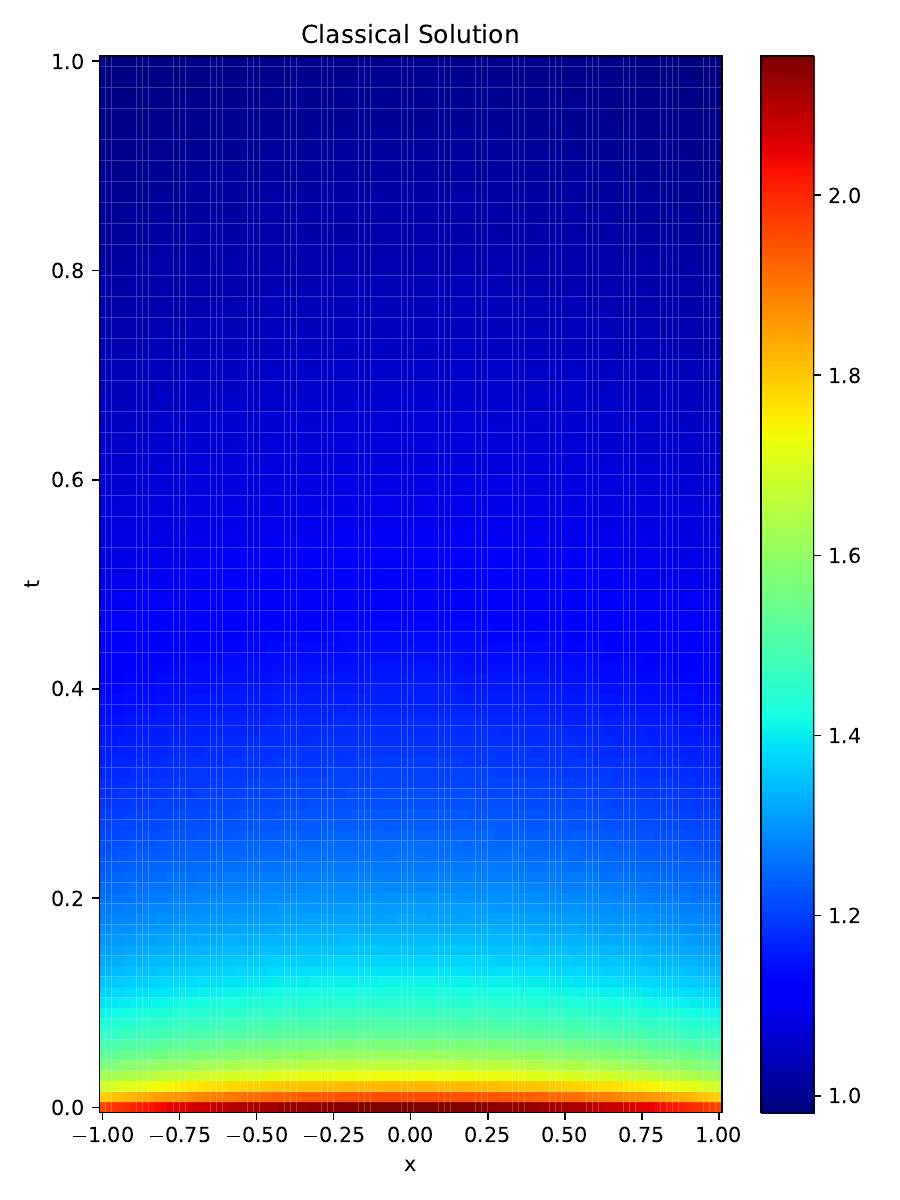}
        \caption{Classical solution}
    \end{subfigure}
\end{figure}

\begin{figure}[H]
    \centering
    \includegraphics[width=0.8\linewidth]{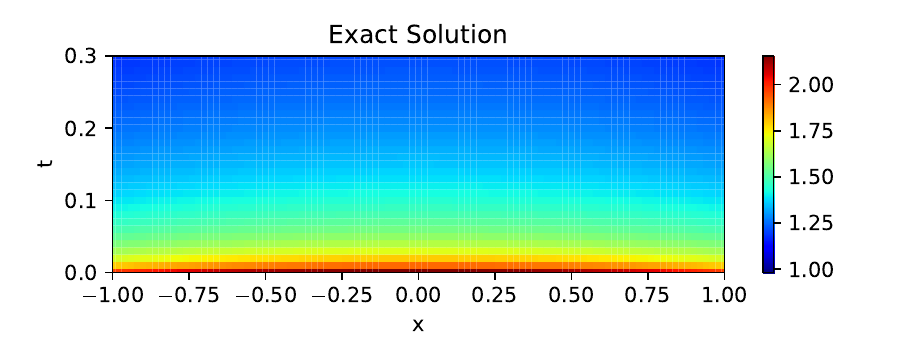}
    \caption{Zoomed-in view of the exact Barenblatt solution for $m = 5 $ for t $\in [0,0.3]$}
\end{figure}

\begin{figure}[H]

\begin{subfigure}{0.50\textwidth}
\includegraphics[width=0.9\linewidth, height=5cm]{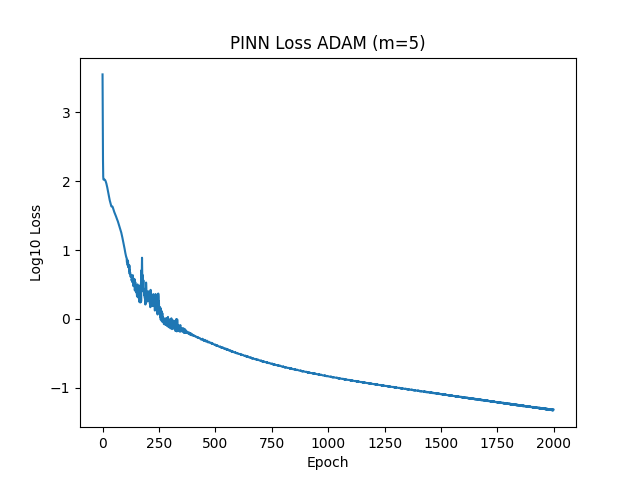} 
\caption{PINN ADAM loss}
\label{fig:subim1}
\end{subfigure}
\begin{subfigure}{0.50\textwidth}
\includegraphics[width=0.9\linewidth, height=5cm]{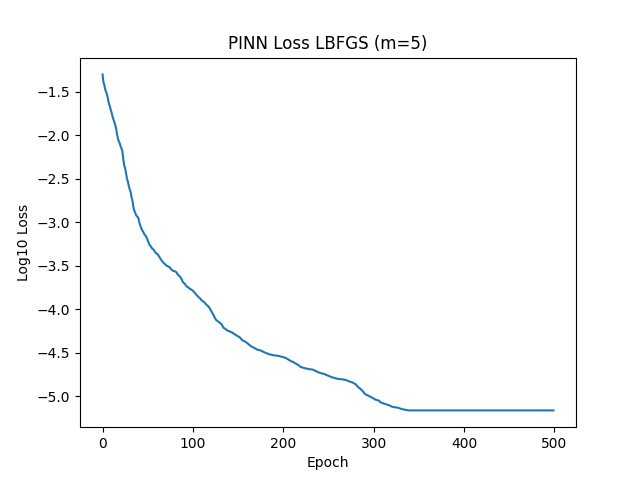}
\caption{PINN LBFGS loss }
\label{fig:subim2}
\end{subfigure}

\end{figure}

\subsection{Polynomial solution}
\label{sec: poly vis }

\vspace{3mm}
\textbf{m = 2}

\begin{figure}[H]
    \centering
    
    \begin{subfigure}[b]{0.32\textwidth}
        \centering
        \includegraphics[width=\textwidth, height = 6.5 cm]{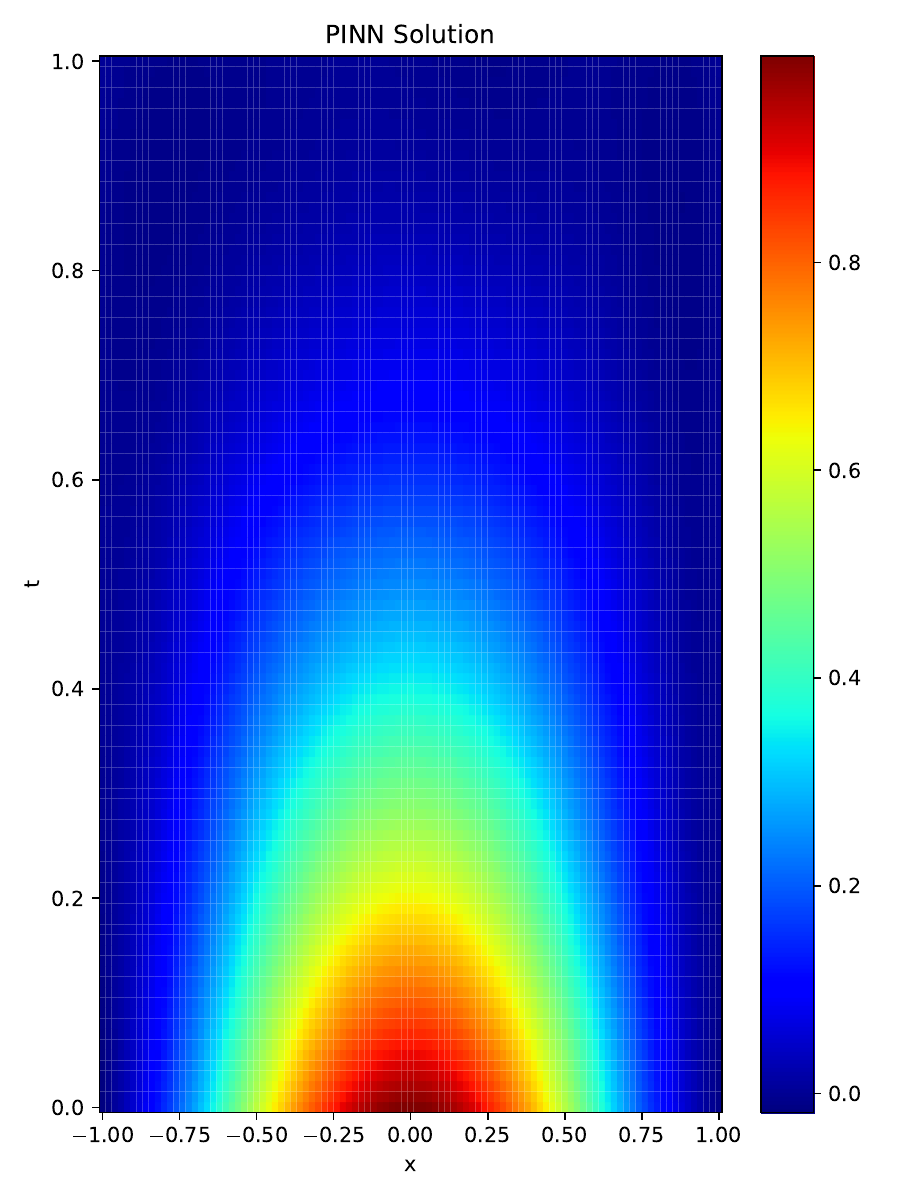}
        \caption{PINN solution}
    \end{subfigure}
    \hfill
    \begin{subfigure}[b]{0.32\textwidth}
        \centering
        \includegraphics[width=\textwidth, height = 6.5 cm]{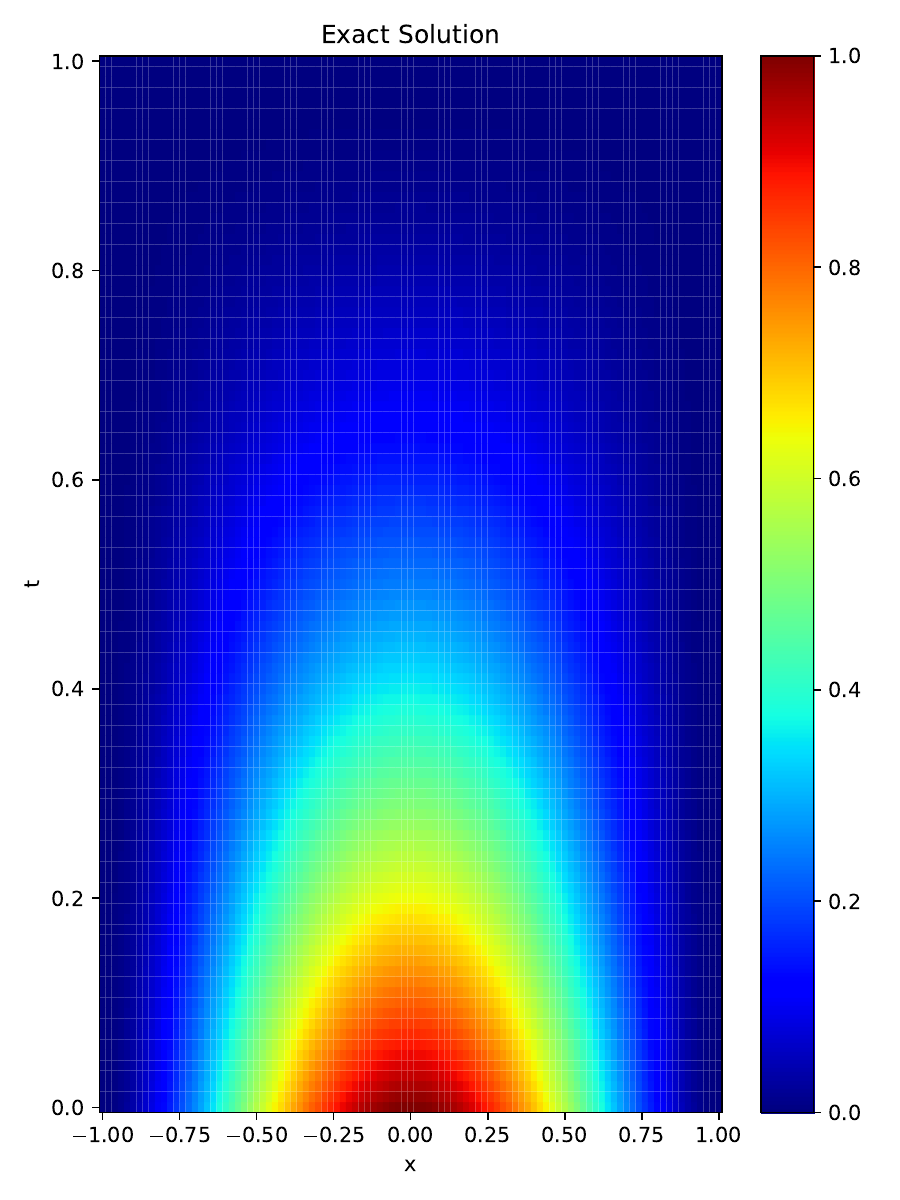}
        \caption{Exact solution}
    \end{subfigure}
    \hfill
    \begin{subfigure}[b]{0.32\textwidth}
        \centering
        \includegraphics[width=\textwidth, height = 6.5 cm]{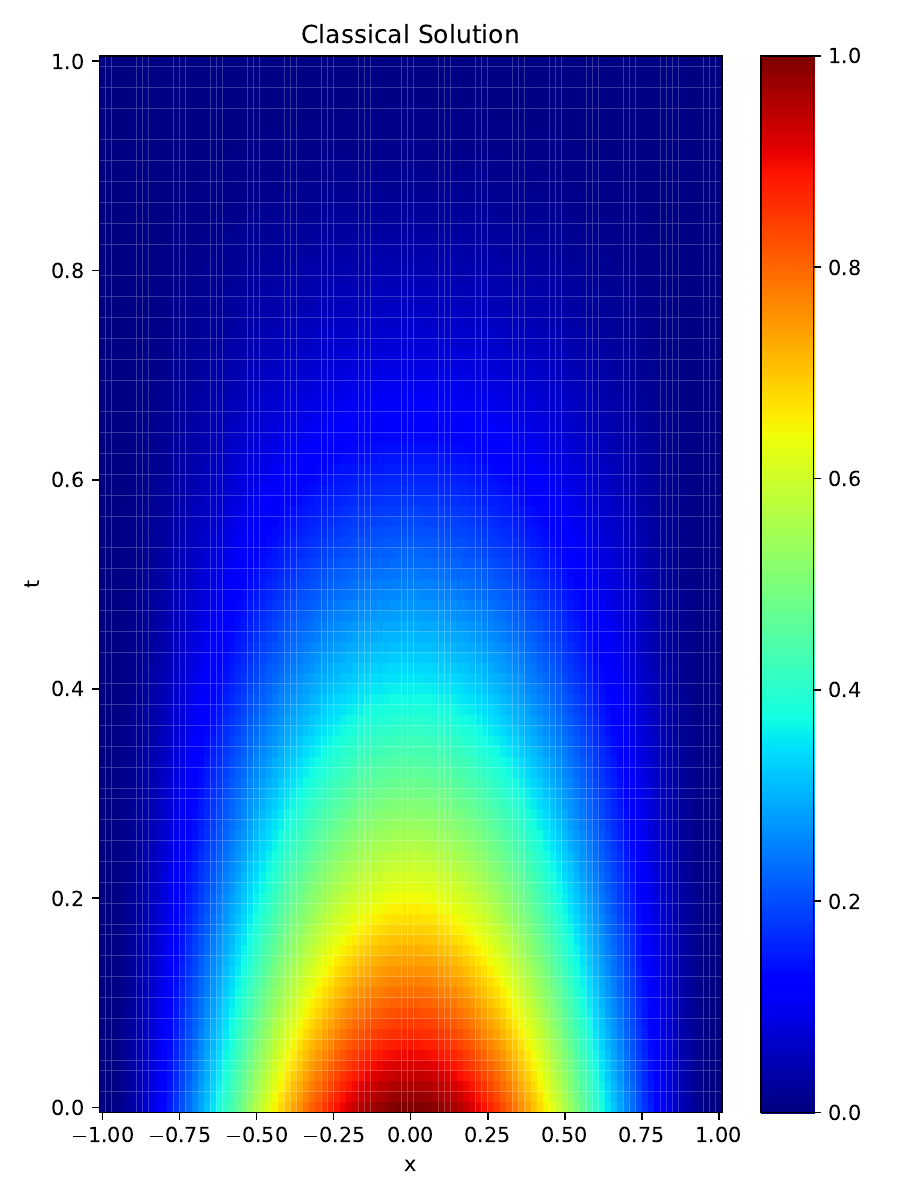}
        \caption{Classical solution}
    \end{subfigure}
\end{figure}

\begin{figure}[H]

\begin{subfigure}{0.50\textwidth}
\includegraphics[width=0.9\linewidth, height=5cm]{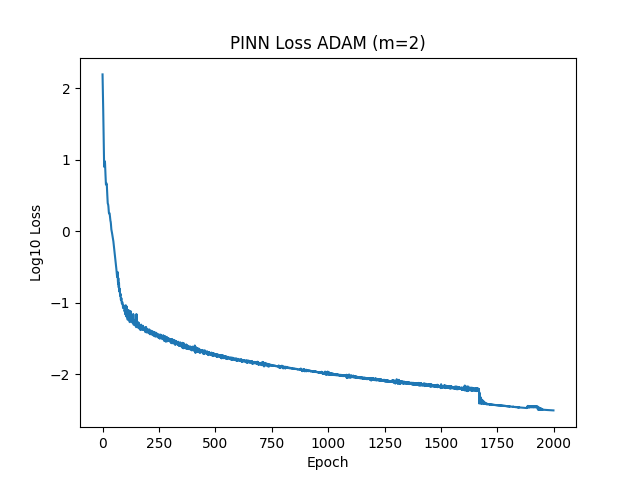} 
\caption{PINN ADAM loss}
\label{fig:subim1}
\end{subfigure}
\begin{subfigure}{0.50\textwidth}
\includegraphics[width=0.9\linewidth, height=5cm]{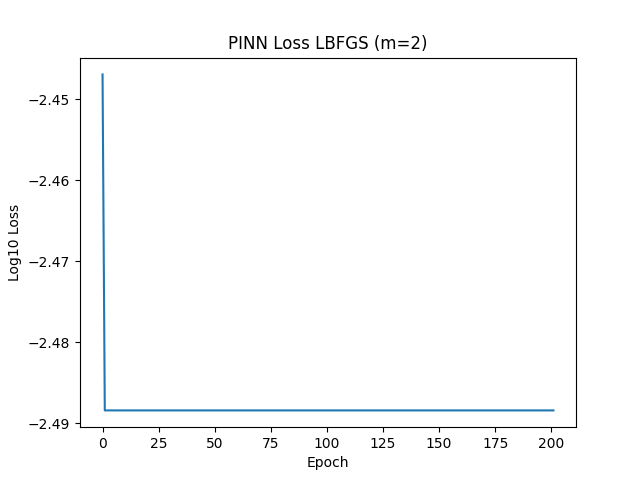}
\caption{PINN LBFGS loss }
\label{fig:subim2}
\end{subfigure}

\end{figure}

\hrule
\vspace{3mm}
\textbf{m = 3}

\begin{figure}[H]
    \centering
    
    \begin{subfigure}[b]{0.32\textwidth}
        \centering
        \includegraphics[width=\textwidth, height = 6.5 cm]{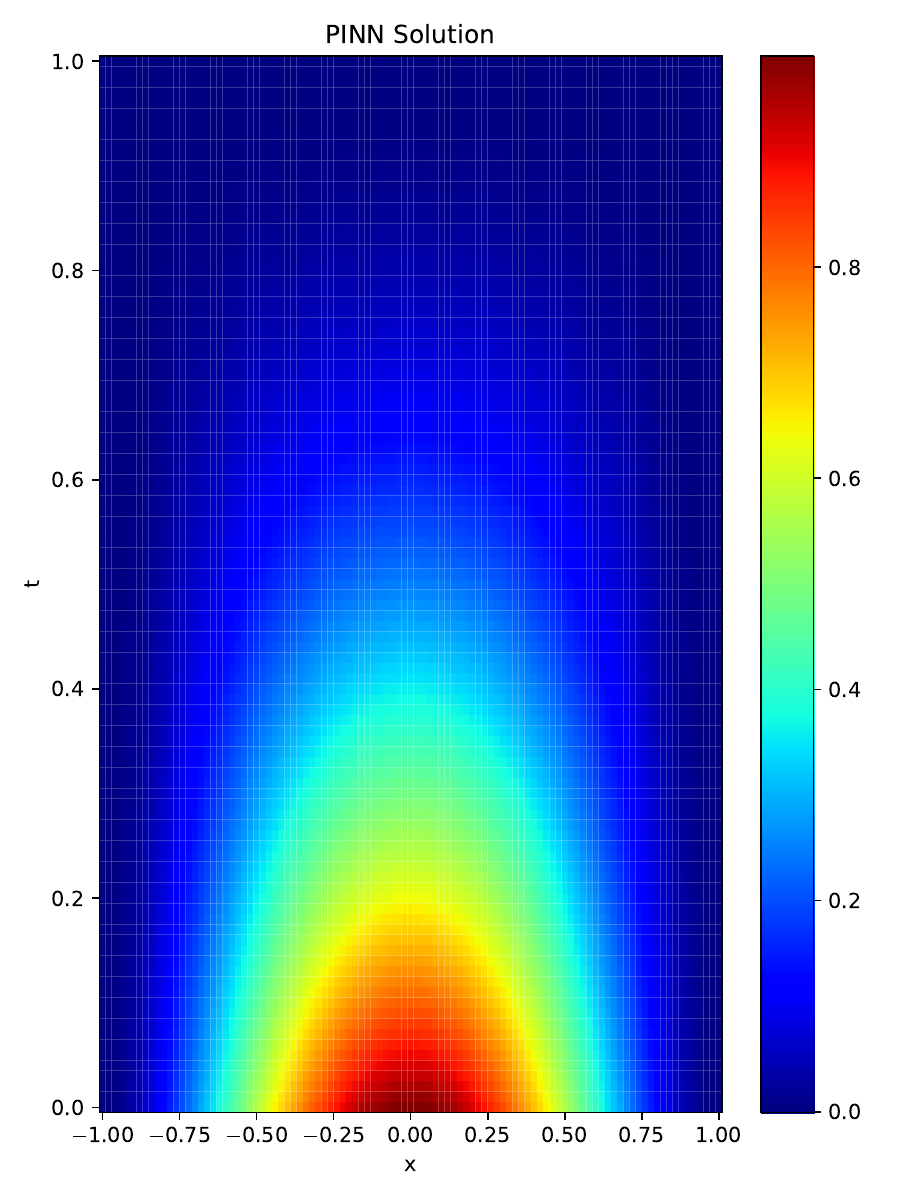}
        \caption{PINN solution}
    \end{subfigure}
    \hfill
    \begin{subfigure}[b]{0.32\textwidth}
        \centering
        \includegraphics[width=\textwidth, height = 6.5 cm]{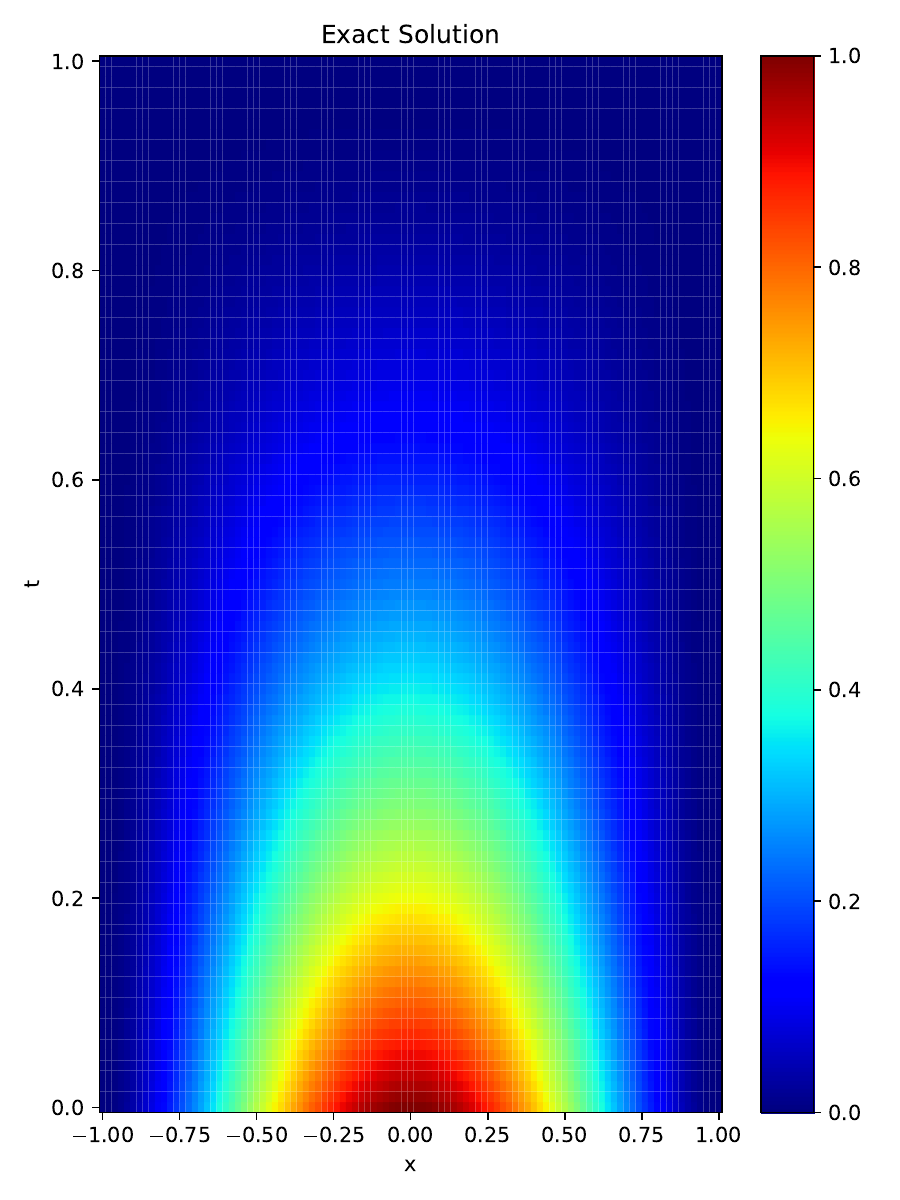}
        \caption{Exact solution}
    \end{subfigure}
    \hfill
    \begin{subfigure}[b]{0.32\textwidth}
        \centering
        \includegraphics[width=\textwidth, height = 6.5 cm]{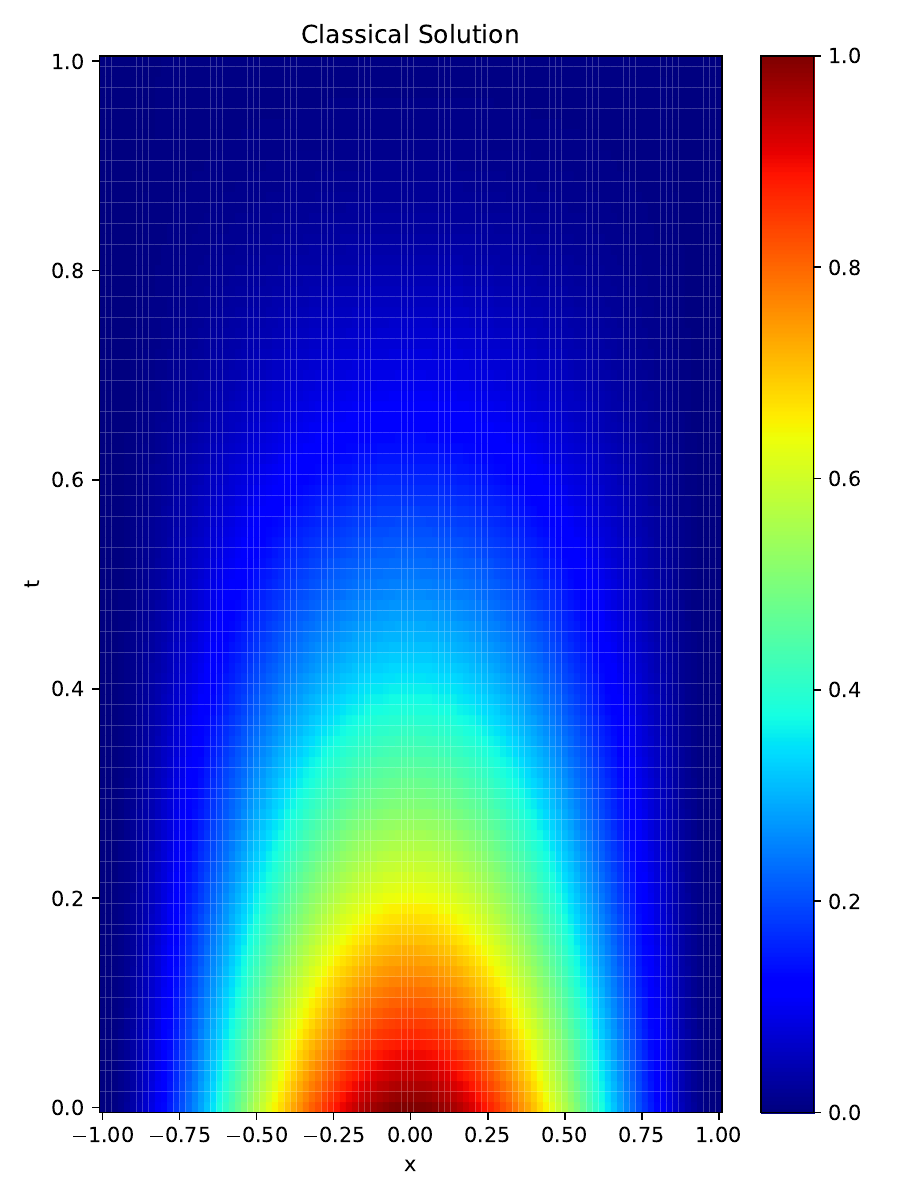}
        \caption{Classical solution}
    \end{subfigure}
\end{figure}

\begin{figure}[H]

\begin{subfigure}{0.50\textwidth}
\includegraphics[width=0.9\linewidth, height=5cm]{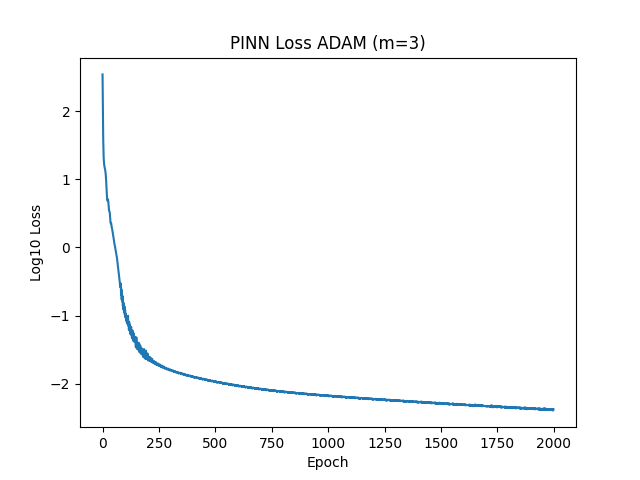} 
\caption{PINN ADAM loss}
\label{fig:subim1}
\end{subfigure}
\begin{subfigure}{0.50\textwidth}
\includegraphics[width=0.9\linewidth, height=5cm]{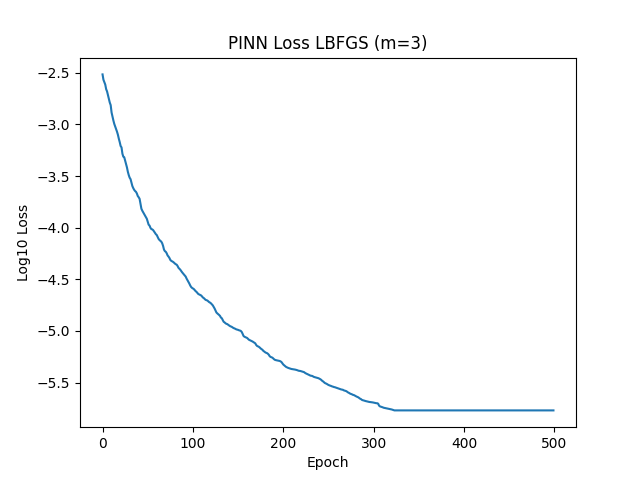}
\caption{PINN LBFGS loss }
\label{fig:subim2}
\end{subfigure}

\end{figure}

\hrule
\vspace{3mm}
\textbf{m = 4}
\begin{figure}[H]
    \centering
    
    \begin{subfigure}[b]{0.32\textwidth}
        \centering
        \includegraphics[width=\textwidth, height =6.5 cm]{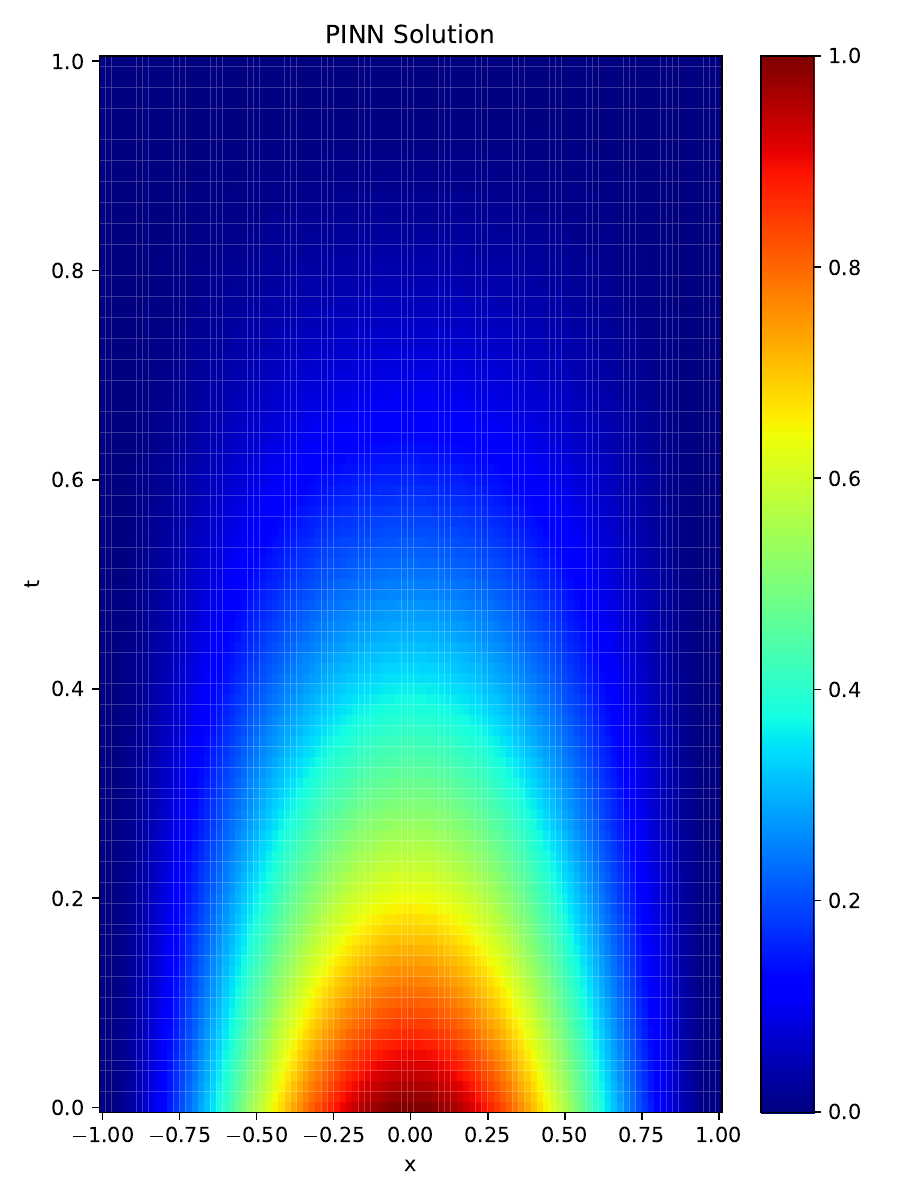}
        \caption{PINN solution}
    \end{subfigure}
    \hfill
    \begin{subfigure}[b]{0.32\textwidth}
        \centering
        \includegraphics[width=\textwidth, height = 6.5 cm]{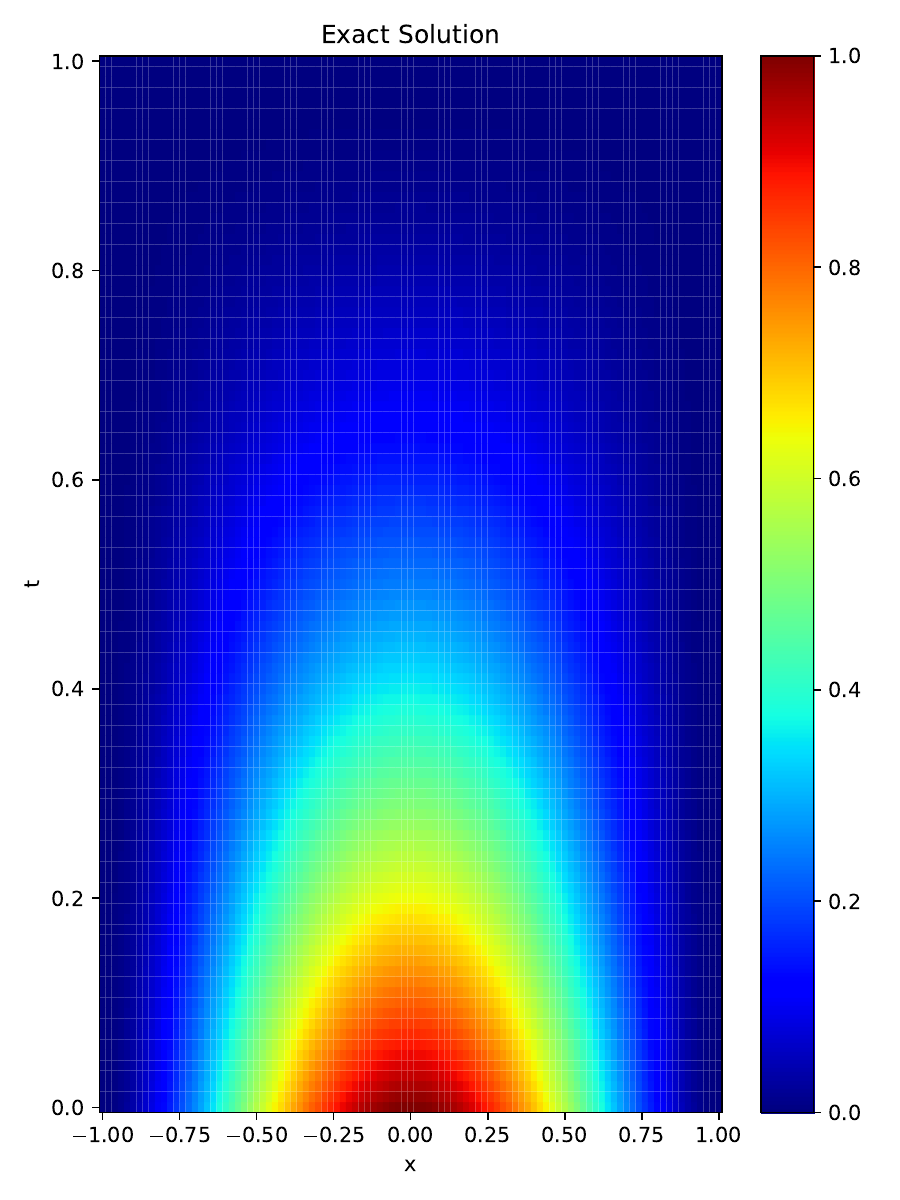}
        \caption{Exact solution}
    \end{subfigure}
    \hfill
    \begin{subfigure}[b]{0.32\textwidth}
        \centering
        \includegraphics[width=\textwidth, height = 6.5 cm ]{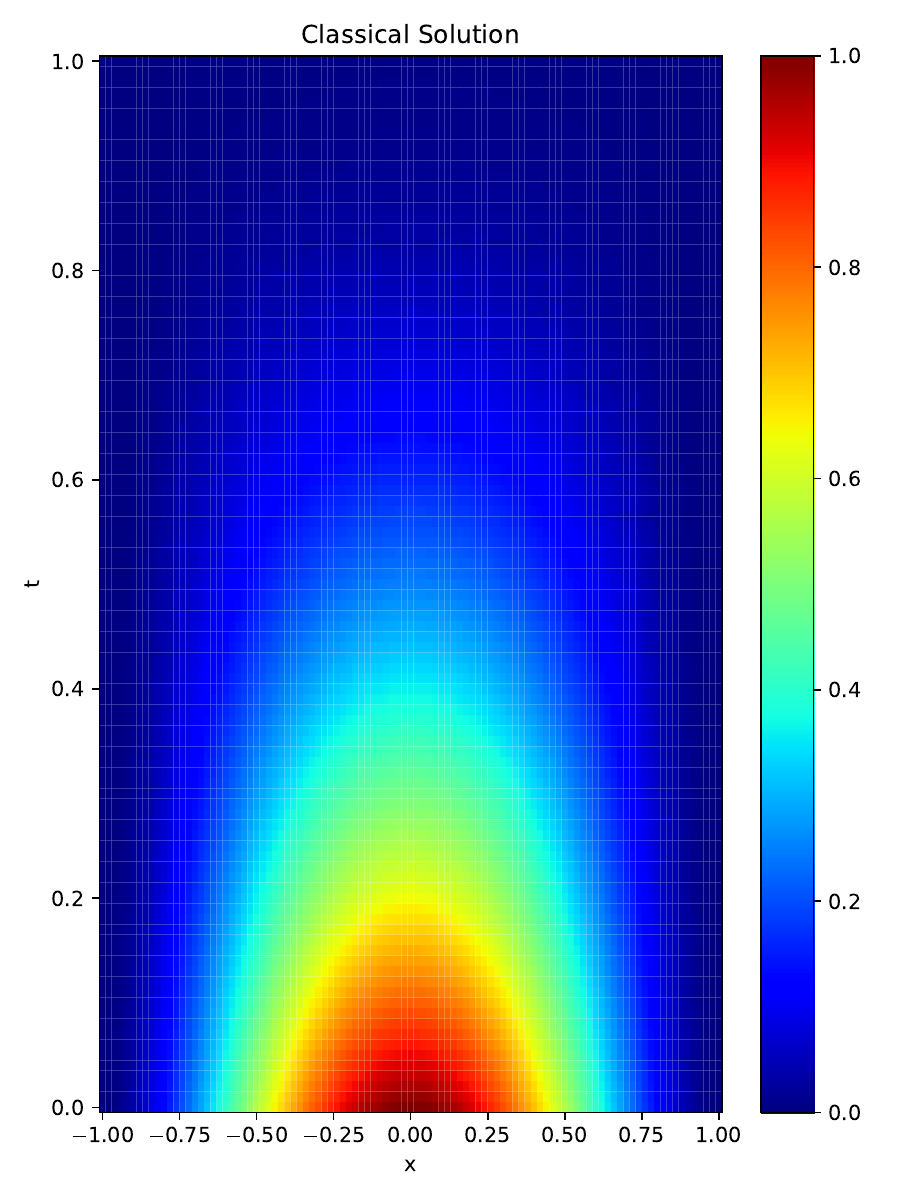}
        \caption{Classical solution}
    \end{subfigure}
\end{figure}

\begin{figure}[H]

\begin{subfigure}{0.50\textwidth}
\includegraphics[width=0.9\linewidth, height=5cm]{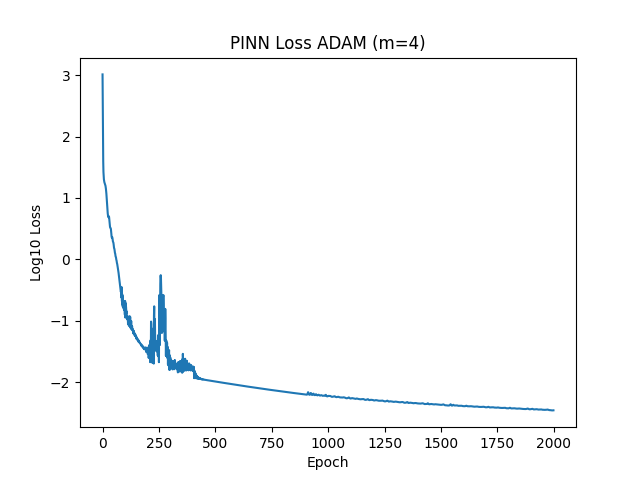} 
\caption{PINN ADAM loss}
\label{fig:subim1}
\end{subfigure}
\begin{subfigure}{0.50\textwidth}
\includegraphics[width=0.9\linewidth, height=5cm]{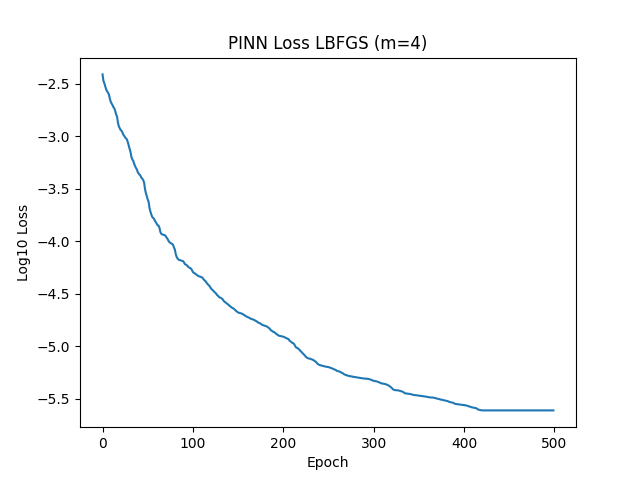}
\caption{PINN LBFGS loss }
\label{fig:subim2}
\end{subfigure}

\end{figure}

\newpage

\textbf{m = 5}

\begin{figure}[H]
    \centering
    
    \begin{subfigure}[b]{0.32\textwidth}
        \centering
        \includegraphics[width=\textwidth, height =6.5 cm]{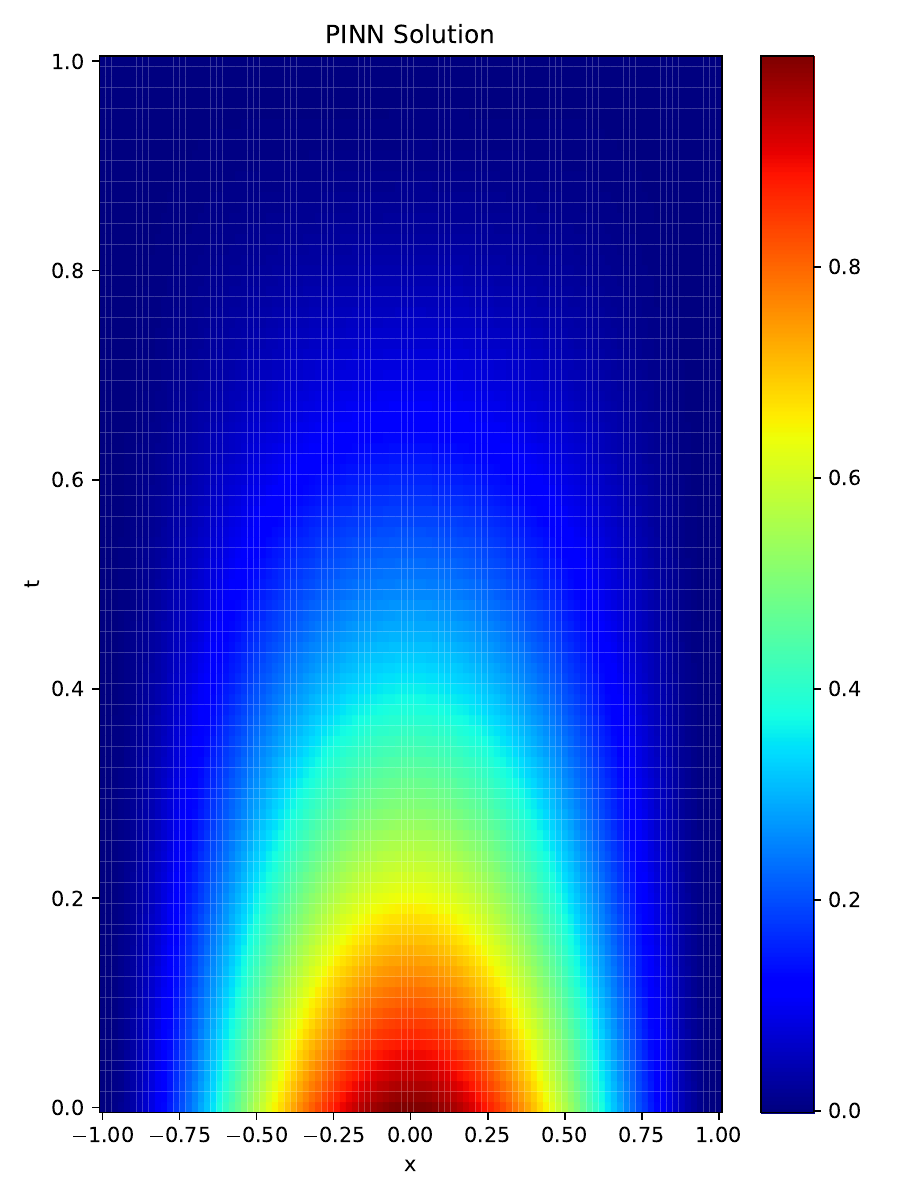}
        \caption{PINN solution}
    \end{subfigure}
    \hfill
    \begin{subfigure}[b]{0.32\textwidth}
        \centering
        \includegraphics[width=\textwidth, height = 6.5 cm]{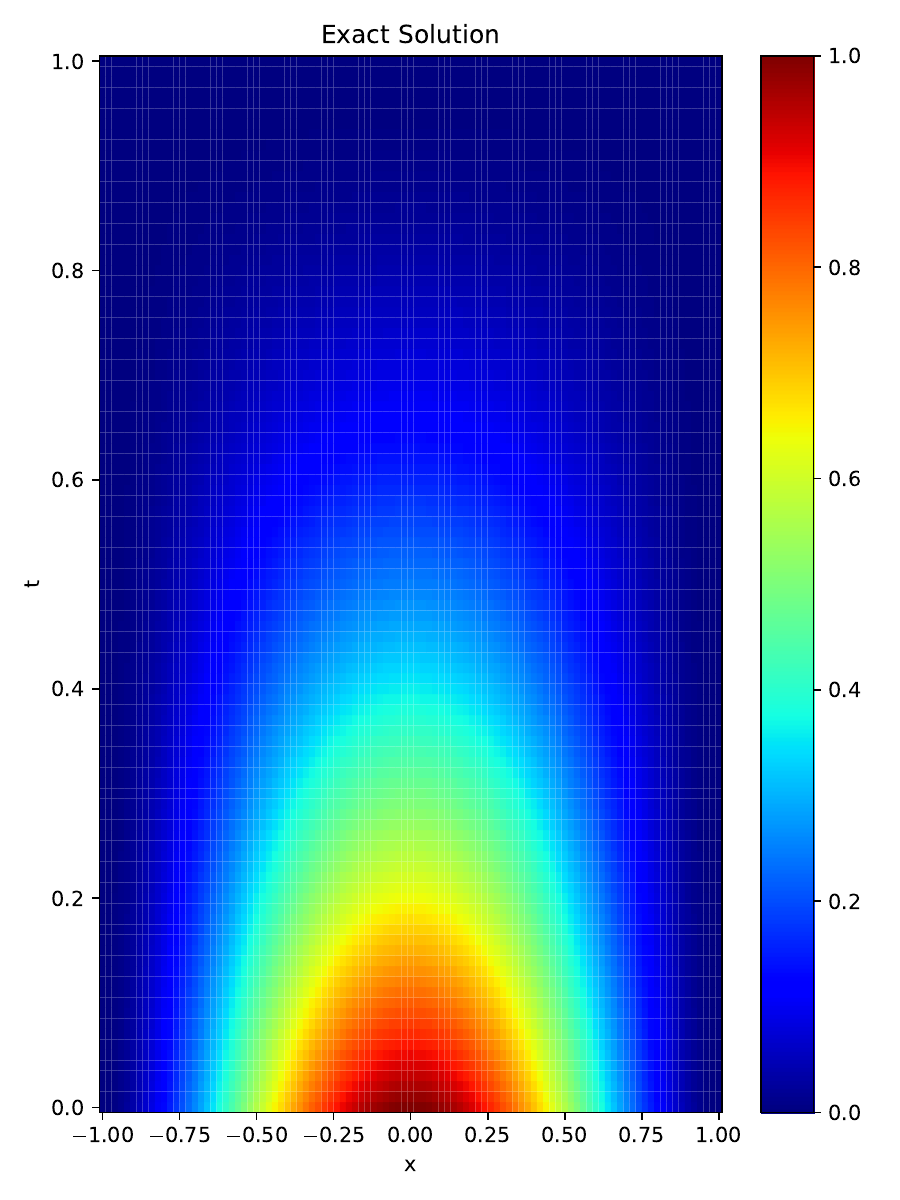}
        \caption{Exact solution}
    \end{subfigure}
    \hfill
    \begin{subfigure}[b]{0.32\textwidth}
        \centering
        \includegraphics[width=\textwidth, height = 6.5 cm ]{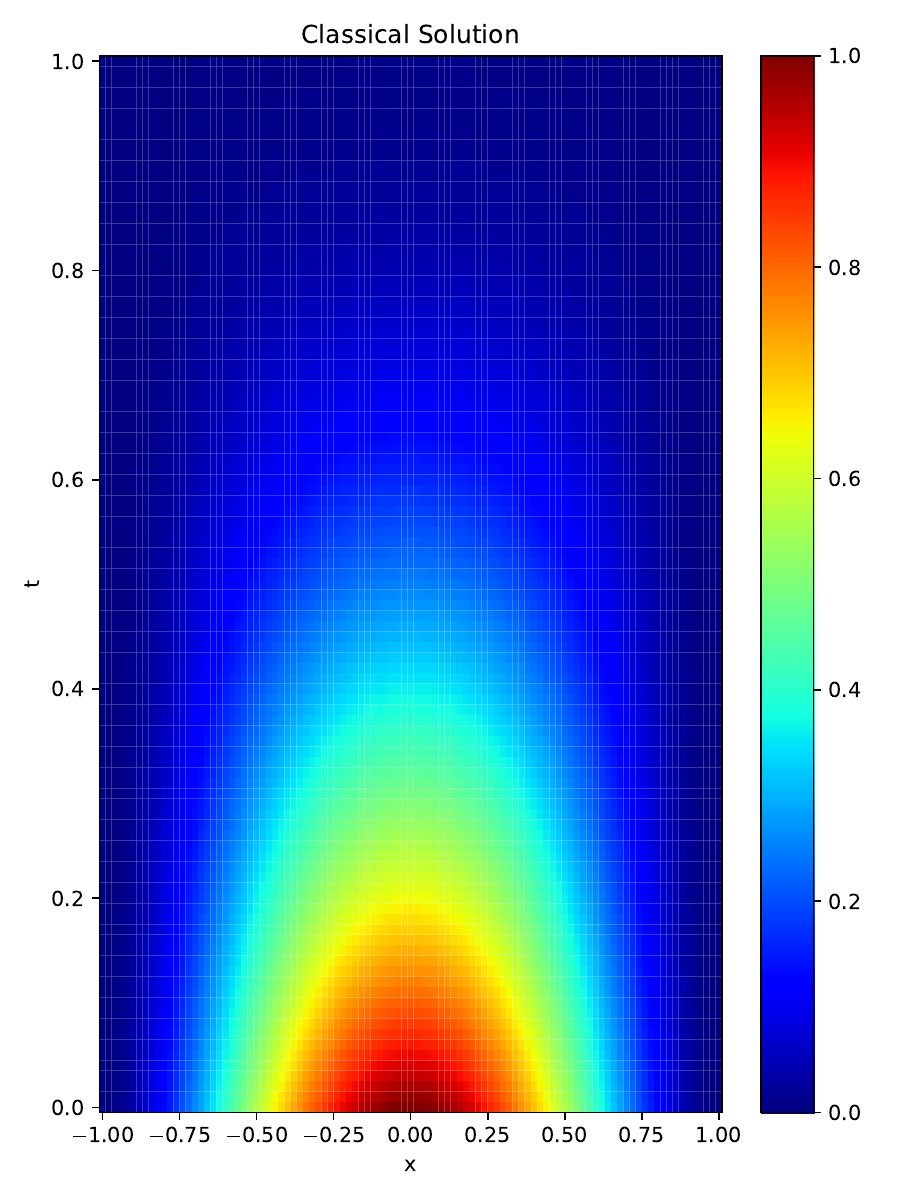}
        \caption{Classical solution}
    \end{subfigure}
\end{figure}\vspace{-4mm}
\begin{figure}[H]
\begin{subfigure}{0.50\textwidth}
\includegraphics[width=0.9\linewidth, height=5cm]{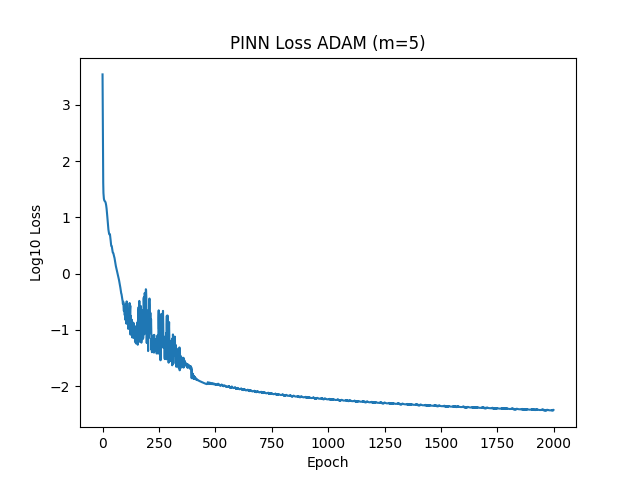} 
\caption{PINN ADAM loss}
\label{fig:subim1}
\end{subfigure}
\begin{subfigure}{0.50\textwidth}
\includegraphics[width=0.9\linewidth, height=5cm]{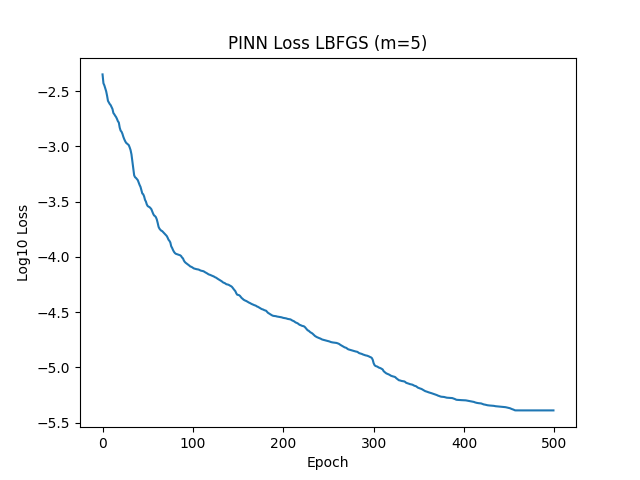}
\caption{PINN LBFGS loss }
\label{fig:subim2}
\end{subfigure}

\end{figure}

\subsection{Damped Harmonic Solution}

\vspace{3mm}
\textbf{m = 2}

\begin{figure}[H]
    \centering
    
    \begin{subfigure}[b]{0.32\textwidth}
        \centering
        \includegraphics[width=\textwidth, height = 6.5 cm]{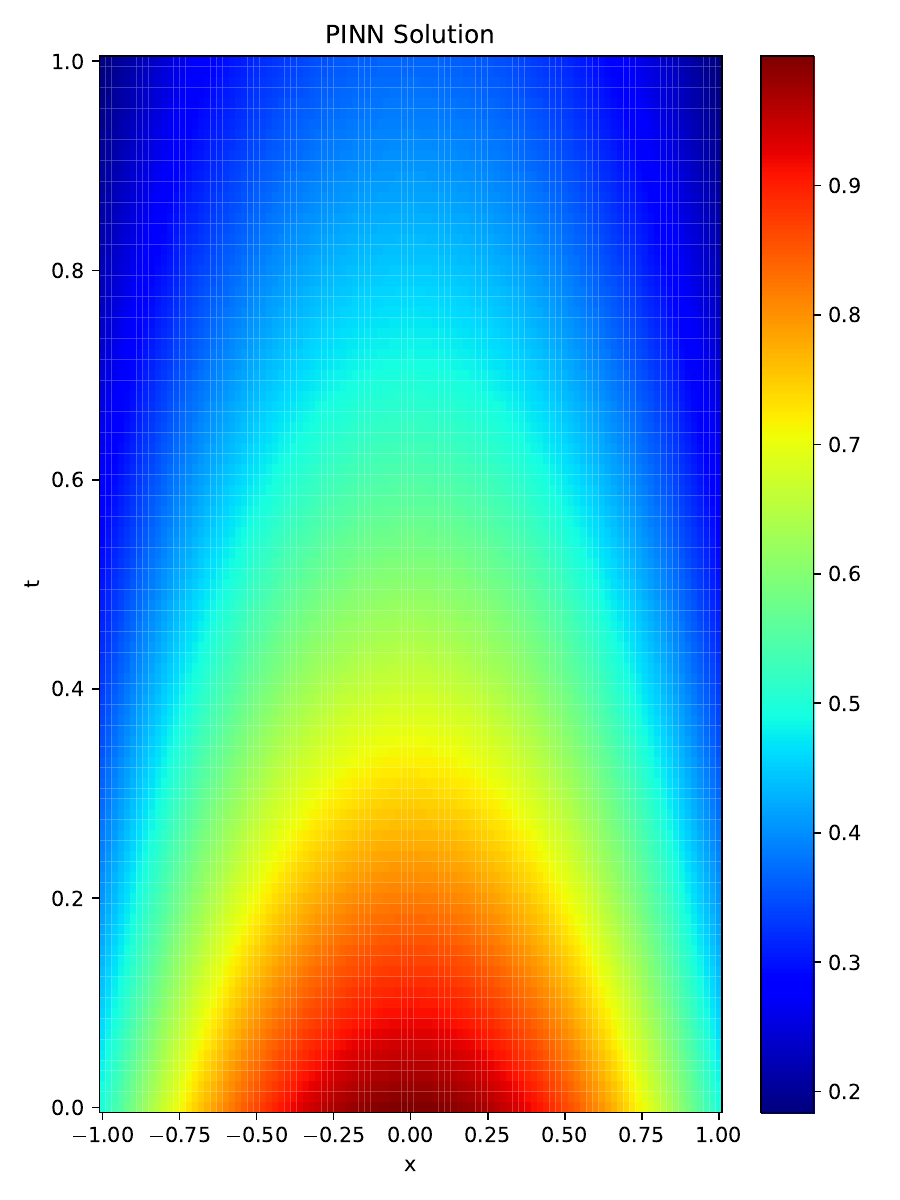}
        \caption{PINN solution}
    \end{subfigure}
    \hfill
    \begin{subfigure}[b]{0.32\textwidth}
        \centering
        \includegraphics[width=\textwidth, height = 6.5 cm]{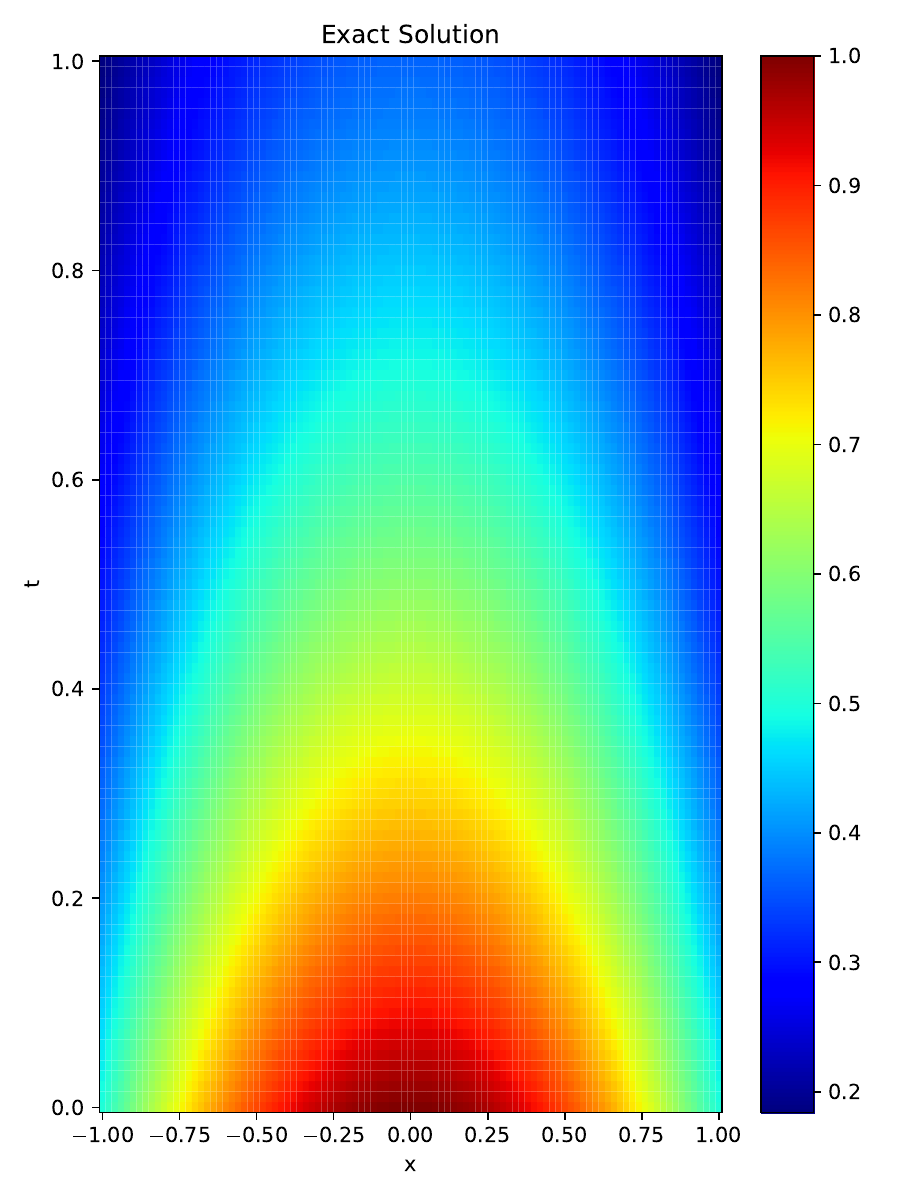}
        \caption{Exact solution}
    \end{subfigure}
    \hfill
    \begin{subfigure}[b]{0.32\textwidth}
        \centering
        \includegraphics[width=\textwidth, height = 6.5 cm]{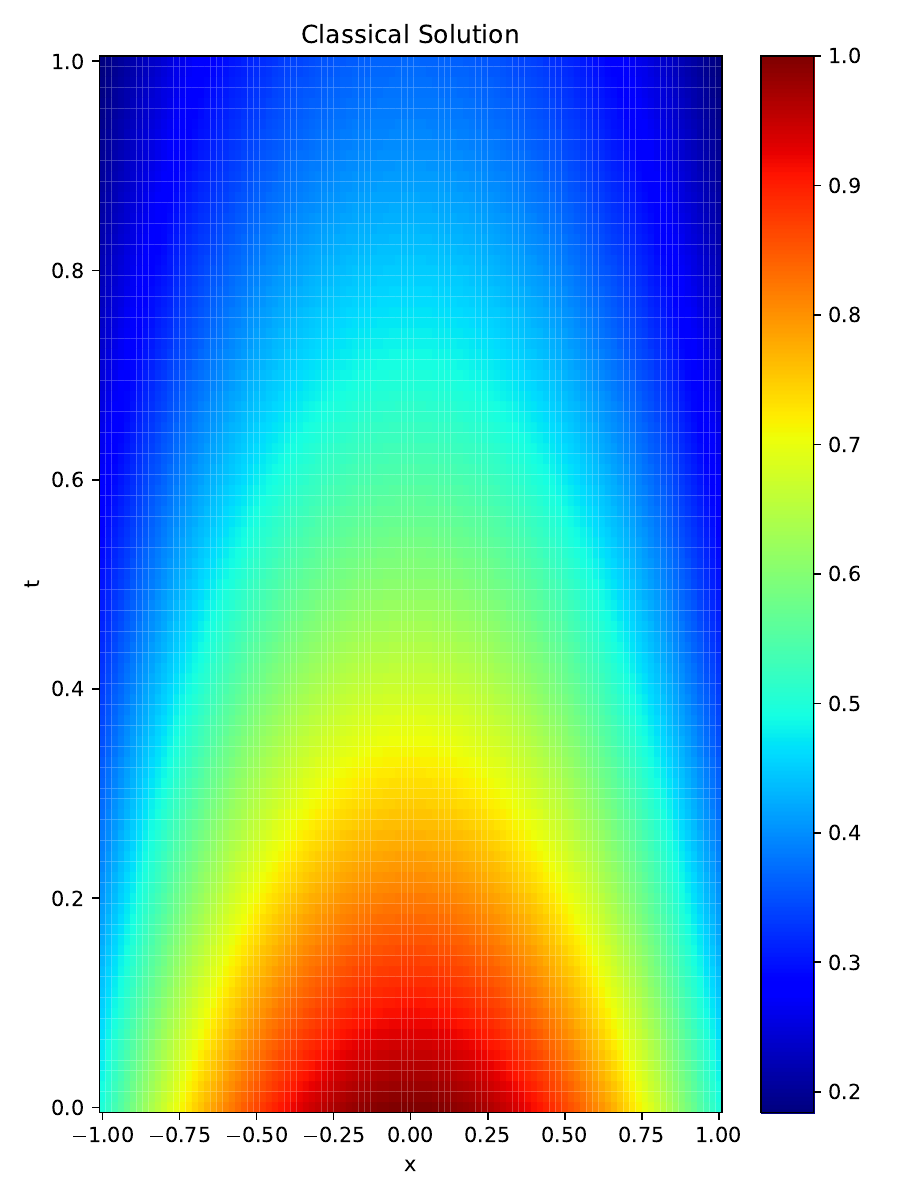}
        \caption{Classical solution}
    \end{subfigure}
\end{figure}

\begin{figure}[H]

\begin{subfigure}{0.50\textwidth}
\includegraphics[width=0.9\linewidth, height=5cm]{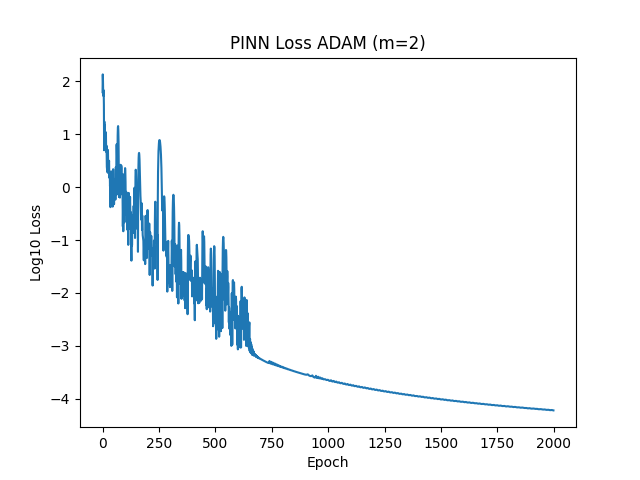} 
\caption{PINN ADAM loss}
\label{fig:subim1}
\end{subfigure}
\begin{subfigure}{0.50\textwidth}
\includegraphics[width=0.9\linewidth, height=5cm]{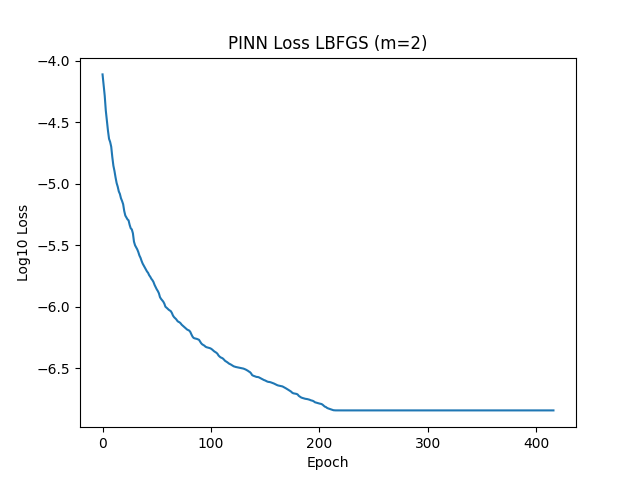}
\caption{PINN LBFGS loss }
\label{fig:subim2}
\end{subfigure}

\end{figure}

\hrule
\vspace{3mm}
\textbf{m = 3}

\begin{figure}[H]
    \centering
    
    \begin{subfigure}[b]{0.32\textwidth}
        \centering
        \includegraphics[width=\textwidth, height = 6.5 cm]{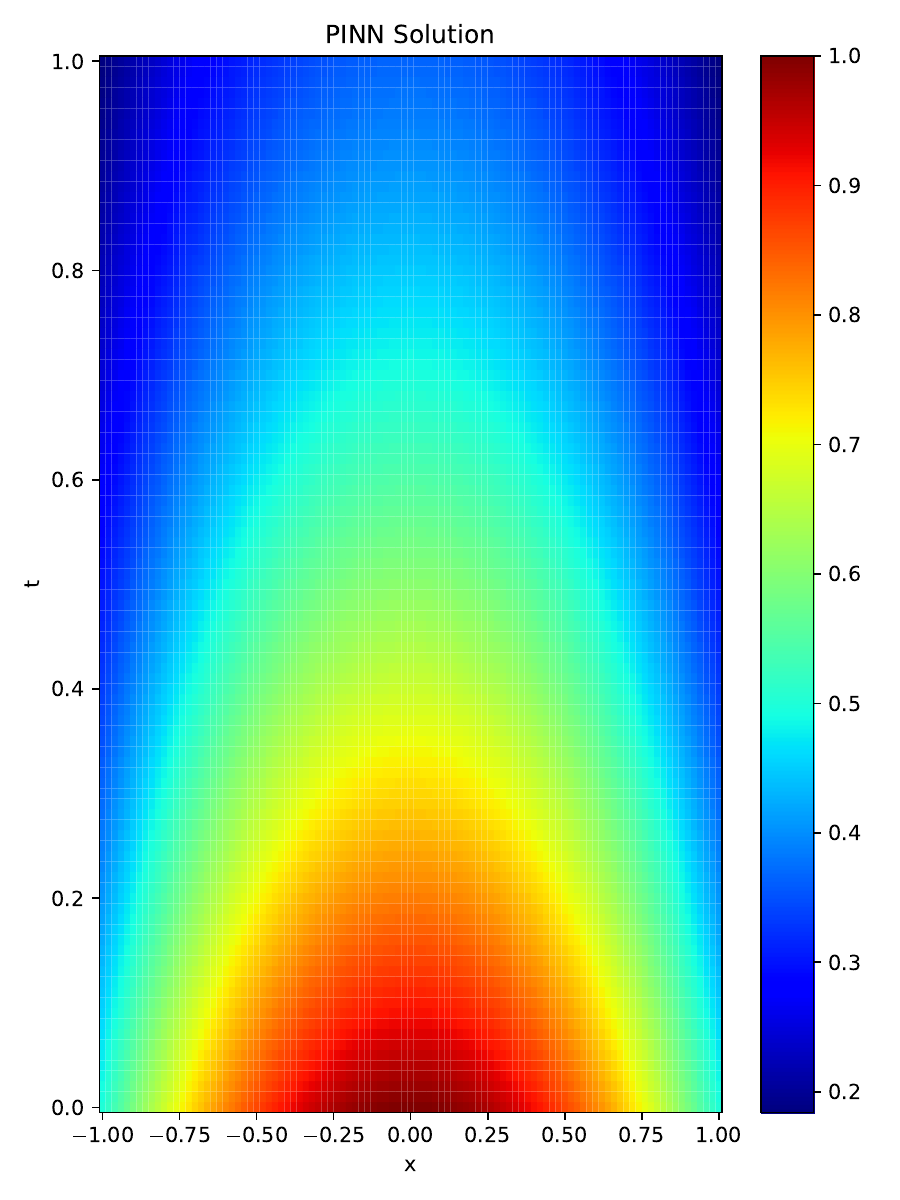}
        \caption{PINN solution}
    \end{subfigure}
    \hfill
    \begin{subfigure}[b]{0.32\textwidth}
        \centering
        \includegraphics[width=\textwidth, height = 6.5 cm]{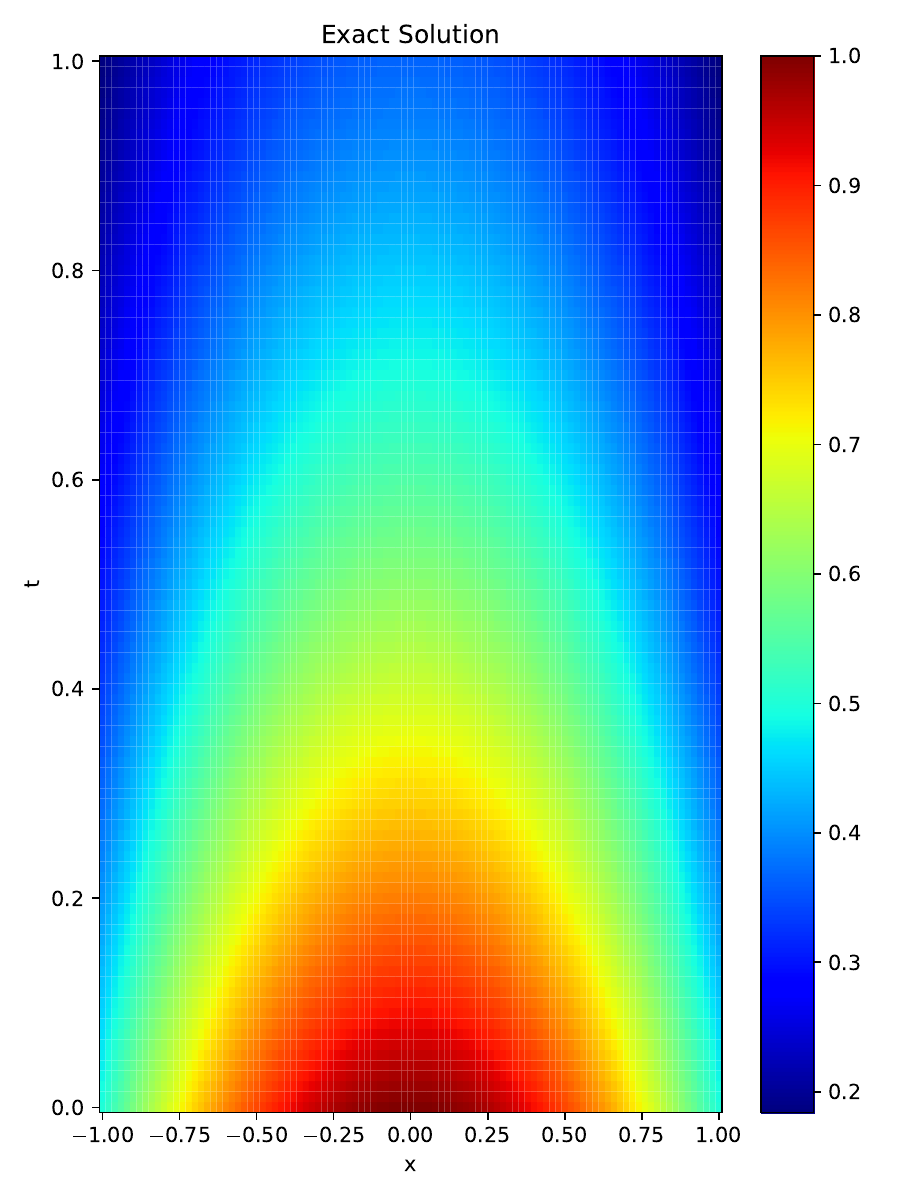}
        \caption{Exact solution}
    \end{subfigure}
    \hfill
    \begin{subfigure}[b]{0.32\textwidth}
        \centering
        \includegraphics[width=\textwidth, height = 6.5 cm]{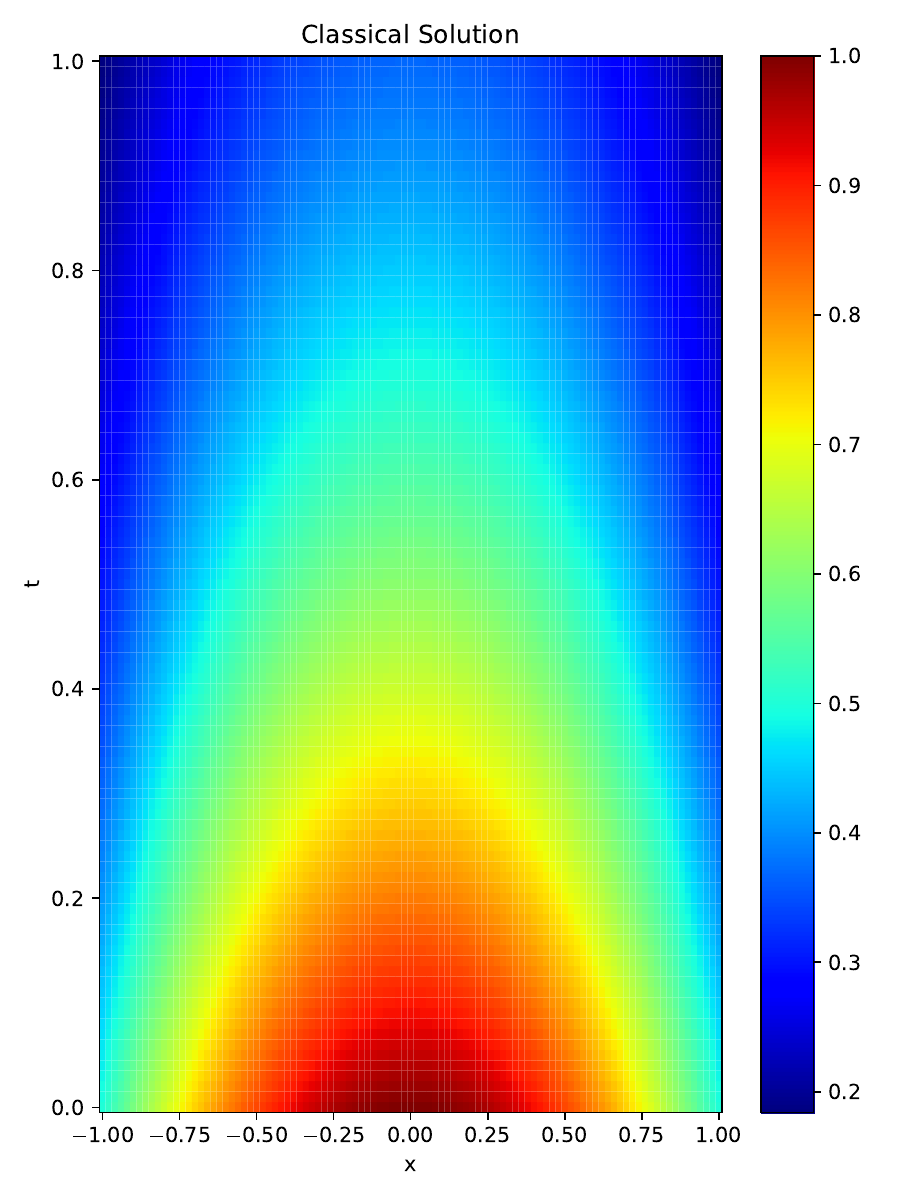}
        \caption{Classical solution}
    \end{subfigure}
\end{figure}

\begin{figure}[H]

\begin{subfigure}{0.50\textwidth}
\includegraphics[width=0.9\linewidth, height=5cm]{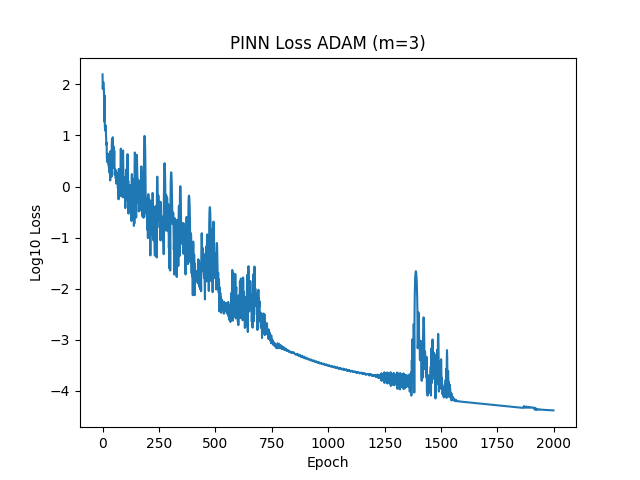} 
\caption{PINN ADAM loss}
\label{fig:subim1}
\end{subfigure}
\begin{subfigure}{0.50\textwidth}
\includegraphics[width=0.9\linewidth, height=5cm]{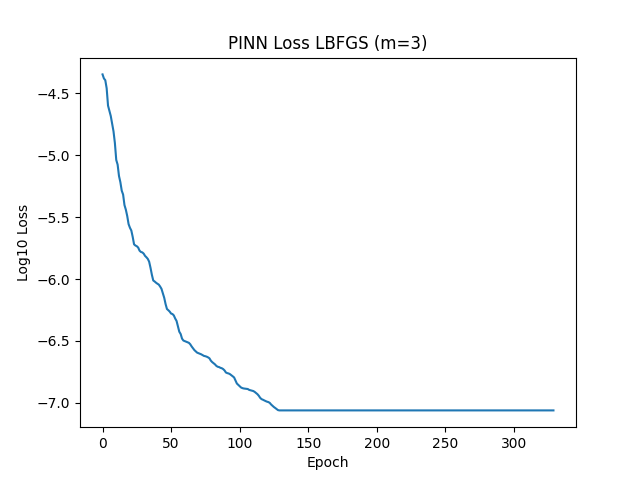}
\caption{PINN LBFGS loss }
\label{fig:subim2}
\end{subfigure}

\end{figure}

\newpage
\textbf{m = 4}
\begin{figure}[H]
    \centering
    
    \begin{subfigure}[b]{0.32\textwidth}
        \centering
        \includegraphics[width=\textwidth, height =6.5 cm]{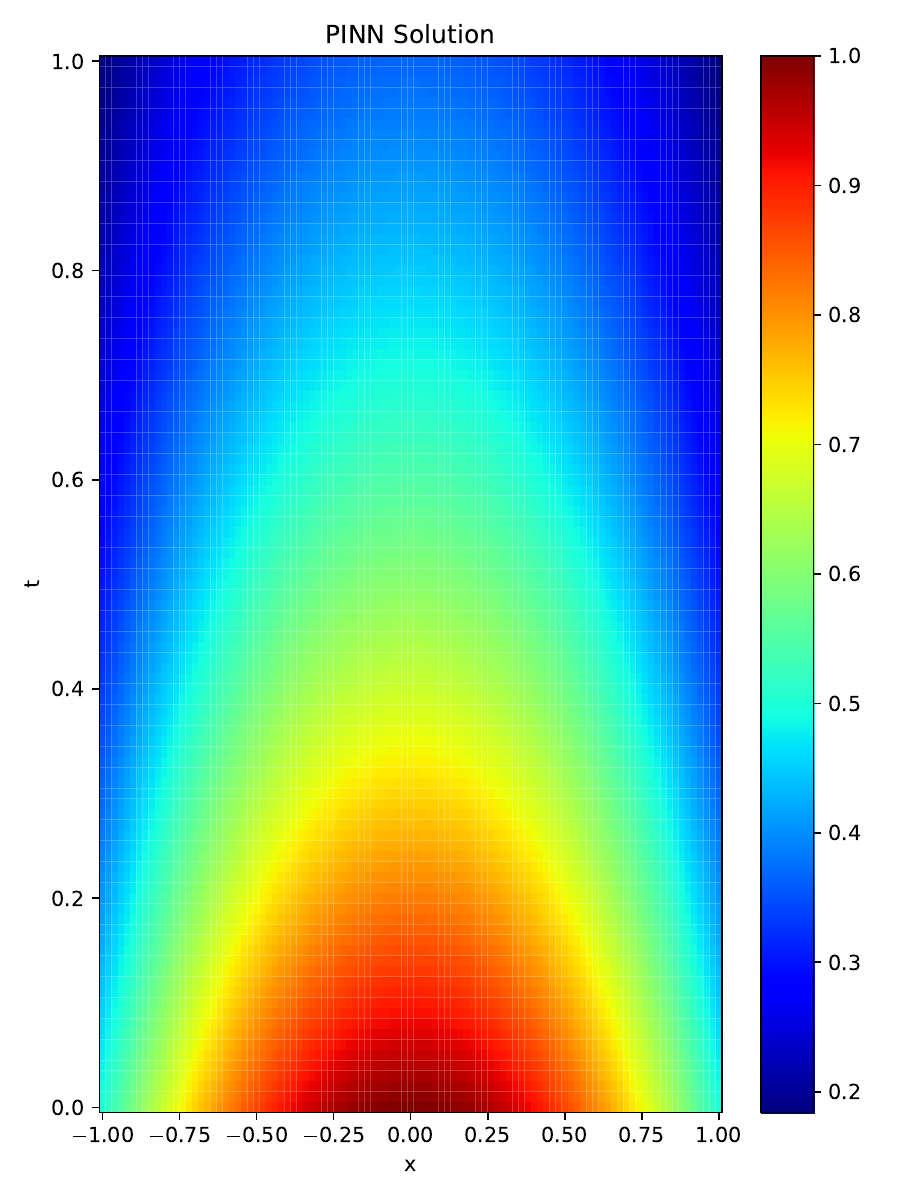}
        \caption{PINN solution}
    \end{subfigure}
    \hfill
    \begin{subfigure}[b]{0.32\textwidth}
        \centering
        \includegraphics[width=\textwidth, height = 6.5 cm]{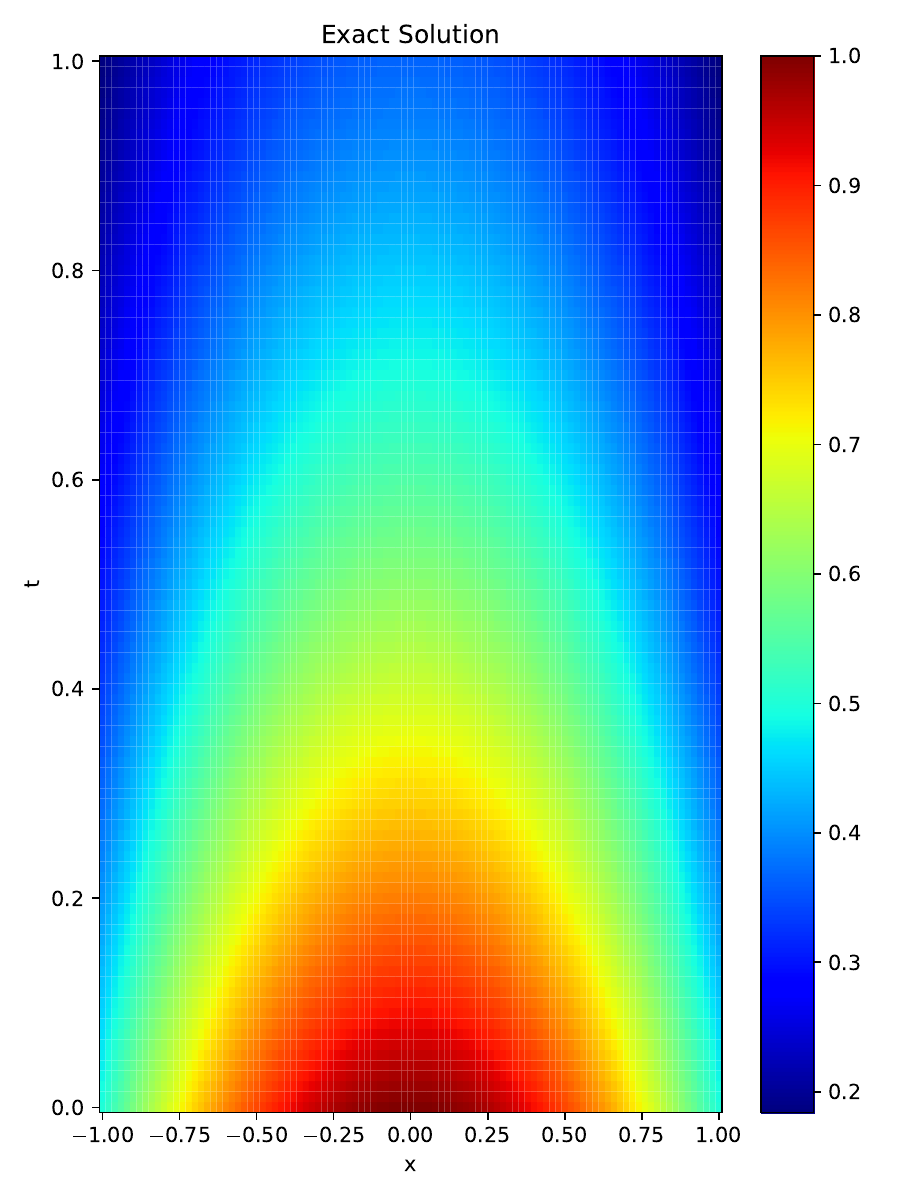}
        \caption{Exact solution}
    \end{subfigure}
    \hfill
    \begin{subfigure}[b]{0.32\textwidth}
        \centering
        \includegraphics[width=\textwidth, height = 6.5 cm ]{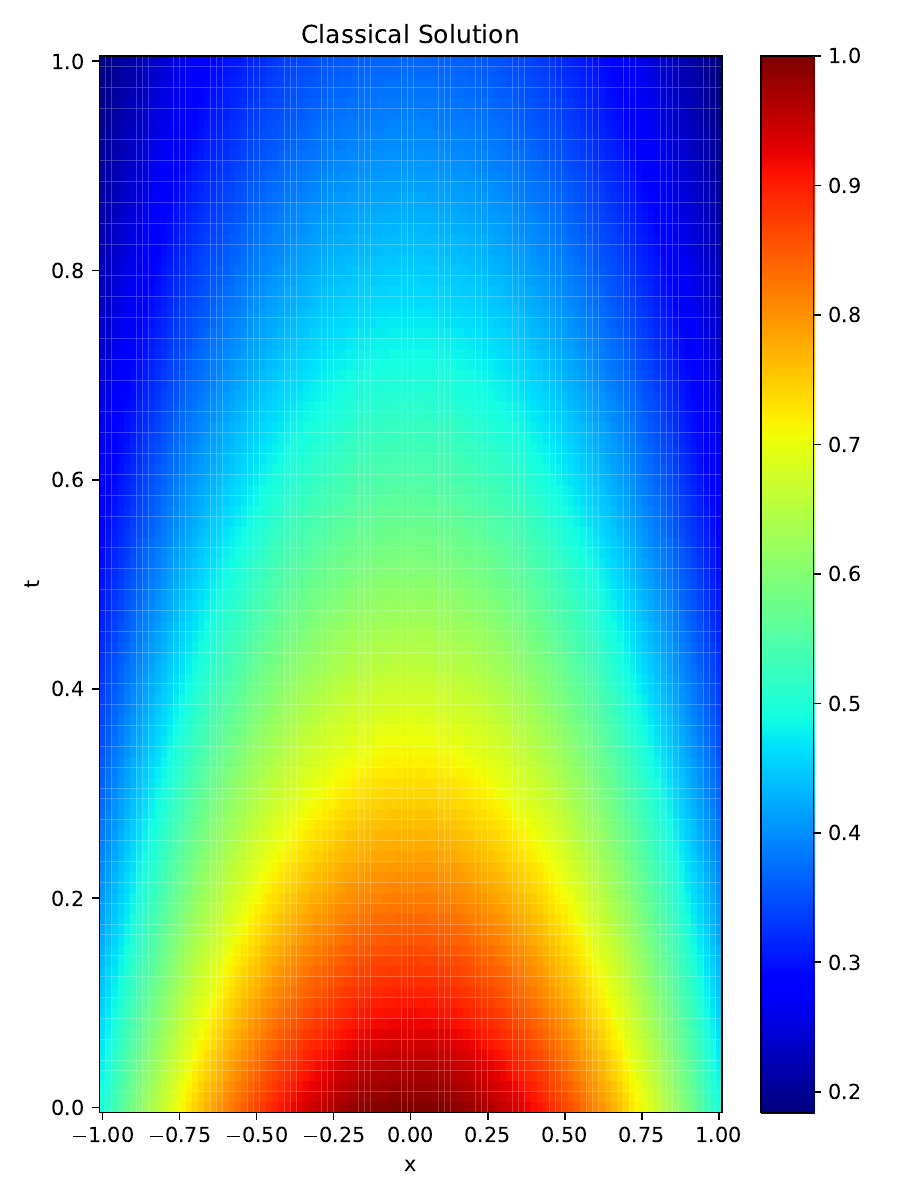}
        \caption{Classical solution}
    \end{subfigure}
\end{figure}

\begin{figure}[H]

\begin{subfigure}{0.50\textwidth}
\includegraphics[width=0.9\linewidth, height=5cm]{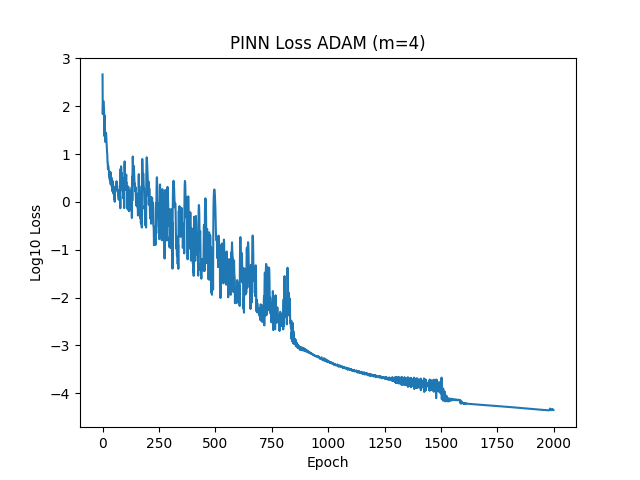} 
\caption{PINN ADAM loss}
\label{fig:subim1}
\end{subfigure}
\begin{subfigure}{0.50\textwidth}
\includegraphics[width=0.9\linewidth, height=5cm]{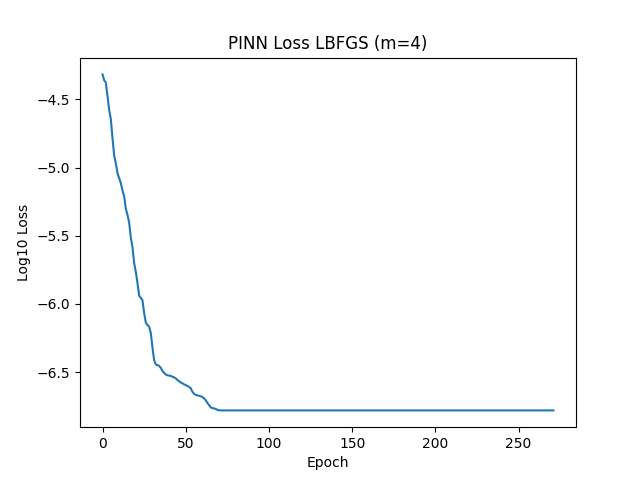}
\caption{PINN LBFGS loss }
\label{fig:subim2}
\end{subfigure}

\end{figure}

\hrule
\vspace{3mm}

\textbf{m = 5}

\begin{figure}[H]
    \centering
    
    \begin{subfigure}[b]{0.32\textwidth}
        \centering
        \includegraphics[width=\textwidth, height =6.5 cm]{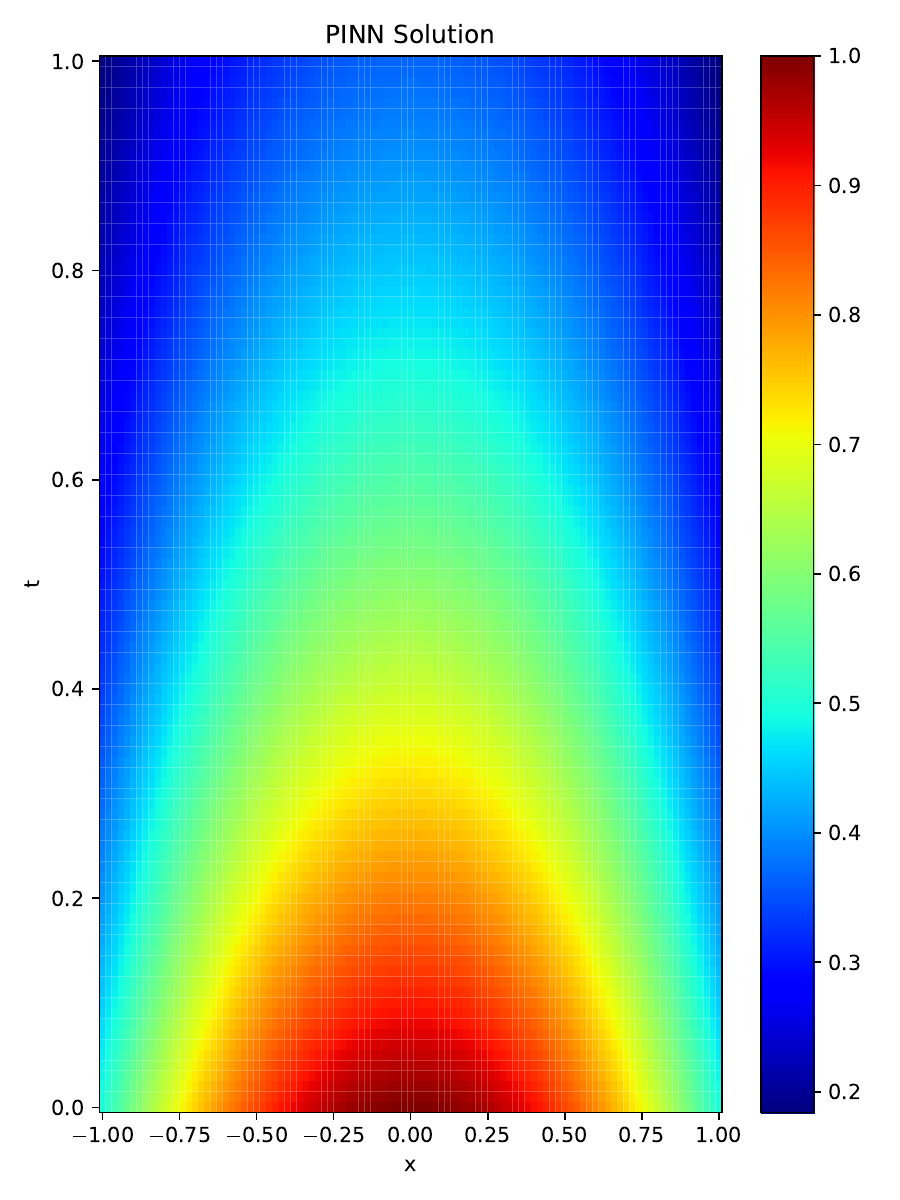}
        \caption{PINN solution}
    \end{subfigure}
    \hfill
    \begin{subfigure}[b]{0.32\textwidth}
        \centering
        \includegraphics[width=\textwidth, height = 6.5 cm]{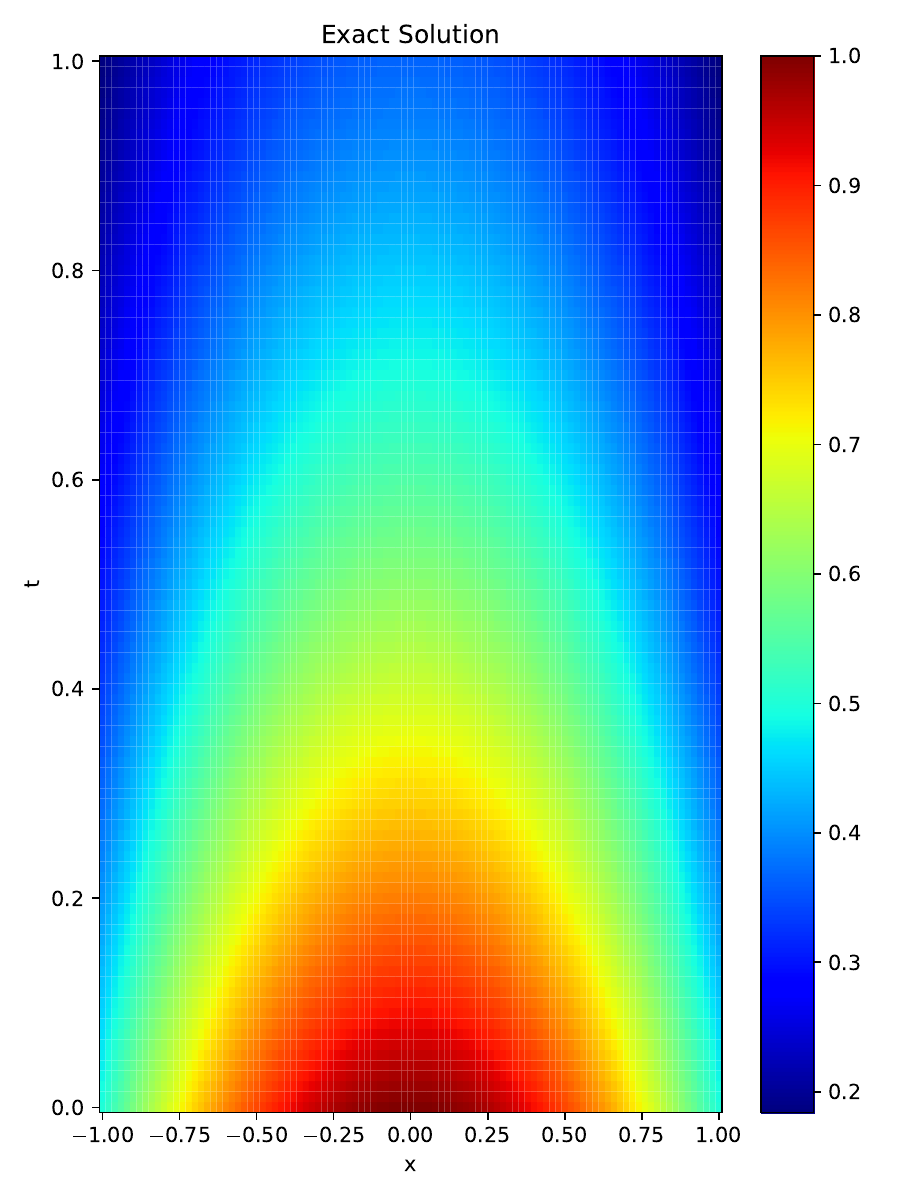}
        \caption{Exact solution}
    \end{subfigure}
    \hfill
    \begin{subfigure}[b]{0.32\textwidth}
        \centering
        \includegraphics[width=\textwidth, height = 6.5 cm ]{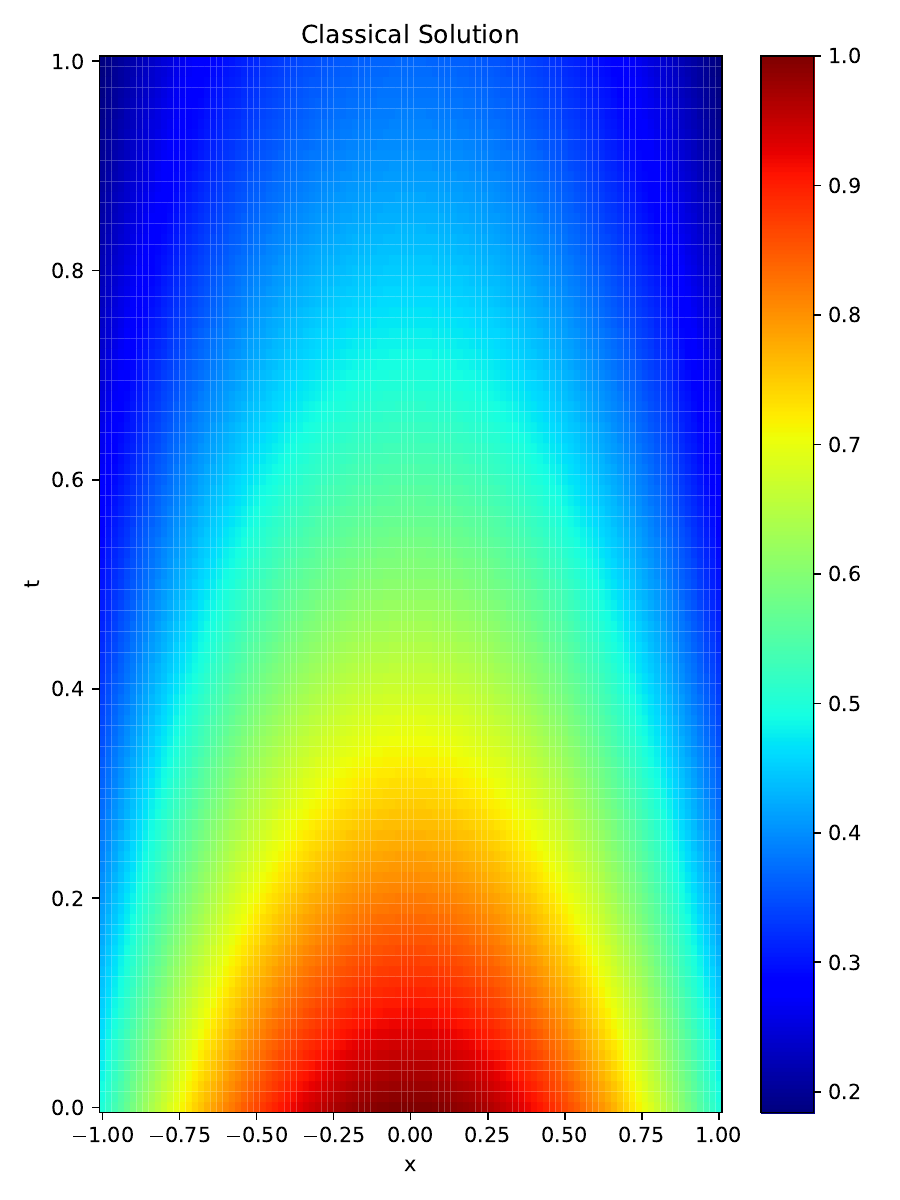}
        \caption{Classical solution}
    \end{subfigure}
\end{figure}

\begin{figure}[H]

\begin{subfigure}{0.50\textwidth}
\includegraphics[width=0.9\linewidth, height=5cm]{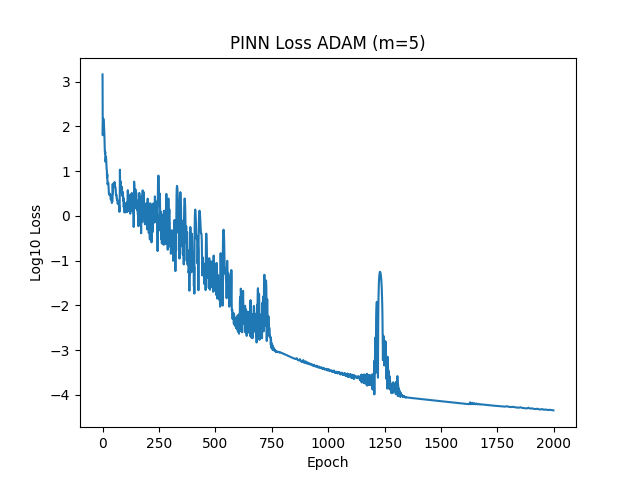} 
\caption{PINN ADAM loss}
\label{fig:subim1}
\end{subfigure}
\begin{subfigure}{0.50\textwidth}
\includegraphics[width=0.9\linewidth, height=5cm]{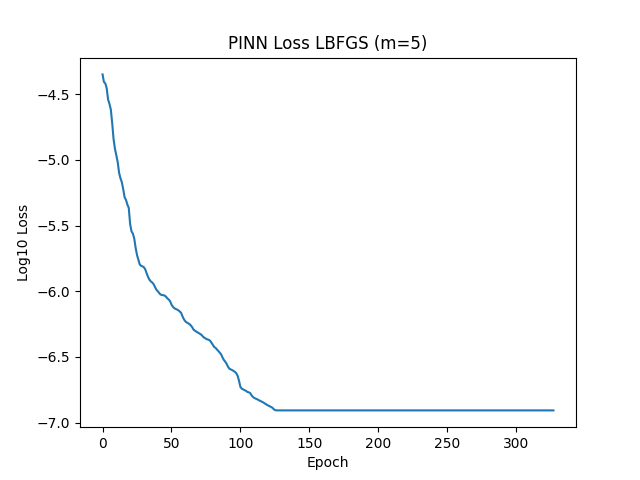}
\caption{PINN LBFGS loss }
\label{fig:subim2}
\end{subfigure}

\end{figure}

\subsection{m-dependent solution}
\vspace{3mm}
\textbf{m = 2}

\begin{figure}[H]
    \centering
    
    \begin{subfigure}[b]{0.32\textwidth}
        \centering
        \includegraphics[width=\textwidth, height = 6.5 cm]{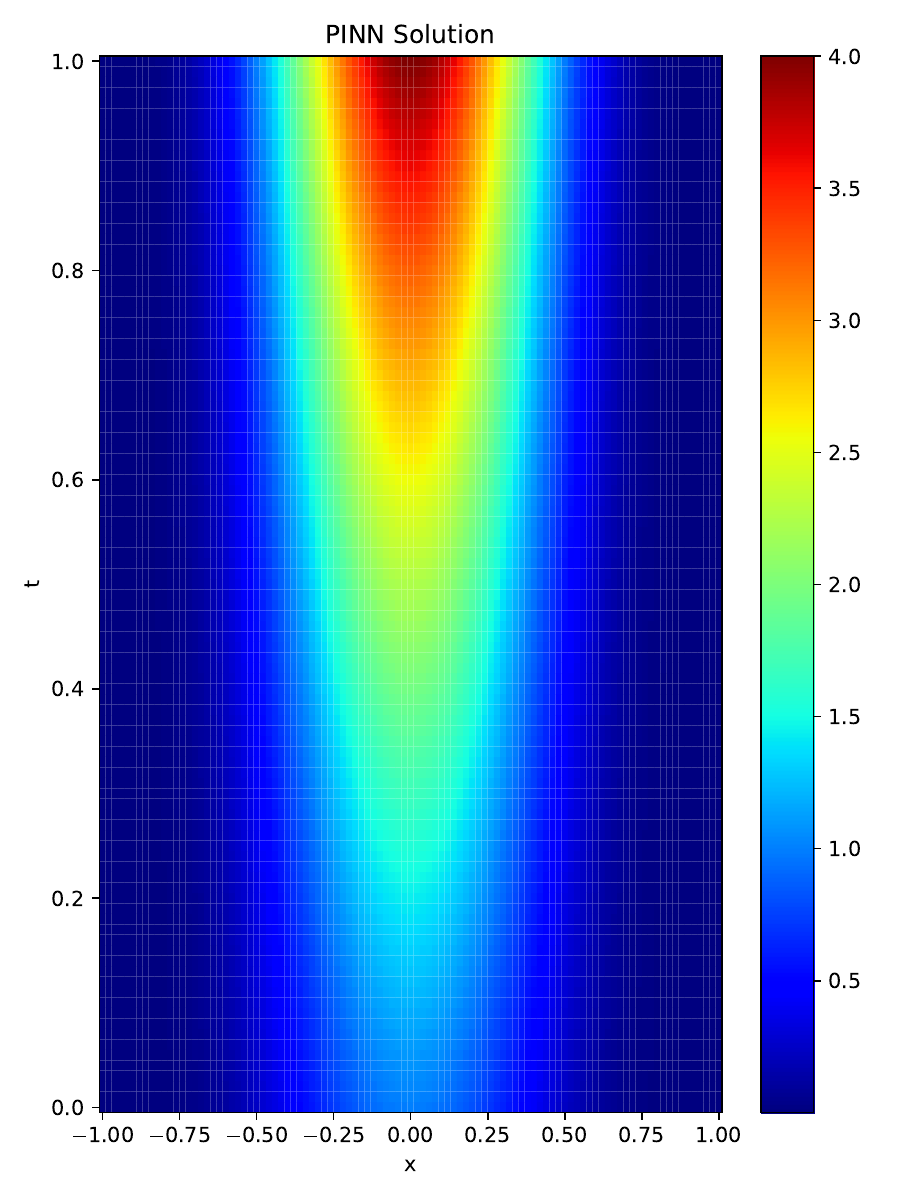}
        \caption{PINN solution}
    \end{subfigure}
    \hfill
    \begin{subfigure}[b]{0.32\textwidth}
        \centering
        \includegraphics[width=\textwidth, height = 6.5 cm]{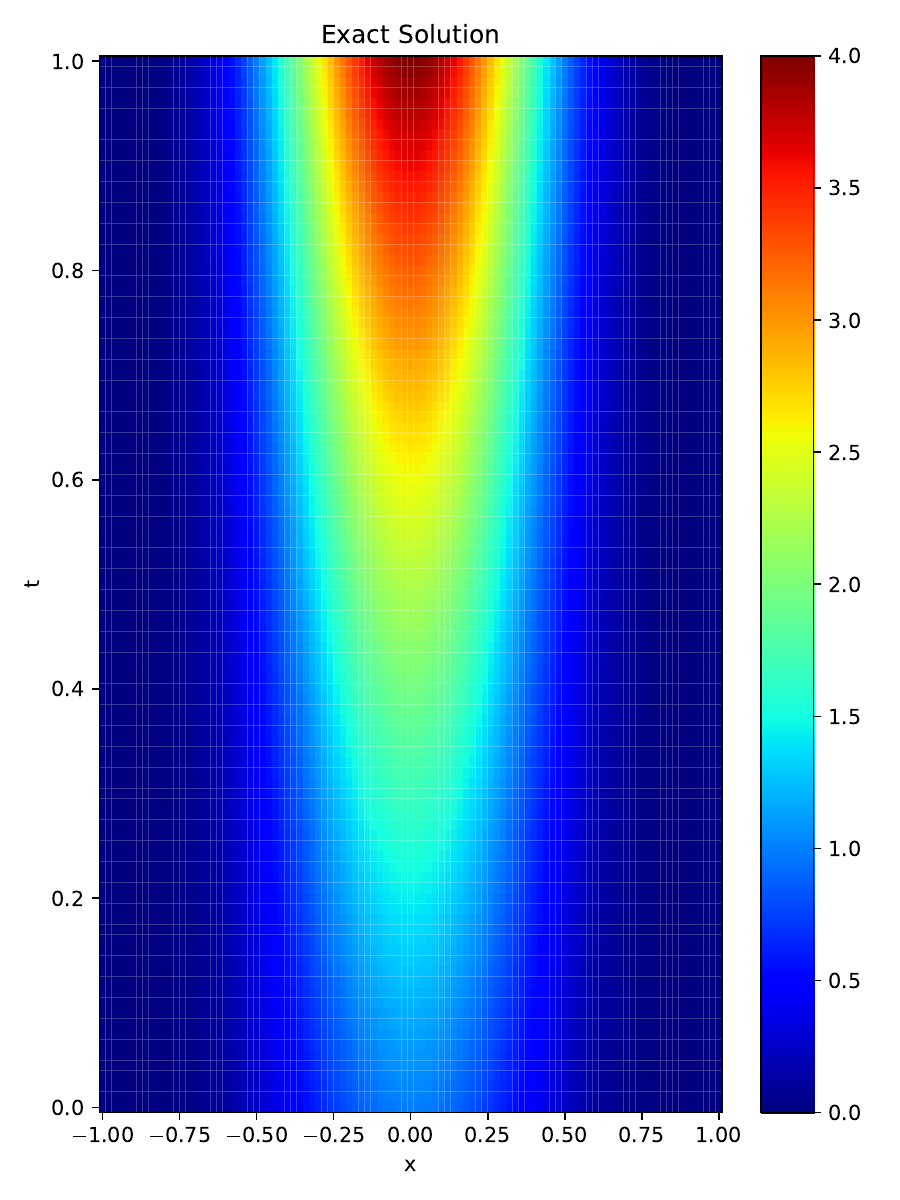}
        \caption{Exact solution}
    \end{subfigure}
    \hfill
    \begin{subfigure}[b]{0.32\textwidth}
        \centering
        \includegraphics[width=\textwidth, height = 6.5 cm]{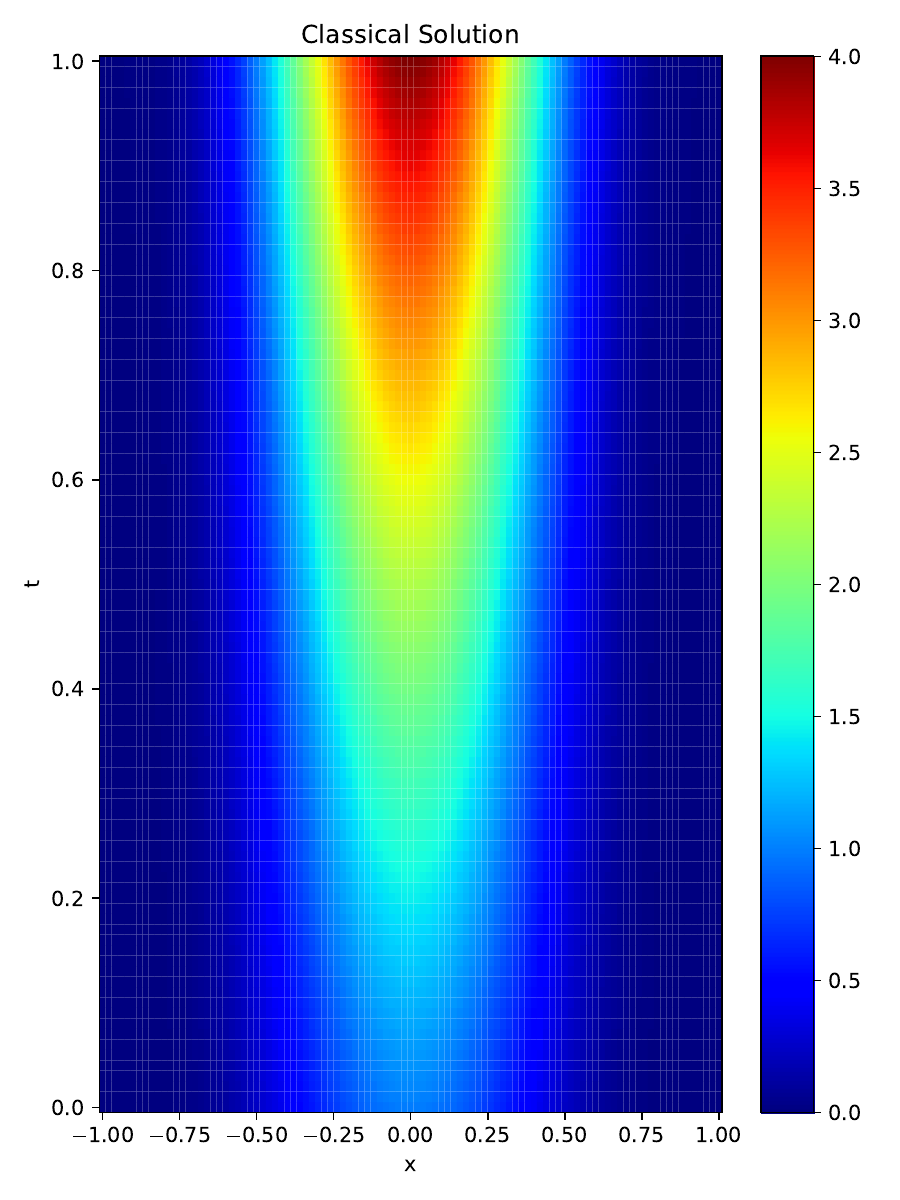}
        \caption{Classical solution}
    \end{subfigure}
\end{figure}

\begin{figure}[H]

\begin{subfigure}{0.50\textwidth}
\includegraphics[width=0.9\linewidth, height=5cm]{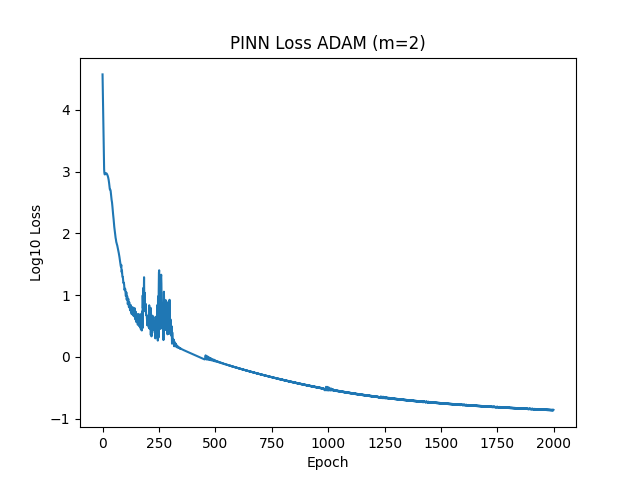} 
\caption{PINN ADAM loss}
\label{fig:subim1}
\end{subfigure}
\begin{subfigure}{0.50\textwidth}
\includegraphics[width=0.9\linewidth, height=5cm]{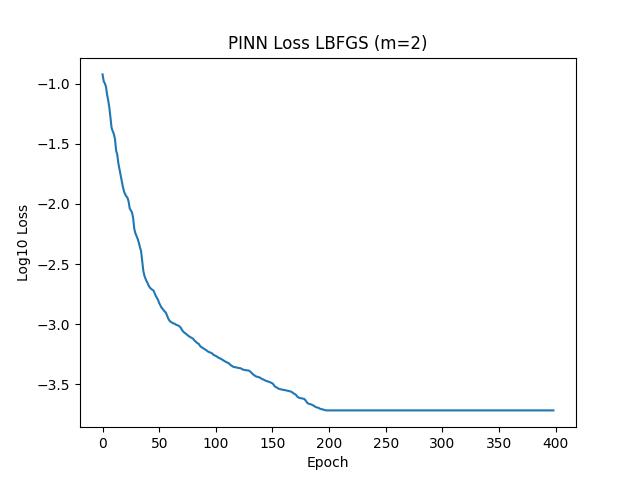}
\caption{PINN LBFGS loss }
\label{fig:subim2}
\end{subfigure}

\end{figure}

\newpage
\textbf{m = 3}

\begin{figure}[H]
    \centering
    
    \begin{subfigure}[b]{0.32\textwidth}
        \centering
        \includegraphics[width=\textwidth, height = 6.5 cm]{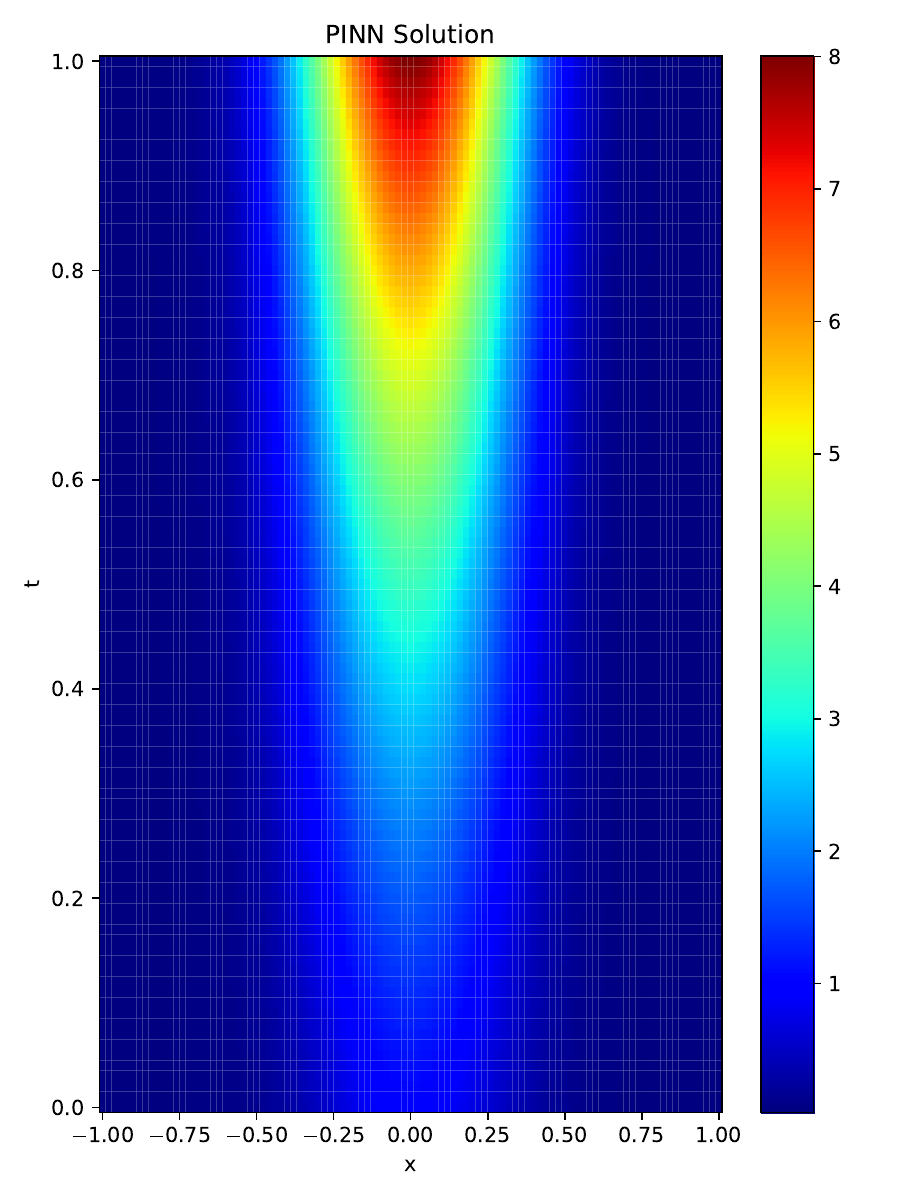}
        \caption{PINN solution}
    \end{subfigure}
    \hfill
    \begin{subfigure}[b]{0.32\textwidth}
        \centering
        \includegraphics[width=\textwidth, height = 6.5 cm]{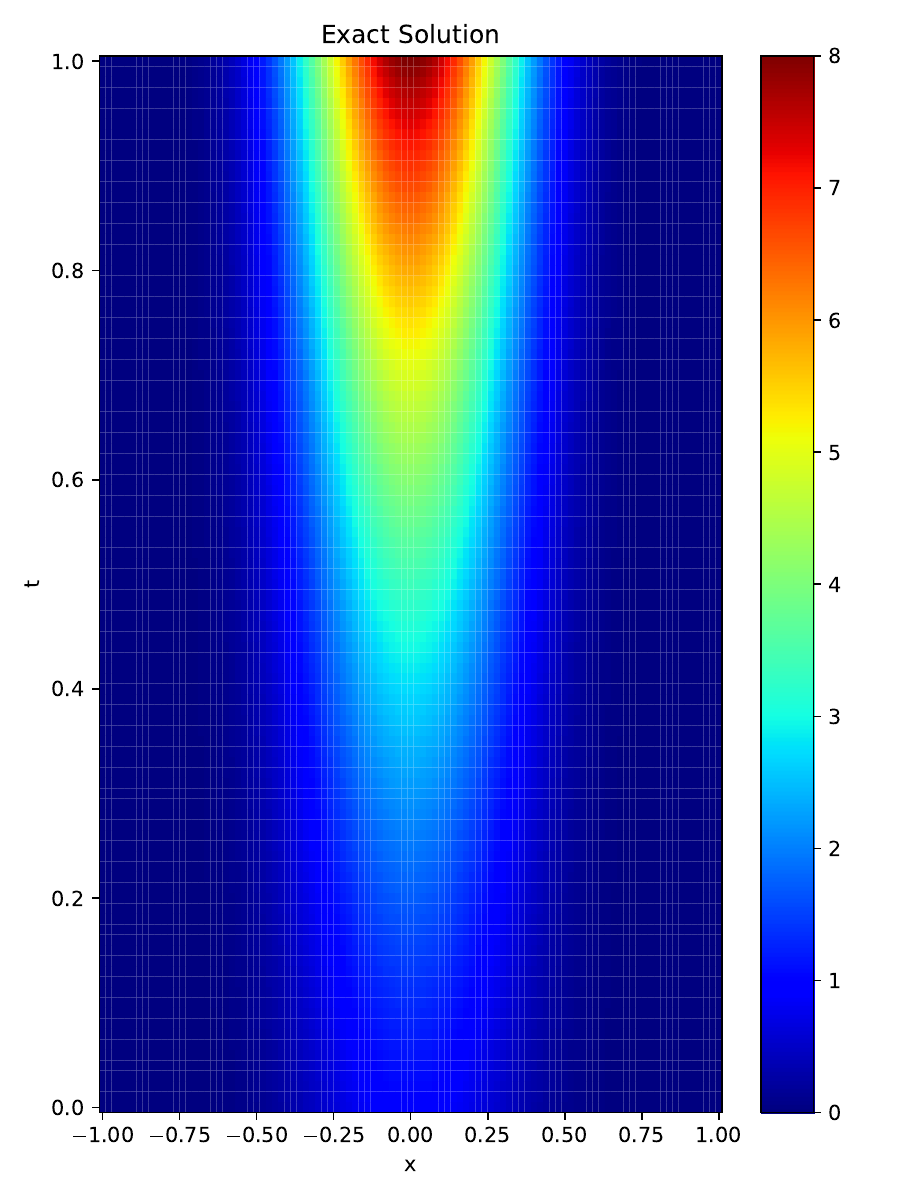}
        \caption{Exact solution}
    \end{subfigure}
    \hfill
    \begin{subfigure}[b]{0.32\textwidth}
        \centering
        \includegraphics[width=\textwidth, height = 6.5 cm]{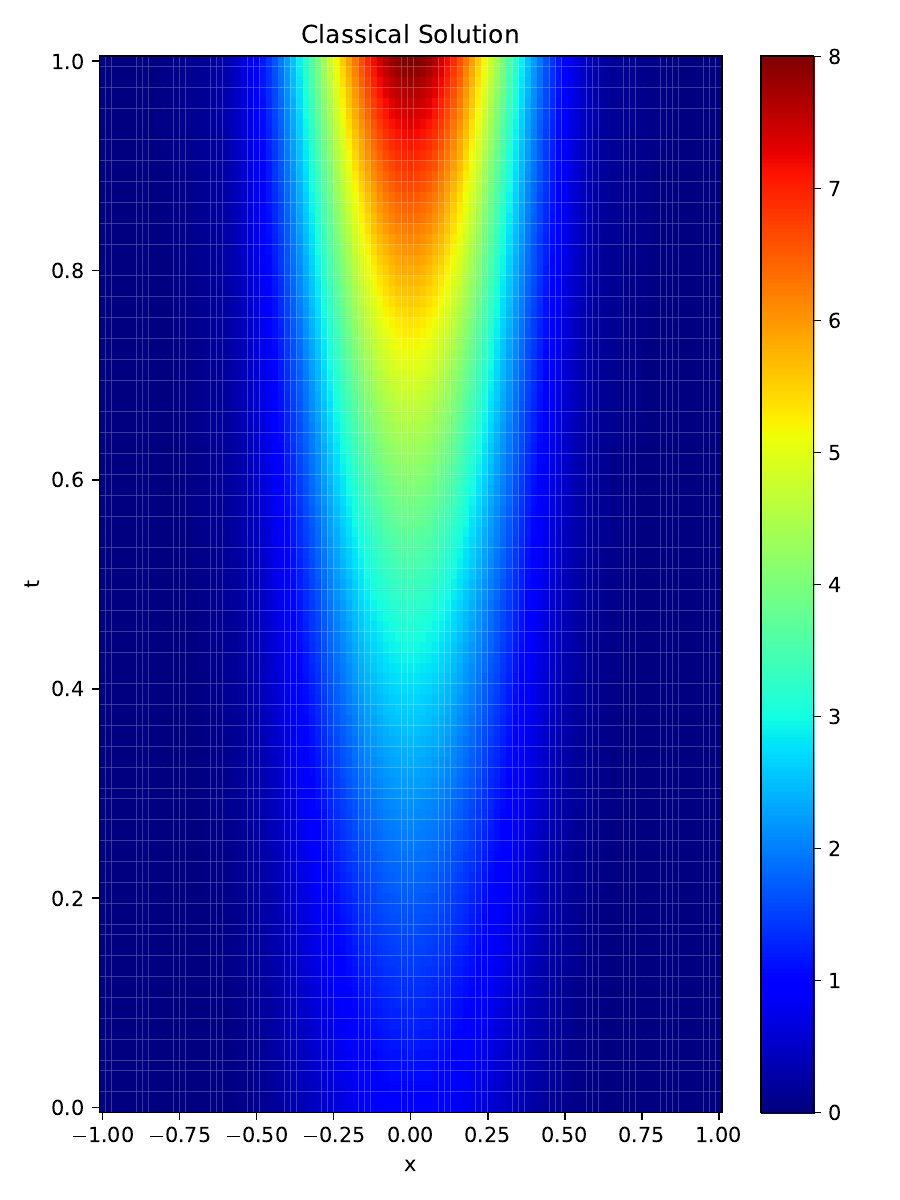}
        \caption{Classical solution}
    \end{subfigure}
\end{figure}

\begin{figure}[H]

\begin{subfigure}{0.50\textwidth}
\includegraphics[width=0.9\linewidth, height=5cm]{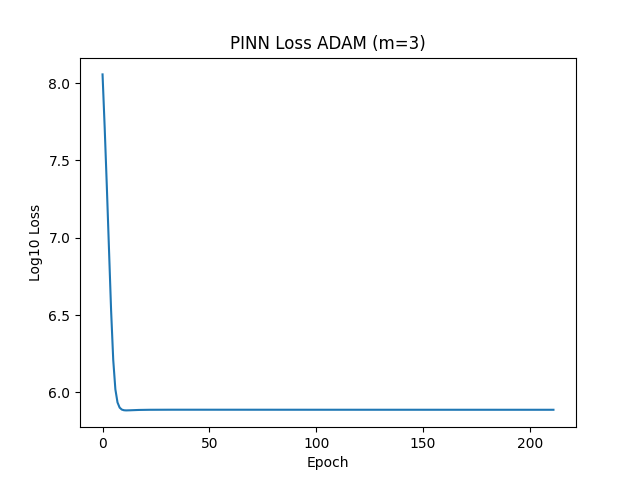} 
\caption{PINN ADAM loss}
\label{fig:subim1}
\end{subfigure}
\begin{subfigure}{0.50\textwidth}
\includegraphics[width=0.9\linewidth, height=5cm]{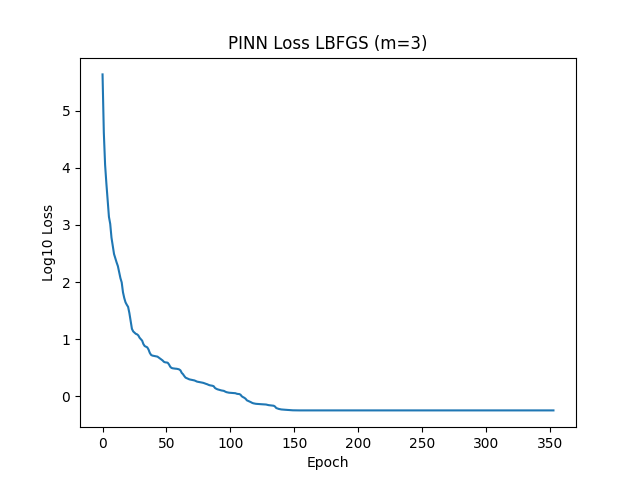}
\caption{PINN LBFGS loss }
\label{fig:subim2}
\end{subfigure}

\end{figure}

\hrule
\vspace{3mm}
\textbf{m = 4}
\begin{figure}[H]
    \centering
    
    \begin{subfigure}[b]{0.32\textwidth}
        \centering
        \includegraphics[width=\textwidth, height =6.5 cm]{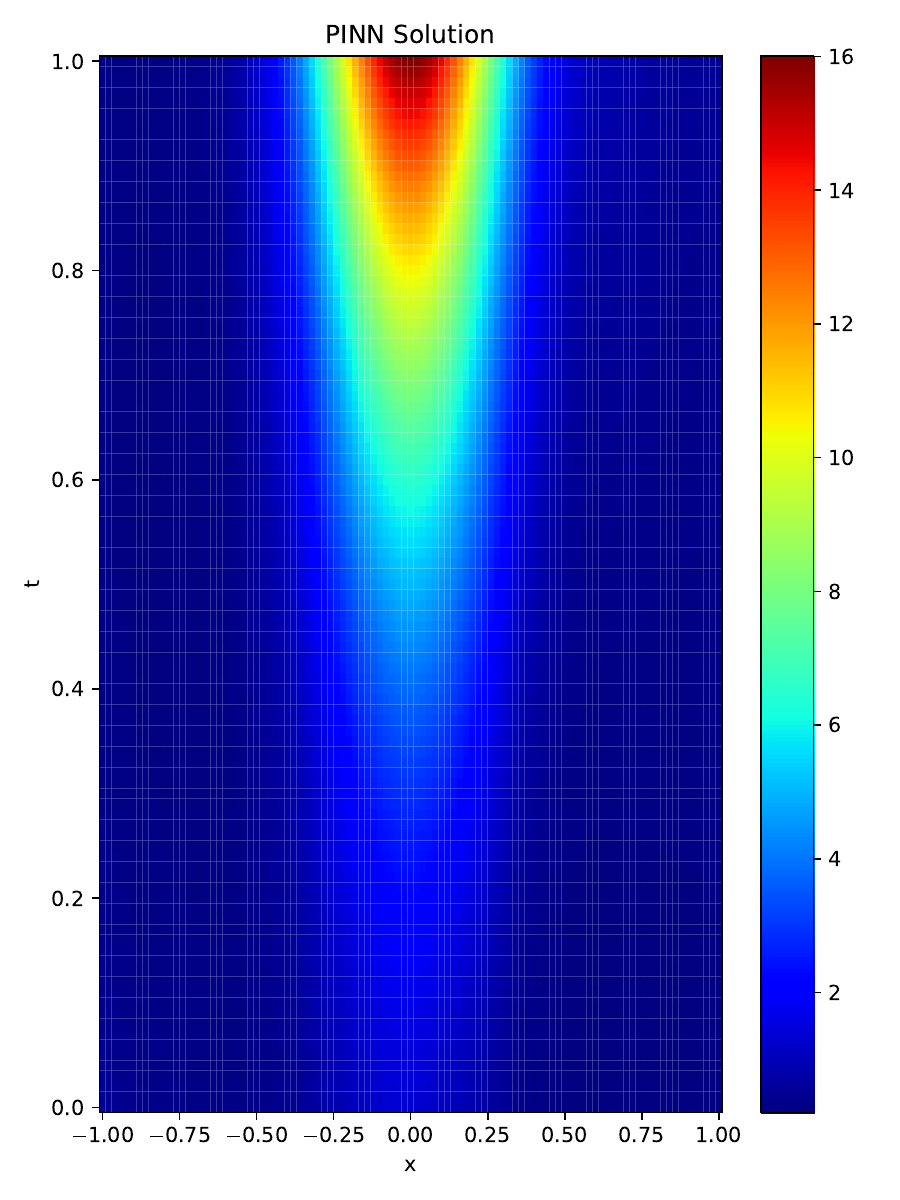}
        \caption{PINN solution}
    \end{subfigure}
    \hfill
    \begin{subfigure}[b]{0.32\textwidth}
        \centering
        \includegraphics[width=\textwidth, height = 6.5 cm]{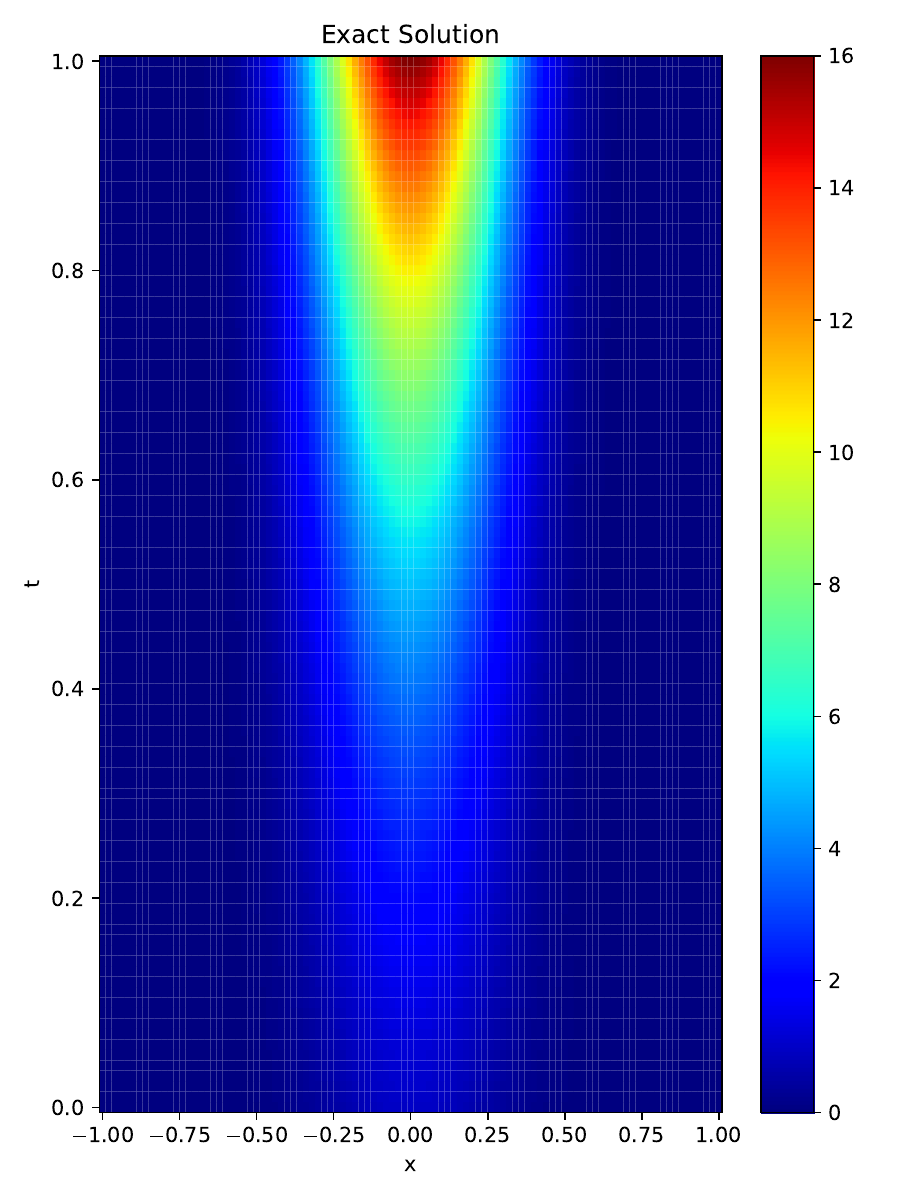}
        \caption{Exact solution}
    \end{subfigure}
    \hfill
    \begin{subfigure}[b]{0.32\textwidth}
        \centering
        \includegraphics[width=\textwidth, height = 6.5 cm ]{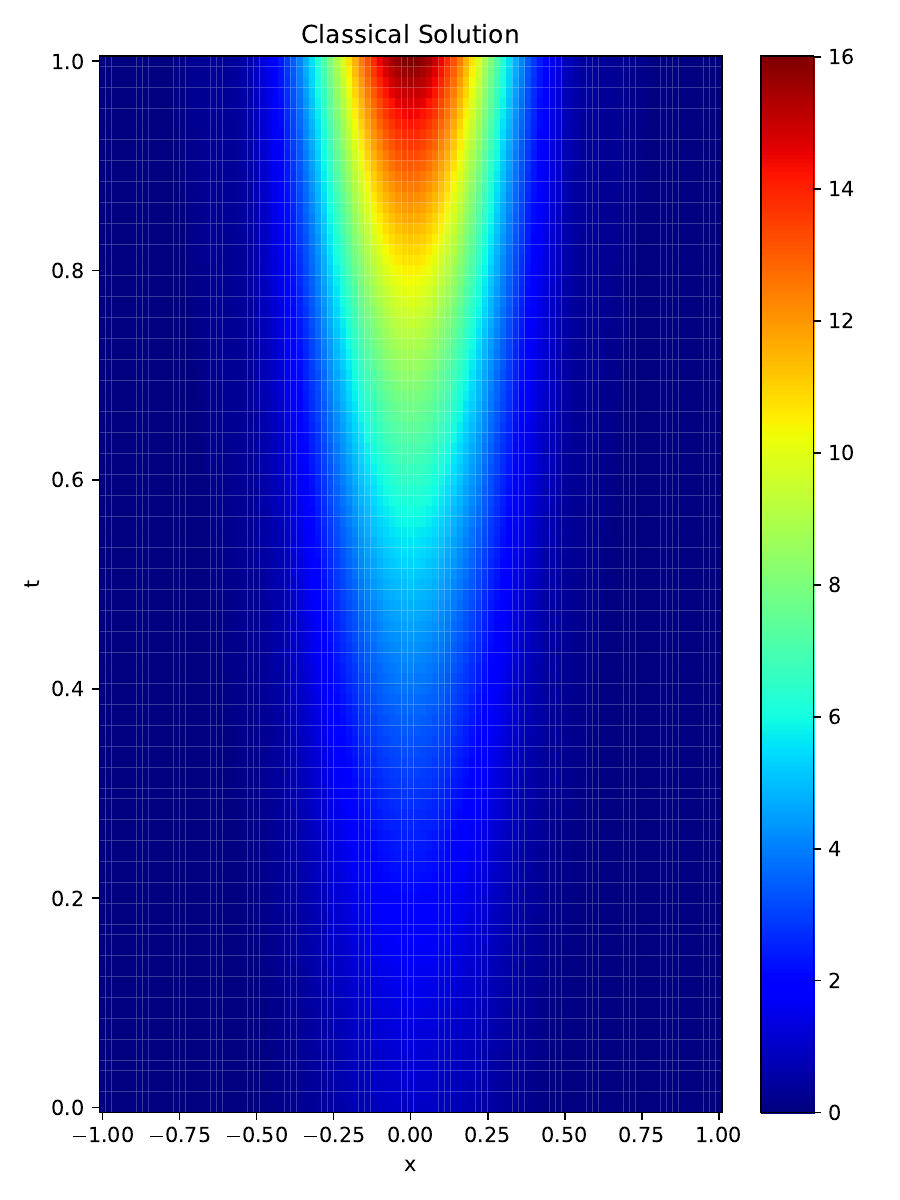}
        \caption{Classical solution}
    \end{subfigure}
\end{figure}

\begin{figure}[H]

\begin{subfigure}{0.50\textwidth}
\includegraphics[width=0.9\linewidth, height=5cm]{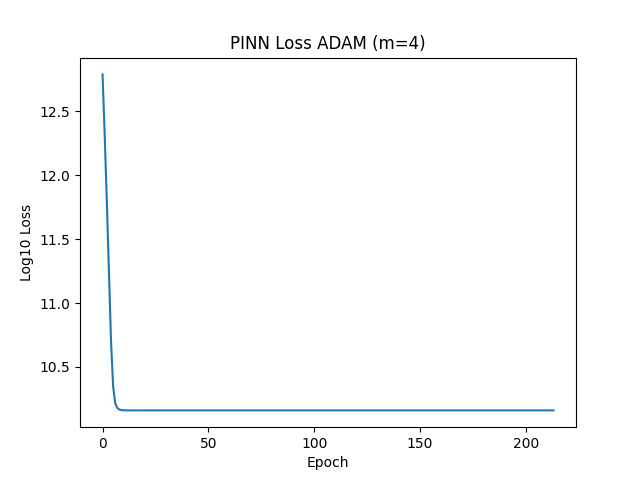} 
\caption{PINN ADAM loss}
\label{fig:subim1}
\end{subfigure}
\begin{subfigure}{0.50\textwidth}
\includegraphics[width=0.9\linewidth, height=5cm]{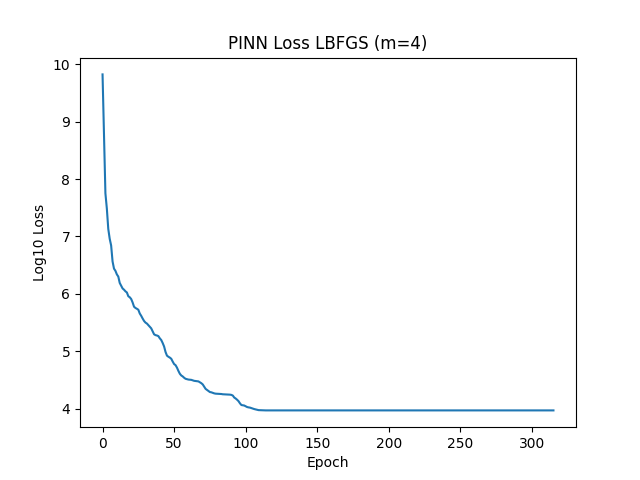}
\caption{PINN LBFGS loss }
\label{fig:subim2}
\end{subfigure}

\end{figure}

\hrule
\vspace{3mm}

\textbf{m = 5}

\begin{figure}[H]
    \centering
    
    \begin{subfigure}[b]{0.32\textwidth}
        \centering
        \includegraphics[width=\textwidth, height =6.5 cm]{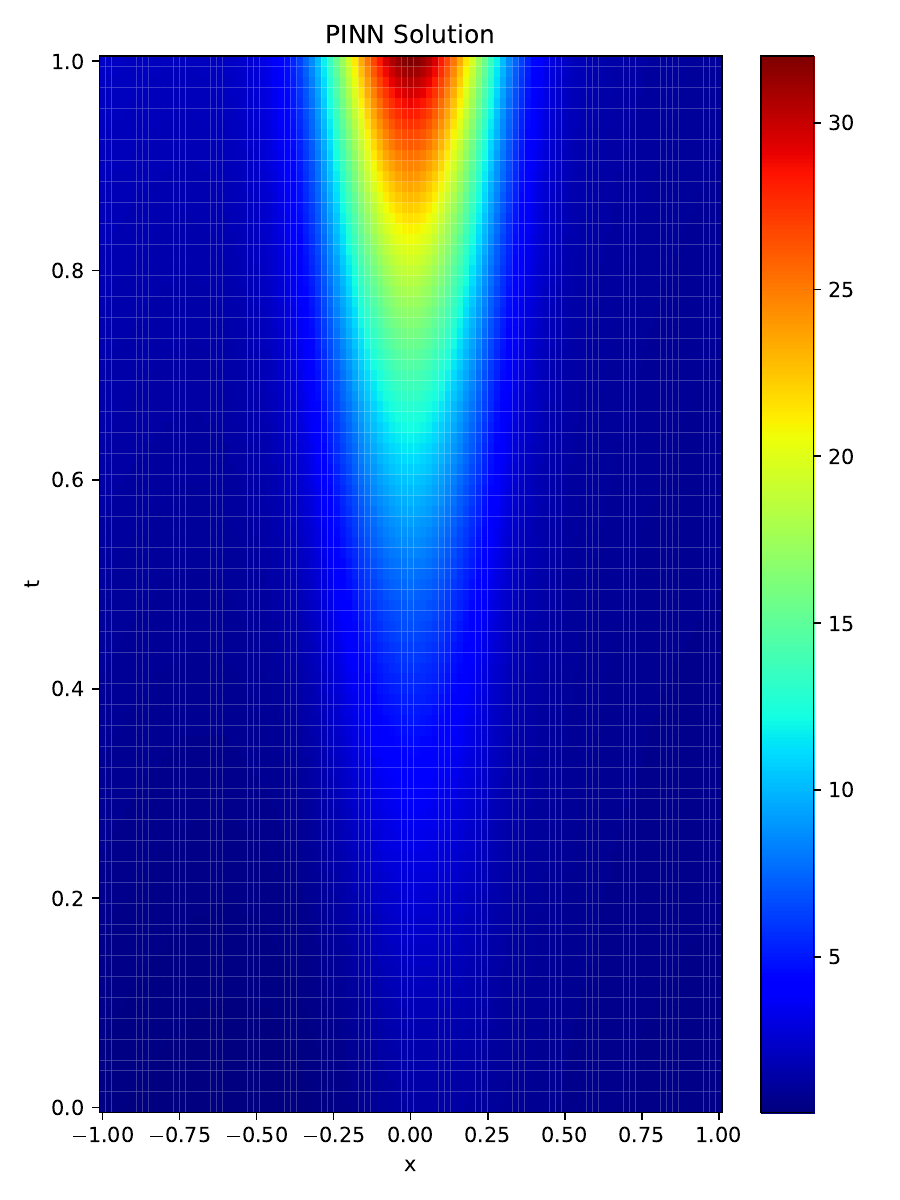}
        \caption{PINN solution}
    \end{subfigure}
    \hfill
    \begin{subfigure}[b]{0.32\textwidth}
        \centering
        \includegraphics[width=\textwidth, height = 6.5 cm]{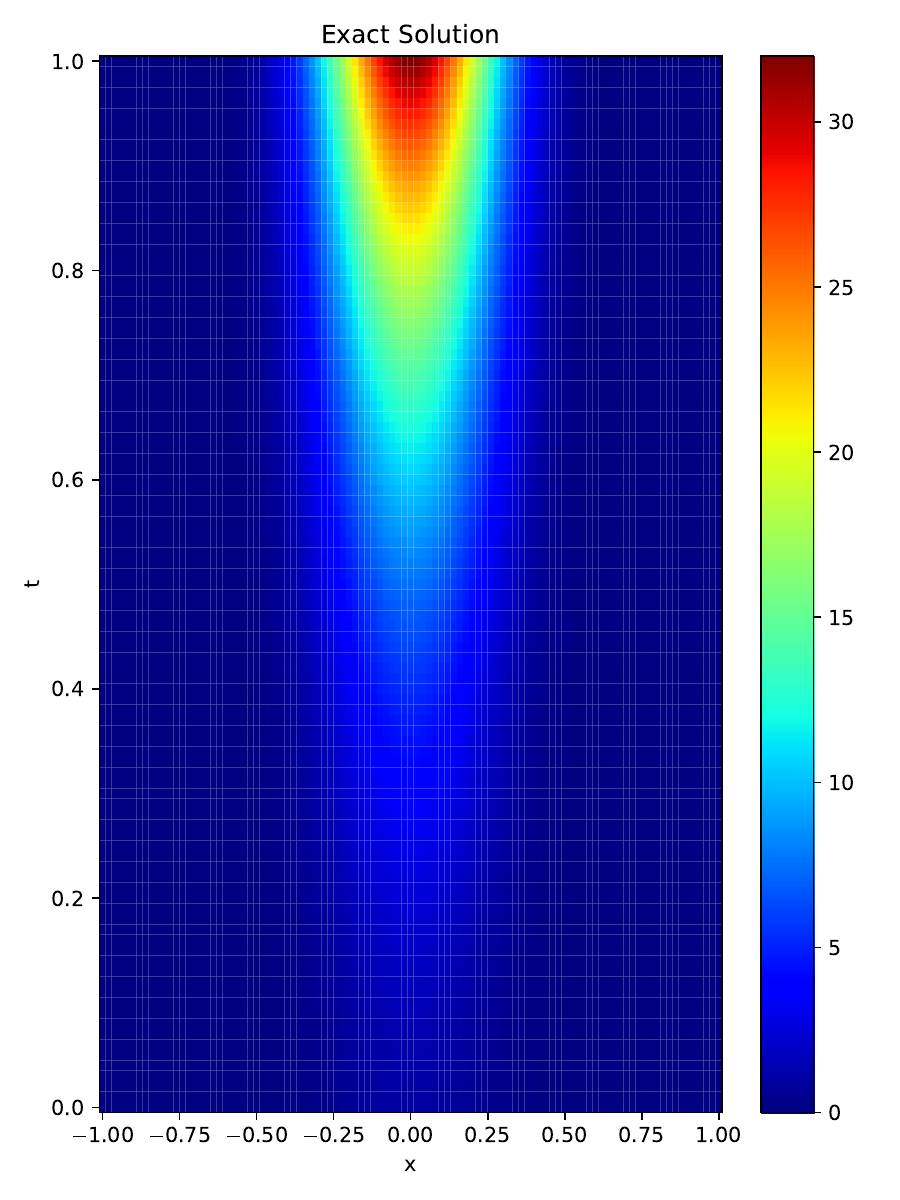}
        \caption{Exact solution}
    \end{subfigure}
    \hfill
    \begin{subfigure}[b]{0.32\textwidth}
        \centering
        \includegraphics[width=\textwidth, height = 6.5 cm ]{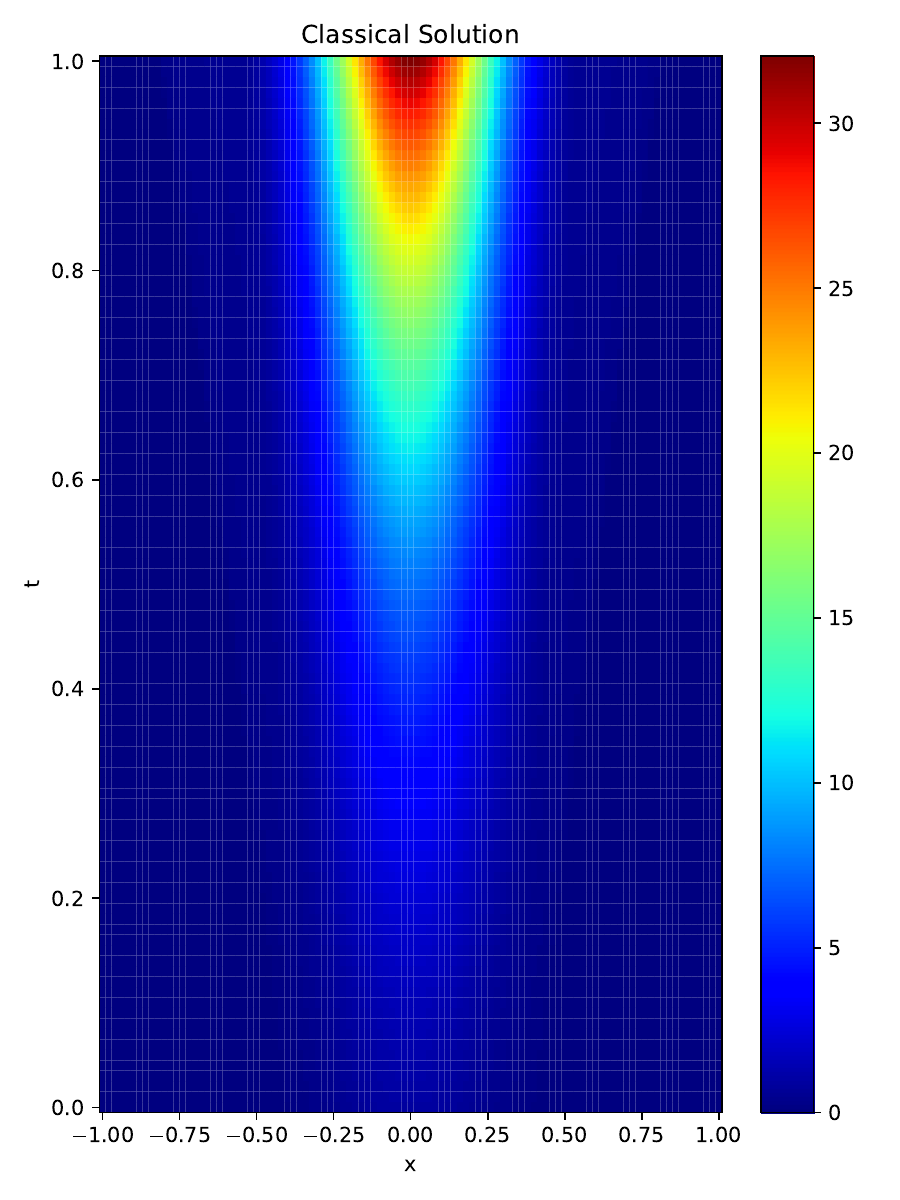}
        \caption{Classical solution}
    \end{subfigure}
\end{figure}

\begin{figure}[H]

\begin{subfigure}{0.50\textwidth}
\includegraphics[width=0.9\linewidth, height=5cm]{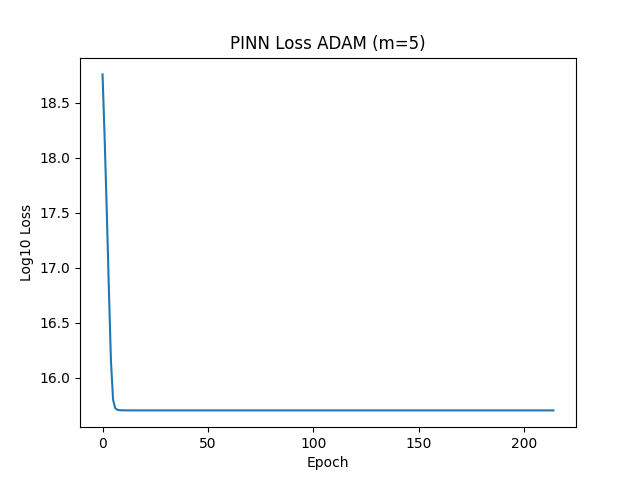} 
\caption{PINN ADAM loss}
\label{fig:subim1}
\end{subfigure}
\begin{subfigure}{0.50\textwidth}
\includegraphics[width=0.9\linewidth, height=5cm]{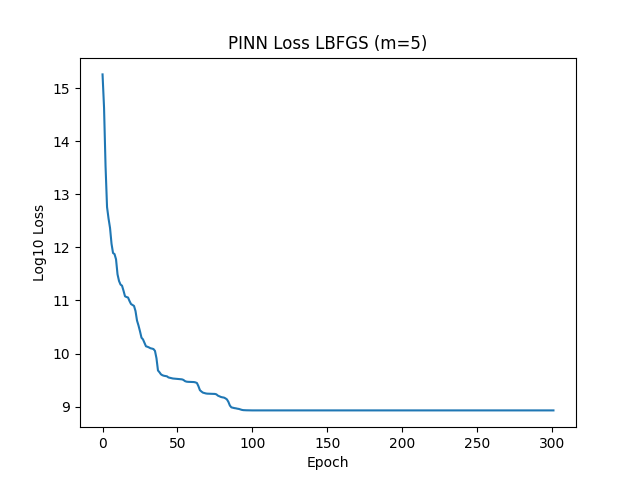}
\caption{PINN LBFGS loss }
\label{fig:subim2}
\end{subfigure}

\end{figure}

\section{Inverse PINN's training loss of the m-dependent solution} 
\label{appendix B}

For the inverse PINN of the m-dependent solution, for each \( m \) and Initial Guess (IG), we show the evolution of the m parameter and the training loss over the epochs. 

\begin{figure}[H]

\begin{subfigure}{0.50\textwidth}
\includegraphics[width=0.9\linewidth, height=5cm]{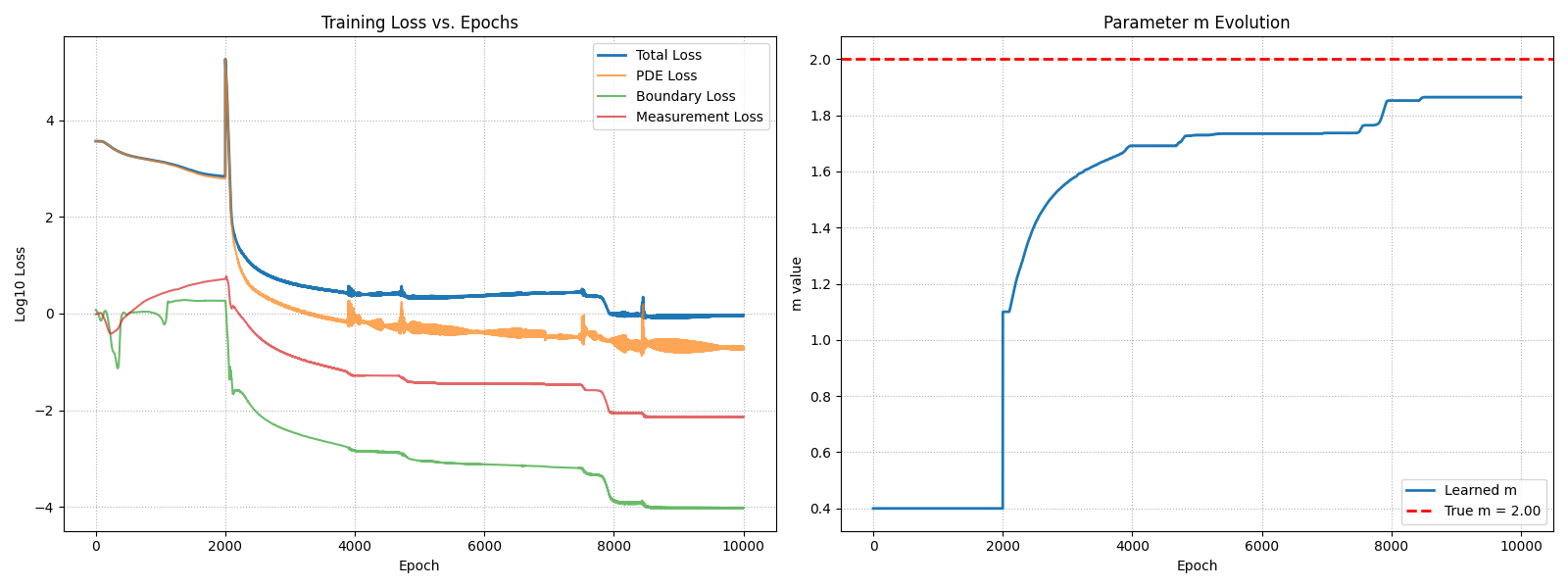} 
\caption*{m = 2 and IG = 0.40}
\label{fig:subim1}
\end{subfigure}
\begin{subfigure}{0.50\textwidth}
\includegraphics[width=0.9\linewidth, height=5cm]{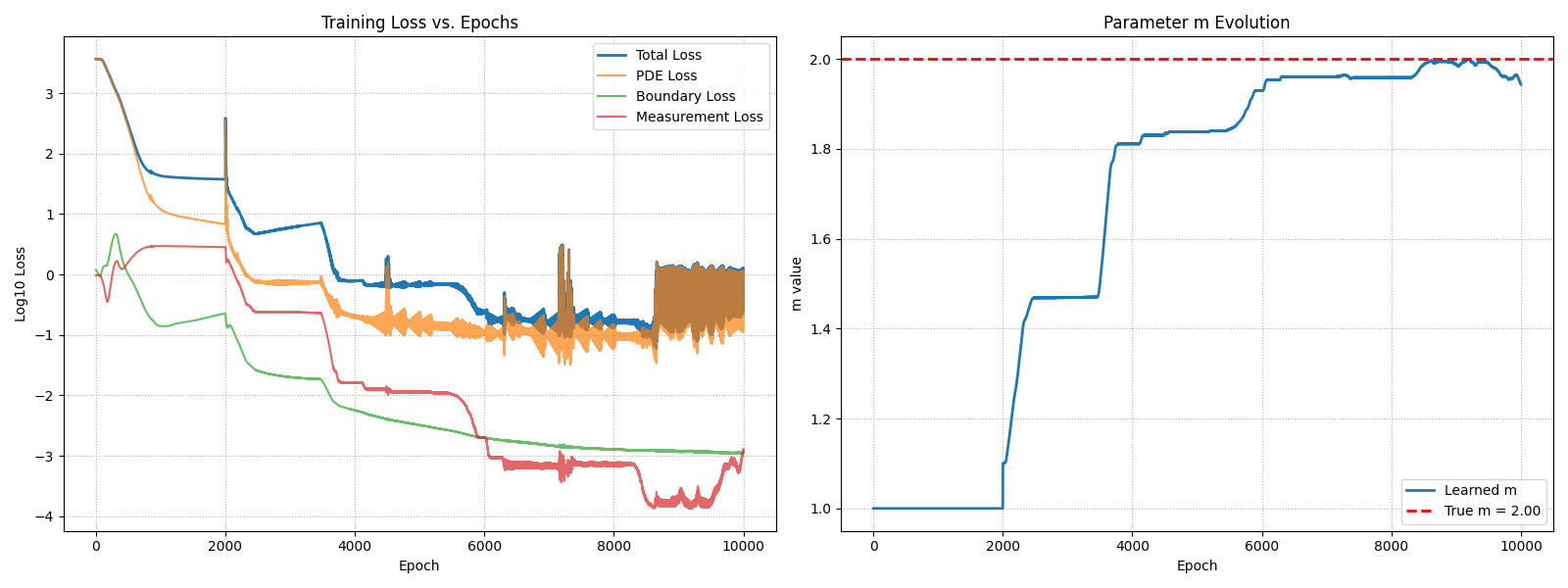}
\caption*{m = 2 and IG = 1.00}
\label{fig:subim2}
\end{subfigure}\vspace{4mm}

\begin{subfigure}{0.50\textwidth}
\includegraphics[width=0.9\linewidth, height=5cm]{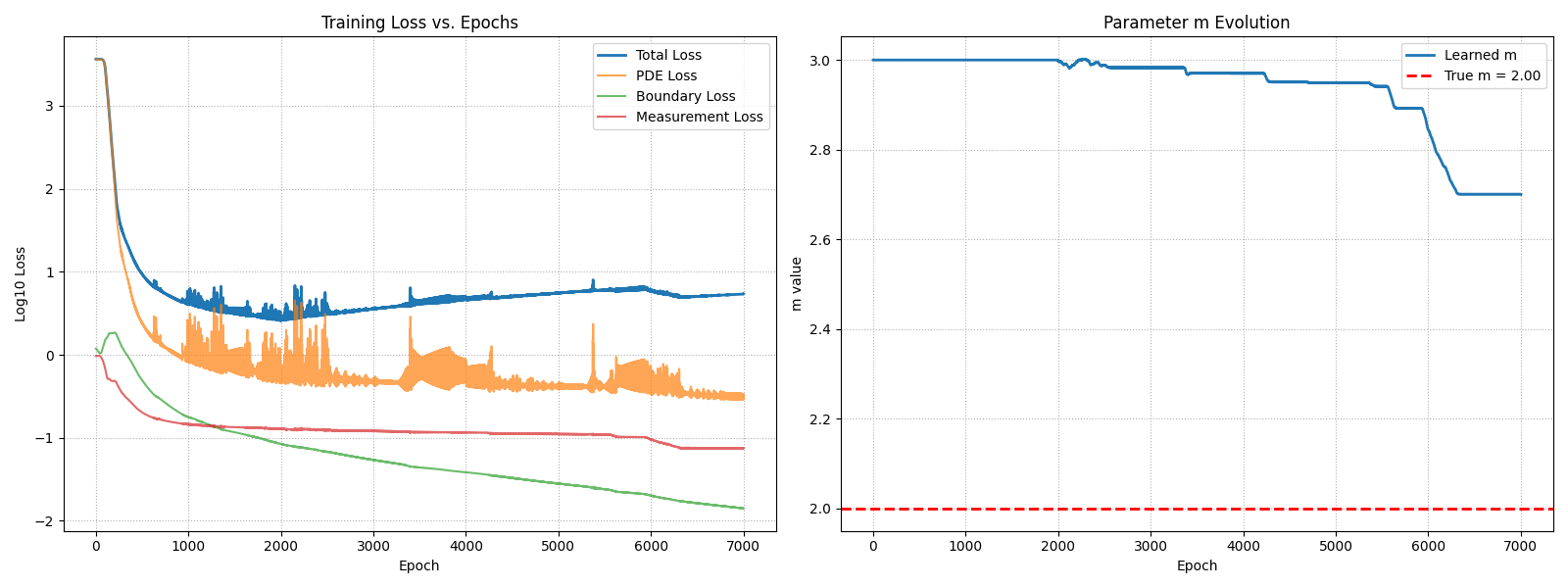} 
\caption*{m = 2 and IG = 3.00}
\label{fig:subim1}
\end{subfigure}
\begin{subfigure}{0.50\textwidth}
\includegraphics[width=0.9\linewidth, height=5cm]{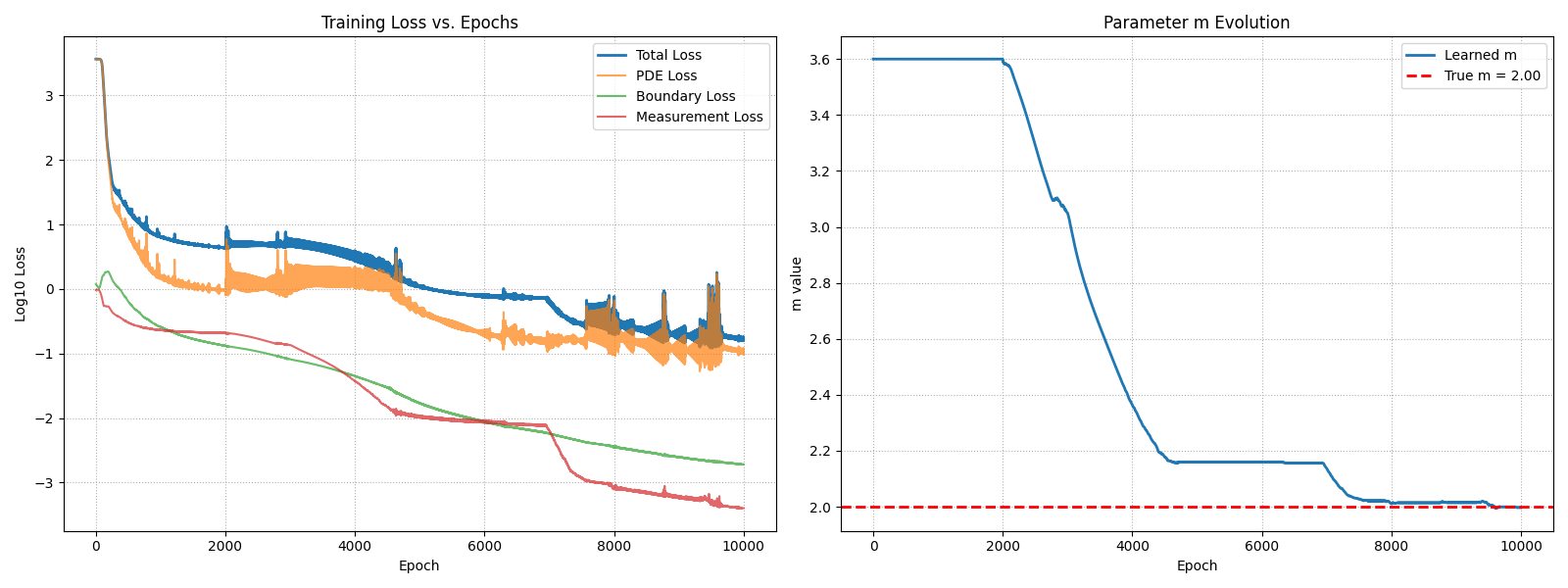}
\caption*{m = 2 and IG = 3.60}
\label{fig:subim2}
\end{subfigure}\vspace{4mm}

\begin{subfigure}{0.50\textwidth}
\includegraphics[width=0.9\linewidth, height=5cm]{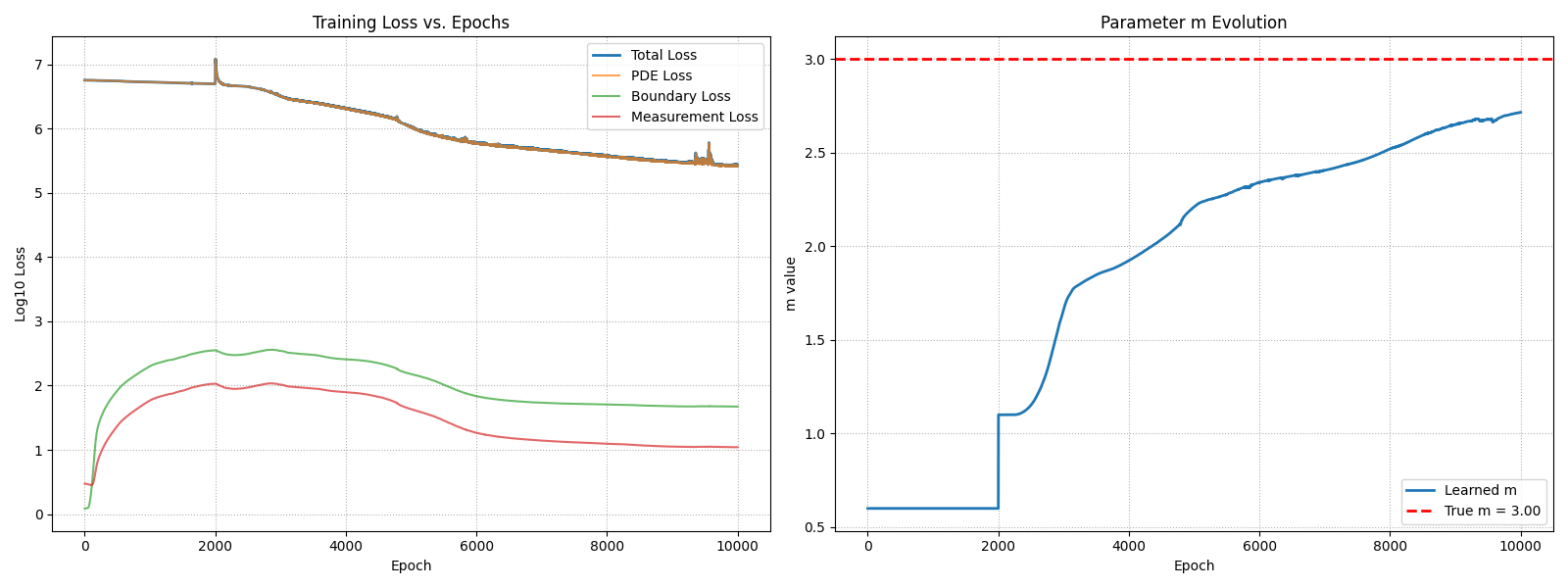} 
\caption*{m = 3 and IG = 0.60}
\label{fig:subim1}
\end{subfigure}
\begin{subfigure}{0.50\textwidth}
\includegraphics[width=0.9\linewidth, height=5cm]{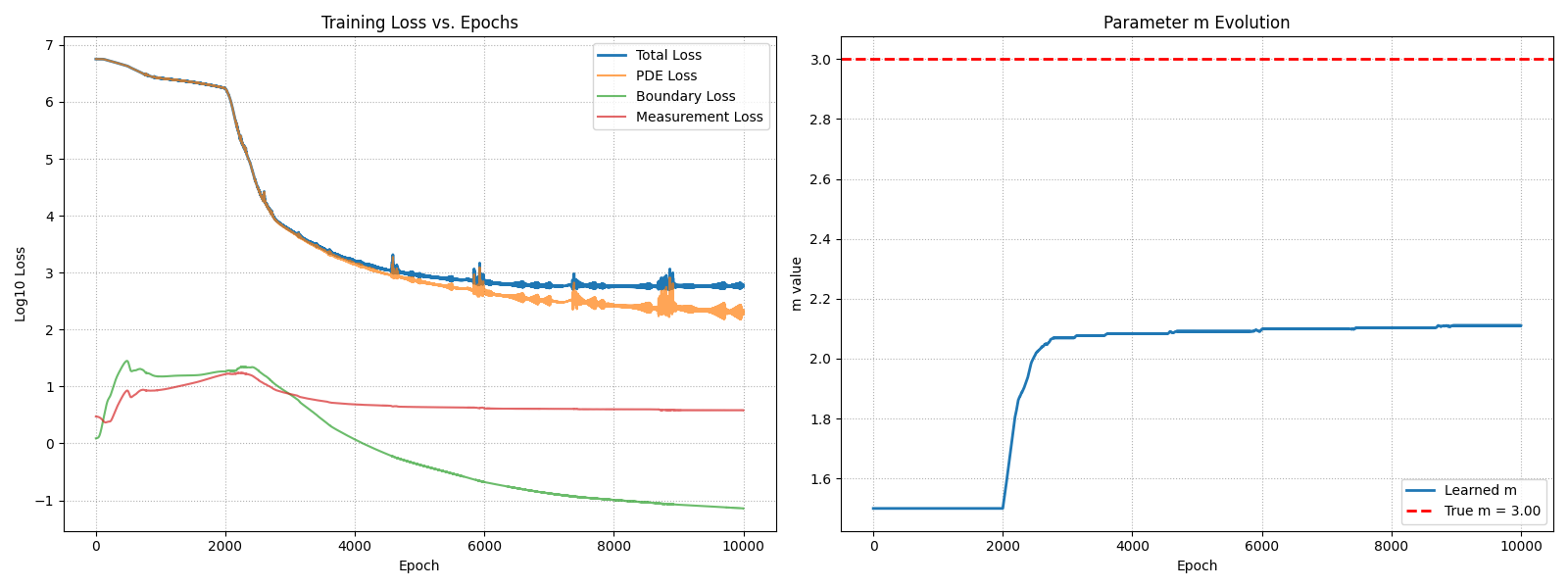}
\caption*{m = 3 and IG = 1.50}
\label{fig:subim2}
\end{subfigure} \vspace{4mm}

\begin{subfigure}{0.50\textwidth}
\includegraphics[width=0.9\linewidth, height=5cm]{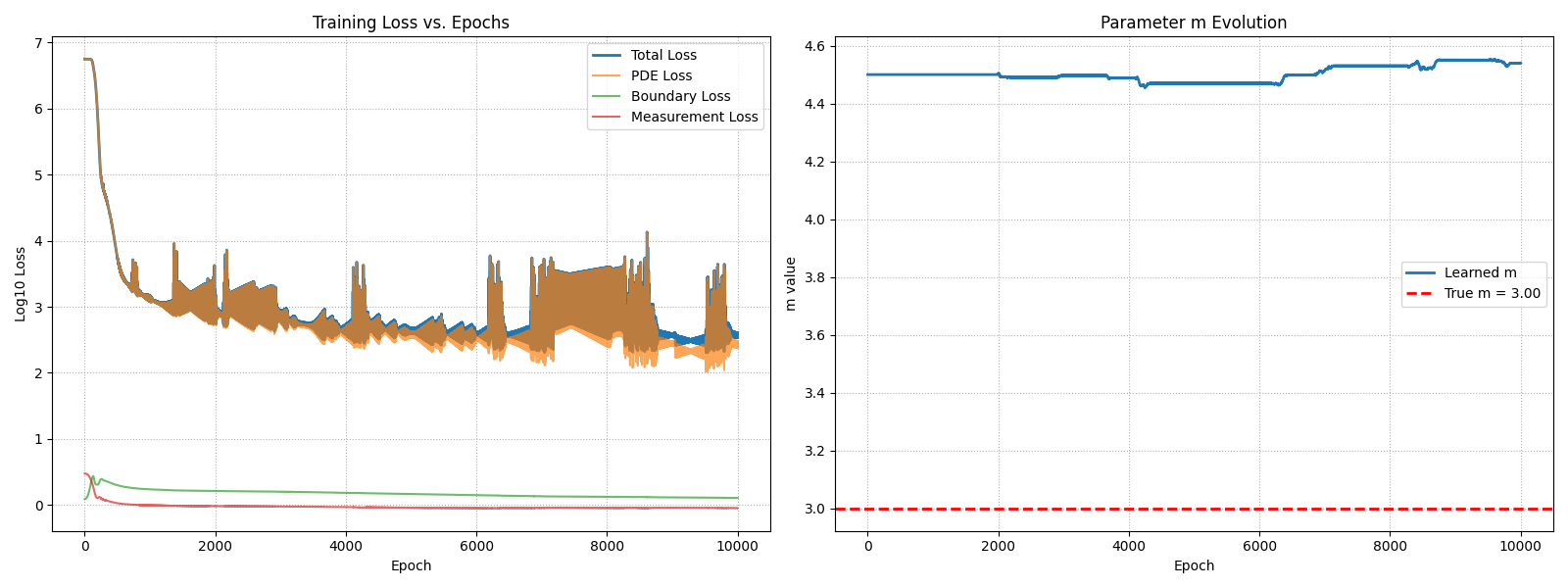} 
\caption*{m = 3 and IG = 4.50}
\label{fig:subim1}
\end{subfigure}
\begin{subfigure}{0.50\textwidth}
\includegraphics[width=0.9\linewidth, height=5cm]{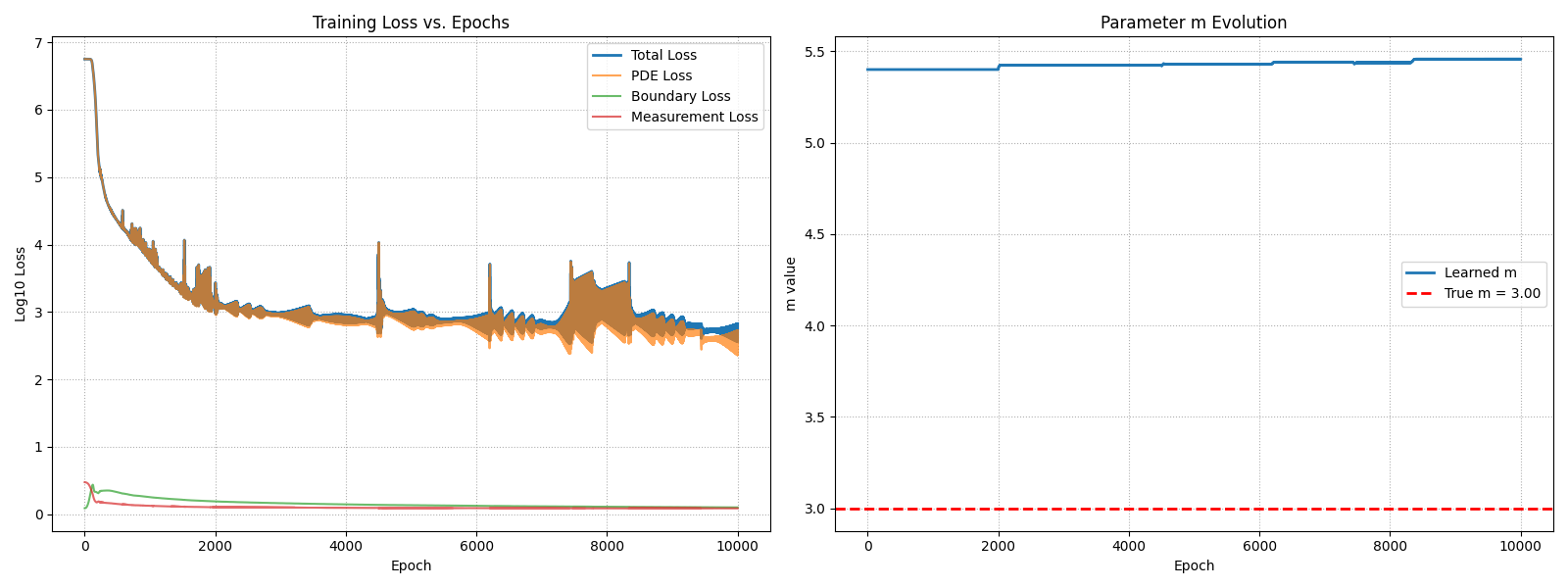}
\caption*{m = 3 and IG = 5.40}
\label{fig:subim2}
\end{subfigure}
\end{figure}

\begin{figure}[H]

\begin{subfigure}{0.50\textwidth}
\includegraphics[width=0.9\linewidth, height=5cm]{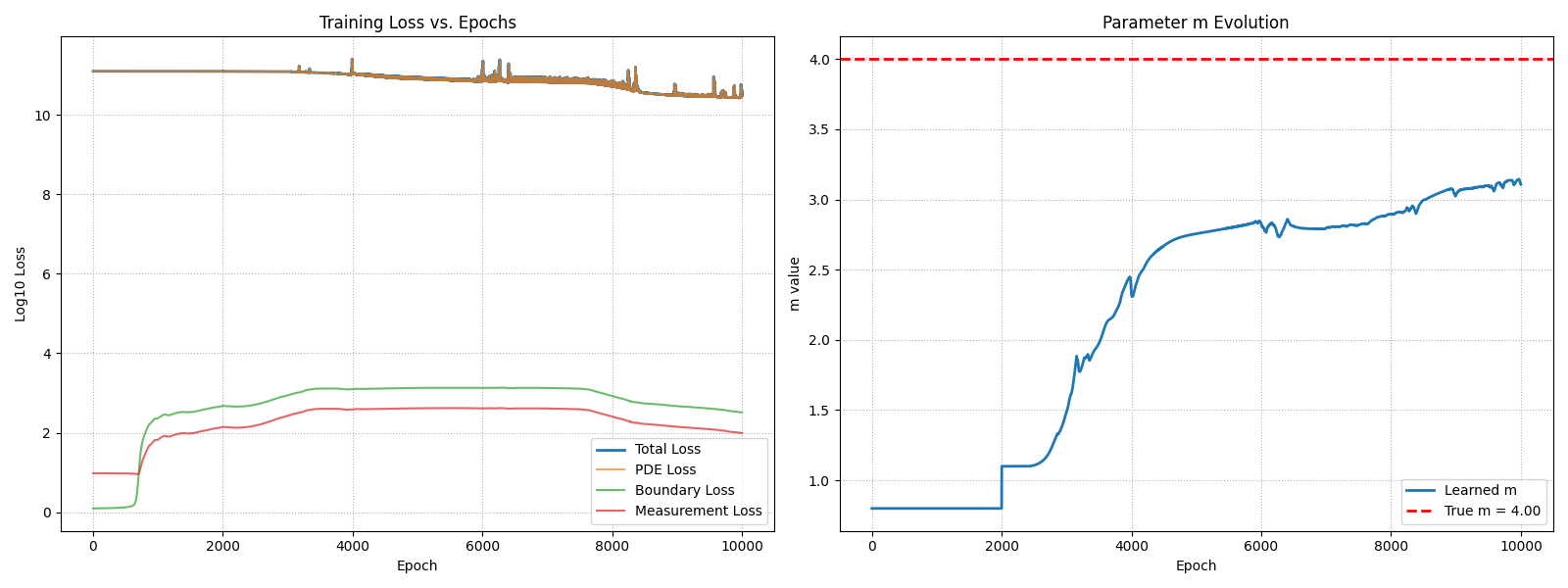} 
\caption*{m = 4 and IG = 0.80}
\label{fig:subim1}
\end{subfigure}
\begin{subfigure}{0.50\textwidth}
\includegraphics[width=0.9\linewidth, height=5cm]{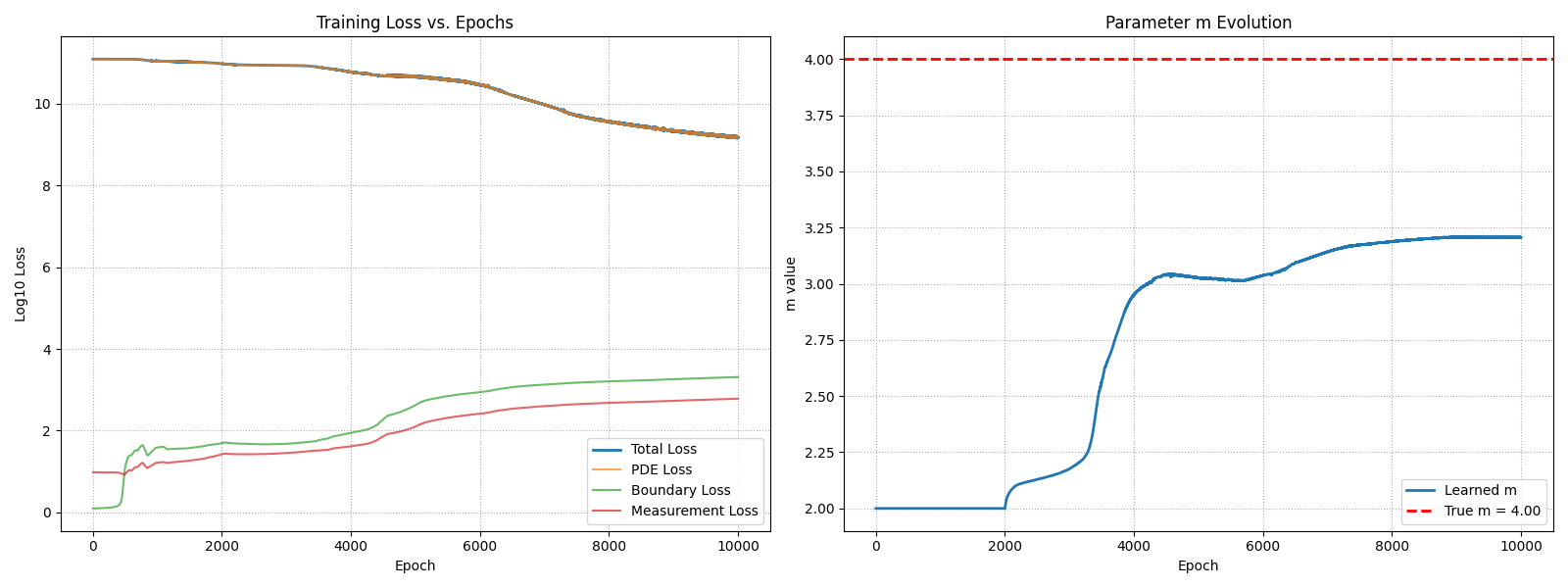}
\caption*{m = 4 and IG = 2.00}
\label{fig:subim2}
\end{subfigure}
\vspace{4mm}

\begin{subfigure}{0.50\textwidth}
\includegraphics[width=0.9\linewidth, height=5cm]{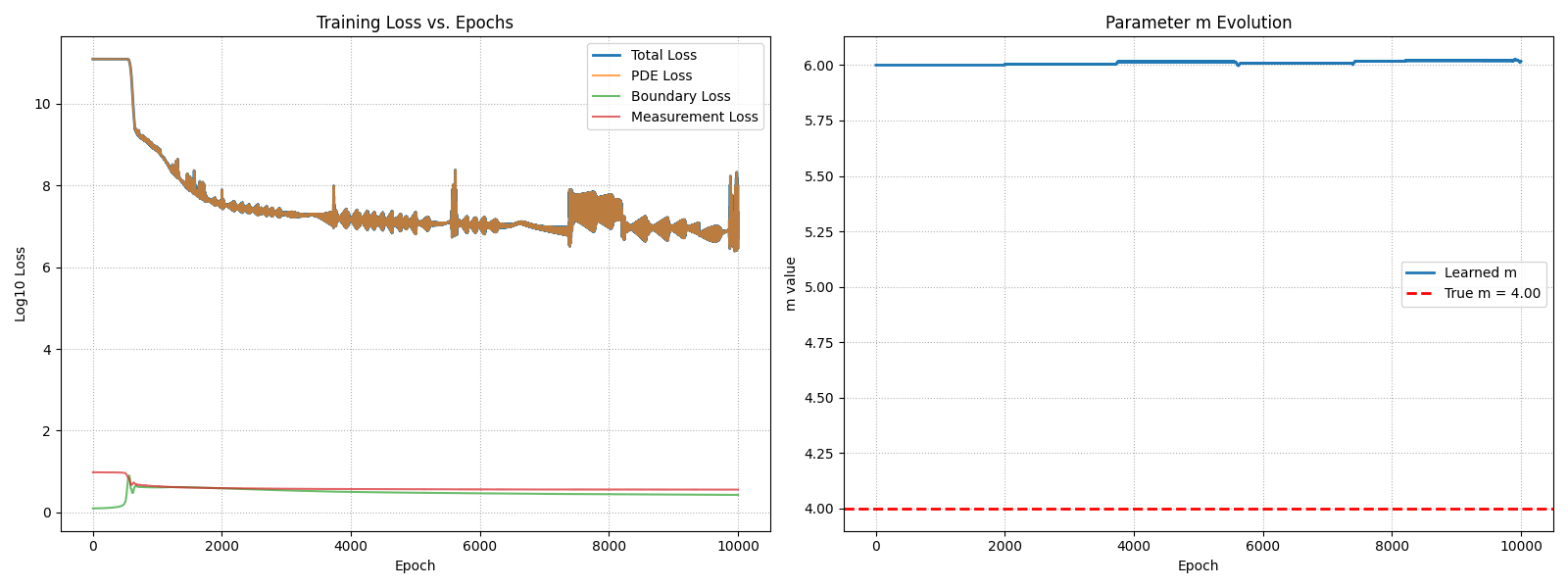} 
\caption*{m = 4 and IG = 6.00}
\label{fig:subim1}
\end{subfigure}
\begin{subfigure}{0.50\textwidth}
\includegraphics[width=0.9\linewidth, height=5cm]{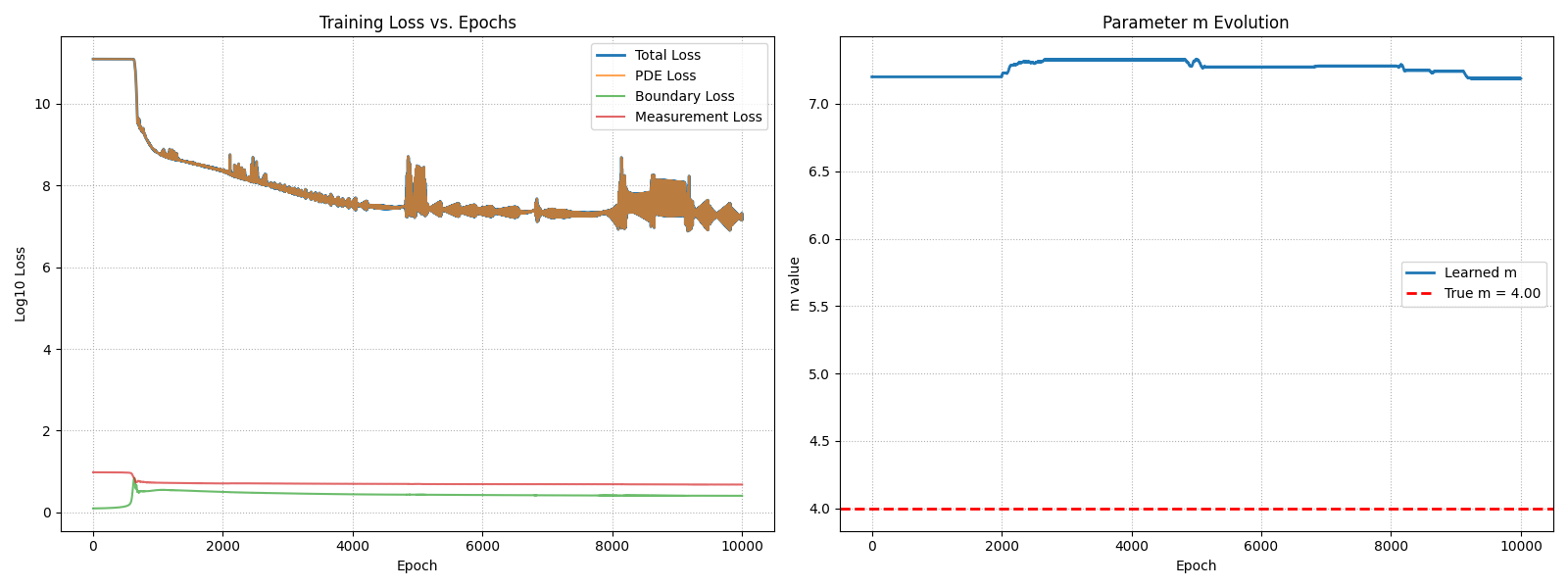}
\caption*{m = 4 and IG = 7.20}
\label{fig:subim2}
\end{subfigure}
\vspace{4mm}

\begin{subfigure}{0.50\textwidth}
\includegraphics[width=0.9\linewidth, height=5cm]{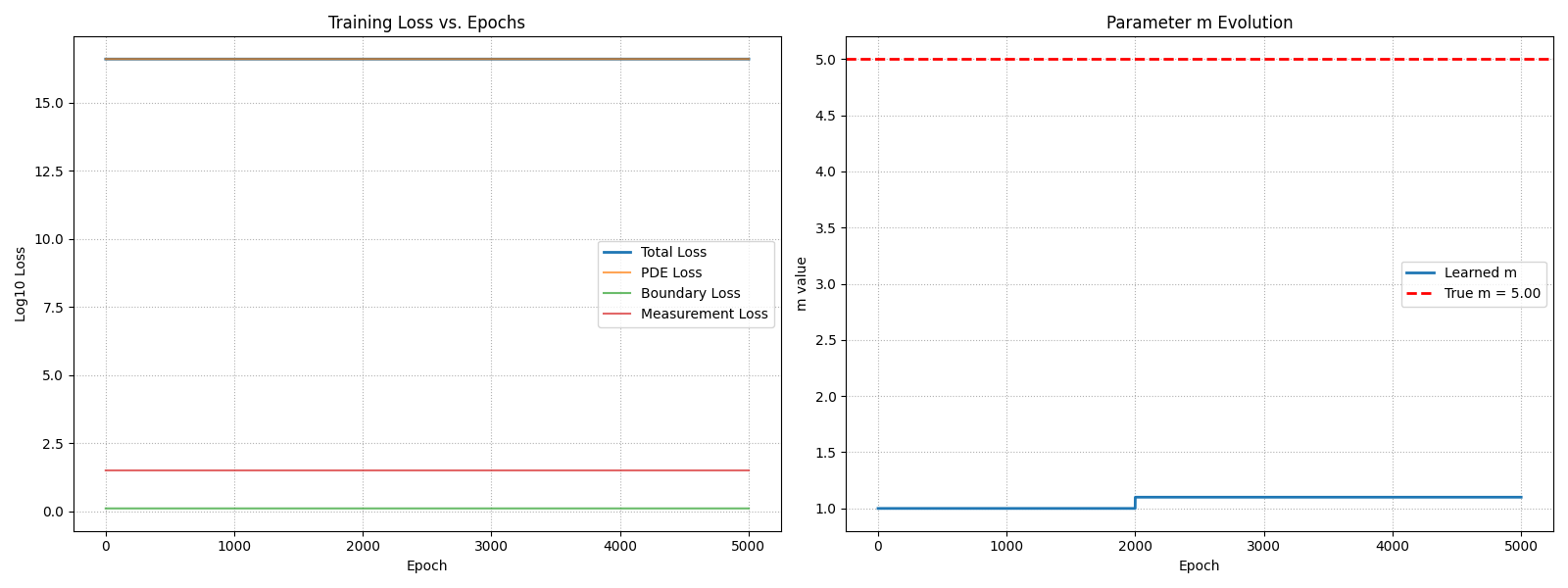} 
\caption*{m = 5 and IG = 1.00}
\label{fig:subim1}
\end{subfigure}
\begin{subfigure}{0.50\textwidth}
\includegraphics[width=0.9\linewidth, height=5cm]{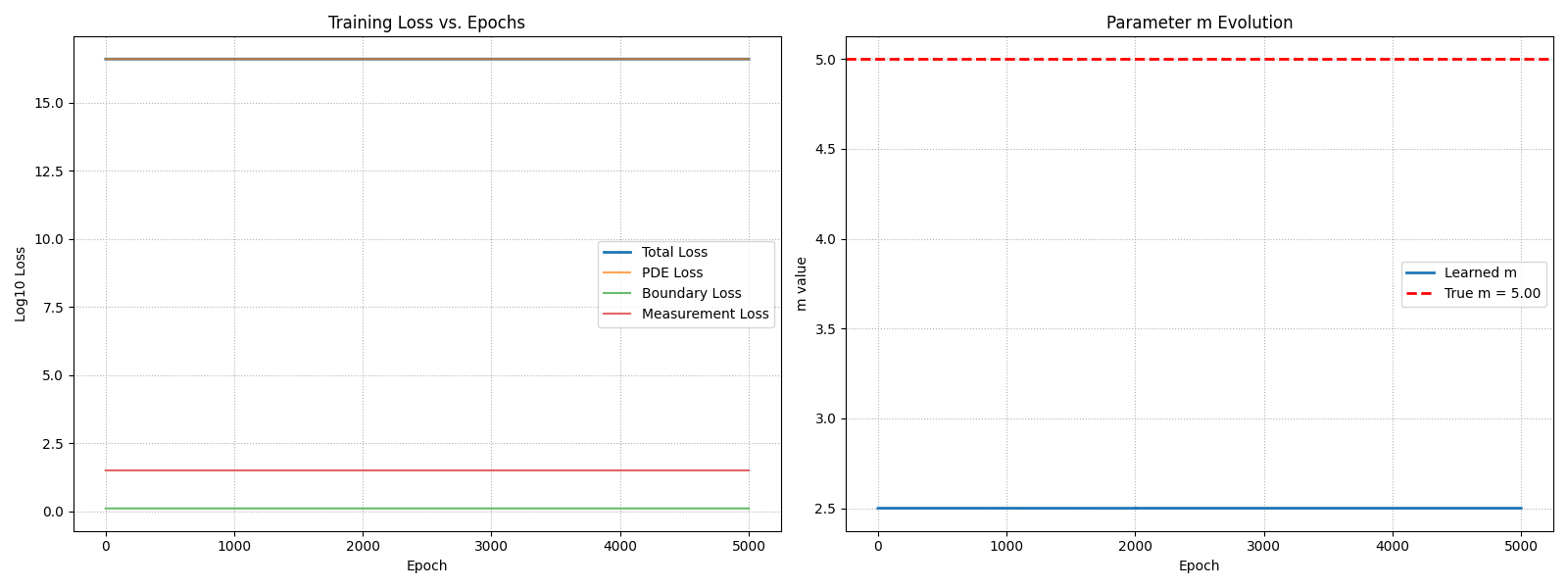}
\caption*{m = 5 and IG = 2.50}
\label{fig:subim2}
\end{subfigure}

\vspace{4mm}

\begin{subfigure}{0.50\textwidth}
\includegraphics[width=0.9\linewidth, height=5cm]{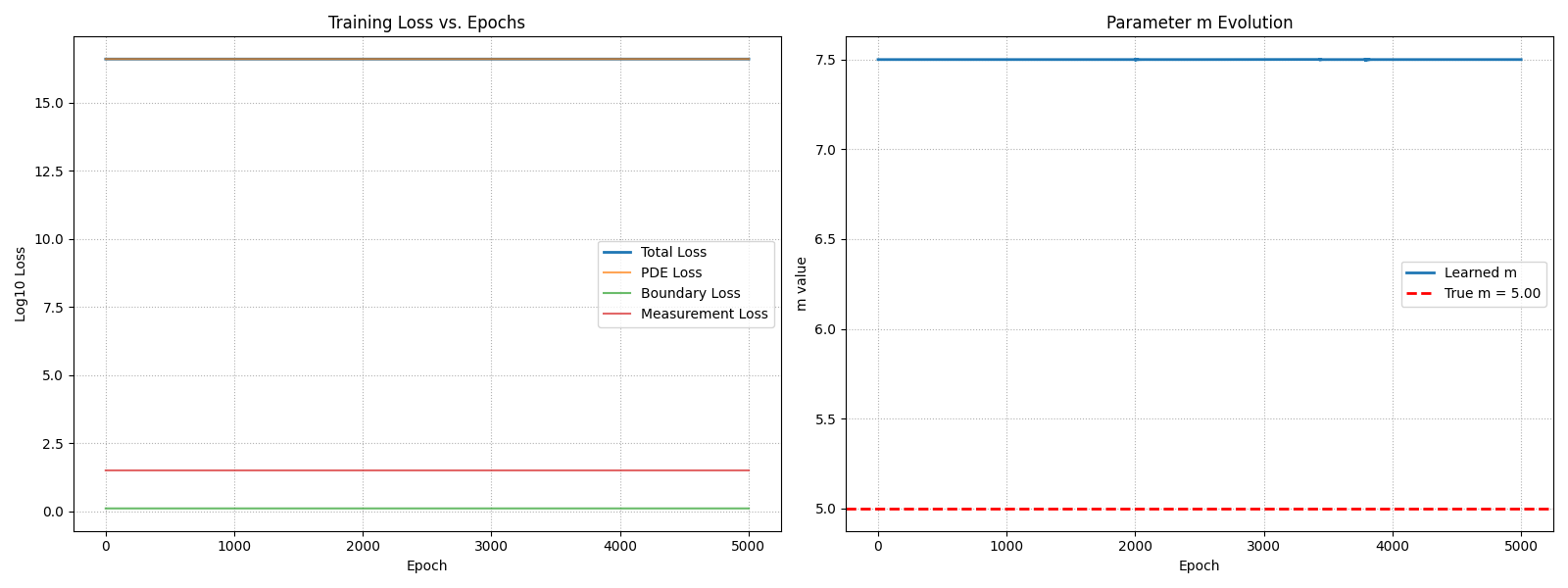} 
\caption*{m = 5 and IG = 7.50}
\label{fig:subim1}
\end{subfigure}
\begin{subfigure}{0.50\textwidth}
\includegraphics[width=0.9\linewidth, height=5cm]{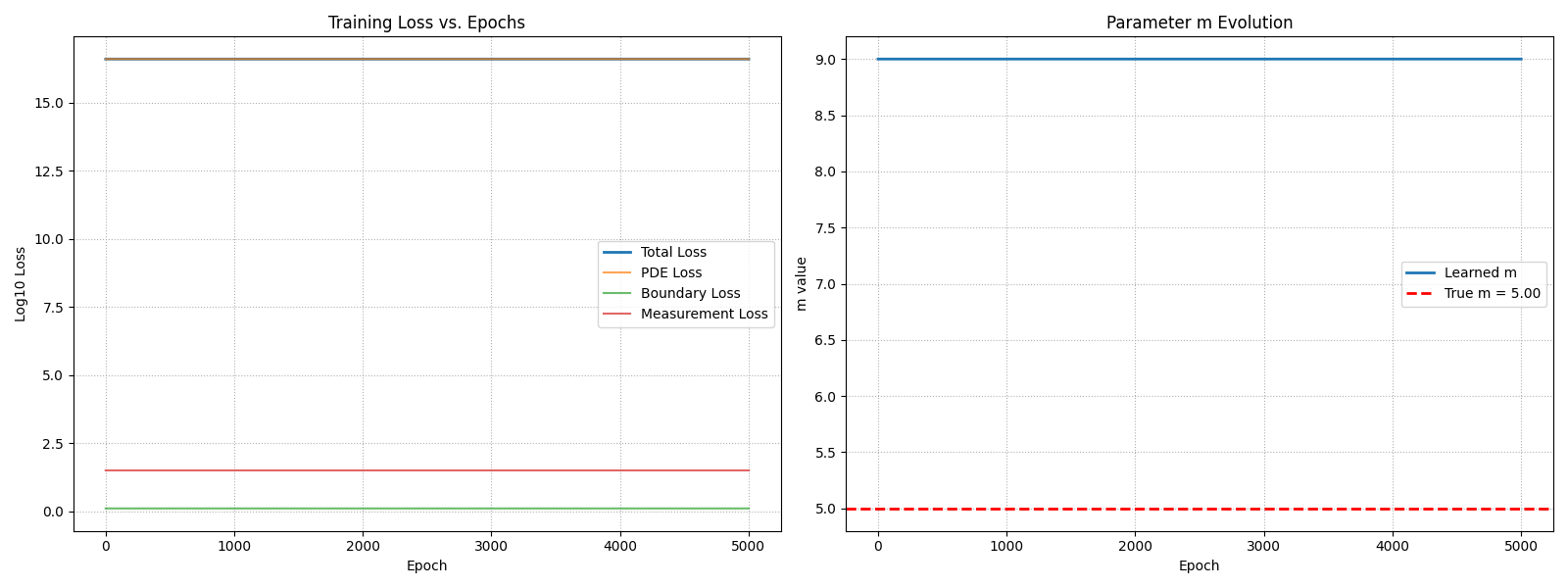}
\caption*{m = 5 and IG = 9.00}
\label{fig:subim2}
\end{subfigure}

\end{figure}

\section{Inverse PINN results before the 2-stage PINN model}
\label{app:C}

In this section, the results of the inverse PINN problem after the standard modifications (i.e., before the 2-stage PINN model) are shown here for the different reference solutions from section \ref{sec: ref sol} following the same plot layout of figure \ref{fig: 2stage}.

\begin{figure}[H]

\begin{subfigure}{0.50\textwidth}
\includegraphics[width=0.9\linewidth, height=5cm]{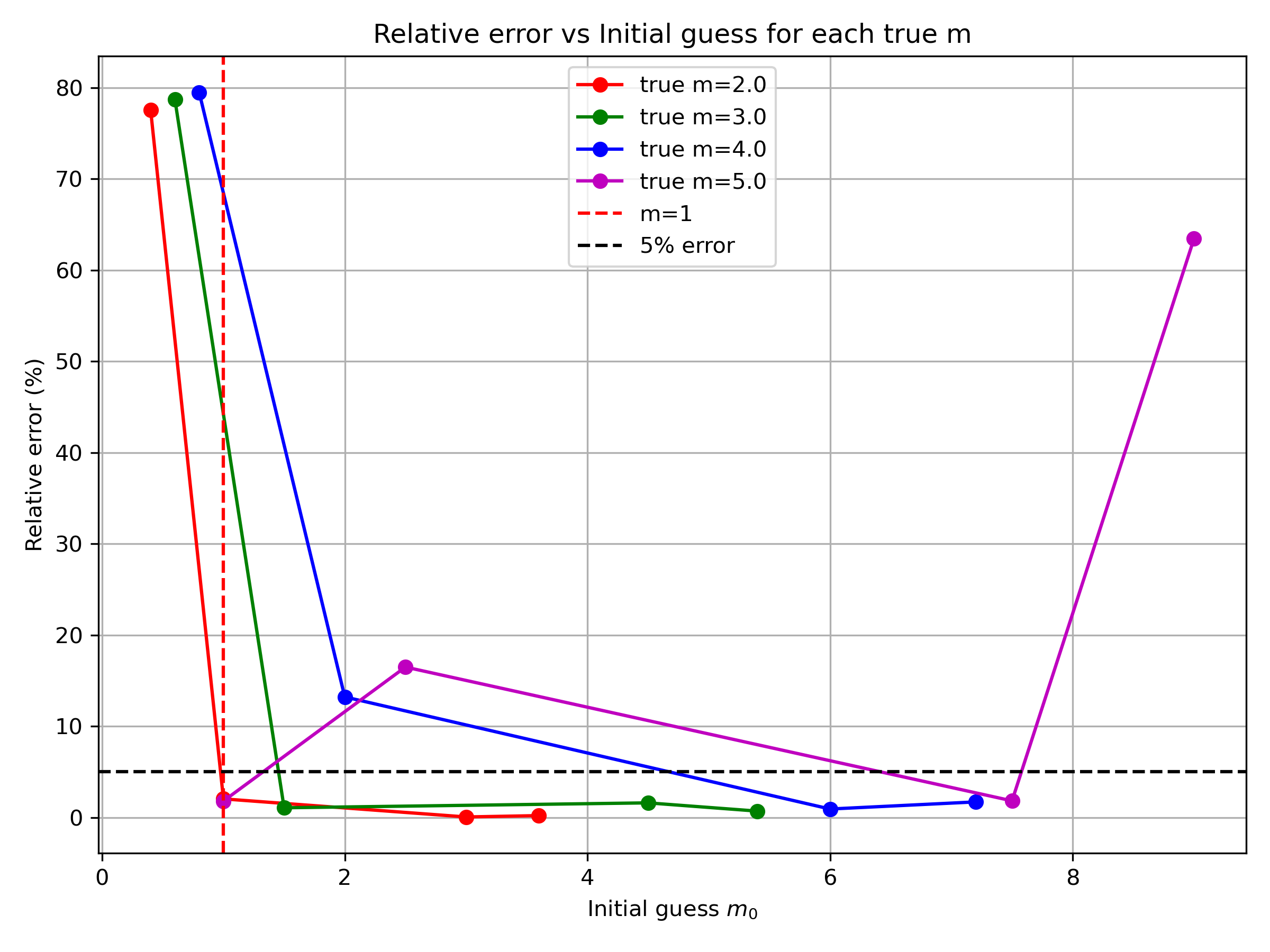} 
\caption*{Barenblatt solution}
\label{fig:subim1}
\end{subfigure}
\begin{subfigure}{0.50\textwidth}
\includegraphics[width=0.9\linewidth, height=5cm]{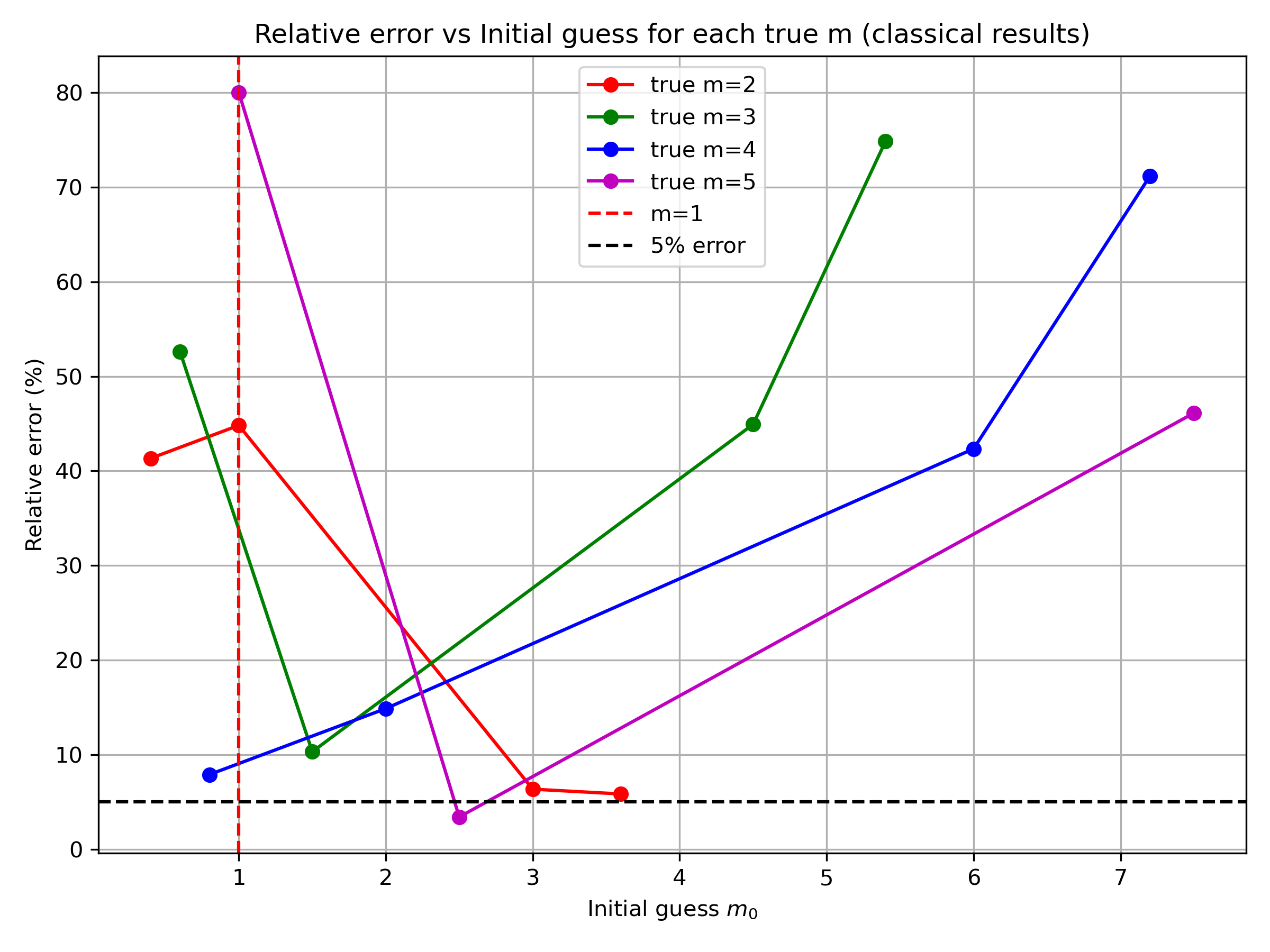}
\caption*{m-dependent solution}
\label{fig:subim2}
\end{subfigure}

\end{figure}

\begin{figure}[H]

\begin{subfigure}{0.50\textwidth}
\includegraphics[width=0.9\linewidth, height=5cm]{PINN_inv_plots_before_2stage/H_IPINN.png} 
\caption*{Harmonic solution}
\label{fig:subim1}
\end{subfigure}
\begin{subfigure}{0.50\textwidth}
\includegraphics[width=0.9\linewidth, height=5cm]{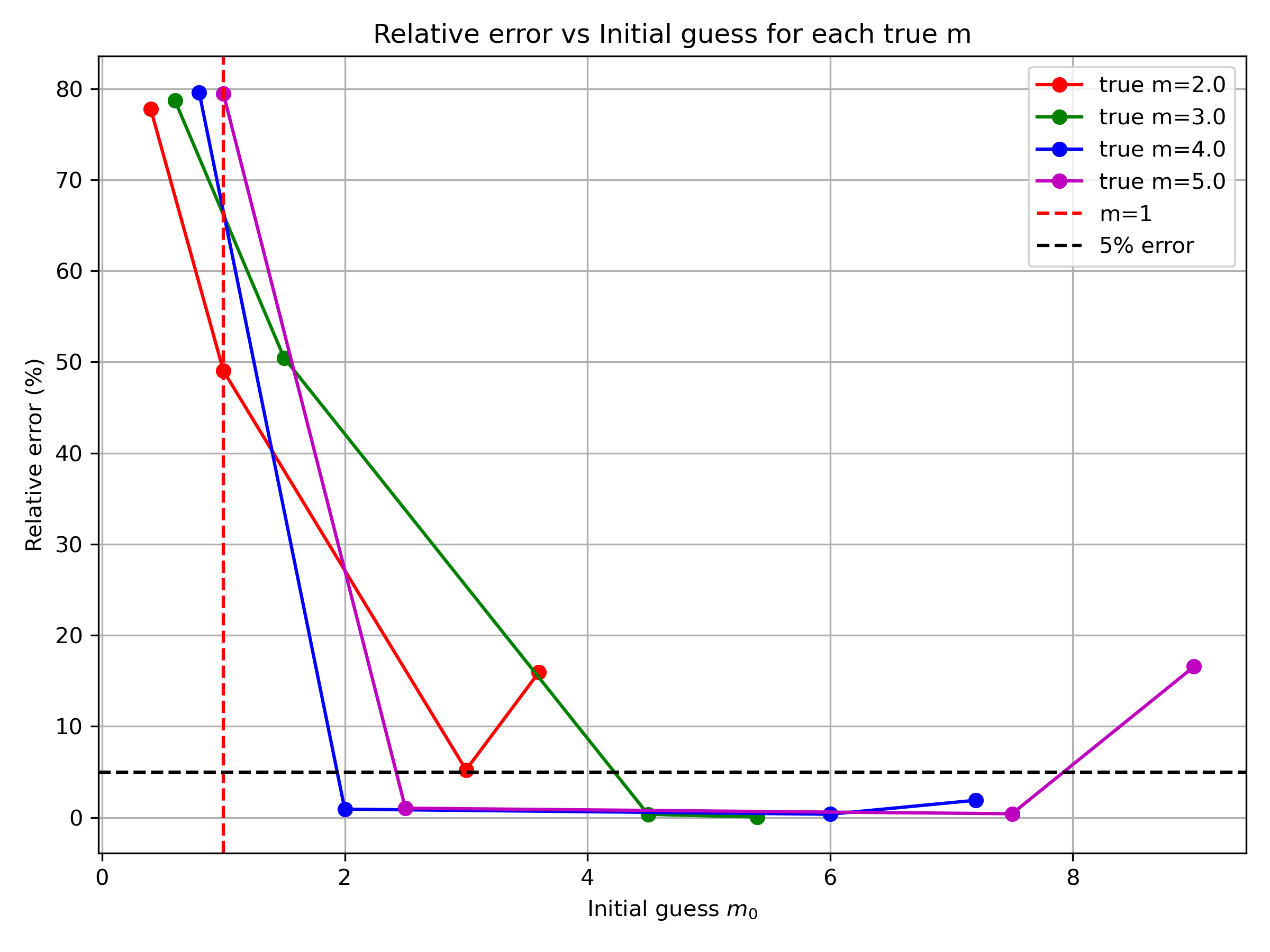}
\caption*{Polynomial solution}
\label{fig:subim2}
\end{subfigure}

\end{figure}

\end{document}